\documentclass[11pt,a4]{book}
\usepackage{indentfirst}
\usepackage{amssymb}
\usepackage{amscd}
\newtheorem{theo}{THEOREM}[section]
\newtheorem{p}[theo]{PROPOSITION}
\newtheorem{lem}[theo]{LEMMA}
\newtheorem{de}[theo]{DEFINITION}
\newtheorem{co}[theo]{COROLLARY}
\newtheorem{e}[theo]{EXAMPLE}
\newtheorem{ax}[theo]{AXIOM}
\newcommand{\fr}[1]{\ensuremath{\mathfrak{#1}}}
\newcommand{\frc}{\ensuremath{\mathfrak{C}_E\,}}
\newcommand{\frm}{\ensuremath{\mathfrak{M}_E\,}}
\newcommand{\frcc}{\ensuremath{\mathfrak{M}_{\bc}\,}}
\newcommand{\h}[1]{Theorem {\ref{#1}}}
\newcommand{\pr}[1]{Proposition {\ref{#1}}}
\newcommand{\dd}[1]{Definition {\ref{#1}}}
\newcommand{\cor}[1]{Corollary {\ref{#1}}}
\newcommand{\lm}[1]{Lemma {\ref{#1}}}
\newcommand{\ee}[1]{Example \ensuremath{\ref{#1}}}
\newcommand{\axi}[1]{Axiom {\ref{#1}}}
\newcommand{\qedd}{\eqno{\rule{3mm}{3mm}}}
\newcommand{\qed}{\hfill {\rule{3mm}{3mm}}}
\newcommand{\mae}[5]{\ensuremath{#1:#2\longrightarrow #3\,, \quad #4 \longmapsto  #5}}
\newcommand{\mad}[4]{\ensuremath{#1\longrightarrow #2, \quad #3 \longmapsto   #4}}
\newcommand{\mac}[3]{\ensuremath{#1:#2\longrightarrow #3}}
\newcommand{\n}[1]{\ensuremath{\left\| #1 \right\|}}
\newcommand{\me}[2]{\ensuremath{\left\{\,\left. #1\;\right|\; #2 \, \right\}}}
\newcommand{\si}[1]{\ensuremath{\sum\limits_{#1}}}
\newcommand{\sii}[2]{\ensuremath{\sum\limits_{#1}^{#2}}}
\newcommand{\pro}[1]{\ensuremath{\prod\limits_{#1}}}
\newcommand{\proo}[2]{\ensuremath{\prod\limits_{#1}^{#2}}}
\newcommand{\prom}{\ensuremath{\prod\limits_{i=1}^{m}}}
\newcommand{\ssa}[1]{\ensuremath{\mathcal S(#1)}}
\newcommand{\f}[2]{\ensuremath{\mathcal F(#1,#2)}}

\newcommand{\ccb}[2]{\ensuremath{\mathcal C\left(#1,#2\right)}}
\newcommand{\cbb}[2]{\ensuremath{\mathcal C_0\left(#1,#2\right)}}
\newcommand{\ccc}[1]{\ensuremath{\mathcal{#1}}}

\newcommand{\bk}{\ensuremath{\mathrm{I\! K}}}
\newcommand{\bm}{\ensuremath{\mathrm{I\! M}}}
\newcommand{\bn}{\ensuremath{\mathrm{I\! N}}}
\newcommand{\bbn}{\ensuremath{n\in \mathrm{I\!N}}}
\newcommand{\bnn}[1]{\ensuremath{\mathrm{I\! N_{#1}}}}
\newcommand{\br}{\ensuremath{\mathrm{I\! R}}}
\newcommand{\bc}{\ensuremath{\mathrm{I\!\!\! C}}}
\newcommand{\bz}{\ensuremath{\mathrm{\,Z\hspace{-0.65em} Z\,}}}
\newcommand{\bzz}[1]{\ensuremath{\mathrm{\,Z\hspace{-0.65em} Z_{#1}\,}}}
\newcommand{\bt}{\ensuremath{\mathrm{\,I\hspace{-0.62em}T}}}
\newcommand{\btt}{\ensuremath{\mathrm{\,I\hspace{-0.5em}T}}}
\newcommand{\bbb}{\ensuremath{\mathrm{I\! B}}}
\newcommand{\bp}{\ensuremath{\mathrm{I\! P}}}
\newcommand{\bs}{\ensuremath{\mathrm{\,S\hspace{-0.65em} S\,}}}
\newcommand{\z}[1]{\ensuremath{\mathrm{\{{#1}\}}}}
\newcommand{\ab}[4]{\ensuremath{ \left\{ \begin{array}{c@{\quad \emph{if} \quad}c} #1&#2 \\
#3&#4 \end{array}} \right.}
\newcommand{\abb}[4]{\ensuremath{ \left\{ \begin{array}{c@{\quad \mbox{if} \quad}c} #1&#2 \\
#3&#4 \end{array}} \right.}

\newcommand{\acc}[6]{\ensuremath{ \left\{ \begin{array}{c@{\quad \mbox{if} \quad}c} #1&#2 \\
#3&#4 \\ #5&#6 \end{array}} \right.}
\newcommand{\ad}[8]{\ensuremath{ \left\{ \begin{array}{c@{\quad \emph{if} \quad}c} #1&#2 \\
#3&#4 \\ #5&#6 \\ #7&#8 \end{array}} \right.}
\renewcommand{\labelenumi}{\alph{enumi})}
\renewcommand{\labelenumii}{$\alph{enumi}_{\arabic{enumii}}$)}
\newcommand{\pp}[1]{\ensuremath{Pr\;#1}}

\newcommand{\unn}[1]{\ensuremath{Un\;#1}}
\newcommand{\unm}[1]{\ensuremath{Un_{\,0}\;#1}}

\newcommand{\unhh}[2]{\ensuremath{Un_{E_{#2}}\;#1}}
\newcommand{\zn}[1]{\ensuremath{un\;#1}}

\newcommand{\zh}[1]{\ensuremath{un_{E}\,#1}}
\newcommand{\oa}[3]{\ensuremath{#1\stackrel{#2}{\longrightarrow }#3}}
\newcommand{\oaa}[3]{\ensuremath{#1\stackrel{#2}{\rightarrow }#3}}
\newcommand{\oag}[4]{\ensuremath{#1\stackrel{#2}{\rule[2.5pt]{#4pt}{0.3pt}\!\!\longrightarrow }#3}}
\newcommand{\ob}[5]{\ensuremath{#1\stackrel{#2}{\longrightarrow }#3\stackrel{#4}{\longrightarrow }#5}}
\newcommand{\obb}[5]{\ensuremath{#1\stackrel{#2}{\rightarrow }#3\stackrel{#4}{\rightarrow }#5}}
\newcommand{\oc}[5]{\ensuremath{0\longrightarrow #1\stackrel{#2}{\longrightarrow }#3\stackrel{#4}{\longrightarrow }#5\longrightarrow 0}}
\newcommand{\of}[7]{\ensuremath{0\longrightarrow #1\stackrel{#2}{\rule[2.5pt]{#6pt}{0.3pt}\!\!\longrightarrow }#3\stackrel{#4}{\rule[2.5pt]{#7pt}{0.3pt}\!\!\longrightarrow }#5\longrightarrow 0}}
\newcommand{\occ}[3]{\ensuremath{0\longrightarrow #1\stackrel{#2}{\longrightarrow }#3}}
\newcommand{\ocd}[3]{\ensuremath{#1\stackrel{#2}{\longrightarrow }#3\longrightarrow 0}}
\newcommand{\od}[6]{\ensuremath{0\longrightarrow #1\stackrel{#2}{\longrightarrow }#3{\scriptscriptstyle{\stackrel{#4}{\longrightarrow }\atop\stackrel{#5}{\longleftarrow }}}#6\longrightarrow 0}}
\newcommand{\og}[9]{\ensuremath{0\longrightarrow #1\stackrel{#2}{\rule[2.6pt]{#7pt}{0.3pt}\!\!\!\longrightarrow }#3{\scriptscriptstyle{\stackrel{#4}{\rule[1.8pt]{#8pt}{0.1pt}\!\longrightarrow }\atop\stackrel{#5}{\longleftarrow\!\!\!\!\!\rule[1.8pt]{#9pt}{0.1pt} }}}#6\longrightarrow 0}}

\newcommand{\oddg}[8]{\ensuremath{0\longrightarrow #1\stackrel{#2}{\rule[2.6pt]{#6pt}{0.3pt}\!\!\!\longrightarrow }#3{\scriptscriptstyle{\stackrel{#4}{\rule[1.8pt]{#7pt}{0.1pt}\!\longrightarrow }\atop\stackrel{#5}{\longleftarrow\!\!\!\!\!\rule[1.8pt]{#8pt}{0.1pt} }}}}}
\newcommand{\odeg}[5]{\ensuremath{{\scriptscriptstyle{\stackrel{#1}{\rule[1.8pt]{#4pt}{0.1pt}\longrightarrow }\atop\stackrel{#2}{\longleftarrow\!\!\!\!\!\rule[1.8pt]{#5pt}{0.1pt} }}}#3\longrightarrow 0}}
\newcommand{\eo}{$E$-C*-algebra\,}
\newcommand{\en}{full $E$-C*-algebra\,}
\newcommand{\kk}[2]{\ensuremath{K_{#1}\left({#2}\right)}}
\begin{document}
\frontmatter

\title{AXIOMATIC K-THEORY FOR C*-ALGEBRAS} 
\author{CORNELIU CONSTANTINESCU \\Bodenacherstr. 53 \\CH 8121 Benglen \\ e-mail: constant@math.ethz.ch \\ ETHZ}
\maketitle



\thispagestyle{empty}
\tableofcontents
\vspace{5mm}

\fbox{\parbox{12cm}{Throughout this book $E$ denotes a fixed commutative unital C*-algebra}}

\mainmatter
\thispagestyle{empty}
\begin{center}
\huge{\bf{Preface}}
\end{center}

In Part I we present an axiomatic frame in which many results of the K-theory for C*-algebras can be proved. In Part II we construct an example for this axiomatic theory, which generalizes the classical theory for C*-algebras. This last theory starts by associating to each C*-algebra $F$ the C*-algebras of square matrices with entries in $F$. Every such C*-algebras of square matrices can be obtained as the projective representation of a certain group with respect to a Schur function for this group with values in $\bc$ (\dd{703}). The above mentioned generalization consists in replacing this Schur function by an arbitrary Schur function which satisfies some axiomatic conditions. Moreover this Schur function can take its values in a commutative unital C*-algebra $E$ instead of  $\bc$. In this case this K-theory does not apply to the category of C*-algebras, but to the category of $E$-C*-algebras (\dd{10.3'}), which are C*-algebras endowed with a supplementary  structure (every C*-algebra can be endowed with such a supplementary structure (\pr{16.10})). Up to some definitions and notation Part II is independent of Part I.

\newpage

In general we use the notation and the terminology of [C1]. In the sequel we give a list of notation used in this book.
\renewcommand{\labelenumi}{\arabic{enumi})} 
\begin{enumerate}
\item $\bc$ (respectively $\br$) denotes the field of complex (respectively real) numbers, $\bn$ denotes the set of natural numbers ($0\not\in \bn$), $\bn^*:=\bn\cup \{0\}$, 
$\bz\;$  denotes the group of integers, and for every $n\in \bn^*$ we put $\bnn{n}:=\me{k\in \bn}{k\leq n}$ and $\;\bzz{n}:=\bz\,/(n\bz)$.  
\item For every set $A$, $Card\; A$ denotes the cardinal number of $A$ and $id_A$ denotes the identity map of $A$. If $x$ is a map defined on $A$ and $B$ is a subset of $A$ then $x|B$ denotes the restriction of $x$ to $B$. 
\item Let $(\Omega _j)_{j\in J}$ be a family of topological spaces and let $\Omega $ be the disjoint union of this family. The topological sum of the family $(\Omega _j)_{j\in J}$ is the topological space obtained by endowing $\Omega $ with the topology $\me{U\subset \Omega }{j\in J\Rightarrow U\cap \Omega _j\;\mbox{is\,an\,open\,set\,of}\;\Omega _j}$. 
\item If $\Omega $ is a topological space and $G$ is a C*-algebra then $\ccc{C}(\Omega ,G)$ denotes the C*-algebra of continuous bounded maps of $\Omega $ into $G$ (endowed with the supremum norm). If $\Omega $ is a locally compact space then $\ccc{C}_0(\Omega ,G)$ denotes the C*-algebra of continuous maps of $\Omega $ into $G$ vanishing at the infinity.
\item $\odot $ denotes the algebraic tensor product of vector spaces.
\item $\approx $ means isomorphic.
\end{enumerate}

\renewcommand{\labelenumi}{\alph{enumi})}

\begin{center}
\part{Axiomatic K-theory}
\end{center}

\begin{center}
\fbox{\parbox{9cm}{Throughout Part I we endow $\z{0,1}$ with the structure of a group by identifying it with $\bz_2$ and take $i\in \z{0,1}$}}
\end{center}

\begin{center}
\chapter{The axiomatic theory}
\end{center}

{\center{\section{$E$-C*-algebras}}}

\begin{de}\label{10.3'}
In this book we call {\bf{$E$-C*-algebra}} a C*-algebra $F$ endowed with a bilinear map (exterior multiplication)
$$\mad{E\times F}{F}{(\alpha ,x)}{\alpha x}$$
such that for all $\alpha ,\beta \in E$ and $x,y\in F$,
$$(\alpha +\beta )x=\alpha x+\beta y\,,\quad (\alpha \beta )x=\alpha (\beta x)\,,\quad (\alpha x)^*=\alpha ^*x^*\,,\quad \n{\alpha x}\leq \n{\alpha }\n{x}\,,$$
$$\alpha (x+y)=\alpha x+\alpha y\,,\qquad\alpha (xy)=(\alpha x)y=x(\alpha y)\,,\qquad 1_Ex=x\;.$$
An {\bf $E$-C*-subalgebra ($E$-ideal)} of $F$ is a C*-subalgebra (a closed ideal) $G$ of $F$ such that 
$$(\alpha ,x)\in E\times G\Longrightarrow \alpha x\in G\;.$$
If $F,G$ are $E$-C*-algebras then a C*-homomorphism $\mac{\varphi }{F}{G}$ is called {\bf{$E$-linear}} or an {\bf{$E$-C*-homomorphism}} if for all $(\alpha ,x)\in E\times F$, $\varphi (\alpha x)=\alpha \varphi x$. A bijective $E$-C*-homomorphism is called {\bf{$E$-C*-isomorphism}}.  We denote by $0$ the $E$-C*-algebra having a unique element. We denote by \frm the category of $E$-C*-algebras for which the morphisms are the $E$-linear C*-homomorphisms. In particular $\fr{M}_{\bc}$ is the category of all C*-algebras.

If $G$ is an $E$-ideal of the \eo $F$ then the C*-algebra $F/G$ has a natural structure of an \eo and 
$$\oc{G}{\varphi }{F}{\psi }{F/G}$$
is an exact sequence in \frm, where $\varphi $ denotes the inclusion map and $\psi $ the quotient map. Conversely, if
$$\oc{F}{\varphi }{G}{\psi }{H}$$
is an exact sequence in \frm then $F$ is an $E$-ideal of $G$ and $H\approx G/F$.
\end{de}

\begin{de}\label{8.6}
If $(F_j)_{j\in J}$ is a finite family of $E$-C*-algebras then we denote by $\pro{j\in J}F_j$ the $E$-C*-algebra obtained by
 endowing the corresponding C*-algebra $\pro{j\in J}F_j$ with the bilinear map
$$\mad{E\times \pro{j\in J}F_j}{\pro{j\in J}F_j}{(\alpha ,(x_j)_{j\in J})}{(\alpha x_j)_{j\in J}}\;.$$
\end{de}

\begin{p}\label{16.10}
Every C*-algebra can be endowed with the structure of an \eo. 
\end{p}

Let $F$ be a C*-algebra. Let $\Omega $ be the spectrum of $E$ and $\omega \in \Omega $ and put
$$\mad{E\times F}{F}{(\alpha ,x)}{\alpha (\omega )x}\;.$$
It is easy to see that $F$ endowed with this exterior multiplication is an \eo.\qed

\begin{e}\label{16.10a}
Let $\Omega $ be a finite set and $E:=\ccb{\Omega }{\bc}$.
\begin{enumerate}
\item Let $(F_\omega )_{\omega \in \Omega }$ be a finite family of C*-algebras and $F:=\pro{\omega \in \Omega }F_\omega $. If we put for all $(\alpha ,x)\in E\times F$,
$$\mae{\alpha x}{\Omega }{F}{\omega }{\alpha (\omega )x_\omega }$$
then $F$ endowed with the exterior multiplication
$$\mad{E\times F}{F}{(\alpha ,x)}{\alpha x}$$
is an \eo. 
\item Let $F$ be an \eo and for every $\omega\in \Omega  $ put 
$$\mae{e_\omega }{\Omega }{\bc}{\omega '}{\ab{1}{\omega '=\omega }{0}{\omega '\not=\omega }}\,,$$
$$F_\omega :=\me{e_\omega x}{x\in F}\;.$$
Then $F_\omega $ is a C*-algebra for all $\omega \in \Omega $ and $F\approx \pro{\omega \in \Omega }F_\omega $, with the meaning of a).\qed
\end{enumerate}
\end{e}

\begin{e}\label{16.10b}
Let $\Omega $ be a discrete locally compact space, $\Omega ^*$ a compactification of $\Omega $, $E:=\ccb{\Omega ^*}{\bc}$, $(F_\omega )_{\omega \in \Omega }$ a family of C*-algebras, and $F:=\pro{\omega \in \Omega }F_\omega$ $\left(resp.\; F:=\me{x\in \pro{\omega \in \Omega }F_\omega }{\lim_{\omega \rightarrow \infty }\n{x_\omega }=0}\right)$. If we put for all $(\alpha ,x)\in E\times F$
$$\mae{\alpha x}{\Omega }{F}{\omega }{\alpha (\omega )x_\omega }$$
then $\alpha x\in F$ for all $(\alpha ,x)\in E\times F$ and $F$ endowed with the exterior multiplication
$$\mad{E\times F}{F}{(\alpha ,x)}{\alpha x}$$
is an \eo.\qed
\end{e}

{\center{\section{The axioms}}}

\begin{de}\label{28.9'c}
We denote by $K_0$ and $K_1$ two covariant functors from the category \frm to the category of additive groups. We denote by $0$ the group which has a unique element and call {\bf{K-null}} an $E$-C*-algebra $F$ for which $\kk{i}{F}=0$. Let $\oaa{F}{\varphi }{G}$ be a morphism  in \frm. We say that $\varphi $ is {\bf{K-null}}  if $\kk{i}{\varphi }=0$. We say that $\varphi $ {\bf{factorizes through null}} if there are morphisms $\oaa{F}{\varphi '}{H}$ and $\oaa{H}{\varphi ''}{G}$ in \frm such that $\varphi =\varphi ''\circ \varphi '$ and such that $H$ is K-null.
\end{de}

 We have $\kk{i}{id_F}=id_{\kk{i}{F}}$ for every $E$-C*-algebra $F$. Every morphism which factorizes through null is K-null. 

\begin{ax}[Null-axiom]\label{27.9'}
$\kk{i}{0}=0$.
\end{ax}

\begin{ax}[Split exact axiom]\label{27.9'a}
If 
$$\od{F}{\varphi }{G}{\psi }{\lambda }{H}$$
is a split exact sequence in \frm then
$$\og{\kk{i}{F}}{\kk{i}{\varphi }}{\kk{i}{G}}{\kk{i}{\psi }}{\kk{i}{\lambda }}{\kk{i}{H}}{20}{20}{20}$$
is a split exact sequence in the category of additive groups.
\end{ax}

It follows that the map
$$\mad{\kk{i}{F}\times \kk{i}{H}}{\kk{i}{G}}{(a,b)}{\kk{i}{\varphi }a+\kk{i}{\lambda }b}$$
is a group isomorphism.

\begin{de}\label{939}
Let $\mac{\varphi ,\psi }{F}{G}$ be morphisms in \frm. We say that $\varphi $ and $\psi $ are {\bf homotopic} if there is a path
$$\mac{\phi  _s}{F}{G},\qquad\qquad s\in [0,1]$$
of morphisms in \frm \,such that $\phi _0=\varphi ,\,\phi _1=\psi $, and the map
$$\mad{[0,1]}{G}{s}{\phi _sx}$$
is continuous for every $x\in F$. 

We say that a pair $\oa{F}{\varphi }{G}$, $\oa{G}{\psi }{F}$ of morphisms in \frm\, is a {\bf homotopy} if $\psi \circ \varphi $ is homotopic to $id_F$ and $\varphi \circ \psi $ is homotopic to $id_G$. In this case we say that $F$ and $G$ are {\bf homotopic}. $F$ is called {\bf null-homotopic} if it is homotopic to the $E$-C*-algebra $0$.
\end{de}

\begin{ax}[Homotopy axiom]\label{27.9'b}
If $\mac{\varphi ,\psi }{F}{G}$ are homotopic morphisms in \frm then $\kk{i}{\varphi }=\kk{i}{\psi }$.
\end{ax}

\begin{de}\label{27.9'c}
We associate to every exact sequence 
$$\oc{F}{\varphi }{G}{\psi }{H}$$
 in \frm two group homomorphisms (called \bf{index maps})
$$\mac{\delta _i}{\kk{i}{H}}{\kk{{i+1}}{F}}\;.$$
\end{de}

\begin{ax}[Six-term axiom]\label{27.9'd}
For every exact sequence in \frm
$$\oc{F}{\varphi }{G}{\psi }{H}$$
the six-term sequence
$$\begin{CD}
K_0(F)@>K_0(\varphi )>>K_0(G)@>K_0(\psi )>>K_0(H)\\
@A\delta _1AA           @.                   @VV\delta _0V\\
K_1(H)@<<K_1(\psi )<K_1(G)@<<K_1(\varphi )<K_1(F)
\end{CD}$$
is exact.
\end{ax}

\begin{ax}[Commutativity of the index maps]\label{27.9'e}
 If the diagram in \frm
$$\begin{CD}
0@>>>F@>\varphi >>G@>\psi >>H@>>>0\\
@.      @V\phi _1 VV  @V\phi _2VV @VV\phi _3V@.\\
0@>>>F'@>>\varphi' >G'@>>\psi' >H'@>>>0
\end{CD}$$
is commutative and has exact rows then the diagram
$$\begin{CD}
K_i(H)@>\delta _i>>K_{i+1}(F)\\
@VK_i(\phi _3)VV    @VVK_{i+1}(\phi _1)V\\
K_1(H')@>>\delta '_1>K_{i+1}(F')
\end{CD}$$
is commutative, where $\delta _i$ and $\delta '_i$ denote the index maps associated to the upper and the lower row of the above diagram, respectively.
\end{ax}

{\it Remark.} The above axioms are fulfilled if $\kk{i}{F}=0$ for all $E$-C*-algebras $F$ .

{\center{\section{Some elementary results}}}

\begin{p}\label{21.10'}
If
$$\od{F}{\varphi }{G}{\psi }{\lambda }{H}$$
is a split exact sequence in \frm then its index maps are $0$.
\end{p}

By the split exact axiom (\axi{27.9'a}),
$$\og{\kk{i}{F}}{\kk{i}{\varphi }}{\kk{i}{G}}{\kk{i}{\psi }}{\kk{i}{\lambda }}{\kk{i}{H}}{20}{20}{20}$$
is a split exact sequence in the category of additive groups and the assertion follows from the six-term axiom (\axi{27.9'd}).\qed

\begin{de}\label{22.10'}
Let $(F_j)_{j\in J}$ be a finite family of $E$-C*-algebras, $F:=\pro{j\in J}F_j$ and for every $j\in J$ let 
$\mac{\varphi _j}{F_j}{F}$ be the canonical inclusion and $\mac{\psi _j}{F}{F_j}$ the canonical projection. We define
$$\mae{\Phi _{(F_j)_{j\in J},i}}{\pro{j\in J }\kk{i}{F_j}}{\kk{i}{F}}{(a_j)_{j\in J}}{\si{j\in J}\kk{i}{\varphi _j}a_j}\,,$$
$$\mae{\Psi _{(F_j)_{j\in J},i}}{\kk{i}{F}}{\pro{j\in J}F_j}{a}{(\kk{i}{\psi _j}a)_{j\in J}}\;.$$
\end{de}

\begin{p}\label{28.9'a}
If $(F_j)_{j\in J}$ is a finite family of $E$-C*-algebras then the map
$$\mac{\Phi _{(F_j)_{j\in J},i}}{\pro{j\in J }\kk{i}{F_j}}{\kk{i}{\pro{j\in J}F_j}}$$
is a group isomorphism and
$$\mac{\Psi _{(F_j)_{j\in J},i}}{\kk{i}{\pro{j\in J}F_j}}{\pro{j\in J}\kk{i}{F_j}}$$
is its inverse.
\end{p}

If $J=\emptyset $ then the assertion follows from the null-axiom (\axi{27.9'}). The assertion is trivial for $Card\,J=1$. We prove the general case by induction with respect to $Card\,J$. Let $j_0\in J$ and assume the assertion holds for $J':=J\setminus \z{j_0}$. We denote by 
$$\mac{\varphi }{F_{j_0}}{\pro{j\in J}F_j}\,,\qquad\qquad \mac{\lambda }{\pro{j\in J'}F_j}{\pro{j\in J}F_j}$$
the canonical inclusion maps and by
$$\mac{\psi }{\pro{j\in J}F_j}{\pro{j\in J'}F_j}$$
the canonical projection. Then
$$\od{F_{j_0}}{\varphi }{\pro{j\in J}F_j}{\psi }{\lambda }{\pro{j\in J'}F_j}$$
is a split exact sequence in \frm. By the split exact axiom (\axi{27.9'a}) the map
$$\mae{\Psi_i }{\kk{i}{F_{j_0}}\times \kk{i}{\pro{j\in J'}F_j}}{\kk{i}{\pro{j\in J}F_j}}{(a,b)}{\kk{i}{\varphi }a+\kk{i}{\lambda }b}$$
is a group isomorphism. Since
$$\Psi_i \circ \left(id_{\kk{i}{F_{j_0}}}\times \Phi _{(F_j)_{j\in J'},i}\right)=\Phi _{(F_j)_{j\in J},i}$$
it follows from the induction hypothesis that $\Phi _{(F_j)_{j\in J},i}$ is a group isomorphism.

The last assertion follows from $\psi _j\circ \varphi _j=id_{F_j}$ for every $j\in J$ and
$$\si{j\in J}\varphi _j\circ \psi _j=id_{\pro{j\in J}F_j}\;.\qedd$$

\begin{p}\label{21.10'a}
Let $(\oaa{F_j}{\phi _j}{F'_j})_{j\in J}$ be a finite family of morphisms in \frm,
$$F:=\pro{j\in J}F_j\,,\qquad\qquad F':=\pro{j\in J}F'_j\,,$$
and for every $j\in J$ let
$$\mac{\varphi _j}{F_j}{F}\,,\qquad\qquad \mac{\varphi '_j}{F'_j}{F'}$$
be the inclusion maps. Then the diagram
$$\begin{CD}
\pro{j\in J}\kk{i}{F_j}@>\si{j\in J}\kk{i}{\varphi _j}>>K_{i}(F)\\
@V\pro{j\in J}\kk{i}{\phi _j}VV    @VVK_{i}\left(\pro{j\in J}\phi _j\right)V\\
\pro{j\in J}\kk{i}{F'_j}@>>\si{j\in J}\kk{i}{\varphi '_j}>K_{i}(F')
\end{CD}$$
is commutative.
\end{p}

For every $j\in J$ the diagram
$$\begin{CD}
F_j@>\varphi _j>>F\\
@V\phi _jVV    @VV\pro{j\in J}\phi _jV\\
F'_j@>>\varphi '_j>F'
\end{CD}$$
is commutative so the diagram
$$\begin{CD}
K_i(F_j)@>\kk{i}{\varphi _j}>>K_{i}(F)\\
@VK_i(\phi _j)VV    @VVK_{i}\left(\pro{j\in J}\phi _j\right)V\\
K_1(F'_j)@>>\kk{i}{\varphi '_j}>K_{i}(F')
\end{CD}$$
is also commutative. For $(a_j)_{j\in J}\in \pro{j\in J}\kk{i}{F_j}$, by the above,
$$\kk{i}{\pro{j\in J}\phi _j}\circ \left(\si{j\in J}\kk{i}{\varphi _j}\right)(a_j)_{j\in J}=\kk{i}{\pro{j\in J}\phi _j}\si{j\in J}\kk{i}{\varphi _j}a_j=$$
$$=\si{j\in J}\kk{i}{\pro{k\in J}\phi _k}\kk{i}{\varphi _j}a_j=\si{j\in J}\kk{i}{\varphi '_j}\kk{i}{\phi _j}a_j=$$
$$=\left(\si{j\in J}\kk{i}{\varphi '_j}\right)(\kk{i}{\phi _j}a_j)_{j\in J}=\left(\si{j\in J}\kk{i}{\varphi '_j}\right)\kk{i}{\pro{j\in J}\phi _j}(a_j)_{j\in J}\,,$$
which proves the  assertion.\qed

\begin{p}\label{28.9'b}
\rule{0pt}{0pt}
\begin{enumerate}
\item If $\oaa{F}{\varphi }{G},\,\oaa{G}{\psi }{F}$ is a homotopy in \frm\, then
$$K_i(\varphi )\circ K_i(\psi )=id_{K_i(G)}\,,\qquad\qquad K_i(\psi )\circ K_i(\varphi )=id_{K_i(F)}\;.$$
\item If $F$ and $G$ are homotopic $E$-C*-algebras then $K_i(F)$ and $K_i(G)$ are isomorphic.
\item If the $E$-C*-algebra $F$ is null-homotopic then it is K-null.
\end{enumerate}
\end{p}

a) follows from the homotopy axiom (\axi{27.9'b}).

b) follows from a).

c) follows from b) and from the null-axiom (\axi{27.9'}).\qed

\begin{p}\label{28.9'd}
Let
$$\oc{F}{\varphi }{G}{\psi }{H}$$
be an exact sequence in \frm.
\begin{enumerate}
\item If $F$ (resp. $H$) is K-null then 
$$\oag{K_i(G)}{K_i(\psi )}{K_i(H)}{10}\qquad (\emph{resp.}\; \oag{K_i(F)}{K_i(\varphi) }{K_i(G)}{10})$$
 is a group isomorphism.
\item If $G$ is K-null then
$$\oa{K_i(H)}{\delta _i}{K_{i+1}(F)}$$
is a group isomorphism.
\item If $\varphi $ is K-null then the sequences
$$\of{K_i(G)}{K_i(\psi )}{K_i(H)}{\delta _i}{K_{i+1}(F)}{20}{0}$$
is exact.
\item If $\psi $ is K-null then the sequences
$$\of{K_i(H)}{\delta _i}{K_{i+1}(F)}{K_{i+1}(\varphi )}{K_{i+1}(G)}{0}{20}$$
is exact.
\item The index maps of a split exact sequence are equal to $0$.
\end{enumerate}
\end{p}

a), b), c), and d) follow from the six-term axiom (\axi{27.9'd}).

e) follows from the six-term axiom (\axi{27.9'd}) and from the split exact axiom (\axi{27.9'a}).\qed

\begin{p}\label{28.9'e}
An {\bf\frm -triple } is a triple $(F_1,F_2,F_3)$ such that $F_1$ is an \eo, $F_2$ is an $E$-ideal of $F_1$, and $F_3$ is an $E$-ideal of $F_1$ and of $F_2$. We denote for all $j,k\in \bnn{3}$, $j<k$, by $\mac{\varphi _{j,k}}{F_k}{F_j}$ the inclusion map, by $\mac{\psi _{j,k}}{F_j}{F_j/F_k}$ the quotient map, and by $\mac{\delta _{j,k,i}}{K_i(F_j/F_k)}{F_k}$ the index maps associated to the exact sequence in \frm
$$\oc{F_k}{\varphi _{j,k}}{F_j}{\psi _{j,k}}{F_j/F_k}\;.$$
\begin{enumerate}
\item There is a unique morphism $\oag{F_2/F_3}{\varphi _{1,2}/F_3}{F_1/F_3}{15}$ in \frm such that 
$$\psi _{1,3}\circ \varphi _{1,2}=(\varphi _{1,2}/F_3)\circ \psi _{2,3}\;.$$
\item The diagram
$$\begin{CD}
K_i(F_3)@>K_i(\varphi _{1,3})>>K_i(F_1)@>K_i(\psi _{1,3})>>K_i(F_1/F_3)@>\delta _{1,3,i}>>K_{i+1}(F_3)\\
@A=AA@AK_i(\varphi _{1,2})AA         @AAK_i(\varphi _{1,2}/F_3)A  @AA=A\\
K_i(F_3)@>>K_i(\varphi _{2,3})>    K_i(F_2)@>>K_i(\psi _{2,3})>K_i(F_2/F_3)@>>\delta _{2,3,i}>K_{i+1}(F_3)\\
\end{CD}$$
is commutative.
\end{enumerate}
\end{p}

a) is easy to see.

b) follows from a), $\varphi _{1,2}\circ \varphi _{2,3}=\varphi _{1,3}$, and from the axiom of commutativity of the index maps (\axi{27.9'e}).\qed

\begin{theo}[The triple theorem]\label{28.9'f}
Let $(F_1,F_2,F_3)$ be an \frm -triple.
\begin{enumerate}
\item Assume $F_2$ K-null.
\begin{enumerate}
\item $\mac{\delta _{2,3,i}}{K_i(F_2/F_3)}{K_{i+1}(F_3)}$ is a group isomorphism.
\item $\delta _{2,3,i}=\delta _{1,3,i}\circ K_i(\varphi_{1,2}/F_3 )$.
\item $\varphi_{1,3} $ is K-null.
\item If we put $\Phi _i:=K_i(\varphi_{1,2}/F_3 )\circ (\delta _{2,3,i})^{-1}$ then
$$\og{K_i(F_1)}{K_i(\psi _{1,3})}{K_i(F_1/F_3)}{\delta _{1,3,i}}{\Phi _i}
{K_{i+1}(F_3)}{20}{10}{5}$$
is a split exact sequence and the map
$$\mad{K_i(F_1)\times K_{i+1}(F_3)}{K_i(F_1/F_3)}
{(a,b)}{K_i(\psi _{1,3})a+\Phi _ib}$$
is a group isomorphism.
\end{enumerate}
\item Assume $F_1/F_3$ K-null.
\begin{enumerate}
\item $\delta _{2,3,i}=0$ and the sequence
$$\of{\kk{i}{F_3}}{\kk{i}{\varphi _{2,3}}}{\kk{i}{F_2}}{\kk{i}{\psi _{2,3}}}{\kk{i}{F_2/F_3}}{20}{20}$$
is exact.
\item $\mac{K_i(\varphi _{1,3})}{K_i(F_3)}{K_i(F_1)}$ is a group isomorphism.
\item If we put $\Phi _i:=K_i(\varphi _{1,3})^{-1}\circ K_i(\varphi _{1,2})$ then the map
$$\mae{\Psi }{K_{i}(F_2)}{K_i(F_3)\times K_{i}(F_2/F_3)}{b}{(\Phi _ib,K_i(\psi _{2,3})b)}$$
is a group isomorphism.
\item If $\psi _{1,2}$ is K-null and if we put $\Phi '_i:=K_i(\varphi _{2,3})\circ K_i(\varphi _{1,3})^{-1}$ then
$$\og{K_{i+1}(F_1/F_2)}{\delta _{1,2,(i+1)}}{K_i(F_2)}{K_i(\varphi _{1,2})}{\Phi '_i}{K_i(F_1)}{20}{20}{20}$$
is a split exact sequence and the map
$$\mad{K_i(F_1)\times K_{i+1}(F_1/F_2)}{K_i(F_2)}{(a,b)}{\Phi '_ia+\delta _{1,2,(i+1)}b}$$
is a group isomorphism. 
\end{enumerate}
\item Assume $F_1$ K-null and denote by $\psi $ the canonical map $F_1/F_3\rightarrow F_1/F_2$.
\begin{enumerate}
\item $\delta _{1,2,i}$ and $\delta _{1,3,i}$ are group isomorphisms.
\item $K_i(\varphi _{2,3})\circ \delta _{1,3,(i+1)}=\delta _{1,2,(i+1)}\circ K_{i+1}(\psi )$.
\item Let $\mac{\varphi}{F_1/F_2}{F_1/F_3}$ be a morphism in \frm such that 
$$K_i(\psi \circ \varphi)=id_{K_i(F_1/F_2)}\;.$$
If we put 
$$\Phi _i:=\delta _{1,3,(i+1)}\circ K_{i+1}(\varphi)\circ (\delta _{1,2,(i+1)})^{-1}$$
then $K_i(\varphi _{2,3})\circ \Phi _i=id_{K_i(F_2)}$. If in addition $\psi _{2,3}$ is K-null then
$$\og{K_{i+1}(F_2/F_3)}{\delta _{2,3,(i+1)}}{K_i(F_3)}{K_i(\varphi _{2,3})}{\Phi _i}{K_i(F_2)}{20}{20}{20}$$
is a split exact sequence and the map
$$\mad{K_{i+1}(F_2/F_3)\times K_i(F_2)}{K_i(F_3)}{(a,b)}{\delta _{2,3,(i+1)}a+\Phi _ib}$$
is a group isomorphism.
\end{enumerate}
\end{enumerate}
\end{theo}

$a_1)$ follows from \pr{28.9'd} b).

$a_2)$ follows from \pr{28.9'e} b).

$a_3)$ $\varphi _{1,3}$ factorizes through null and so it is K-null.

$a_4)$  By $a_2)$,
$$\delta _{1,3,i}\circ \Phi _i=\delta _{1,3,i}\circ K_i(\varphi_{1,2}/F_3 )\circ (\delta _{2,3,i})^{-1}=\delta _{2,3,i}\circ (\delta _{2,3,i})^{-1}=id_{K_i(F_3)}$$
and this implies the assertion.

$b_1)$ By \pr{28.9'e} b), $\delta _{2,3,i}$ factorizes through null and so it is K-null. By the six-term axiom (\axi{27.9'd}) the sequence 
$$\of{\kk{i}{F_3}}{\kk{i}{\varphi _{2,3}}}{\kk{i}{F_2}}{\kk{i}{\psi _{2,3}}}{\kk{i}{F_2/F_3}}{20}{20}$$
is exact.

$b_2)$ follows from \pr{28.9'd} a).

$b_3)$
\begin{center}
Step 1 $\Phi _i\circ K_i(\varphi _{2,3})=id_{K_i(F_3)}$
\end{center}

Since $\varphi _{1,3}=\varphi _{1,2}\circ \varphi _{2,3}$,
$$\Phi _i\circ K_i(\varphi _{2,3})=K_i(\varphi _{1,3})^{-1}\circ K_i(\varphi _{1,2})\circ K_i(\varphi _{2,3})=$$
$$=K_i(\varphi _{1,3})^{-1}\circ K_i(\varphi _{1,3})=id_{K_i(F_3)}\;.$$

\begin{center}
Step 2 $\Psi $ is injective
\end{center}

Let $b\in K_i(F_2)$ with $\Psi b=0$. Then $K_i(\psi _{2,3})b=0$ so by $b_1)$,
$$b\in Ker\,K_i(\psi _{2,3})=Im\,K_i(\varphi _{2,3})$$
and there is an $a\in K_i(F_3)$ with $b=K_i(\varphi _{2,3})a$. By Step 1, 
$$a=\Phi _iK_i(\varphi _{2,3})a=\Phi _ib=0\,,$$
so $b=0$ and $\Psi $ is injective.

\begin{center}
Step 3 $\Psi $ is surjective
\end{center}

Let $(a,c)\in K_i(F_3)\times K_i(F_2/F_3)$. Put $b':=K_i(\varphi _{2,3})a$. By $b_1)$,
$$K_i(\psi _{2,3})b'=K_i(\psi _{2,3})K_i(\varphi _{2,3})a=0$$
and by Step 1, $\Phi _ib'=\Phi _iK_i(\varphi _{2,3})a=a$. By $b_1)$, there is a $b''\in K_i(F_2)$ with $c=K_i(\psi _{2,3})b''$. By Step 1, 
$$\Phi _i(b''-K_i(\varphi _{2,3})\Phi _ib'')=\Phi _ib''-\Phi _iK_i(\varphi _{2,3})\Phi _ib''=\Phi _ib''-\Phi _ib''=0\;.$$
Thus by $b_1)$,
$$\Psi (b'+b''-K_i(\varphi _{2,3})\Phi _ib'')=$$
$$=(\Phi _ib',K_i(\psi _{2,3})b''-K_i(\psi _{2,3})K_i(\varphi _{2,3})\Phi _ib'')=(a,c)$$
and $\Psi $ is surjective.

$b_4)$ Since $\varphi _{1,3}=\varphi _{1,2}\circ \varphi _{2,3}$,
$$K_i(\varphi _{1,2})\circ \Phi '_i=K_i(\varphi _{1,2})\circ K_i(\varphi _{2,3})\circ K_i(\varphi _{1,3})^{-1}=$$
$$=K_i(\varphi _{1,3})\circ K_i(\varphi _{1,3})^{-1}=id_{K_i(F_1)}$$
and the assertion follows.

$c_1)$ follows from \pr{28.9'd} b)).

$c_2)$ follows from the commutativity of the index maps (\axi{27.9'e}).

$c_3)$ By $c_2)$,
$$K_i(\varphi _{2,3})\circ \Phi _i=K_i(\varphi _{2,3})\circ \delta _{1,3,(i+1)}\circ K_{i+1}(\varphi)\circ (\delta _{1,2,(i+1)})^{-1}=$$
$$=\delta _{1,2,(i+1)}\circ K_{i+1}(\psi )\circ K_{i+1}(\varphi)\circ (\delta _{1,2,(i+1)})^{-1}=$$
$$=\delta _{1,2,(i+1)}\circ K_{i+1}(\psi \circ \varphi)\circ (\delta _{1,2,(i+1)})^{-1}=\delta _{1,2,(i+1)}\circ (\delta _{1,2,(i+1)})^{-1}=id_{K_i(F_2)}\;.$$

The last assertion follows from the first one.\qed

{\it Remark.} a) still holds with the weaker assumption that $F_2$ is only an $E$-C*-subalgebra of $F_1$.

{\center{\section{Tensor products}}}

\begin{center}
\fbox{\parbox{8.8cm}{Throughout this section $F$ denotes an $E$-C*-algebra}}
\end{center}

\begin{de}\label{26.3'}
Let $G$ be a C*-algebra. We denote by $F\otimes G$ the spatial tensor product of $F$ and $G$ endowed with the structure of an \eo by using the exterior multiplication
$$\mad{E\times (F\otimes G)}{F\otimes G}{(\alpha ,x\otimes y)}{(\alpha x)\otimes y}$$
\emph{([W] Proposition T.5.14 and T.5.17 Remark)}. If $\oaa{F}{\varphi }{F'}$ is a morphism in \frm and $\oaa{G}{\psi }{G'}$ a morphism in $\fr{M}_{\bc}$ then  $\oa{F\otimes G}{\varphi\otimes \psi  }{F'\otimes G'}$ denotes the morphism in \frm defined by
$$\mae{\varphi\otimes \psi }{F\otimes G}{F'\otimes G'}{x\otimes y}{\varphi x\otimes \psi y}\;.$$
If $(G_j)_{j\in J}$ is a family of C*-algebras then we put
$$\bigotimes_{j\in \emptyset }G_j:=\bc\;. $$
\end{de}

We have $F\otimes \bc\approx F$ and $id_F\otimes id_G=id_{F\otimes G}$. If $\ob{F}{\varphi }{F'}{\varphi '}{F''}$ are morphisms in \frm and $\ob{G}{\psi }{G'}{\psi '}{G''}$ are morphisms in $\fr{M}_{\bc}$ then
$$(\varphi \otimes \psi )\circ (\varphi '\otimes \psi ')=(\varphi \circ \varphi ')\otimes (\psi \circ \psi ')\;.$$
If $G$ and $H$ are C*-algebras then 
$$F\otimes (G\times H)\approx (F\otimes G)\times (F\otimes H)\,,\qquad F\otimes (G\otimes H)\approx (F\otimes G)\otimes H\;.$$
If $G$ is a C*-algebra and $F_1,F_2$ are $E$-C*-algebras then
$$(F_1\times F_2)\otimes G\approx (F_1\otimes G)\times (F_2\otimes G)\;.$$ 

\begin{p}\label{29.3'}
Let $G,H$ be C*-algebras.
\begin{enumerate}
\item If  $\mac{\varphi_0 \,,\varphi _1}{G}{H}$ are homotopic C*-homomorphisms then $id_F\otimes \varphi _0$ and $id_F\otimes \varphi _1$ are also homotopic.
\item If $\oaa{G}{\varphi }{H}$, $\oaa{H}{\psi }{G}$ is a homotopy in \frcc then
$$\oa{F\otimes G}{id_F\otimes \varphi }{F\otimes H}\,,\qquad\oa{F\otimes H}{id_F\otimes \psi }{F\otimes G}$$
is a homotopy in \frm.
\item If $G$ is homotopic to $0$ then $F\otimes G$ is also homotopic to $0$ and so K-null.
\end{enumerate}
\end{p}

a) Let $[0,1]\longrightarrow \varphi _s$ be a pointwise continuous map of C*-homomorphisms $G\rightarrow H$. Let $z\in F\odot G$. There are finite  families $(x_j)_{j\in J}$ in $F$ and $(y_j)_{j\in J}$ in $G$ such that
$$z=\si{j\in J}x_j\otimes y_j\;.$$
For $s\in [0,1]$,
$$(id_F\otimes \varphi _s)z=\si{j\in J}x_j\otimes \varphi _sy_j$$
so the map
$$\mad{[0,1]}{F\otimes H}{s}{(id_F\otimes \varphi _s)}z$$
is continuous.

Let now $z\in F\otimes G$, $s_0\in [0,1]$, and $\varepsilon >0$. There is a $z'\in F\odot G$ such that $\n{z-z'}<\frac{\varepsilon }{3}$. By the above, there is a $\delta >0$ such that
$$\n{(id_F\otimes \varphi _s)z'-(id_F\otimes \varphi _{s_0})z'}<\frac{\varepsilon }{3}$$
for all $s\in [0,1]$, $|s-s_0|<\delta $. It follows
$$\n{(id_F\otimes \varphi _s)z-(id_F\otimes \varphi _{s_0})z}\leq \n{(id_F\otimes \varphi _s)(z-z')}+$$
$$+\n{(id_F\otimes \varphi _s)z'-(id_F\otimes \varphi _{s_0})z'}+\n{(id_F\otimes \varphi _{s_0})(z-z')}<\varepsilon \,,$$
which proves the assertion.

b) follows from a).

c) follows from b) and \pr{28.9'b} c)).\qed

\begin{p}\label{26.3'a}
Let
$$\od{G_1}{\varphi }{G_2}{\psi }{\lambda }{G_3}$$
be a split exact sequence in $\fr{M}_{\bc}$.
\begin{enumerate}
\item The sequence in \frm
 $$\og{F\otimes G_1}{id_F\otimes \varphi }{F\otimes G_2}{id_F\otimes \psi }{id_F\otimes \lambda }{F\otimes G_3}{20}{20}{20}$$
is split exact.
\item The sequence
$$\og{K_i(F\otimes G_1)}{K_i(id_F\otimes \varphi) }{K_i(F\otimes G_2)}{K_i(id_F\otimes \psi )}{K_i(id_F\otimes \lambda )}{K_i(F\otimes G_3)}{30}{30}{30}$$
is split exact and the map
$$K_i(F\otimes G_1)\times K_i(F\otimes G_3)\longrightarrow  K_i(F\otimes G_2)\,,$$
$$(a,b)\longmapsto K_i(id_F\otimes \varphi)a+K_i(id_F\otimes \lambda )b$$
is a group isomorphism.
\end{enumerate}
\end{p}

a) By [W] Corollary T.5.19, $id_F\otimes \varphi $ is injective. We have
$$(id_F\otimes \psi )\circ (id_F\otimes \lambda )=id_F\otimes (\psi \circ \lambda )=id_F\otimes id_{G_3}=id_{F\otimes G_3}\,,$$
$$(id_F\otimes \psi )\circ (id_F\circ \varphi )=id_F\otimes (\psi \circ \varphi )=0\,,$$
so 
$$Im\,(id_F\otimes \varphi )\subset Ker\,(id_F\otimes \psi )\;.$$

Let $z\in (F\odot G_2)\cap Ker\,(id_F\otimes \psi )$. There is a linearly independent finite family $(x_j)_{j\in J}$ in $F$ and a family $(y_j)_{j\in J}$ in $G_2$ such that
$$z=\si{j\in J}x_j\otimes y_j\;.$$
From $$0=(id_F\otimes \psi )z=\si{j\in J}x_j\otimes \psi y_j$$
we get $\psi y_j=0$ for all $j\in J$. Thus for every $j\in J$ there is a $y'_j\in G_1$ with $\varphi y'_j=y_j$. It follows
$$z=\si{j\in J}x_j\otimes \varphi y'_j=(id_F\otimes \varphi )\si{j\in J}x_j\otimes y'_j\in Im\,(id_F\otimes \varphi )\;.$$

Let $z\in Ker\,(id_F\otimes \psi )$. Then
$$(id_F\otimes (\lambda \circ \psi ))z=(id_F\otimes \lambda )(id_F\otimes \psi )z=0\;.$$
Let $(z_n)_{\bbn}$ be a sequence in $F\odot G_2$ converging to $z$. For \bbn, by the above,
$$(id_F\otimes \psi )(z_n-(id_F\otimes (\lambda \circ \psi ))z_n)=$$
$$=(id_F\otimes \psi )z_n-(id_F\otimes \psi )(id_F\otimes \lambda )(id_F\otimes \psi )z_n=$$
$$=(id_F\otimes \psi )z_n-(id_F\otimes \psi )z_n=0\,,$$
$$z_n-(id_F\otimes (\lambda \circ \psi))z_n\in Im\,(id_F\otimes \psi )\;.$$
Since $Im\,(id_F\otimes \varphi )$ is closed,
$$z=z-(id_F\otimes (\lambda \circ \psi))z=\lim_{n\rightarrow \infty }(z_n-(id_F\otimes (\lambda \circ \psi))z_n)\in Im\,(id_F\otimes \varphi )\,,$$
which proves the Proposition.

b) follows from a) and the split exact axiom (\axi{27.9'a}).\qed

\begin{de}\label{3.10'}
\item We denote for every C*-algebra $G$ by $\tilde{G} $ its unitization \emph{(see e.g. [R] Exercise 1.3)} and by
$$\od{G}{\iota_G}{\tilde{G} }{\pi_ G}{\lambda_ G}{\bc}$$
its associated split exact sequence. If $G$ and $H$ are C*-algebras and $\mac{\varphi }{G}{H}$ is a C*-homomorphism then $\mac{\tilde{\varphi } }{\tilde{G} }{\tilde{H} }$ denotes the unitization of $\varphi $.
\end{de}

\begin{co}\label{26.3'b}
Let $G$ be a C*-algebra. 
\begin{enumerate}
\item The sequence in \frm
$$\og{F\otimes G}{id_F\otimes \iota_G}{F\otimes \tilde{G} }{id_F\otimes \pi_G}{id_F\otimes \lambda_G}{F}{20}{20}{20}$$
is split exact.
\item The sequence
$$\og{K_i(F\otimes G)}{K_i\left(id_F\otimes \iota_G\right)}{K_i\left(F\otimes \tilde{G}\right) }{K_i\left(id_F\otimes \pi_G\right)}{K_i\left(id_F\otimes \lambda_G\right)}{K_i(F)}{30}{35}{35}$$
is split exact and the map
$$K_i(F)\times K_i(F\otimes G)\longrightarrow K_i\left(F\otimes \tilde{G}\right)\,,$$
$$(a,b)\longmapsto K_i\left(id_F\otimes \lambda_G\right)a+K_i\left(id_F\otimes \iota_G\right)b$$
is a group isomorphism.
\item Let $\oaa{F}{\varphi }{F'}$ be a morphism in \frm and $\oaa{G}{\psi }{G'}$ a morphism in $\fr{M}_{\bc}$. If we identify the isomorphic groups of b) then
$$\mac{\kk{i}{\varphi \otimes \tilde{\psi} }}{\kk{i}{F\otimes \tilde{G} }}{\kk{i}{F'\otimes \tilde{G'} }},$$
$$(a,b)\longmapsto (\kk{i}{\varphi }a,\kk{i}{\varphi \otimes \psi }b)$$
is a group isomorphism.
\item Let $\mac{\varphi }{G}{G'}$ be a morphism in \frcc. If we denote by $\Psi _i$ and $\Psi '_i$ the group isomorphisms of b) associated to $G$ and $G'$, respectively, then
$$\kk{i}{id_F\otimes \tilde{\varphi } }\circ \Psi _i=\Psi '_i\circ \left(id_{\kk{i}{F}}\times \kk{i}{id_F\otimes \varphi }\right)\;.$$
\end{enumerate}
\end{co}

a) and b) follow from \pr{26.3'a} a),b).

c) follows from b) and the commutativity of the following diagram:
$$\begin{CD}
F\otimes G@>id_F\otimes \iota_G>>F\otimes \tilde{G} @<id_F\otimes \lambda_G<<F\otimes \bc\\
@V\varphi \otimes \psi VV@V\varphi \otimes \tilde{\psi }VV         @VV\varphi \otimes id_{\bc}V\\
F'\otimes G'@>>id_{F'}\otimes \iota_{G'}>    F'\otimes \tilde{G'} @<<id_{F'}\otimes \lambda_{G'}<F'\otimes \bc\\
\end{CD}\hspace{1cm}.$$

d) For $(a,b)\in \kk{i}{F}\times \kk{i}{F\otimes G}$, since $\tilde{\varphi }\circ \lambda _G=\lambda _{G'} $ and $\iota _{G'}\circ \varphi =\tilde{\varphi }\circ \iota _G $,
$$\kk{i}{id_F\otimes \tilde{\varphi } }\Psi _i(a,b)=\kk{i}{id_F\otimes \tilde{\varphi } }(\kk{i}{id_F\otimes \lambda _G}a+\kk{i}{id_F\otimes \iota _G}b)=$$
$$=\kk{i}{id_F\otimes \tilde{\varphi } }\kk{i}{id_F\otimes \lambda _G}a+\kk{i}{id_F\otimes \tilde{\varphi } }\kk{i}{id_F\otimes \iota _G}b=$$
$$=\kk{i}{id_F\otimes (\tilde{\varphi }\circ \lambda _G )}a+\kk{i}{id_F\otimes (\tilde{\varphi }\circ \iota _G )}b=$$
$$=\kk{i}{id_F\otimes \lambda _{G'}}a+\kk{i}{id_F\otimes (\iota _{G'}\circ \varphi )}b=$$
$$=\kk{i}{id_F\otimes \lambda _{G'}}a+\kk{i}{id_F\otimes \iota _{G'}}\kk{i}{id_F\otimes \varphi }b=$$
$$=\Psi '_i(a,\kk{i}{id_F\otimes \varphi }b)=\Psi '_i(id_{\kk{i}{F}}\times \kk{i}{id_F\otimes \varphi })(a,b)\,,$$
so
$$\kk{i}{id_F\otimes \tilde{\varphi } }\circ \Psi _i=\Psi '_i\circ \left(id_{\kk{i}{F}}\times \kk{i}{id_F\otimes \varphi }\right)\;.\qedd$$

\begin{p}\label{27.3'}
If $(G_j)_{j\in J}$ is a finite family of C*-algebras then
$$\kk{i}{F\otimes \left(\bigotimes_{j\in J}\tilde{G_j} \right)}\approx \pro{I\subset J}\kk{i}{F\otimes \left(\bigotimes _{j\in I}G_j\right)}\;.$$
\end{p}

We prove the assertion by induction with respect to $Card\,J$. The assertion is trivial for $Card\,J=0$ (\dd{26.3'} and Null-axiom (\axi{27.3'})). Let $j_0\in J$, $J':=J\setminus \z{j_0}$, and assume the assertion holds for $J'$. By \cor{26.3'b} b),
$$\kk{i}{F\otimes \left(\bigotimes_{j\in J}\tilde{G_j}\right)}\approx \kk{i}{\left(F\otimes \left(\bigotimes_{j\in J'}\tilde{G_j}  \right)\right)\otimes \tilde{G_{j_0}} }\approx $$
$$\approx \kk{i}{F\otimes \left(\bigotimes_{j\in J'}\tilde{G_j}\right)}\times \kk{i}{\left(F\otimes \left(\bigotimes_{j\in J'}\tilde{G_j}  \right)\right)\otimes G_{j_0}}\approx $$
$$\approx \kk{i}{F\otimes \left(\bigotimes_{j\in J'}\tilde{G_j}\right)}\times \kk{i}{(F\otimes G_{j_0})\otimes \left(\bigotimes_{j\in J'}\tilde{G_j}  \right)}\approx $$
$$\approx \pro{I\subset J'}\kk{i}{F\otimes \left(\bigotimes_{j\in I}G_j \right)}\times \pro{I\subset J'}\kk{i}{F\otimes \left(\bigotimes_{j\in I\cup \z{j_0}}G_j \right)}\approx $$
$$\approx \pro{I\subset J}\kk{i}{F\otimes \left(\bigotimes _{j\in I}G_j\right)}\;.\qedd$$

\begin{co}\label{27.3'a}
If $G$ is a C*-algebra then for all $n\in \bn^*$
$$\kk{i}{F\otimes \left(\bigotimes_{j\in \bnn{n}}\tilde{G}\right)}\approx \proo{k=0}{n}\kk{i}{F\otimes \left(\bigotimes_{j\in \bnn{k}}G \right)}^{n\choose k}\;.\qedd$$
\end{co}

\begin{p}\label{23.4'}
Let $G$ be a C*-algebra and
$$\od{F_1}{\varphi }{F_2}{\psi }{\lambda }{F_3}$$
a split exact sequence in \frm.
\begin{enumerate}
\item The sequence in \frm
$$\og{F_1\otimes G}{\varphi \otimes id_G}{F_2\otimes G}{\psi \otimes id_G}{\lambda \otimes id_G}{F_3\otimes G}{20}{20}{20}$$
is split exact.
\item The sequence
$$\og{K_i(F_1\otimes G)}{K_i(\varphi \otimes id_G)}{K_i(F_2\otimes G)}{K_i(\psi \otimes id_G)}{K_i(\lambda \otimes id_G )}{K_i(F_3\otimes G)}{30}{30}{30}$$
is split exact and the map
$$K_i(F_1\otimes G)\times K_i(F_3\otimes G)\longrightarrow  K_i(F_2\otimes G)\,,$$
$$(a,b)\longmapsto K_i(\varphi \otimes id_G)a+K_i(\lambda \otimes id_G)b$$
is a group isomorphism.
\end{enumerate}
\end{p}

The proof is similar to the proof of \pr{26.3'a}.\qed

\begin{p}\label{29.3'a}
Let
$$\oc{G_1}{\varphi }{G_2}{\psi }{G_3}$$
be an exact sequence in $\fr{M}_{\bc}$. If $F$ or $G_3$ is nuclear then the sequence in \frm
$$\of{F\otimes G_1}{id_F\otimes \varphi }{F\otimes G_2}{id_F\otimes \psi }{F\otimes G_3}{20}{20}$$
is exact and so
$$\frac{F\otimes G_2}{F\otimes G_1}\approx F\otimes \frac{G_2}{G_1}\;.$$
\end{p}

[W] Theorem T.6.26.\qed

\begin{p}\label{23.4'b}
Let G be a C*-algebra and
$$\oc{F_1}{\phi_1 }{F_2}{\phi_2 }{F_3}$$
an exact sequence in $\fr{M}_E$.
If $F_3$ or $G$ is nuclear then
$$\of{F_1\otimes G}{\phi_1 \otimes id_G }{F_2\otimes G}{\phi_2 \otimes id_G}{F_3\otimes G}{20}{20}$$
is exact. 
\end{p}

[W] Theorem T.6.26.\qed

\begin{de}\label{26.8'}
Let
$$\oc{F_1}{\phi _1}{F_2}{\phi _2}{F_3}$$
be an exact sequence in \frm and $G$ a C*-algebra. If $\delta _i$ denotes the index maps associated to the above exact sequence in \frm and if the sequence in \frm
$$\of{F_1\otimes G}{\phi _1\otimes id_G}{F_1\otimes G}{\phi _2\otimes id_G}{F_3\otimes G}{20}{20}$$
is exact \emph{(e.g. $F_3$ or $G$ is nuclear  ([W] T.6.26))}
 then we denote by $\delta _{G,i}$ the index maps associated to this last exact sequence in \frm. 
\end{de}

In this case the six-term sequence
$$\begin{CD}
K_0(F_1\otimes G)@>K_0(\phi_1\otimes id_G )>>K_0(F_2\otimes G)@>K_0(\phi _2\otimes id_G)>>K_0(F_3\otimes G)\\
@A\delta _{G,1}AA           @.                   @VV\delta _{G,0}V\\
K_1(F_3\otimes G)@<<K_1(\phi _2\otimes id_G)<K_1(F_2\otimes G)@<<K_1(\phi_1\otimes id_G )<K_1(F_1\otimes G)
\end{CD}$$
is exact (by the six-term axiom (\axi{27.9'd})).

\begin{co}\label{23.4'c}
Let $G$ be a unital C*-algebra,
$$\oc{F_1}{\phi_1 }{F_2}{\phi_2 }{F_3}$$
an exact sequence in $\fr{M}_E$, and $\delta _i$ its index maps. We assume that $F_3$ or $G$ is nuclear and put for every $j\in \z{1,2,3}$
$$\mae{\varphi _j}{F_j}{F_j\otimes G}{x}{x\otimes 1_G}\;.$$
Then $\delta _{G,i}\circ \kk{i}{\varphi _3}=\kk{i+1}{\varphi _1}\circ \delta _i$.
\end{co}

The diagram
$$\begin{CD}
F_1@>\phi_1 >>F_2 @>\phi_2 >>F_3\\
@V\varphi_1 VV@V\varphi_2  VV         @VV\varphi _3 V\\
F_1\otimes G@>>\phi_1 \otimes id_G>    F_2\otimes G  @>>\phi_2 \otimes id_{G}>F_3\otimes G\\
\end{CD}$$
is commutative and the assertion follows from \pr{23.4'b} and the commutativity of the index maps (\axi{27.9'e}).\qed

{\center{\section{The class $\Upsilon $  }}}

\begin {center}
\fbox{\parbox{8.8cm}{Throughout this section $F$ denotes an $E$-C*-algebra}}
\end{center}

\begin{de}\label{5.7'}
Let $\Upsilon $ be the class of those C*-algebras $G$ for which there are $p(G),q(G)\in \bn^*$ and group isomorphisms
$$\mac{\Phi _{i,G,F}}{\kk{i}{F}^{p(G)}\times \kk{i+1}{F}^{q(G)}}{\kk{i}{F\otimes G}}$$
such that for every morphism $\oaa{F}{\phi }{F'}$ in \frm the diagram
$$\begin{CD}
\kk{i}{F}^{p(G)}\times \kk{i+1}{F}^{q(G)}@>\Phi _{i,G,F} >>\kk{i}{F\otimes G}\\
@V\kk{i}{\phi }^{p(G)}\times \kk{i+1}{\phi }^{q(G)}VV@VV\kk{i}{\phi \otimes id_G}V \\        
\kk{i}{F'}^{p(G)}\times \kk{i+1}{F'}^{q(G)}@>>\Phi _{i,G,F'}>\kk{i}{F'\otimes G}\\
\end{CD}$$
is commutative. We denote by $\vec{G} $ the class of group isomorphisms
$$\mac{\Phi _{i,G,F}}{\kk{i}{F}^{p(G)}\times \kk{i+1}{F}^{q(G)}}{\kk{i}{F\otimes G}}$$
having the above property. A C*-algebra $G$ is called {\bf{$\Upsilon $-null}} if $G\in \Upsilon $ and $p(G)=q(G)=0$.
\end{de}

If $G$ is $\Upsilon $-null or if $F$ is K-null and $G\in \Upsilon $  then $F\otimes G$ is K-null. In general we shall use $\Phi _{i,G,F}$ without writing $\z{\Phi _{i,G,F}}\in \vec{G} $.

\begin{p}\label{13.11'}
Let $p.q\in \bn^*$ and let $\Lambda $ be the class of group isomorphisms
$$\mac{\Lambda _{i,F}}{\kk{i}{F}^p\times \kk{i+1}{F}^q}{\kk{i}{F}^p\times \kk{i+1}{F}^q}$$
such that for all morphisms $\oaa{F}{\phi }{F'}$ in \frm the diagram
$$\begin{CD}
\kk{i}{F}^{p}\times \kk{i+1}{F}^{q}@>\Lambda _{i,F} >>\kk{i}{F}^{p}\times \kk{i+1}{F}^{q}\\
@V\kk{i}{\phi }^{p}\times \kk{i+1}{\phi }^{q}VV@VV\kk{i}{\phi }^{p}\times \kk{i+1}{\phi }^{q}V \\        
\kk{i}{F'}^{p}\times \kk{i+1}{F'}^{q}@>>\Lambda _{i,F'}>\kk{i}{\phi }^{p}\times \kk{i+1}{\phi }^{q}\\
\end{CD}$$
is commutative. Let $G\in \Upsilon $ with $p(G)=p$, $q(G)=q$, and let $\{\Phi _{i,G,F}\}\in \vec{G}$.
\begin{enumerate}
\item If ${\Lambda _{i,F}}\in \Lambda $ and if we put
$$\Phi '_{i,G,F}:=\mac{\Phi _{i,G,F}\circ \Lambda _{i,F}}{\kk{i}{\phi }^{p}\times \kk{i+1}{\phi }^{q}}{\kk{i}{F\otimes G}}$$
then $\{\Phi '_{i,G,F}\}\in \vec{G} $.
\item If $\{\Phi '_{i,G,F}\}\in \vec{G} $ and if we put
$$\Lambda _{i,F}:=\mac{\Phi ^{-1}_{i,G,F}\circ \Phi '_{i,G,F}}{\kk{i}{\phi }^{p}\times \kk{i+1}{\phi }^{q}}{\kk{i}{\phi }^{p}\times \kk{i+1}{\phi }^{q}}$$
then $\{\Lambda _{i,F}\}\in \Lambda $.
\item If $\{\Lambda _{i,F}\},\z{\Lambda '_{i,F}}\in \Lambda $ then $\z{\Lambda _{i,F}\circ \Lambda '_{i,F}}\in \Lambda $, $\z{\Lambda ^{-1}_{i,F}}\in \Lambda $.\qed
\end{enumerate}
\end{p}

\begin{de}\label{25.8'}
We denote for every nuclear $G\in \Upsilon $ by $G_\Upsilon $ the class of exact sequences in \frm
$$\oc{F_1}{\phi _1}{F_2}{\phi _2}{F_3}$$
such that if $\delta _i$ denote its index maps then the diagram
$$\begin{CD}
\kk{i}{F_3}^{p(G)}\times \kk{i+1}{F_3}^{q(G)}@>\Phi _{i,G,F_3} >>\kk{i}{F_3\otimes G}\\
@V\delta _i^{p(G)}\times \delta _{i+1}^{q(G)}VV@VV\delta _{G,i}V \\        
\kk{i+1}{F_1}^{p(G)}\times \kk{i}{F_1}^{q(G)}@>>\Phi _{(i+1),G,F_1}>\kk{i+1}{F_1\otimes G}\\
\end{CD}$$
is commutative.
\end{de} 

If $G$ is $\Upsilon $-null then every exact sequence in \frm belongs to $G_\Upsilon $.

\begin{p}\label{5.7'a}
$\rule{0mm}{0mm}$
\begin{enumerate}
\item $0$ is $\Upsilon $-null.
\item $\bc\in \Upsilon $, $p(\bc)=1$, $q(\bc)=0$, $\Phi _{i,\bc,F}=\kk{i}{\phi _{\,\bc,F}}$, where
$$\mae{\phi _{\,\bc,F}}{F}{F\times \bc}{x}{x\otimes 1_{\bc}}\;.$$
Every exact sequence in \frm belongs to $\bc_\Upsilon $.
\item Let $\oa{G}{\varphi }{G'}$, $\oa{G'}{\psi }{G}$ be a homotopy in $\frcc$. 
 If $G\in \Upsilon $ then
$$G'\in \Upsilon\,,\qquad p(G')=p(G)\,,\qquad q(G')=q(G)\,,$$
$$\Phi _{i,G',F}=\kk{i}{id_F\otimes \varphi }\circ \Phi _{i,G,F}\;.$$
If in addition $G$ and $G'$ are nuclear then $G_\Upsilon =G'_\Upsilon $.
\item If $G$ is null-homotopic then $G$ is $\Upsilon $-null.
\end{enumerate}
\end{p}

a) By the null-axiom (\axi{27.9'}), $0$ is $\Upsilon $-null.

b) The first assertion is easy to see. The second one follows from the commutativity of the index maps (\axi{27.9'e}).

c) By \pr{29.3'} b),
$$\oa{F\otimes G}{id_{F}\otimes \varphi }{F\otimes G'}\,,\qquad \oa{F\otimes G'}{id_F\otimes \psi }{F\otimes G}$$
is a homotopy in \frm. By \pr{28.9'e} a),
$$\mac{\kk{i}{id_F\otimes \varphi }}{\kk{i}{F\otimes G}}{\kk{i}{F\otimes G'}}\,,$$
$$\mac{\kk{i}{id_F\otimes \psi  }}{\kk{i}{F\otimes G'}}{\kk{i}{F\otimes G}}$$
are group isomorphisms and $\kk{i}{id_F\otimes \psi }=\kk{i}{id_F\otimes \varphi }^{-1}$. Thus
$$\mac{\kk{i}{id_F\otimes \varphi }\circ \Phi _{i,G,F}}{\kk{i}{F}^{p(G)}\times \kk{i+1}{F}^{q(G)}}{\kk{i}{F\otimes G'}}$$
is a group isomorphism. If $\oaa{F}{\phi}{F'}$ is a morphism in \frm then the diagram
$$\begin{CD}
\kk{i}{F}^{p(G)}\times \kk{i+1}{F}^{q(G)}@>\Phi _{i,G,F} >>\kk{i}{F\otimes G}@>\kk{i}{id_F\otimes \varphi }>>\kk{i}{F\otimes G'}\\
@VV\kk{i}{\phi }^{p(G)}\times \kk{i+1}{\phi }^{q(G)}V@VV\kk{i}{\phi \otimes id_G}V@V\kk{i}{\phi\otimes id_{G'}}VV \\        
\kk{i}{F'}^{p(G)}\times \kk{i+1}{F'}^{q(G)}@>>\Phi _{i,G,F'}>\kk{i}{F'\otimes G}@>>\kk{i}{id_{F'}\otimes \varphi }>\kk{i}{F'\otimes G'}\\
\end{CD}$$
is commutative and the first assertion follows.

Assume now that $G$ and $G'$ are nuclear, let 
$$(\oc{F_1}{\phi _1}{F_2}{\phi _2}{F_3})\in G_\Upsilon \,,$$
and let $\delta _i$ be its associated index maps. By the commutativity of the index maps (\axi{27.9'e} a)) the diagram
$$\begin{CD}
\kk{i}{F_3}^{p(G)}\times \kk{i+1}{F_3}^{q(G)}@>\delta _i^{p(G)}\times \delta _{i+1}^{q(G)}>>\kk{i+1}{F_1}^{p(G)}\times \kk{i}{F_1}^{q(G)}\\
@V\Phi _{i,G,F_3}VV@VV\Phi _{(i+1),G,F_1}V \\        
\kk{i}{F_3\otimes G}@>>\delta _{G,i}>\kk{i+1}{F_1\otimes G}\\
@V\kk{i}{id_{F_3}\otimes \varphi }VV@VV\kk{i+1}{id_{F_1}\otimes \varphi }V\\
\kk{i}{F_3\otimes G'}@>>\delta _{G',i}>\kk{i+1}{F_1\otimes G'}
\end{CD}$$
is commutative. Since the maps of the columns are group isomorphisms, it follows by the above, that the diagram
$$\begin{CD}
\kk{i}{F_3}^{p(G')}\times \kk{i+1}{F_3}^{q(G')}@>\delta _i^{p(G')}\times \delta _{i+1}^{q(G')} >>\kk{i+1}{F_1}^{p(G')}\times \kk{i}{F_1}^{q(G')}\\
@V\Phi _{i,G',F_3}VV@VV\Phi _{(i+1),G',F_1}V \\        
\kk{i}{F_3\otimes G'}@>>\delta _{G',i}>\kk{i+1}{F_1\otimes G'}\\
\end{CD}$$
is also commutative.

d) follows from a) and c).\qed

\begin{p}\label{27.8'a}
Let $G$ be a nuclear C*-algebra belonging to $\Upsilon $.
\begin{enumerate}
\item Every split exact sequence in \frm belongs to $G_\Upsilon $.
\item Every exact sequence in \frm 
$$\oc{F_1}{}{F_2}{}{F_3}$$
for which $F_1$ or $F_3$ is homotopic to $0$ belongs to $G_\Upsilon $ .
\end{enumerate}
\end{p}

a) Let
$$\od{F_1}{\varphi }{F_2}{\psi }{\lambda }{F_3}$$
be a split exact sequence in \frm and let $\delta _i$ be its index maps. By \pr{23.4'} a),
$$\og{F_1\otimes G}{\varphi \otimes id_G}{F_2\otimes G}{\psi \otimes id_G}{\lambda \otimes id_G}{F_3\otimes G}{20}{20}{20}$$
is split exact and so by \pr{21.10'}, $\delta _i=\delta _{G,i}=0$.

b) By \pr{29.3'} c), $F_1\otimes G$ or $F_3\otimes G$ is null-homotopic and so K-null. Thus by the six-term axiom (\axi{27.9'd}), $\delta _i=\delta _{G,i}=0$, where $\delta _i$ denote the index maps associated to
$$\oc{F_1}{}{F_2}{}{F_3}\;.\qedd$$ 

\begin{p}\label{6.7'}
Let 
$$\oc{G_1}{\varphi }{G_2}{\psi }{G_3}$$
be an exact sequence in \frcc such that $G_3$ is nuclear.
\begin{enumerate}
\item Assume $G_1$ is $\Upsilon $-null.
\begin{enumerate}
\item $\mac{\kk{i}{id_F\otimes \psi }}{\kk{i}{F\otimes G_2}}{\kk{i}{F\otimes G_3}}$ is a group isomorphism.
\item If $G_2\in \Upsilon $ or $G_3\in \Upsilon $ then
$$G_2,G_3\in \Upsilon\, ,\qquad 
p(G_2)=p(G_3)\,,\qquad q(G_2)=q(G_3)\,,$$
$$ \Phi _{i,G_3,F}=\kk{i}{id_F\otimes \psi }\circ \Phi _{i,G_2,F}\;.$$
If in addition $G_2$ is nuclear then $(G_2)_\Upsilon =(G_3)_\Upsilon $.
\end{enumerate}
\item Assume $G_2$ is $\Upsilon $-null and let $\delta^F _i$ denote the index maps of the exact sequence in \frm
$$\of{F\otimes G_1}{id_F\otimes \varphi }{F\otimes G_2}{id_F\otimes \psi }{F\otimes G_3}{20}{20}\;.$$
\begin{enumerate}
\item $\mac{\delta^F _i}{\kk{i}{F\otimes G_3}}{\kk{i+1}{F\otimes G_1}}$ is a group isomorphism.
\item If $G_1\in \Upsilon $ or $G_3\in \Upsilon $ then
$$G_1,G_3\in \Upsilon\,,\qquad 
p(G_1)=q(G_3)\,,\qquad q(G_1)=p(G_3)\,,$$
$$ \Phi _{i,G_3,F}=\Phi _{(i+1),G_1,F}\circ \delta^F _i\;.$$
\end{enumerate}
\item Assume $G_3$ is $\Upsilon $-null.
\begin{enumerate}
\item $\mac{\kk{i}{id_F\otimes \varphi }}{\kk{i}{F\otimes G_1}}{\kk{i}{F\otimes G_2}}$ is a group isomorphism.
\item If $G_1\in \Upsilon $ or $G_2\in \Upsilon $ then
$$G_1,G_2\in \Upsilon\, ,\qquad
 p(G_1)=p(G_2)\,,\qquad q(G_1)=q(G_2)\,,$$
$$ \Phi _{i,G_2,F}=\kk{i}{id_F\otimes \varphi }\circ \Phi _{i,G_1,F}\;.$$
If in addition $G_1$ and $G_2$ are nuclear then $(G_1)_\Upsilon =(G_2)_\Upsilon $.
\end{enumerate}
\end{enumerate}
\end{p}

By \pr{29.3'a} a), the sequence in \frm
$$\of{F\otimes G_1}{id_{F}\otimes \varphi }{F\otimes G_2}{id_F\otimes \psi }{F\otimes G_3}{18}{18}$$
is exact. If $G_j$ is $\Upsilon $-null then $F\otimes G_j$ is K-null so $a_1),b_1),c_1)$ follow from \pr{28.9'd} a),b).

$a_2)$ By $a_1)$, it is easy to see that
$$G_2,G_3\in \Upsilon\, ,\qquad 
p(G_2)=p(G_3)\,,\qquad q(G_2)=q(G_3)\,,$$
$$ \Phi _{i,G_3,F}=\kk{i}{id_F\otimes \psi }\circ \Phi _{i,G_2,F}\;.$$

Assume now $G_2 $  nuclear. Let 
$$\oc{F_1}{\phi _1}{F_2}{\phi _2}{F_3}$$
belong to $(G_2)_\Upsilon $ or $(G_3)_\Upsilon $  and let $\delta _i$ be its associated index maps. Consider the diagram
$$\begin{CD}
\kk{i}{F_3}^{p(G_2)}\times \kk{i+1}{F_3}^{q(G_2)}@>\delta _i^{p(G_2)}\times \delta _{i+1}^{q(G_2)}>>\kk{i+1}{F_1}^{p(G_2)}\times \kk{i}{F_1}^{q(G_2)}\\
@V\Phi _{i,G_2,F_3}VV@VV\Phi _{(i+1),G_2,F_1}V \\        
\kk{i}{F_3\otimes G_2}@>\delta _{G_2,i}>>\kk{i+1}{F_1\otimes G_2}\\
@V\kk{i}{id_{F_3}\otimes \psi  }VV@VV\kk{i+1}{id_{F_1}\otimes \psi  }V\\
\kk{i}{F_3\otimes G_3}@>>\delta _{G_3,i}>\kk{i+1}{F_1\otimes G_3}\\
@A\Phi _{i,G_3,F_3}AA@AA\Phi _{(i+1),G_3,F_1}A\\
\kk{i}{F_3}^{p(G_3)}\times \kk{i+1}{F_3}^{q(G_3)}@>>\delta _i^{p(G_3)}\times \delta _{i+1}^{q(G_3)}>\kk{i+1}{F_1}^{p(G_3)}\times \kk{i}{F_1}^{q(G_3)}.
\end{CD}$$
Its upper part or lower part is commutative and the maps of the columns are group isomorphisms. It follows, by the above, that the diagram is  commutative. Thus $(G_2)_\Upsilon =(G_3)_\Upsilon $.

$b_2)$ Let $\oaa{F}{\phi }{F'}$ be a morphism in \frm. Then the diagram
$$\begin{CD}
0@>>>\kk{i}{F\otimes G_1}@>\kk{i}{id_F\otimes \varphi } >>\kk{i}{F\otimes G_2}@>\kk{i}{id_F\otimes \psi } >>\\
@.      @V\kk{i}{\phi \otimes id_{G_1}}VV  @VV\kk{i}{\phi \otimes id_{G_2}}V \\
0@>>>\kk{i}{F'\otimes G_1}@>>\kk{i}{id_{F'}\otimes \varphi }>\kk{i}{F'\otimes G_2}@>>\kk{i}{id_{F'}\otimes \psi }>
\end{CD}$$
$$\begin{CD}
@>\kk{i}{id_F\otimes \varphi } >>\kk{i}{F\otimes G_2}@>\kk{i}{id_F\otimes \psi } >>\kk{i}{F\otimes G_3}@>>>0\\
@.  @V\kk{i}{\phi \otimes id_{G_2}}VV@VV\kk{i}{\phi \otimes id_{G_3}}V\\
@>>\kk{i}{id_{F'}\otimes \varphi }>\kk{i}{F'\otimes G_2}@>>\kk{i}{id_{F'}\otimes \psi }>\kk{i}{F'\otimes G_3}@>>>0
\end{CD}$$
is commutative and has exact rows. By the commutativity of the index maps (\axi{27.9'e}), the diagram
$$\begin{CD}
\kk{i}{F\otimes G_3}@>\delta _i^F>>\kk{i+1}{F\otimes G_1}\\
@V\kk{i}{\phi \otimes id_{G_3}}VV@VV\kk{i+1}{\phi \otimes id_{G_1}}V \\        
\kk{i}{F'\otimes G_3}@>>\delta^{F'} _{i}>\kk{i+1}{F'\otimes G_1}\\
\end{CD}$$
is commutative. By $b_1)$,
 $$G_1,G_3\in \Upsilon\,,\qquad 
p(G_1)=q(G_3)\,,\qquad q(G_1)=p(G_3)\,,$$
$$ \Phi _{i,G_3,F}=\Phi _{(i+1),G_1,F}\circ \delta^F _i\;.$$

$c_2)$ The proof is similar to the proof of $a_2)$.\qed 

\begin{p}\label{6.7'a}
Let
$$\od{G_1}{\varphi }{G_2}{\psi }{\lambda }{G_3}$$
be a split exact sequence in \frcc.
\begin{enumerate}
\item If $G_1,G_3\in \Upsilon $ then
$$G_2\in \Upsilon \,,\qquad p(G_2)=p(G_1)+p(G_3)\,,\qquad q(G_2)=q(G_1)+q(G_3)\,,$$
$$\Phi _{i,G_2,F}=(\kk{i}{id_F\otimes \varphi }\times \kk{i}{id_F\otimes \lambda })\circ (\Phi _{i,G_1,F}\times \Phi _{i,G_3,F})\;.$$
\item If in addition $G_1$,$G_2$,and $G_3$ are nuclear then $(G_1)_\Upsilon \cap (G_3)_\Upsilon \subset (G_2)_\Upsilon  $.
\end{enumerate}  
\end{p}

a) By \pr{26.3'a} b), the sequence
$$\og{\kk{i}{F\otimes G_1}}{\kk{i}{id_F\otimes \varphi }}{\kk{i}{F\otimes G_2}}{\kk{i}{id_F\otimes \psi }}{\kk{i}{id_F\otimes \lambda }}{\kk{i}{F\otimes G_3}}{30}{30}{30}$$
is split exact. Thus the maps
$$\Big(\kk{i}{F}^{p(G_1)}\times \kk{i+1}{F}^{q(G_1)}\Big)\times \Big(\kk{i}{F}^{p(G_3)}\times \kk{i+1}{F}^{q(G_3)}\Big) $$
$$\stackrel{\Phi _{i,G_1,F}\times \Phi _{i,G_3,F}}{\rule[2.3pt]{50pt}{0.3pt}\!\!\longrightarrow  }
\kk{i}{F\otimes G_1}\times \kk{i}{F\otimes G_3} 
\stackrel{\kk{i}{id_F\otimes \varphi }\times \kk{i}{id_F\otimes \lambda }}{\rule[2.3pt]{75pt}{0.3pt}\!\!\longrightarrow }
\kk{i}{F\otimes G_2}$$
are group isomorphisms.

Let $\oa{F}{\phi }{F'}$ be a morphism in \frm. Since the diagram with split exact rows
$$\og{F\otimes G_1}{id_F\otimes \varphi }{F\otimes G_2}{id_F\otimes \psi }{id_F\otimes \lambda }{F\otimes G_3}{20}{20}{20}\,,$$
$$\og{F'\otimes G_1}{id_{F'}\otimes \varphi }{F'\otimes G_2}{id_{F'}\otimes \psi }{id_{F'}\otimes \lambda }{F'\otimes G_3}{20}{20}{20}\,,$$
(\pr{26.3'a} a)) and with columns $\phi \otimes id_{G_1}$, $\phi \otimes id_{G_2}$, and $\phi \otimes id_{G_3}$ is commutative, the assertion follows from \pr{26.3'a} b).

b) Let 
$$(\oc{F_1}{\phi _1}{F_2}{\phi _2}{F_3})\in (G_1)_\Upsilon \cap (G_3)_\Upsilon $$
 and let $\delta _i$ be its associated index maps. 
Consider the diagram (by a))
$$\begin{CD}
\kk{i}{F_3}^{p(G_2)}\times \kk{i+1}{F_3}^{q(G_2)}@>\delta _i^{p(G_2)}\times \delta _{i+1}^{q(G_2)}>>\kk{i+1}{F_1}^{p(G_2)}\times \kk{i}{F_1}^{q(G_2)}\\
@V\Phi _{i,G_1,F_3}\times \Phi _{i,G_3,F_3}VV@V\Phi _{(i+1),G_1,F_1}\times \Phi _{(i+1),G_3,F_1}VV \\        
\kk{i}{F_3\otimes G_1}\times \kk{i}{F_3\otimes G_3}@>>\delta _{G_1,i}\times \delta _{G_3,i}>\kk{i+1}{F_1\otimes G_1}\times \kk{i+1}{F_1\otimes G_3}\\
@VAVV@V\kk{i+1}{id_{F_1}\otimes \varphi }\times \kk{i+1}{id_{F_1}\otimes \lambda }VV\\
\kk{i}{F_3\otimes G_2}@>>\delta _{G_2,i}>\kk{i+1}{F_1\otimes G_2}\\
@A\Phi _{i,G_2,F_3}AA@AA\Phi _{(i+1),G_2,F_1}A\\
\kk{i}{F_3}^{p(G_2)}\times \kk{i+1}{F_3}^{q(G_2)}@>\delta _i^{p(G_2)}\times \delta _{i+1}^{q(G_2)}>>\kk{i+1}{F_1}^{p(G_2)}\times \kk{i}{F_1}^{q(G_2)},
\end{CD}$$
where
$$A:=\kk{i}{id_{F_3}\otimes \varphi }\times \kk{i}{id_{F_3}\otimes \lambda }\;.$$
Its upper part is commutative and the maps of the columns are group isomorphisms. It follows that the lower part of the diagram is also commutative.\qed

\begin{co}\label{6.7'b}
If $G\in \Upsilon $ then $\tilde{G}\in \Upsilon  $, $p(\tilde{G})=p(G)+1$, $q(\tilde{G})=q(G)$. If in addition $G$ and $\tilde{G} $ are nuclear then $G_\Upsilon \subset \tilde{G}_\Upsilon $.\qed
\end{co}

\begin{p}\label{5.7'b}
Let $(G_j)_{j\in J}$ be a finite family in $\Upsilon $.
\begin{enumerate}
\item $$G:=\pro{j\in J}G_j\in \Upsilon \,,\qquad p(G)=\si{j\in J}p(G_j)\,,\qquad q(G)=\si{j\in J}q(G_j)\,,$$
$$\Phi _{i,G,F}=\left(\pro{j\in J}\Phi _{i,G_j,F}\right)\circ \Phi _{(F\otimes G_j)_{j\in J},i}\;.$$
In particular if $G_j$ is $\Upsilon $-null for every $j\in J$ then $G$ is $\Upsilon $-null.
\item If in addition $G$ and all $G_j$, $j\in J$, are nuclear then
$$\bigcap _{j\in J}(G_j)_\Upsilon \subset G_\Upsilon \;.$$
\item $\bc^J\in \Upsilon $, $p(\bc^J)=Card\,J$, $q(\bc^J)=0$, and every exact sequence in \frm belongs to $\left(\bc^J\right)_\Upsilon $.
\end{enumerate}
\end{p}

a) We put 
$$\bar{p}:=\si{j\otimes J}p(G_j)\,,\qquad\qquad \bar{q}:=\si{j\otimes J}q(G_j)\;. $$
Since
$$F\otimes \pro{j\in J}G_j\approx \pro{j\in J}(F\otimes G_j)\,,$$
by \pr{28.9'a}, the maps
$$\pro{j\in J}\Phi _{i,G_j,F}:\kk{i}{F}^{\bar{p} }\times \kk{i+1}{F}^{\bar{q} }=\pro{j\in J}\left(\kk{i}{F}^{p(G_j)}\times \kk{i+1}{F}^{q(G_j)}\right)\longrightarrow $$
$$\longrightarrow {\pro{j\in J}\kk{i}{F\otimes G_j}\stackrel{\Phi _{(F\otimes G_j)_{j\in J},i}}{\rule[2.5pt]{40pt}{0.6pt}\!\!\longrightarrow } \kk{i}{F\otimes G}}$$
are group isomorphisms. Let $\oaa{F}{\phi }{F'}$ be a morphism in \frm. The diagram 
$$\begin{CD}
\kk{i}{F}^{\bar{p} }\times \kk{i+1}{F}^{\bar{q} }@>\pro{j\in J}\Phi _{i,G_j,F}>>\pro{j\in J}\kk{i}{F\otimes G_j}\\
@V\kk{i}{\phi }^{\bar{p} }\times \kk{i+1}{\phi }^{\bar{q} }VV@VV\pro{j\in J}\kk{i}{\phi \otimes id_{G_j}}V\\        
\kk{i}{F'}^{\bar{p}}\times \kk{i+1}{F'}^{\bar{q} }@>>\pro{j\in J}\Phi _{i,G_j,F'}>\pro{j\in J}\kk{i}{F'\otimes G_j}\\
\end{CD}$$
is obviously commutative and by \pr{21.10'a}  the diagram
$$\begin{CD}
\pro{j\in J}\kk{i}{F\otimes G_j}@>\Phi _{(F\otimes G_j)_{j\in J},i}>>\kk{i}{F\otimes G}\\
@V\pro{j\in J}\kk{i}{\phi \otimes id_{G_j}}VV@VV\kk{i}{\phi \otimes id_{G}}V \\        
\pro{j\in J}\kk{i}{F'\otimes G_j}@>>\Phi _{(F'\otimes G_j)_{j\in J},i}>\kk{i}{F'\otimes G}\\
\end{CD}$$
is also commutative and this proves the assertion.

b) follows from \pr{6.7'a} by complete induction.

c) follows from a), b), and \pr{5.7'a} b).\qed

\begin{p}\label{12.7'}
Let $J$ be a finite set and for every $j\in J$ let 
$$\oc{F_{j,1}}{\varphi _j}{F_{j,2}}{\psi _j}{F_{j,3}}$$
be an exact sequence in \frm and $\delta _{j,i}$ its associated index maps. For every $k\in \z{1,2,3}$ put
$$F_k:=\pro{j\in J}F_{j,k}$$
and for every $j\in J$ denote by 
$$\mac{\varphi _{j,k}}{F_{j,k}}{F_k}\,,\qquad \mac{\psi _{j,k}}{F_k}{F_{j,k}}$$
the canonical inclusion and projection, respectively. Then
$$\oc{F_1}{\pro{j\in J}\varphi _j}{F_2}{\pro{j\in J}\psi _j}{F_3}$$
is an exact sequence in \frm and if we denote by $\delta _i$ its index maps then the diagram
$$\begin{CD}
\pro{j\in J}\kk{i}{F_{j,3}}@>\Psi _{3,i}>>\kk{i}{F_3}\\
@V\pro{j\in J}\delta _{j,i}VV@VV\delta _iV \\        
\pro{j\in J}\kk{i+1}{F_{j,1}}@>>\Psi _{1,(i+1)}>\kk{i+1}{F_1}\\
\end{CD}$$
is commutative, where for every $k\in \z{1,3}$,
$$\mae{\Psi _{k,i}}{\pro{j\in J}\kk{i}{F_{j,k}}}{\kk{i}{F_k}}{(a_j)_{j\in J}}{\si{j\in J}\kk{i}{\varphi _{j,k}}a_j}\;.$$
\end{p}

For every $j\in J$ the diagram
$$\begin{CD}
0@>>>F_{j,1}@>\varphi_j >>F_{j,2}@>\psi_j >>F_{j,3}@>>>0\\
@.      @V\varphi _{j,1} VV  @V\varphi _{j,2} VV @VV\varphi _{j,3}V@.\\
0@>>>F_1@>>\pro{j\in J}\varphi_j >F_2@>>\pro{j\in J}\psi_j >F_3@>>>0\,,
\end{CD}$$
is commutative. By the commutativity of the index maps (\axi{27.9'e}), the diagram
$$\begin{CD}
\kk{i}{F_{j,3}}@>\kk{i}{\varphi _{j,3}} >>\kk{i}{F_3}\\
@V\delta _{j,i}VV@VV\delta _iV \\        
\kk{i+1}{F_{j,1}}@>>\kk{i+1}{\varphi _{j,1}}>\kk{i+1}{F_1}\\
\end{CD}$$
is commutative. Let $(a_j)_{j\in J}\in \pro{j\in J}\kk{i}{F_{j,3}}$. Then
$$\delta _i\Psi _{3,i}(a_j)_{j\in J}=\delta _i\si{j\in J}\kk{i}{\varphi _{j,3}}a_j=$$
$$=\si{j\in J}\kk{i+1}{\varphi _{j,1}}\delta _{j,i}a_j=\Psi _{1,(i+1)}\left(\pro{j\in J}\delta _{j,i}\right)(a_j)_{j\in J}\;.$$
Thus the diagram
$$\begin{CD}
\pro{j\in J}\kk{i}{F_{j,3}}@>\Psi _{3,i}>>\kk{i}{F_3}\\
@V\pro{j\in J}\delta _{j,i}VV@VV\delta _iV \\        
\pro{j\in J}\kk{i+1}{F_{j,1}}@>>\Psi _{1,(i+1)}>\kk{i+1}{F_1}\\
\end{CD}$$
is commutative.\qed

\begin{p}\label{27.3'b}
Let $(G_j)_{j\in J}$ be a finite family in $\Upsilon $.
\begin{enumerate}
\item $$G:=\bigotimes_{j\in J}G_j\in \Upsilon \,,$$
$$p(G)=\frac{1}{2}\left(\pro{j\in J}(p(G_j)+q(G_j))+\pro{j\in J}(p(G_j)-q(G_j))\right)\,,$$
$$q(G)=\frac{1}{2}\left(\pro{j\in J}(p(G_j)+q(G_j))-\pro{j\in J}(p(G_j)-q(G_j))\right)\;.$$ 
\item If $G_{j_0}$ is K-null for a $j_0\in J$ then $F\otimes \left(\bigotimes\limits _{j\in J}G_j\right)$ is also K-null.
\item If $p(G_{j_0})=q(G{j_0})$ for a $j_0\in J$ then $p(G)=q(G)$.
\item Let $j_0\in J$, $J':=J\setminus \z{j_0}$, and $G':=\bigotimes\limits_{j\in J'}G_j $.
\begin{enumerate}
\item If $p(G_{j_0})=1,\,q(G_{j_0})=0$ then $p(G')=p(G),\,q(G')=q(G)$.
\item If $p(G_{j_0})=0,\,q(G_{j_0})=1$ then $p(G')=q(G),\,q(G')=p(G)$.
\end{enumerate}
\item If we put 
$$H:=\bigotimes _{j\in J}\tilde{G_j}\qquad \emph{and}\qquad G_I:=\bigotimes_{j\in I}G_j  $$
for every $I\subset J$ then
$$H\in \Upsilon, \qquad p(H)=\si{I\subset J}p(G_I),\qquad q(H)=\si{I\subset J}q(G_I);.$$
\item If in addition $G$ and all $(G_j)_{j\in J}$ are nuclear then 
$$\bigcap_{j\in J}(G_j)_\Upsilon \subset G_\Upsilon \;. $$
\end{enumerate} 
\end{p}

a) Assume first $J=\z{1,2}$. The maps
$$\kk{i}{F}^{p(G_1)p(G_2)+q(G_1)q(G_2)}\times \kk{i+1}{F}^{p(G_1)q(G_2)+p(G_2)q(G_1)}=$$
$$=\Big(\kk{i}{F}^{p(G_1)}\times \kk{i+1}{F}^{q(G_1)}\Big)^{p(G_2)}\times \Big(\kk{i+1}{F}^{p(G_1)}\times \kk{i}{F}^{q(G_1)}\Big)^{q(G_2)}\longrightarrow $$
$$\stackrel{(\Phi _{i,G_1,F})^{p(G_2)}\times (\Phi _{(i+1),G_1,F})^{q(G_2)}}{\rule[3pt]{100pt}{0.3pt}\!\!\longrightarrow }$$
$$\longrightarrow \kk{i}{F\otimes G_1}^{p(G_2)}\times \kk{i+1}{F\otimes G_1}^{q(G_2)}\longrightarrow $$
$$\stackrel{\Phi _{i,G_2,F\otimes G_1}}{\rule[3pt]{40pt}{0.3pt}\!\!\longrightarrow }$$
$$\longrightarrow \kk{i}{(F\otimes G_1)\otimes G_2}\approx \kk{i}{F\otimes (G_1\otimes G_2)}$$
are group isomorphisms and
$$p(G_1\otimes G_2):=p(G_1)p(G_2)+q(G_1)q(G_2)=$$
$$=\frac{1}{2}\lbrack(p(G_1)+q(G_1))(p(G_2)+q(G_2))+(p(G_1)-q(G_1))(p(G_2)-q(G_2))\rbrack\,,$$
$$q(G_1\otimes G_2):=p(G_1)q(G_2)+p(G_2)q(G_1)=$$
$$=\frac{1}{2}\lbrack(p(G_1)+q(G_1))(p(G_2)+q(G_2))-(p(G_1)-q(G_1))(p(G_2)-q(G_2))\rbrack\;.$$
If $\oaa{F}{\phi }{F'}$ is a morphism in \frm then the diagrams
$$\begin{CD}
\kk{i}{F\otimes G_1}^{p(G_2)}\times \kk{i+1}{F\otimes G_1}^{q(G_2)}@>\Phi _{i,G_2,(F\otimes G_1)} >>\kk{i}{(F\otimes G_1)\otimes G_2}\\
@VV\kk{i}{\phi\otimes id_{G_1} }^{p(G_2)}\times \kk{i+1}{\varphi\otimes id_{G_1} }^{q(G_2)}V@VV\kk{i}{(\phi \otimes id_{G_1})\otimes id_{G_2}}V \\        
\kk{i}{F'\otimes G_1}^{p(G_2)}\times \kk{i+1}{F'\otimes G_1}^{q(G_2)}@>>\Phi _{i,G_2,(F'\otimes G_1)}>\kk{i}{(F'\otimes G_1)\otimes G_2}\\
\end{CD}$$
$$\begin{CD}
\kk{i}{(F\otimes G_1)\otimes G_2}@>\approx >>\kk{i}{F\otimes (G_1\otimes G_2)}\\
@VV\kk{i}{(\phi \otimes id_{G_1})\otimes id_{G_2}}V@VV\kk{i}{\phi \otimes id_{(G_1\otimes G_2)}}V \\        
\kk{i}{(F'\otimes G_1)\otimes G_2}@>>\approx >\kk{i}{F'\otimes (G_1\otimes G_2)}\\
\end{CD}$$
are commutative, which proves the assertion in this case.

The general case is obtained now by induction with respect to $Card\,J$. Let $Card\,J>1$, $k\in J$, $J':=J\setminus \z{k}$, $G':=\otimes _{j\in J'}G_j$ and assume the assertion holds for $J'$. By the above,
$$p(G)=\frac{1}{2}[(p(G')+q(G'))(p(G_k)+q(G_k))+(p(G')-q(G'))(p(G_k)-q(G_k))]=$$
$$=\frac{1}{2}\left(\pro{j\in J'}(p(G_j)+q(G_j))(p(G_k)+q(G_k))+\right.$$
$$\left.+\pro{j\in J'}(p(G_j)-q(G_j))(p(G_k)-q(G_k))\right)=$$
$$=\frac{1}{2}\left(\pro{j\in J}(p(G_j)+q(G_j))+\pro{j\in J}(p(G_j)-q(G_j))\right)\,,$$  
$$q(G)=\frac{1}{2}[(p(G')+q(G'))(p(G_k)+q(G_k))-(p(G')-q(G'))(p(G_k)-q(G_k))]=$$
$$=\frac{1}{2}\left(\pro{j\in J'}(p(G_j)+q(G_j))(p(G_k)+q(G_k))-\right.$$
$$\left.-\pro{j\in J'}(p(G_j)-q(G_j))(p(G_k)-q(G_k))\right)=$$
$$=\frac{1}{2}\left(\pro{j\in J}(p(G_j)+q(G_j))-\pro{j\in J}(p(G_j)-q(G_j))\right)\;.$$

b), c), and d) follow directly from a).

e) By \cor{6.7'b}, $\tilde{G_j}\in \Upsilon  $ for every $j\in J$. By a) and \pr{27.3'}, $H\in \Upsilon $,
$$\kk{i}{F\otimes H}\approx \pro{I\subset J}\kk{i}{F\otimes G_I}\approx 
 \pro{I\subset J}\Big(\kk{i}{F}^{p(G_{I})}\times \kk{i+1}{F}^{q(G_I)}\Big)=$$
$$=\kk{i}{F}^{\si{I\subset J}p(G_I)}\times \kk{i+1}{F}^{\si{I\subset J}q(G_I)}\;.$$

f) Assume first $J:=\z{1,2}$, let 
$$(\oc{F_1}{\phi _1}{F_2}{\phi _2}{F_3})\in (G_1)_\Upsilon \cap (G_2)_\Upsilon \,,$$
 and let $\delta _i$ be its index maps.
Then (by a)) the diagram
$$\begin{CD}
\kk{i}{F_3}^{p(G)}\times \kk{i+1}{F_3}^{q(G)}@>\delta _i^{p(G)}\times \delta _{i+1}^{q(G)}>>A\\
@V\Phi _{i,G_1,F_3}^{p(G_2)}\times \Phi _{(i+1),G_1,F_3}^{q(G_2)}VV@V\Phi _{(i+1),G_1,F_1}^{p(G_2)}\times \Phi _{i,G_1,F_1}^{q(G_2)}VV \\ 
\kk{i}{F_3\otimes G_1}^{p(G_2)}\times \kk{i+1}{F_3\otimes G_1}^{q(G_2)}@>\delta _{G_1,i}^{p(G_1)}\times \delta _{G_1,i}^{q(G_2)}>>B\\ 
@V\Phi _{i,G_2,(F_3\otimes G_1)}VV@V\Phi _{(i+1),G_2,(F_1\otimes G_1)}VV\\
\kk{i}{(F_3\otimes G_1)\otimes G_2}\approx \kk{i}{F_3\otimes (G_1\otimes G_2)}@>>\delta_{G_1,(G_2,i)}>C
\end{CD}$$
is commutative, where
$$A:=\kk{i+1}{F_1}^{p(G)}\times \kk{i}{F_1}^{q(G)}\,,$$
$$B:=\kk{i+1}{F_1\otimes G_1}^{p(G_2)}\times \kk{i}{F_1\otimes G_1}^{q(G_2)}\,,$$
$$C:=\kk{i+1}{((F_1\otimes G_1)\otimes G_2)}\approx \kk{i+1}{(F_1\otimes (G_1\otimes G_2))}\;.$$
Thus 
$$(\oc{F_1}{\phi _1}{F_2}{\phi _2}{F_3})\in G_\Upsilon\;.$$

The general case follows by induction with respect to $Card\,J$.\qed

\begin{co}\label{27.3'c}
Let $G\in \Upsilon $, $\bbn$, and $H:=\otimes _{j\in \bnn{n}}G$. Then $H\in \Upsilon $, $G_\Upsilon \subset H_\Upsilon $, and
$$p(H)=\frac{1}{2}\left((p(G)+q(G))^n+(p(G)-q(G))^n\right)\,,$$
$$q(H)=\frac{1}{2}\left((p(G)+q(G))^n-(p(G)-q(G))^n\right)\;.$$
\end{co}

The assertion follows from \pr{27.3'b} a).\qed

\begin{p}\label{8.7'}
Let $(G_1,G_2,G_3)$ be an \frcc -triple such that $G_1/G_3$ and $G_2/G_3$ are nuclear, $G_2$ is $\Upsilon $-null, and $G_1,G_3\in \Upsilon $. We use the notation of the triple theorem \emph{(\h{28.9'f} a))} associated to the \frm -triple 
$$(F\otimes G_1,F\otimes G_2,F\otimes G_3)$$
 \emph{(\pr{29.3'a})}, put $\varphi :=\varphi _{1,2}/(F\otimes G_3)$ \emph{(as in \pr{28.9'e} a))}, and denote by
$$\mac{\Psi _{F,i}}{\kk{i}{F\otimes G_1}\times \kk{i+1}{F\otimes G_3}}{\kk{i}{F\otimes (G_1/G_3)}},$$
$$(a,b)\longmapsto \kk{i}{\psi _{1,3}}a+\Phi _ib$$
the corresponding group isomorphism \emph{(\h{28.9'f} $a_4)$, \pr{29.3'a})}. Then
$$G_1/G_3\in \Upsilon ,\qquad p(G_1/G_3)=p(G_1)+q(G_3),\qquad q(G_1/G_3)=q(G_1)+p(G_3),$$
$$\Phi _{i,(G_1/G_3),F}=\Psi _{F,i}\circ (\Phi _{i,G_1,F}\times \Phi _{(i+1),G_3,F})\;.$$ 
\end{p}

Since $G_1,G_3\in \Upsilon $, the map
$$\Psi _{F,i}\circ (\Phi _{i,G_1,F}\times \Phi _{(i+1),G_3,F}):
\Big(\kk{i}{F}^{p(G_1)}\times \kk{i+1}{F}^{q(G_1)}\Big)\times $$
$$\times \Big(\kk{i+1}{F}^{p(G_3)}\times \kk{i}{F}^{q(G_3)}\Big)\longrightarrow
\kk{i}{F\otimes (G_1/G_3)}$$
is a group isomorphism. We put
$$\bar{p}(G_1/G_3):=p(G_1)+Q(G_3)\,,\qquad\qquad \bar{q}(G_1/G_3):=q(G_1)+p(G_3)\,,  $$
$$\bar{\Phi} _{i,G_1/G_3,F}:=\Psi _{F,i}\circ (\Phi _{i,G_1,F}\times \Phi _{(i+1),G_3,F})\;.$$

Let $\oaa{F}{\phi }{F'}$ be a morphism in \frm. We mark with a prime the above notation associated to $F'$. By the commutativity of the index maps (\axi{27.9'e}),
$$\kk{i+1}{\phi \otimes id_{G_3}}\circ \delta _{2,3,i}=\delta '_{2,3,i}\circ \kk{i}{\phi \otimes id_{(G_2/G_3)}}\;.$$
Moreover
$$\kk{i}{\phi \otimes id_{(G_1/G_3)}}\circ \kk{i}{\varphi }=\kk{i}{\varphi '}\circ \kk{i}{\phi \otimes id_{(G_2/G_3)}}\,,$$
$$\kk{i}{\phi \otimes id_{(G_1/G_3)}}\circ \kk{i}{\psi _{1,3}}=\kk{i}{\psi _{1,3}'}\circ \kk{i}{\phi \otimes id_{G_1}}\;.$$
It follows
$$\kk{i}{\phi \otimes id_{(G_1/G_3)}}\circ \Phi _i=\kk{i}{\phi \otimes id_{(G_1/G_3)}}\circ \kk{i}{\varphi }\circ (\delta _{2,3,i})^{-1}=$$
$$=\kk{i}{\varphi '}\circ \kk{i}{\phi \otimes id_{(G_2/G_3)}}\circ (\delta _{2,3,i})^{-1}=$$
$$=\kk{i}{\varphi '}\circ (\delta '_{2,3,i})^{-1}\circ \kk{i+1}{\phi \otimes id_{G_3}}=\Phi '_i\circ \kk{i+1}{\phi \otimes id_{G_3}}\;.$$

We want to prove that the diagram
$$\begin{CD}
\kk{i}{F\otimes G_1}\times \kk{i+1}{F\otimes G_3}@>\Psi _{F,i}>>\kk{i}{F\otimes (G_1/G_3)}\\
@V\kk{i}{\phi \otimes id_{G_1}}\times \kk{i+1}{\phi \otimes id_{G_3}}VV@VV\kk{i}{\phi \otimes id_{(G_1/G_3)}}V \\        
\kk{i}{F'\otimes G_1}\times \kk{i+1}{F'\otimes G_3}@>>\Psi _{F',i}>\kk{i}{F'\otimes (G_1/G_3)}\\
\end{CD}$$
is commutative. For $(a,b)\in \kk{i}{F\otimes G_1}\times \kk{i+1}{F\otimes G_3}$, by the above,
$$\kk{i}{\phi \otimes id_{(G_1/G_3)}}\Psi_ {F,i}(a,b)=\kk{i}{\phi \otimes id_{(G_1/G_3)}}(\kk{i}{\psi _{1,3}}a+\Phi _ib)=$$
$$=\kk{i}{\phi \otimes id_{(G_1/G_3)}}\kk{i}{\psi _{1,3}}a+\kk{i}{\phi \otimes id_{(G_1/G_3)}}\Phi _ib=$$
$$=\kk{i}{\psi '_{1,3}}\kk{i}{\phi \otimes id_{G_1}}a+\Phi '_i\kk{i+1}{\phi \otimes id_{G_3}}b=$$
$$=\Psi _{F',i}(\kk{i}{\phi \otimes id_{G_1}}a,\kk{i+1}{\phi \otimes id_{G_3}}b)=$$
$$=\Psi _{F',i}(\kk{i}{\phi \otimes id_{G_1}}\times \kk{i+1}{\phi \otimes id_{G_3}})(a,b)\;.$$
Thus the above diagram is commutative. It follows, since $G_1,G_3\in \Upsilon $, that the diagram
$$\begin{CD}
\kk{i}{F}^{p(G_1/G_3)}\times \kk{i+1}{F}^{q(G_1/G_3)}@>\bar{\Phi} _{i,(G_1/G_3),F}>>\kk{i}{F\otimes (G_1/G_3)}\\
@V\kk{i}{\phi}^{p(G_1/G_3)}\times \kk{i+1}{\phi}^{q(G_1/G_3)}VV@VV\kk{i}{\phi \otimes id_{(G_1/G_3)}}V \\        
\kk{i}{F'}^{p(G_1/G_3)}\times \kk{i+1}{F'}^{q(G_1/G_3)}@>>\bar{\Phi} _{i,(G_1/G_3),F'}>\kk{i}{F'\otimes (G_1/G_3)}\\
\end{CD}$$
is commutative. Hence 
$$G_1/G_3\in \Upsilon ,\qquad p(G_1/G_3)=p(G_1)+q(G_3),\qquad q(G_1/G_3)=q(G_1)+p(G_3),$$
$$\Phi _{i,(G_1/G_3),F}=\Psi _{F,i}\circ (\Phi _{i,G_1,F}\times \Phi _{(i+1),G_3,F})\;.\qedd$$

\begin{p}\label{9.7'}
Let $(G_1,G_2,G_3)$ be an \frcc -triple such that $G_1/G_2$ and $G_1/G_3$ are nuclear, $G_1/G_3$ is $\Upsilon $-null, and $G_1,G_1/G_2\in \Upsilon $. We use the notation of the triple theorem \emph{(\h{28.9'f} b))} associated to the \frm -triple 
$$(F\otimes G_1,F\otimes G_2,F\otimes G_3)$$
 \emph{(\pr{29.3'a})}, assume $\psi _{12}$ K-null for all $E$-C*-algebras $F$, and denote by
$$\mac{\Psi _{F,i}}{\kk{i}{F\otimes G_1}\times \kk{i+1}{F\otimes (G_1/G_2)}}{\kk{i}{F\otimes G_2}},$$
$$(a,b)\longmapsto \Phi '_ia+\delta _{1,2,(i+1)}b$$
the corresponding group isomorphism \emph{(\h{28.9'f} $b_4)$, \pr{29.3'a})}. Then
$$G_2\in \Upsilon ,\qquad p(G_2)=p(G_1)+q(G_1/G_2),\qquad q(G_2)=q(G_1)+p(G_1/G_2),$$
$$\Phi _{i,G_2,F}=\Psi _{F,i}\circ (\Phi _{i,G_1,F}\times \Phi _{(i+1),(G_1/G_2),F})\;.$$ 
\end{p}

Since $G_1,G_1/G_2\in \Upsilon $, the map
$$\Psi _{F,i}\circ (\Phi _{i,G_1,F}\times \Phi _{(i+1),(G_1/G_2),F}):
\Big(\kk{i}{F}^{p(G_1)}\times \kk{i+1}{F}^{q(G_1)}\Big)\times $$
$$\times \Big(\kk{i+1}{F}^{p(G_1/G_2)}\times \kk{i}{F}^{q(G_1/G_2)}\Big)\longrightarrow
\kk{i}{F\otimes G_2}$$
is a group isomorphism. We put
$$\tilde{p}(G_2):=p(G_1)+q(G_1/G_2),\qquad\qquad \tilde{q}(G_2):=q(G_1)+p(G_1/G_2),$$
$$\tilde{\Phi} _{i,G_2,F}:=\Psi _{F,i}\circ (\Phi _{i,G_1,F}\times \Phi _{(i+1),(G_1/G_2),F})\;.$$

Let $\oaa{F}{\phi }{\bar{F}}$ be a morphism in \frm. We mark with a bar the above notation associated to $\bar{F}$. By the commutativity of the index maps (\axi{27.9'e}),
$$\kk{i}{\phi \otimes id_{G_2}}\circ \delta _{1,2,(i+1)}=\bar{\delta}_{1,2,(i+1)}\circ \kk{i+1}{\phi \otimes id_{(G_1/G_2)}}\;.$$
Moreover
$$\kk{i}{\phi \otimes id_{G_1}}\circ \kk{i}{\varphi _{1,3}}=\kk{i}{\bar{\varphi }_{1,3} }\circ \kk{i}{\phi \otimes id_{G_3}}\,,$$
$$\kk{i}{\phi \otimes id_{G_2}}\circ \kk{i}{\varphi _{2,3}}=\kk{i}{\bar{\varphi }_{2,3} }\circ \kk{i}{\phi \otimes id_{G_3}}\;.$$
It follows
$$\kk{i}{\phi \otimes id_{G_2}}\circ \Phi' _i=\kk{i}{\phi \otimes id_{G_2}}\circ \kk{i}{\varphi_{2,3} }\circ \kk{i}{\varphi _{1,3}}^{-1}=$$
$$=\kk{i}{\bar{\varphi}_{2,3}}\circ \kk{i}{\phi \otimes id_{G_3}}\circ \kk{i}{\varphi _{1,3}}^{-1}=$$
$$=\kk{i}{\bar{\varphi}_{2,3}}\circ \kk{i}{\bar{\varphi }_{1,3} }^{-1}\circ \kk{i}{\phi \otimes id_{G_1}}=\bar{\Phi}'_i\circ \kk{i}{\phi \otimes id_{G_1}}\;.$$

We want to prove that the diagram
$$\begin{CD}
\kk{i}{F\otimes G_1}\times \kk{i+1}{F\otimes (G_1/G_2)}@>\Psi _{F,i}>>\kk{i}{F\otimes G_2}\\
@V\kk{i}{\phi \otimes id_{G_1}}\times \kk{i+1}{\phi \otimes id_{(G_1/G_2)}}VV@VV\kk{i}{\phi \otimes id_{G_2}}V \\        
\kk{i}{\bar{F}\otimes G_1}\times \kk{i+1}{\bar{F}\otimes (G_1/G_2)}@>>\Psi _{\bar{F},i}>\kk{i}{\bar{F}\otimes G_2}\\
\end{CD}$$
is commutative. For $(a,b)\in \kk{i}{F\otimes G_1}\times \kk{i+1}{F\otimes (G_1/G_2)}$, by the above,
$$\kk{i}{\phi \otimes id_{G_2}}\Psi_ {F,i}(a,b)=\kk{i}{\phi \otimes id_{G_2}}(\Phi '_ia+\delta _{1,2,(i+1)}b)=$$
$$=\kk{i}{\phi \otimes id_{G_2}}\Phi '_ia+\kk{i}{\phi \otimes id_{G_2}}\delta _{1,2,(i+1)}b=$$
$$=\bar{\Phi }'_i\kk{i}{\phi \otimes id_{G_1}}a+\bar{\delta }_{1,2,(i+1)}\kk{i+1}{\phi \otimes id_{(G_1/G_2)}}b=$$
$$=\Psi _{\bar{F},i}(\kk{i}{\phi \otimes id_{G_1}}a,\kk{i+1}{\phi \otimes id_{(G_1/G_2)}}b)=$$
$$=\Psi _{\bar{F},i}\left(\kk{i}{\phi \otimes id_{G_1}}\times \kk{i+1}{\phi \otimes id_{(G_1/G_2)}}\right)(a,b)\;.$$
Thus the above diagram is commutative. Since $G_1,G_1/G_2\in \Upsilon $, It follows that the diagram
$$\begin{CD}
\kk{i}{F}^{\tilde{p}(G_2)}\times \kk{i+1}{F}^{\tilde{q}(G_2)}@>\tilde{\Phi} _{i,G_2,F}>>\kk{i}{F\otimes G_2}\\
@V\kk{i}{\phi}^{p(G_2)}\times \kk{i+1}{\phi}^{q(G_2)}VV@VV\kk{i}{\phi \otimes id_{G_2}}V \\        
\kk{i}{\bar{F}}^{\tilde{p}(G_2)}\times \kk{i+1}{\bar{F}}^{\tilde{q}(G_2)}@>>\tilde{\Phi} _{i,G_2,\bar{F}}>\kk{i}{\bar{F}\otimes G_2}\\
\end{CD}$$
is commutative. Hence
$$G_2\in \Upsilon ,\qquad p(G_2)=p(G_1)+q(G_1/G_2),\qquad q(G_2)=q(G_1)+p(G_1/G_2),$$
$$\Phi _{i,G_2,F}=\Psi _{F,i}\circ (\Phi _{i,G_1,F}\times \Phi _{(i+1),(G_1/G_2),F})\;.\qedd$$

\begin{p}\label{27.8'}
Let
$$\oc{G}{\varphi }{H}{\psi }{\bc}$$
be an exact sequence in \frcc with $G$ nuclear and $H$ $\Upsilon $-null and let $\delta^F _i$ denote the index maps associated to the exact sequence in \frm
$$\of{F\otimes G}{id_F\otimes \varphi }{F\otimes H}{id_F\otimes \psi }{F}{17}{17}\;.$$
Then
$$G\in \Upsilon ,\qquad p(G)=0,\qquad q(G)=1,\qquad \Phi _{i,G,F=\delta^F _{i+1}}\,,$$
$$(\of{F\otimes G}{id_F\otimes \varphi }{F\otimes H}{id_F\otimes \psi }{F}{16}{16})\in G_\Upsilon \;.$$
\end{p}

By \pr{6.7'} b) and \pr{5.7'a} b),
$$G\in \Upsilon ,\qquad p(G)=0,\qquad q(G)=1,\qquad \Phi _{i,G,F=\delta^F _{i+1}}\;.$$
Since the diagram
$$\begin{CD}
K_{i+1}(F)@>\delta^F _{i+1}>>K_i(F\otimes G)\\
@V\Phi _{i,G,F}=\delta^F _{i+1}VV    @VV\Phi _{(i+1),G,(F\otimes G)}=\delta^F _{G,i}V\\
K_{i}(F\otimes G)@>>\delta^F _{G,i}>K_{i+1}((F\otimes G)\otimes G)
\end{CD}$$
is obviously commutative,
$$(\of{F\otimes G}{id_F\otimes \varphi }{F\otimes H}{id_F\otimes \psi }{F}{16}{16})\in G_\Upsilon \;.\qedd$$

{\center{\section{The class $\Upsilon_1 $  }}}

\begin {center}
\fbox{\parbox{8.8cm}{Throughout this section $F$ denotes an $E$-C*-algebra}}
\end{center}

\begin{de}\label{10.10'}
We denote by $\Upsilon _1$ the class of unital C*-algebras $G$ belonging to $\Upsilon $ such that
$$p(G)=1\,,\qquad q(G)=0\,,\qquad \Phi _{i,G,F}=\kk{i}{\phi _{G,F}}\,,$$
where 
$$\mae{\phi _{G,F}}{F}{F\otimes G}{x}{x\otimes 1_G}\;.$$
\end{de}

\begin{p}\label{10.10'd}
$\bc\in \Upsilon _1$.
\end{p}

In fact
$$\mae{\phi _{\bc,F} }{F}{F\otimes \bc}{x}{x\otimes 1_{\bc}}$$
is an isomorphism.\qed

\begin{p}\label{10.10'a}
Let $G\in \Upsilon _1$ and let $\oaa{F}{\phi }{F'}$ be a morphism in \frm. If we identify $\kk{i}{F}$ with $\kk{i}{F\otimes G}$ for all $E$-C*-algebras $F$ using the group isomorphisms $\Phi _{i,G,F}$ then $\kk{i}{\phi \otimes id_G}$ is identified with $\kk{i}{\phi }$.
\end{p}
The assertion follows from the commutativity of the diagram
$$\begin{CD}
K_{i}(F)@>\kk{i}{\phi }>>K_i(F')\\
@V\Phi _{i,G,F}VV    @VV\Phi _{i,G,F'}V\\
K_{i}(F\otimes G)@>>\kk{i}{\phi \otimes id_G}>K_{i}(F'\otimes G)
\end{CD}\qedd$$

\begin{p}\label{10.10'c}
Let $G,H$ be C*-algebras and $\mac{\varphi }{G}{H}$ and $\mac{\psi }{H}{G}$ be a homotopy such that $\varphi $ and $\psi $ are unital. If $G\in \Upsilon _1$ then $H\in \Upsilon _1$.
\end{p}

By \pr{5.7'a} c),
$$H\in \Upsilon \,,\qquad p(H)=1\,,\qquad q(H)=0\,,$$
$$\Phi _{i,H,F}=\kk{i}{id_F\otimes \varphi }\circ \Phi _{i,G,F}=\kk{i}{id_F\otimes \varphi }\circ \kk{i}{\phi _{G,F}}=\kk{i}{\phi _{H,F}}\;.\qedd$$

\begin{p}\label{10.10'b}
If $(G_j)_{j\in J}$ is a finite family in $\Upsilon _1$, $J\not=\emptyset $, then 
$\bigotimes\limits_{j\in J}G_j\in \Upsilon _1 $.
\end{p}

$\bigotimes\limits_{j\in J}G_j$ is unital and by \pr{27.3'c} a), $\bigotimes\limits_{j\in J}G_j\in \Upsilon $. Assume $J=\z{1,2}$ and let $\oaa{F}{\phi }{F'}$ be a morphism in \frm. Then the diagram
$$\begin{CD}
\kk{i}{F}@>\kk{i}{\phi _{G_1,F}}>>\kk{i}{F\otimes G_1}@>\kk{i}{\phi _{G_2,(F\otimes G_1)}}>>\kk{i}{F\otimes G_1\otimes G_2}\\
@V\kk{i}{\phi }VV @VV\kk{i}{\phi \otimes id_{G_1}}V @VV\kk{i}{\phi \otimes id_{(G_1\otimes G_2)}}V\\
\kk{i}{F'}@>>\kk{i}{\phi _{G_1,F'}}>\kk{i}{F'\otimes G_1}@>>\kk{i}{\phi _{G_2,(F'\otimes G_1)}}>\kk{i}{F'\otimes G_1\otimes G_2}
\end{CD}$$
is commutative. Since 
$$\phi _{(G_1\otimes G_2),F}=\phi _{G_2,(F\otimes G_1)}\circ \phi _{G_1,F}\,,\qquad \phi _{(G_1\otimes G_2),F'}=\phi _{G_2,(F'\otimes G_1)}\circ \phi _{G_1,F'}\,,$$
the diagram
$$\begin{CD}
\kk{i}{F}@>\kk{i}{\phi _{(G_1\otimes G_2),F}}>>\kk{i}{F\otimes G_1\otimes G_2}\\
@V\kk{i}{\phi }VV @VV\kk{i}{\phi \otimes id_{(G_1\otimes G_2)}}V \\
\kk{i}{F'}@>>\kk{i}{\phi _{(G_1\otimes G_2),F'}}>\kk{i}{F'\otimes G_1\otimes G_2}
\end{CD}$$
is commutative, which proves the assertion in this case. The general case follows now by induction with respect to $Card\,J$.\qed

\begin{p}\label{11.10'a}
If $G\in \Upsilon _1$ is nuclear then every exact sequence in \frm belongs to $G_\Upsilon $.
\end{p}

Let
$$\oc{F_1}{\phi _1}{F_2}{\phi _2}{F_3}$$
be an exact sequence in \frm. Then the diagram
$$\begin{CD}
0@>>>F_1@>\phi_1 >>F_2@>\phi_2 >>F_3@>>>0\\
@.      @V\phi _{G,F_1} VV  @V\phi _{G,F_2}VV @VV\phi _{G,F_3}V@.\\
0@>>>F_1\otimes G@>>\phi_1\otimes id_G >F_2\otimes G@>>\phi_2\otimes id_G >F_3\otimes G@>>>0
\end{CD}$$
is commutative and has exact rows. By the commutativity of the index maps (\axi{27.9'e}) the diagram
$$\begin{CD}
K_i(F_3)@>\delta _i>>K_{i+1}(F_1)\\
@V\Phi _{i,G,F_3}=K_i(\phi _{G,F_3})VV    @VV\Phi _{(i+1),G,F_1}=K_{i+1}(\phi _{G,F_1})V\\
K_i(F_3\otimes G)@>>\delta _{G,i}>K_{i+1}(F_1\otimes G)
\end{CD}$$
is commutative, where $\delta _i$ denotes the index maps of the exact sequence
$$\oc{F_1}{\phi _1}{F_2}{\phi _2}{F_3}\;.\qedd$$

\begin{p}\label{12.10'}
Let $G$ be a C*-algebra. 
\begin{enumerate}
\item $\phi _{\tilde{G},F }=(id_F\otimes \lambda _G)\circ \phi _{\bc,F}$.
\item $G$ is $\Upsilon $-null iff $\tilde{G}\in \Upsilon _1 $.
\item If $G$ is $\Upsilon $-null and $\mac{\varphi }{G}{G'}$, $\mac{\psi }{G}{G'}$ are C*-homomorphisms then $\kk{i}{id_F\otimes \tilde{\varphi } }=\kk{i}{id_F\otimes \tilde{\psi  } }$. In particular if $G=G'$ then 
$$\kk{i}{id_F\otimes \tilde{\varphi } }=id_{\kk{i}{F\otimes \tilde{G} }}\approx id_{\kk{i}{F}}\;.$$
\end{enumerate}
\end{p}
 
a) is easy to see.

b) By \cor{26.3'b} b), the sequence
$$\og{K_i(F\otimes G)}{K_i\left(id_F\otimes \iota_G\right)}{K_i\left(F\otimes \tilde{G}\right) }{K_i\left(id_F\otimes \pi_G\right)}{K_i\left(id_F\otimes \lambda_G\right)}{K_i(F)}{30}{35}{35}$$
is split exact. By a) and \pr{10.10'd},
$$\kk{i}{\phi _{\tilde{G},F }}=\kk{i}{id_F\otimes \lambda _G}\circ \Phi _{i,\bc,F}\;.$$
If $\tilde{G}\in \Upsilon _1 $ then
$$\Phi _{i,\tilde{G},F }=\kk{i}{\phi _{\tilde{G},F }}=\kk{i}{id_F\otimes \lambda _G}\circ \Phi _{i,\bc,F}\,,$$
so by \pr{10.10'd}, $\kk{i}{id_F\otimes \lambda _G}$ is an isomorphism, $\kk{i}{id_F\otimes \iota _G}=0$, $\kk{i}{F\otimes G}=0$, and $G$ is $\Upsilon $-null. If $G$ is $\Upsilon $-null then $\kk{i}{id_F\otimes \lambda _G}$ is an isomorphism so 
$$\mac{\kk{i}{\phi _{\tilde{G},F }}}{\kk{i}{F}}{\kk{i}{F\otimes G}}$$
is an isomorphism and $\tilde{G}\in \Upsilon _1 $.

c) Since $\tilde{\varphi }\circ \lambda _G=\tilde{\psi }\circ \lambda _G  $,
$$\kk{i}{id_F\otimes \tilde{\varphi } }\circ \kk{i}{id_F\otimes \lambda _G}=\kk{i}{id_F\otimes \tilde{\psi } }\circ \kk{i}{id_F\otimes \lambda _G}\;.$$
By b), $\tilde{G}\in \Upsilon _1 $ and so $\kk{i}{id_F\otimes \lambda _G}$ is an isomorphism. Thus $\kk{i}{id_F\otimes \tilde{\varphi } }=\kk{i}{id_F\otimes \tilde{\psi  } }$.\qed

\begin{co}\label{12.10'a}
If $(G_j)_{j\in J}$ is a finite family of $\Upsilon $-null C*-algebras and $G:=\pro{j\in J}G_j$ then $\tilde{G}\in \Upsilon _1 $.
\end{co}

By \pr{5.7'b} a), $G$ is $\Upsilon $-null and by \pr{12.10'} b), $\tilde{G}\in \Upsilon _1 $.\qed

\begin{p}\label{19.10'}
Let
$$\oc{G_1}{\varphi }{G_2}{\psi }{G_3}$$
be an exact sequence in \frcc such that $G_1$ is $\Upsilon $-null, $G_3$ is nuclear, and $G_2,\,G_3$ are unital. Then $G_2\in \Upsilon _1$ iff $G_3\in \Upsilon _1$.
\end{p}

Since $G_2$ and $G_3$ are unital and $\psi $ is surjective, $\psi (1_{G_2})=1_{G_3}$. It follows
$$\phi _{G_3,F}=(id_F\otimes \psi )\circ \phi _{G_2,F}\,,\qquad\qquad \kk{i}{\phi _{G_3,F}}=\kk{i}{id_F\otimes \psi }\circ \kk{i}{\phi _{G_2,F}}\;.$$
By \pr{6.7'} a), $\kk{i}{id_F\otimes \psi }$ is a group isomorphism,
$$G_2,G_3\in \Upsilon \,,\qquad p(G_2)=p(G_3)=1\,,\qquad q(G_2)=q(G_3)=0\,,$$
$$\Phi _{i,G_3,F}=\kk{i}{id_F\otimes \psi }\circ \Phi _{i,G_2,F}\;.$$
If $G_2\in \Upsilon _1$ then by the above,
$$\Phi _{i,G_3,F}=\kk{i}{id_F\otimes \psi }\circ \kk{i}{\phi _{G_2,F}}=\kk{i}{\phi _{G_3,F}}\,,$$
so $G_3\in \Upsilon _1$. If $G_3\in \Upsilon _1$ then by the above, 
$$\kk{i}{id_F\otimes \psi }\circ \kk{i}{\phi _{G_2,F}}=\kk{i}{\phi _{G_3,F}}=\Phi _{i,G_3,F}=\kk{i}{id_F\otimes \psi }\circ \Phi _{i,G_2,F}\,,$$
so $\Phi _{i,G_2,F}=\kk{i}{\phi _{G_2,F}}$ and $G_2\in \Upsilon _1$.\qed 

{\center{\chapter{Locally compact spaces}}}

{\center{\section{Tietze's Theorem}}}

\begin{de}\label{28.3'a}
Let $\Omega $ be a topological space and $F$ an $E$-C*-algebra. We endow canonically the C*-algebra $\ccb{\Omega }{F}$ with the structure of an $E$-C*-algebra by putting
$$\mae{\alpha x}{\Omega }{F}{\omega }{\alpha x(\omega )}$$
for all $(\alpha ,x)\in E\times F$. If $\Omega $ is a locally compact space then we endow $\cbb{\Omega }{F}$ with the structure of on $\eo$ in a similar way. If $\Omega '$ is an open set of a locally compact space $\Omega $ then we identify $\cbb{\Omega '}{F}$ with the $E$-ideal $\me{x\in \cbb{\Omega }{F}}{x|(\Omega \setminus \Omega ')=0}$ of $\cbb{\Omega }{F}$.
\end{de}

\begin{de}\label{2.9'}
Let $\Omega $ be a locally compact space with $\cbb{\Omega }{\bc}\in \Upsilon $. We put
$$\Omega \in \Upsilon\, ,\qquad p(\Omega) :=p(\cbb{\Omega }{\bc})\,,\qquad q(\Omega ):=q(\cbb{\Omega }{\bc})\,,$$
$$ \Phi _{i,\Omega ,F}:=\Phi _{i,\cbb{\Omega }{\bc},F}\,,\quad \Omega _\Upsilon :=\cbb{\Omega }{\bc}_\Upsilon \,,\quad \Omega \in \Upsilon _1:\Longleftrightarrow \cbb{\Omega }{\bc}\in \Upsilon _1\;.$$
We say that $\Omega $ is {\bf{$\Upsilon$-null}} if $\cbb{\Omega }{\bc}$ is $\Upsilon $-null. We say that $\Omega $ is {\bf{null-homotopic}} if $\cbb{\Omega }{\bc}$ is null-homotopic.
\end{de}

\begin{p}\label{12.10'c}
If $\Omega $ is a locally compact space and if $\Omega ^*$ denotes its Alexandroff compactification then $\Omega $ is $\Upsilon $-null iff $\Omega ^*\in \Upsilon _1$.
\end{p}

The Proposition is a particular case of \pr{12.10'}.\qed

\begin{lem}\label{14.9'}
Let $\Omega $ be a locally compact space.
\begin{enumerate}
\item $\cbb{\Omega }{\bc}$ is nuclear.
\item $\cbb{\Omega }{F}\approx F\otimes \cbb{\Omega }{\bc}$.
\item If $\Omega $ is a finite compact space then $\Omega \in \Upsilon $, $p(\Omega) =Card\,\Omega $, $q(\Omega) =0$, and every exact sequence in \frm belongs to $\Omega _\Upsilon $.
\end{enumerate} 
\end{lem}

a) [W] Theorem T.6.20.

b) [W] Proposition T.5.11,

c) follows from \pr{5.7'b} c).\qed 

\begin{co}[Tietze's Theorem]\label{23.11a}
Let $\Omega $ be a locally compact space, $\Gamma $ a closed set of $\Omega $, $\mac{\varphi }{\cbb{\Omega \setminus \Gamma }{F}}{\cbb{\Omega }{F}}$ the inclusion map, and
$$\mae{\psi }{\cbb{\Omega }{F}}{\cbb{\Gamma }{F}}{x}{x|\Gamma }\;.$$
Then
$$\oc{\cbb{\Omega \setminus \Gamma }{F}}{\varphi }{\cbb{\Omega }{F}}{\psi }{\cbb{\Gamma }{F}}$$
is an exact sequence in \frm.
\end{co}

By \lm{14.9'} a),b), the assertion follows from \pr{29.3'a}.\qed

\begin{co}\label{15.9'}
If 
$$\oc{F_1}{\phi_1 }{F_2}{\phi_2 }{F_3}$$
is an exact sequence in $\fr{M}_E$ and $\Omega $ a locally compact space then 
$$\of{\cbb{\Omega }{F_1}}{\phi_1 \otimes id_G}{\cbb{\Omega }{F_2}}{\phi_2 \otimes id_G}{\cbb{\Omega }{F_3}}{16}{16}$$
is an exact sequence in \frm.
\end{co}

By \lm{14.9'} a),b), the assertion follows from  \pr{23.4'b}.\qed

\begin{p}\label{14.6'a}
Let
$$\oc{F_1}{\phi_1 }{F_2}{\phi_2 }{F_3}$$
be an exact sequence in \frm, $\Omega $ a locally compact space, $\Gamma $ a closed set of $\Omega $, $\mac{\varphi }{\cbb{\Omega \setminus \Gamma }{\bc}}{\cbb{\Omega }{\bc}}$ the inclusion map, and
$$\mae{\psi }{\cbb{\Omega }{\bc}}{\cbb{\Gamma }{\bc}}{x}{x|\Gamma }\;.$$
\begin{enumerate}
\item $G:=\me{x\in \cbb{\Omega }{F_2}}{x|\Gamma \in \cbb{\Gamma }{F_1}}$ is a closed $E$-ideal of $\cbb{\Omega }{F_2}$; we denote by 
$\mac{\varphi '}{G}{\cbb{\Omega }{F_2}}$
 the inclusion map.
\item The sequence in \frm
$$\of{G}{\varphi '}{\cbb{\Omega }{F_2}}{\phi_2 \otimes \psi }{\cbb{\Gamma }{F_3}}{0}{16}$$
is exact.
\end{enumerate}
\end{p}

a) is easy to see.

b) We put 
$$G_1:=\cbb{\Omega \setminus \Gamma }{\bc}\,,\qquad G_2:=\cbb{\Omega }{\bc}\,,\qquad G_3:=\cbb{\Gamma }{\bc}\;.$$
Let us consider the following commutative diagram.
$$\begin{CD}
\rule{0mm}{0mm}@.0@.0@.0@.@.\\
@.@VVV@VVV@VVV@.\\
0@>>>F_1\otimes G_1@>\phi_1 \otimes id_{G_1}>>F_2\otimes G_1@>\phi_2 \otimes id_{G_1}>>F_3\otimes G_1@>>>0\\
@.@Vid_{F_1}\otimes \varphi VV@Vid_{F_2}\otimes \varphi VV         @Vid_{F_3}\otimes \varphi VV  @.\\
0@>>>F_1\otimes G_2@>\phi_1 \otimes id_{G_2}>>F_2\otimes G_2@>\phi_2 \otimes id_{G_2}>>F_3\otimes G_2@>>>0\\
@.@Vid_{F_1}\otimes \psi VV@Vid_{F_2}\otimes \psi VV         @Vid_{F_3}\otimes \psi VV @. \\
0@>>>F_1\otimes G_3@>\phi_1 \otimes id_{G_3}>>F_2\otimes G_3@>\phi_2 \otimes id_{G_3}>>F_3\otimes G_3@>>>0\\
@.@VVV@VVV@VVV@.\\
\rule{0mm}{0mm}@.0@.0@.0@.@.\\
\end{CD}$$
By \lm{14.9'} a), \pr{29.3'a}, and \pr{23.4'b}, its columns and rows are exact. It follows that $\phi_2 \otimes \psi $ is surjective. Let $x\in Ker\,(\phi_2 \otimes \psi )$. Then 
$$(id_{F_3}\otimes \psi )(\phi_2 \otimes id_{G_2})x=(\phi_2 \otimes \psi )x=0\,,$$
so there is a $y\in F_2\otimes G_1$ with
$$(\phi_2 \otimes \varphi )y=(id_{F_3}\otimes \varphi )(\phi_2 \otimes id_{G_1})y=(\phi_2 \otimes id_{G_2})x\;.$$
Then
$$(\phi_2 \otimes id_{G_2})(x-(id_{F_2}\otimes \varphi )y)=(\phi_2 \otimes id_{G_2})x-(\phi_2 \otimes \varphi )y=0\,,$$
so there is a $z\in F_1\otimes G_2$ with
$$(\phi_1 \otimes id_{G_2})z=x-(id_{F_2}\otimes \varphi )y\;.$$
Thus
$$x=(id_{F_2}\otimes \varphi )y+(\phi_1 \otimes id_{G_2})z\in G\,,\qquad Ker\,(\phi_2 \otimes \psi )\subset G\;. $$

Let now $x\in G$. By \pr{29.3'a}, there is a $y\in \cbb{\Omega }{F_1}=F_1\otimes G_2$ with $x|\Gamma =y|\Gamma $. There is a $z\in \cbb{\Omega \setminus \Gamma }{F_2}=F_2\otimes G_1$ with 
$$(id_{F_2}\otimes \varphi )z=x-(\phi_1 \otimes id_{G_2})y\;.$$
 We get
$$(\phi_2 \otimes \psi )x=(\phi_2 \otimes \psi )(\phi_1 \otimes id_{G_2})y+(\phi_2 \otimes \psi )(id_{F_2}\otimes \varphi )z=$$
$$=((\phi_2 \circ \phi_1 )\otimes \psi )y+(\phi_2 \otimes (\psi \circ \varphi ))z=0\,,$$
$G\subset Ker\,(\phi_2 \otimes \psi )$.\qed

{\it Remark.} If we put $F_1:=0$ and  $F_2=F_3$ in the above Proposition then we obtain Tietze's Theorem (\cor{23.11a}).

\begin{p}[Topological six-term sequence]\label{24.11}
Let $\Omega $ be a locally compact space, $\Gamma $ a closed set of $\Omega $, $\mac{\varphi }{\cbb{\Omega \setminus \Gamma }{F}}{\cbb{\Omega }{F}}$ the inclusion map, 
$$\mae{\psi }{\cbb{\Omega }{F}}{\cbb{\Gamma }{F}}{x}{x|\Gamma }\,,$$
and $\delta _i$ the 
 index maps associated to the exact sequence in \frm \emph{(Tietze's Theorem (\cor{23.11a}))}
$$\oc{\cbb{\Omega \setminus \Gamma }{F}}{\varphi }{\cbb{\Omega }{F}}{\psi }{\cbb{\Gamma }{F}}\;.$$
\begin{enumerate}
\item Assume $\Omega \setminus \Gamma $ is $\Upsilon $-null.
\begin{enumerate}
\item $\mac{K_i(\psi )}{K_i(\cbb{\Omega }{F})}{K_i(\cbb{\Gamma }{F})}$
is a group isomorphism.
\item If $\Omega \in \Upsilon $ or $\Gamma \in \Upsilon $ then
$$\Omega ,\Gamma \in \Upsilon\, ,\qquad 
p(\Omega) =p(\Gamma) \,,\qquad q(\Omega )=q(\Gamma )\,,$$
$$ \Phi _{i,\Gamma ,F}=\kk{i}{id_F\otimes \psi }\circ \Phi _{i,\Omega ,F}\,,\qquad 
 \Omega _\Upsilon =\Gamma _\Upsilon\;. $$
\end{enumerate}
\item Assume $\Omega $ is $\Upsilon $-null.
\begin{enumerate}
\item $\mac{\delta _i}{K_i(\cbb{\Gamma }{F})}{K_{i+1}(\cbb{\Omega \setminus \Gamma }{F})}$
is a group isomorphism.
\item If $\Omega \setminus \Gamma \in \Upsilon $ or $\Gamma \in \Upsilon $ then
$$\Omega \setminus \Gamma ,\Gamma \in \Upsilon\,,\qquad 
p(\Omega \setminus \Gamma )=q(\Gamma) \,,\qquad q(\Omega \setminus \Gamma) =p(\Gamma) \,,$$
$$\Phi _{i,\Gamma ,F}=\Phi _{(i+1),(\Omega \setminus \Gamma ),F}\circ \delta _i\;.$$
\end{enumerate}
\item Assume $\Gamma $ is $\Upsilon $-null.
\begin{enumerate}
\item $\mac{K_i(\varphi )}{K_i(\cbb{\Omega\setminus \Gamma }{F})}{K_i(\cbb{\Omega }{F})}$
is a group isomorphism.
\item If $\Omega \setminus \Gamma \in \Upsilon $ or $\Omega \in \Upsilon $ then
$$\Omega \setminus \Gamma ,\Omega \in \Upsilon\,,\qquad 
p(\Omega \setminus \Gamma) =p(\Omega) \,,\qquad q(\Omega \setminus \Gamma) =q(\Omega) \,,$$
$$ \Phi _{i,\Omega ,F}=\kk{i}{id_F\otimes \varphi }\circ \Phi _{i,(\Omega \setminus \Gamma ),F}\,,\qquad (\Omega \setminus \Gamma )_\Upsilon =\Omega _\Upsilon\;.$$
\end{enumerate}
\end{enumerate}
\end{p}

The assertions follow from \lm{14.9'} a),b) and \pr{6.7'}.\qed

\begin{co}\label{5.9'}
Let $\Omega $ be a locally compact space, $\omega \in \Omega $ such that $\Omega \setminus \z{\omega} $ is $\Upsilon $-null, $\Gamma $ a closed set of $\Omega $,
$$\Omega ':=(\Omega \setminus \z{\omega })\setminus \Gamma \,,\qquad \Gamma ':=\Gamma \setminus \z{\omega }\,,$$
$\mac{\varphi }{\cbb{\Omega '}{F}}{\cbb{\Omega \setminus \z{\omega }}{F}}$ the inclusion map,
$$\mae{\psi }{\cbb{\Omega \setminus \z{\omega }}{F}}{\cbb{\Gamma '}{F}}{x}{x|\Gamma '}\,,$$
and $\delta _i$ the index maps of the exact sequence in \frm \emph{(Tietze's Theorem (\cor{23.11a}))}
$$\oc{\cbb{\Omega '}{F}}{\varphi }{\cbb{\Omega \setminus \z{\omega }}{F}}{\psi }{\cbb{\Gamma '}{F}}\;.$$ 
\begin{enumerate}
\item $\mac{\delta _i}{\kk{i}{\cbb{\Gamma '}{F}}}{\kk{i+1}{\cbb{\Omega '}{F}}}$ is a group isomorphism.
\item If $\Omega '\in \Upsilon $ or $\Gamma '\in \Upsilon $ then
$$\Omega ',\Gamma '\in \Upsilon \,,\qquad p(\Omega ')=q(\Gamma ')\,,\qquad q(\Omega ')=p(\Gamma ')\,,$$
$$\Phi _{i,\Gamma ',F}=\Phi _{(i+1),\Omega ',F}\circ \delta _i\;. $$
\item If $\Gamma $ is finite then
$$\Omega '\in \Upsilon \,,\qquad p(\Omega ')=0\,,\qquad q(\Omega ')=Card\,\Gamma '\;.$$
\end{enumerate}
\end{co}

a) and b) follow from the Topological six-term sequence (\pr{24.11} b)).

c) follows from b) and \lm{14.9'} c).\qed

\begin{co}\label{28.2'}
Let $\Omega ,\Omega '$ be locally compact spaces, $\omega \in \Omega $, and $\omega '\in \Omega'$ such that $\Omega '\setminus \z{\omega '}$ is null-homotopic. 
\begin{enumerate}
\item $K_i\left(\cbb{(\Omega \times \Omega ')\setminus \z{(\omega ,\omega ')}}{F}\right)\approx K_i\left(\cbb{(\Omega \setminus \z{\omega })\times \Omega '}{F}\right)\;.$
\item If also $\Omega \setminus \z{\omega }$ is null-homotopic then $\cbb{(\Omega \times \Omega ')\setminus \z{(\omega ,\omega ')}}{F}$ is K-null.
\end{enumerate}
\end{co}

a) The sequence in \frm (with obvious notation)
$$\occ{\cbb{(\Omega\setminus \z{\omega }) \times \Omega '}{F}}{\varphi }{\cbb{(\Omega \times \Omega ')\setminus \z{(\omega ,\omega ')}}{F}}$$
$$\ocd{\cbb{(\Omega \times \Omega ')\setminus \z{(\omega ,\omega ')}}{F}}{\psi }{\cbb{\z{\omega }\times (\Omega '\setminus \z{\omega '})}{F}}$$
is exact and the assertion follows from the Topological six-term sequence (\pr{24.11} $c_1)$).

b) By \pr{29.3'} c) and \lm{14.9'} b), $(\Omega \setminus \z{\omega })\times \Omega '$ is null-homotopic and so K-null (\pr{5.7'a} a)). By a), 
$$K_i\left(\cbb{(\Omega \times \Omega ')\setminus \z{(\omega ,\omega ')}}{F}\right)$$
 is K-null.\qed

\begin{p}[Topological triple]\label{3.12}
Let $\Omega _1$ be a locally compact space, $\Omega _2$ an open set of $\Omega _1$, $\Omega _3$ an open set of $\Omega _2$, and $\mac{\varphi }{\cbb{\Omega _2\setminus \Omega _3}{F}}{\cbb{\Omega _1\setminus \Omega _3}{F}}$ the inclusion map. For all $j,k\in \z{1,2,3}$, $j<k$, put
$$\mae{\psi _{j,k}}{\cbb{\Omega _j}{F}}{\cbb{\Omega _j\setminus \Omega _k}{F}}{x}{x|(\Omega _j\setminus \Omega _k)}$$
and denote by $\mac{\varphi _{j,k}}{\cbb{\Omega _k}{F}}{\cbb{\Omega _j}{F}}$ the inclusion map and by $\delta _{j,k,i}$ the index maps associated to the exact sequence in \frm
$$\oc{\cbb{\Omega _k}{F}}{\varphi _{j,k}}{\cbb{\Omega _j}{F}}{\psi _{j,k}}{\cbb{\Omega _j\setminus \Omega _k}{F}}\;.$$
\begin{enumerate}
\item Assume $\cbb{\Omega _2}{F}$ K-null.
\begin{enumerate}
\item $\mac{\delta _{2,3,i}}{K_i(\cbb{\Omega _2\setminus \Omega _3}{F})}{K_{i+1}(\cbb{\Omega _3}{F})}$ is a group isomorphism.
\item $\delta _{2,3,i}=\delta _{1,3,i}\circ K_i(\varphi )$.
\item $\varphi_{1,3} $ is K-null.
\item If we put $\Phi _i:=K_i(\varphi )\circ (\delta _{2,3,i})^{-1}$ then
$$\oddg{K_i(\cbb{\Omega _1}{F})}{K_i(\psi _{1,3})}{K_i(\cbb{\Omega _1\setminus \Omega _3}{F})}{\delta _{1,3,i}}{\Phi _i}{25}{10}{10}$$
$$\odeg{\delta _{1,3,i}}{\Phi _i}{K_{i+1}(\cbb{\Omega _3}{F})}{10}{10}$$
is a split exact sequence and the map
$$K_i(\cbb{\Omega _1}{F})\times K_{i+1}(\cbb{\Omega _3}{F})\longrightarrow K_i(\cbb{\Omega _1\setminus \Omega _3}{F})\,,$$
$$(a,b)\longmapsto K_i(\psi _{1,3})a+\Phi _ib$$
is a group isomorphism.
\item If $\Omega _2$ is $\Upsilon $-null and $\Omega _1,\Omega _3\in \Upsilon $ then
$$\Omega _1\setminus \Omega _3\in \Upsilon \,,\; p(\Omega _1\setminus \Omega _3)=p(\Omega _1)+q(\Omega _3)\,,\; q(\Omega _1\setminus \Omega _3)=q(\Omega _1)+p(\Omega _3)\,,$$
and (with the notation of \emph{\pr{8.7'}})
$$\Phi _{i,(\Omega _1\setminus \Omega _3),F}=\Psi _{F,i}\circ (\Phi _{i,\Omega _1,F}\times \Phi _{(i+1),\Omega _3,F})\;.$$
\end{enumerate}
\item Assume $\cbb{\Omega _1\setminus \Omega _3}{F}$ K-null.
\begin{enumerate}
\item $\delta _{2,3,i}=0$.
\item $\mac{K_i(\varphi _{1,3})}{K_i(\cbb{\Omega _3}{F})}{K_i(\cbb{\Omega _1}{F})}$ is a group isomorphism.
\item If we put $\Phi _i:=K_i(\varphi _{1,3})^{-1}\circ K_i(\varphi _{1,2})$ then the map
$$\mac{\Psi }{K_{i}(\cbb{\Omega _2}{F})}{K_i(\cbb{\Omega _3}{F})\times K_{i}(\cbb{\Omega _2\setminus \Omega _3}{F})}\,,$$
$$b\longmapsto (\Phi _ib,K_i(\psi _{2,3})b)$$
is a group isomorphism.
\item If $\psi _{1,2}$ is K-null and if we put $\Phi '_i:=K_i(\varphi _{2,3})\circ K_i(\varphi _{1,3})^{-1}$ \emph{(by $c_2)$)} then
$$\oddg{K_{i+1}(\cbb{\Omega _1\setminus \Omega _2}{F})}{\delta _{1,2,(i+1)}}{K_i(\cbb{\Omega _2}{F})}{K_i(\varphi _{1,2})}{\Phi '_i}{20}{20}{20}$$
$$\odeg{K_i(\varphi _{1,2})}{\Phi '_i}{K_i(\cbb{\Omega _1}{F})}{20}{20}$$
is a split exact sequence and the map
$$K_i(\cbb{\Omega _1}{F})\times K_{i+1}(\cbb{\Omega _1\setminus 
\Omega _2}{F})\longrightarrow K_i(\cbb{\Omega _2}{F}),$$
$$(a,b)\longmapsto \Phi '_ia+\delta _{1,2,(i+1)}b$$
is a group isomorphism. 
\item If $\Omega _1\setminus \Omega _3$ is $\Upsilon $-null, $\Omega _1,\Omega _1\setminus \Omega _2\in \Upsilon $, and $\psi _{1,2}$ is K-null then
$$\Omega _2\in \Upsilon \,,\; p(\Omega _2)=p(\Omega _1)+q(\Omega _1\setminus \Omega _2)\,,\; q(\Omega _2)=q(\Omega _1)+p(\Omega _1\setminus \Omega _2)\;.$$
\end{enumerate}
\item Assume $\cbb{\Omega _1}{F}$ K-null and put
$$\mae{\psi }{\cbb{\Omega _1\setminus \Omega _3}{F}}{\cbb{\Omega _1\setminus \Omega _2}{F}}{x}{x|(\Omega _1\setminus \Omega _2)}\;.$$
\begin{enumerate}
\item $\delta _{1,2,i}$ and $\delta _{1,3,i}$ are group isomorphisms.
\item $K_i(\varphi _{2,3})\circ \delta _{1,3,(i+1)}=\delta _{1,2,(i+1)}\circ K_{i+1}(\psi )$.
\item Let $\mac{\varphi '}{\cbb{\Omega _1\setminus \Omega _2}{F}}{\cbb{\Omega _1\setminus \Omega _3}{F}}$ be a morphism in \frm such that 
$$K_i(\psi \circ \varphi ')=id_{K_i(\cbb{\Omega _1\setminus \Omega _2}{F})}\;.$$ 
If we put 
$$\Phi _i:=\delta _{1,3,(i+1)}\circ K_{i+1}(\varphi ')\circ (\delta _{1,2,(i+1)})^{-1}$$
then $K_i(\varphi _{2,3})\circ \Phi _i=id_{K_i(\cbb{\Omega _2}{F})}$. If in addition $\psi _{2,3}$ is K-null then
$$\oddg{K_{i+1}(\cbb{\Omega _2\setminus \Omega _3}{F})}{\delta _{2,3,(i+1)}}{K_i(\cbb{\Omega _3}{F})}{K_i(\varphi _{2,3})}{\Phi _i}{20}{20}{20}$$
$$\odeg{K_i(\varphi _{2,3})}{\Phi _i}{K_i(\cbb{\Omega _2}{F})}{20}{20}$$
is a split exact sequence and the map
$$K_{i+1}((\cbb{\Omega _2\setminus \Omega _3}{F})\times K_i(\cbb{\Omega _2}{F})\longrightarrow K_i(\cbb{\Omega _3}{F})\,,$$
$$(a,b)\longmapsto \delta _{2,3,(i+1)}a+\Phi _ib$$
is a group isomorphism.
\end{enumerate}
\end{enumerate}
\end{p}

Up to $a_5)$ and $b_5)$ the Proposition follows from Tietze's Theorem (\cor{23.11a}) and from the triple theorem (\h{28.9'f}) (and \lm{14.9'} a),b)). $a_5)$ follows from \pr{8.7'} and $b_5)$ follows from \pr{9.7'}.\qed

\begin{co}\label{1.1'}
Let $\oaa{F}{\phi }{F'}$ be a morphism in \frm. We use the notation and hypotheses of \emph{\pr{3.12}} and the hypothesis that $\cbb{\Omega _2}{F}$ and $\cbb{\Omega _2}{F'}$ are K-null, and mark with an accent those notation associated to $F'$. We put for all $j\in \z{1,2,3}$ and for all $j,k\in \z{1,2,3}$, $j<k$,
$$\mae{\phi _j}{\cbb{\Omega _j}{F}}{\cbb{\Omega _j}{F'}}{x}{\phi\circ x}\,,$$
$$\mae{\phi _{j,k}}{\cbb{\Omega _j\setminus \Omega _k}{F}}{\cbb{\Omega _j\setminus \Omega _k}{F'}}{x}{\phi\circ x}\;.$$
\begin{enumerate}
\item $\Phi '_i\circ K_{i+1}(\phi _3)=K_i(\phi _{1,3})\circ \Phi _i$.
\item If we identify $K_i(\cbb{\Omega _1\setminus \Omega _3}{F})$ with $K_i(\cbb{\Omega _1}{F})\times K_{i+1}(\cbb{\Omega _3}{F})$ and $K_i(\cbb{\Omega _1\setminus \Omega _3}{F'})$ with $K_i(\cbb{\Omega _1}{F'})\times K_{i+1}(\cbb{\Omega _3}{F'})$ using the isomorphisms of \emph{\pr{3.12} $a_4)$} then
$$K_i(\phi _{1,3}):K_i(\cbb{\Omega _1\setminus \Omega _3}{F})\longrightarrow K_i(\cbb{\Omega _1\setminus \Omega _3}{F'}),$$
$$(a,b)\longrightarrow (K_i(\phi _1)a,K_{i+1}(\phi _3)b)$$
is a group isomorphism.
\end{enumerate}
\end{co}

a) The diagram
$$\begin{CD}
0\longrightarrow \cbb{\Omega _3}{F}@>\varphi_{2,3}>>\cbb{\Omega _2}{F}@>\psi _{2,3}>>\cbb{\Omega _2\setminus \Omega _3}{F}\longrightarrow 0\\
@VV\phi _3V         @VV\phi_2V@VV\phi _{2,3} V      \\
0\longrightarrow \cbb{\Omega _3}{F'})@>>\varphi'_{2,3}>\cbb{\Omega _2}{F'}@>>\psi'_{2,3}>\cbb{\Omega _2\setminus \Omega _3}{F'}\longrightarrow 0
\end{CD}$$
is obviously commutative and has exact rows. By the commutativity of the index maps (\axi{27.9'e}), 
$$K_{i+1}(\phi _3)\circ \delta _{2,3,i}=\left(\delta _{2,3,i}\right)'\circ K_i(\phi _{2,3})\,,$$
$$\left(\left(\delta _{2,3,i}\right)'\right)^{-1}\circ K_{i+1}(\phi _3)=K_i(\phi _{2,3})\circ \left(\delta _{2,3,i}\right)^{-1}\;.$$
By the above, since $\phi _{1,3}\circ \varphi =\varphi '\circ \phi _{2,3}$, 
$$K_i(\phi _{1,3})\circ \Phi _i=K_i(\phi _{1,3})\circ K_i(\varphi )\circ \left(\delta _{2,3,i}\right)^{-1}=K_i(\varphi ')\circ K_i(\phi _{2,3})\circ \left(\delta _{2,3,i}\right)^{-1}=$$
$$=K_i(\varphi ')\circ \left(\left(\delta _{2,3,i}\right)'\right)^{-1}\circ K_{i+1}(\phi _3)=\Phi '_i\circ K_{i+1}(\phi _3)\;.$$

b) follows from a) and \pr{3.12} $a_4)$.\qed

{\center{\section{Alexandroff compactification}}}

\begin{theo}[Alexandroff K-theorem]\label{20.4}
Let $\Omega $ be a locally compact space and $\Omega ^*$ its Alexandroff compactification. We denote by 
$$\mac{\varphi }{\cbb{\Omega }{F}}{\ccb{\Omega ^*}{F}}$$
the inclusion map and put
$$\mae{\lambda }{F}{\ccb{\Omega ^*}{F}}{y}{y1_{\ccb{\Omega ^*}{\bc}}}$$
\begin{enumerate}
\item The map
$$\mad{K_i(\cbb{\Omega }{F})\times K_i(F)}{K_i(\ccb{\Omega ^*}{F})}{(a,b)}{K_i(\varphi )a+K_i(\lambda )b}$$
is a group isomorphism.
\item If $\Omega \in \Upsilon $ then
$$\Omega ^*\in \Upsilon \,,\qquad p(\Omega ^*)=p(\Omega) +1\,,\qquad q(\Omega ^*)=q(\Omega) \,\qquad \Omega _\Upsilon \subset \Omega^* _{\Upsilon }\;.$$
\item $\Omega $ is $\Upsilon $-null iff $\Omega ^*\in \Upsilon _1$.
\end{enumerate}
\end{theo}

$\ccb{\Omega ^*}{\bc}$ is the unitization of $\cbb{\Omega }{\bc}$.

a) Since
$$\cbb{\Omega }{F}\approx F\otimes \cbb{\Omega }{\bc}\,,\qquad \cbb{\Omega^* }{F}\approx F\otimes \cbb{\Omega^* }{\bc}$$
(\lm{14.9'} b)), the assertion follows from \cor{26.3'b} b).

b) follows from \cor{6.7'b}.

c) follows from \pr{12.10'} b).\qed

\begin{co}\label{30.12}
Let $\Omega _1$ and $\Omega _2$ be locally compact spaces, $\Omega _1^*,\,\Omega _2^*$ their Alexandroff compactification, respectively, $\mac{\vartheta }{\Omega _1}{\Omega _2}$ a proper continuous map, $\mac{\vartheta ^*}{\Omega _1^*}{\Omega _2^*}$ its continuous extension, and
$$\mae{\phi }{\cbb{\Omega _2}{F}}{\cbb{\Omega _1}{F}}{x}{x\circ \vartheta }\,,$$
$$\mae{\phi^* }{\ccb{\Omega _2^*}{F}}{\ccb{\Omega _1^*}{F}}{x}{x\circ \vartheta^* }\;.$$
\begin{enumerate}
\item If we identify $K_i\left(\ccb{\Omega _j^*}{F}\right)$ with $K_i(\cbb{\Omega _j}{F})\times K_i(F)$ for every $j\in \z{1,2}$ using the group isomorphisms of the Alexandroff K-theorem \emph{(\h{20.4} a))} then
$$\mae{K_i(\phi ^*)}{K_i(\ccb{\Omega _2^*}{F})}{K_i(\ccb{\Omega _1^*}{F})}{(a,b)}{(K_i(\phi )a,b)}\;.$$
\item Let $\mac{\vartheta '}{\Omega _1}{\Omega _2}$ be a proper continuous map and let $\phi ',\phi '^*$ be the above maps associated to $\vartheta '$. If $\Omega _2$ is $\Upsilon $-null then $\kk{i}{id_F\otimes \phi ^*}=\kk{i}{id_F\otimes \phi '^*}$. In particular if $\Omega _1=\Omega _2$ then $$\kk{i}{id_F\otimes \phi ^*}=id_{\kk{i}{\ccb{\Omega _1^*}{F}}}\;.$$
\end{enumerate}
 \end{co}

a) follows from \cor{26.3'b} c).

b) follows from \pr{12.10'} c).\qed

\begin{co}\label{20.2'}
Let $\oaa{F}{\phi }{F'}$ be a morphism in \frm. We use the notation of the Alexandroff K-theorem \emph{(\h{20.4})} and put
$$\mae{\phi _\Omega }{\cbb{\Omega }{F}}{\cbb{\Omega }{F'}}{x}{\phi \circ x}\,,$$
$$\mae{\phi _{\Omega^*} }{\ccb{\Omega^* }{F}}{\ccb{\Omega^* }{F'}}{x}{\phi \circ x}\;.$$
If we identify $K_i(\ccb{\Omega ^*}{F})$ with $K_i(\cbb{\Omega }{F})\times K_i(F)$ and $K_i(\ccb{\Omega ^*}{F'})$ with $K_i(\cbb{\Omega }{F'})\times K_i(F')$ using the group isomorphism of the Alexandroff K-theorem \emph{(\h{20.4} a))} then
$$\mae{K_i(\phi _{\Omega ^*})}{K_i(\ccb{\Omega ^*}{F})}{K_i\left(\ccb{\Omega ^*}{F'}\right)}{(a,b)}{(K_i(\phi _\Omega )a,K_i(\phi )b)}\;.$$
\end{co}

The assertion follows from \cor{26.3'b} c).\qed

\begin{co}\label{20.2'a}
We use the notation of the Alexandroff K-theorem \emph{(\h{20.4} a))} and denote by $\omega _\infty $ the Alexandroff point of $\Omega $. Let $\Omega '$ be a locally compact space,
$$\mac{\varphi '}{\cbb{\Omega \times \Omega '}{F}}{\cbb{\Omega ^*\times \Omega '}{F}}$$
the inclusion map, and
$$\mae{\lambda '}{\cbb{\Omega '}{F}}{\cbb{\Omega ^*\times \Omega '}{F}}{x}{\tilde{x} }\,,$$
where
$$\mae{\tilde{x} }{\Omega ^*\times \Omega '}{F}{(\omega ,\omega ')}{x(\omega ')}\;.$$
Then the map
$$K_i\left(\cbb{\Omega \times \Omega '}{F}\right)\times K_i(F)\longrightarrow K_i\left(\cbb{\Omega ^*\times \Omega '}{F}\right),$$
$$(a,b)\longmapsto K_i(\varphi ')a+K_i(\lambda ')b$$
is a group isomorphism.
\end{co}

If we put
$$\mae{\psi '}{\cbb{\Omega ^*\times \Omega '}{F}}{\cbb{\Omega '}{F}}{x}{x(\omega _\infty ,\,\cdot\,) }$$
then
$$\od{\cbb{\Omega \times \Omega '}{F}}{\varphi '}{\cbb{\Omega ^*\times \Omega '}{F}}{\psi '}{\lambda '}{\cbb{\Omega '}{F}}$$
is a split exact sequence in \frm and the assertion follows from the split exact axiom (\axi{27.9'a}).\qed

{\center{\section{Topological sums of locally compact spaces}}}

\begin{p}[Product Theorem]\label{14.11}
Let $(\Omega _j)_{j\in J}$ be a finite family of locally compact spaces, $\Omega $ its topological sum, and for every $j\in J$ let $\mac{\varphi _j}{\cbb{\Omega _j}{F}}{\cbb{\Omega }{F}}$ be the inclusion map and
$$\mae{\psi _j}{\cbb{\Omega }{F}}{\cbb{\Omega _j}{F}}{x}{x|\Omega _j}\;.$$
\begin{enumerate}
\item 
$$\mae{\Phi_i }{\pro{j\in J}K_i(\cbb{\Omega _j}{F})}{K_i(\cbb{\Omega }{F})}{(a_j)_{j\in J}}{\si{j\in J}}K_i(\varphi _j)a_j$$
is a group isomorphism and
$$\mae{\Psi_i }{K_i(\cbb{\Omega }{F})}{\pro{j\in J}K_i(\cbb{\Omega _j}{F})}{a}{(K_i(\psi _j)a)_{j\in J}}$$
is its inverse.
\item If all $\Omega _j$, $j\in J$, belong to $\Upsilon $ then
$$\Omega \in \Upsilon \,,\qquad p(\Omega) =\si{j\in J}p(\Omega _j),\qquad q(\Omega) =\si{j\in J}q(\Omega _j)\,,$$
$$\Phi _{i,\Omega ,F}=\pro{j\in J}\Phi _{i,\Omega _j,F}\,,\qquad \bigcap_{j\in J}(\Omega _j)_\Upsilon \subset \Omega _\Upsilon \;. $$
\item If $\Omega _j$ is $\Upsilon $-null for every $j\in J$ then $\Omega $ is also $\Upsilon $-null and $\Omega ^*\in \Upsilon _1$, where $\Omega ^*$ denotes the Alexandroff compactification of $\Omega $.
\end{enumerate}
\end{p}

a) follows from \pr{28.9'a}.

b) follows from \pr{5.7'b}.

c) By b), $\Omega $ is $\Upsilon $-null and by Alexandroff's K-theorem (\h{20.4} a)), $\Omega ^*\in \Upsilon _1$.\qed

\begin{co}\label{30.1'}
Let $\Omega $ be a locally compact space, $\Gamma $ a closed set of $\Omega $, and $(\Omega _j)_{j\in J}$ a finite family of pairwise disjoint open sets of $\Omega $ such that $\bigcup\limits_{j\in J}\Omega _j=\Omega \setminus \Gamma  $. We denote for every $j\in J$ by $\mac{\varphi _j}{\cbb{\Omega _j}{F}}{\cbb{\Omega }{F}}$ the inclusion map and assume that the maps
$$\mac{K_i(\varphi _j)}{K_i(\cbb{\Omega _j}{F})}{K_i(\cbb{\Omega }{F})}$$
are group isomorphisms. If $\mac{\varphi }{\cbb{\Omega \setminus \Gamma }{F}}{\cbb{\Omega }{F}}$ denotes the inclusion map and if we identify the above groups then $K_i(\cbb{\Omega \setminus \Gamma }{F})\approx K_i(\cbb{\Omega }{F})^J$ and
$$\mae{K_i(\varphi )}{K_i(\cbb{\Omega \setminus \Gamma }{F})}{K_i(\cbb{\Omega }{F})}{(a_j)_{j\in J}}{\si{j\in J}}a_j\;.\qedd$$
\end{co}

\begin{co}\label{5.11}
Let $\Omega $ be a locally compact space such that $\cbb{\Omega }{F}$ is K-null and $\Gamma $  a closed set of $\Omega $.
\begin{enumerate}
\item $K_i(\cbb{\Omega\setminus \Gamma  }{F})\approx K_{i+1}(\ccb{\Gamma }{F})\;.$
\item Assume $\Gamma $ finite and $\Omega $ $\Upsilon $-null, put
$$\mae{\psi }{\cbb{\Omega }{F}}{\ccb{\Gamma }{F}}{x}{x|\Gamma }\,,$$
and denote by $\mac{\varphi }{\cbb{\Omega \setminus \Gamma }{F}}{\cbb{\Omega }{F}}$ the inclusion map and by $\delta _i$ the index maps associated to the exact sequence in \frm
$$\oc{\cbb{\Omega \setminus \Gamma }{F}}{\varphi }{\cbb{\Omega }{F}}{\psi }{\ccb{\Gamma }{F}}\;.$$
Then 
$$K_i(\cbb{\Omega \setminus \Gamma }{F})\approx K_{i+1}(F)^\Gamma \,,$$
$$\Omega \setminus \Gamma \in \Upsilon \,,\qquad p(\Omega \setminus \Gamma )=0\,,\qquad q(\Omega \setminus \Gamma )=Card\,\Gamma \,,\qquad \Phi _{i,(\Omega \setminus \Gamma) ,F}=\delta _{i+1}\;.$$
\end{enumerate}
\end{co} 

a) Since $\cbb{\Omega }{F}$ is K-null, the assertion follows from the six-term axiom (\axi{27.9'd}).

b) follows from a), \lm{14.9'} c), and the Product Theorem (\pr{14.11}).\qed

\begin{co}\label{5.5a}
Let $(\Omega _j)_{j\in J}$ be a finite family of locally compact spaces, $\Omega $ its topological sum, and $\Omega ^*$ the Alexandroff compactification of $\Omega $.
\begin{enumerate}
\item $$K_i(\ccc{C}(\Omega ,F))\approx \pro{j\in J}K_i(\ccc{C}_0(\Omega _j,F))\,,$$
$$ K_i(\ccc{C}(\Omega^* ,F))\approx K_i(F)\times \pro{j\in J}K_i(\ccc{C}_0(\Omega _j,F))\;.$$
\item If all $\Omega _j$, $j\in J$, belong to $\Upsilon $ then
$$\Omega ^*\in \Upsilon \,,\qquad p(\Omega ^*)=1+\si{j\in J}p(\Omega _j)\,,\qquad q(\Omega ^*)=\si{j\in J}q(\Omega _j)\;.$$
\end{enumerate}
\end{co}

The assertion follows immediately from the Product Theorem (\pr{14.11} a)) and the Alexandroff K-theorem (\h{20.4} a)).\qed

\begin{co}\label{26.11}
Let $(\Omega _j)_{j\in J}$ be a finite family of locally compact spaces such that $\cbb{\Omega _j}{F}$ is K-null for every $j\in J$ and let $\Gamma _j$ be a closed set of $\Omega _j$ for every $j\in J$. We denote by $\Omega $ the Alexandroff compactification of the topological sum of the family $(\Omega _j\setminus \Gamma _j)_{j\in J}$.
\begin{enumerate}
\item $K_i(\ccb{\Omega }{F})\approx K_i(F)\times \pro{j\in J}K_{i+1}(\cbb{\Gamma _j}{F})$.
\item If for every $j\in J$, $\Omega _j$ is $\Upsilon $-null and $\Gamma _j$ is finite then
$$\Omega \in \Upsilon \,,\qquad p(\Omega )=1\,,\qquad q(\Omega) =\si{j\in J}Card\,\Gamma _j\;.$$
\end{enumerate}
\end{co}

a) By \cor{5.11} a), $K_i(\cbb{\Omega _j\setminus \Gamma _j}{F})\approx K_{i+1}(\cbb{\Gamma _j}{F})$ for every $j\in J$ so by \cor{5.5a} a), 
$$K_i(\ccb{\Omega }{F})\approx K_i(F)\times \pro{j\in J}K_{i+1}(\cbb{\Gamma _j}{F})\;.$$

b) By \cor{5.11} b), for every $j\in J$,
$$\Omega _j\setminus \Gamma _j\in \Upsilon \,,\qquad p(\Omega _j\setminus \Gamma _j)=0\,,\qquad q(\Omega _j\setminus \Gamma _j)=Card\,\Gamma _j\;.$$
Thus by \cor{5.5a} b),
$$\Omega \in \Upsilon \,,\qquad p(\Omega) =1\,,\qquad q(\Omega) =\si{j\in J}Card\,\Gamma _j\;.\qedd$$

\begin{p}\label{12.10'd}
Let $\Omega $ be a compact space belonging to $\Upsilon _1$,
 $\Gamma $ a closed set of $\Omega $, $\omega_0 \in \Gamma $, and $\Gamma ':=\Gamma \setminus \z{\omega_0 }$. We use the notation of \emph{ the Topological triple (\pr{3.12})} and put there 
$$\Omega _1:=\Omega \,,\qquad\Omega _2:=\Omega \setminus \z{\omega_0 }\,,\qquad\Omega _3:=\Omega \setminus \Gamma\;. $$
\begin{enumerate}
\item $\Omega \setminus \z{\omega _0}$ is $\Upsilon $-null.
\item $\kk{i}{\cbb{\Omega \setminus \Gamma }{F}}\approx \kk{i+1}{\cbb{\Gamma '}{F}}$.
\item  
$$\oddg{K_i(\ccb{\Omega }{F})}{K_i(\psi _{1,3})}{K_i(\ccb{\Gamma }{F})}{\delta _{1,3,i}}{\Phi _i}{20}{10}{10}$$
$$\odeg{\delta _{1,3,i}}{\Phi _i}{K_{i+1}(\cbb{\Omega \setminus \Gamma }{F})}{10}{10}$$
is a split exact sequence, and the maps
$$K_i(\ccb{\Omega }{F})\times K_{i+1}(\cbb{\Omega \setminus \Gamma }{F})\longrightarrow K_i(\ccb{\Gamma }{F}),$$
$$(a,b)\longmapsto K_i(\psi _{1,3})a+\Phi _ib\,,$$
$$\mac{\delta _{2,3,i}}{K_{i}\left(\cbb{\Gamma'}{F}\right)}{K_{i+1}(\cbb{\Omega \setminus \Gamma }{F})}$$
are group isomorphisms.
\item If $\Omega \setminus \Gamma \in \Upsilon $ or $\Gamma '\in \Upsilon $ then with the notation of \emph{\cor{5.9'}}
$$\mac{\delta _i}{\kk{i}{\cbb{\Gamma '}{F}}}{\kk{i+1}{\cbb{\Omega \setminus \z{\omega _0}}{F}}}$$
is a group isomorphism and
$$\Omega \setminus \Gamma ,\Gamma '\in \Upsilon \,,\qquad p(\Omega \setminus \Gamma )=q(\Gamma ')\,,\qquad q(\Omega \setminus \Gamma )=p(\Gamma ')\,,$$
$$ \Phi _{i,(\Omega \setminus \Gamma ),F}=\delta _{i+1}\circ \Phi _{(i+1),\Gamma ',F}\;.$$
\item Assume $\Gamma $ finite.
\begin{enumerate}
\item 
$$\mac{(\delta _{2,3,i})^{-1}}{K_{i+1}(\cbb{\Omega\setminus \Gamma }{F})}{K_{i}(F)^{\Gamma'}}$$ 
is a group isomorphism.
\item 
$$\Omega\setminus \Gamma \in \Upsilon \,,\qquad p(\Omega\setminus \Gamma )=0\,,\qquad q(\Omega\setminus \Gamma) =Card\,\Gamma '\;.$$
\end{enumerate}
\end{enumerate}
\end{p}

a) follows from Alexandroff's K-theorem (\h{20.4} c)).

b) follows from \cor{5.11} a).

c) By a), $\Omega \setminus \z{\omega }$ is K-null and the assertion follows from the Topological triple (\pr{3.12} a)). 

d) follows from \cor{5.9'}.

$e_1)$ follows from c) and the Product Theorem (\pr{14.11} $a_4)$).

$e_2)$ follows from a) and \cor{5.9'} c).\qed

\begin{p}\label{17.11}
Let $\Omega $ be a locally compact space, $\Gamma $ a closed set of $\Omega $, 
$\mac{\varphi }{\cbb{\Omega\setminus \Gamma }{F}}{\cbb{\Omega }{F}}$ the inclusion map,
$$\mae{\psi }{\cbb{\Omega }{F}}{\cbb{\Gamma }{F}}{x}{x|\Gamma }\,,$$
and $\delta _i$ the index maps associated to the exact sequence in \frm
$$\oc{\cbb{\Omega \setminus \Gamma }{F}}{\varphi }{\cbb{\Omega }{F}}{\psi }{\cbb{\Gamma }{F}}\;.$$
Let $(\Omega _j)_{j\in J}$ be a finite family of pairwise disjoint open sets of $\Omega $ the union of which is $\Omega \setminus \Gamma $ and for every $j\in J$ put
$$\mae{\psi_j }{\cbb{\bar{\Omega}_j }{F}}{\cbb{\bar{\Omega }_j\setminus \Omega _j }{F}}{x}{x|(\bar{\Omega }_j\setminus \Omega _j )}\,,$$
$$\mae{\psi'_j }{\cbb{\Omega\setminus \Gamma  }{F}}{\cbb{\Omega _j}{F}}{x}{x|\Omega _j}\,,$$
$$\mae{\psi''_j }{\cbb{\Gamma }{F}}{\cbb{\bar{\Omega }_j\setminus \Omega _j }{F}}{x}{x|(\bar{\Omega }_j\setminus \Omega _j )}$$
and denote by
$$\mac{\varphi _j}{\cbb{\Omega _j}{F}}{\cbb{\bar{\Omega }_j }{F}}\,,$$
$$\mac{\varphi' _j}{\cbb{\Omega _j}{F}}{\cbb{\Omega \setminus \Gamma }{F}}\,,$$
$$\mac{\varphi'' _j}{\cbb{\Omega _j}{F}}{\cbb{\Omega }{F}}$$
the inclusion maps and by $\delta _{j,i}$ the index maps associated to the exact sequence in \frm
$$\oc{\cbb{\Omega _j}{F}}{\varphi _j}{\cbb{\bar{\Omega }_j }{F}}{\psi _j}{\cbb{\bar{\Omega }_j\setminus \Omega _j }{F}}\;.$$
\begin{enumerate}
\item For every $j\in J$,
$$\delta _{j,i}\circ K_i(\psi ''_j)=K_{i+1}(\psi '_j)\circ \delta _i$$
and
$$\delta _i=\si{j\in J}K_{i+1}(\varphi '_j)\circ \delta _{j,i}\circ K_i(\psi ''_j)\;.$$
\item 
$K_i(\varphi )=\si{j\in J}K_i(\varphi ''_j)\circ K_i(\psi '_j)$.
\item Let $j_0\in J$ such that $\cbb{\Omega \setminus \Omega _{j_0}}{F}$ is K-null.
\begin{enumerate}
\item $K_i(\varphi ''_{j_0})$ is a group isomorphism.
\item Assume $\psi $ K-null. If we put
$$\mac{\Phi_i :=K_i(\varphi '_{j_0})\circ K_i(\varphi ''_{j_0})^{-1}}{K_i(\cbb{\Omega }{F})}{K_i(\cbb{\Omega \setminus \Gamma }{F})}$$
then
$$\oddg{K_{i+1}(\cbb{\Gamma }{F})}{\delta _{i+1}}{K_i(\cbb{\Omega \setminus \Gamma }{F})}{K_i(\varphi )}{\Phi_i }{0}{20}{20}$$
$$\odeg{K_i(\varphi )}{\Phi_i }{K_i(\cbb{\Omega }{F})}{20}{20}$$
is a split exact sequence and the map
$$K_{i+1}(\cbb{\Gamma }{F})\times K_i(\cbb{\Omega }{F})\longrightarrow K_i(\cbb{\Omega \setminus \Gamma }{F}),$$
$$(a,b)\longmapsto {\delta _{i+1}a+\Phi_i b}$$
is a group isomorphism.
\end{enumerate}
\end{enumerate}
\end{p}

a) By the commutativity of the index maps (\axi{27.9'e}),
$$\delta _{j,i}\circ K_i(\psi ''_j)=K_{i+1}(\psi '_j)\circ \delta _i\;.$$
Since $\si{j\in J}\varphi '_j\circ \psi '_j$ is the identity map of $\cbb{\Omega \setminus \Gamma }{F}$,
$$\si{j\in J}K_{i+1}(\varphi '_j)\circ \delta _{j,i}\circ K_i(\psi ''_j) =\si{j\in J}K_{i+1}(\varphi '_j)\circ K_{i+1}(\psi '_j)\circ \delta _i=$$
$$=K_{i+1}\left(\si{j\in J}\varphi '_j\circ \psi '_j\right)\circ \delta _i=\delta _i\;.$$

b) We have $\varphi ''_j=\varphi \circ \varphi '_j$ for every $j\in J$. Since $\si{j\in J}\varphi '_j\circ \psi '_j$ is the identity map of $\cbb{\Omega \setminus \Gamma }{F}$,
$$K_i(\varphi) =K_i(\varphi )\circ K_i\left(\si{j\in J}\varphi '_j\circ \psi '_j\right)=$$
$$=\si{j\in J}K_i(\varphi )\circ K_i(\varphi '_j)\circ K_i(\psi '_j)=\si{j\in J}K_i(\varphi ''_j)\circ K_i(\psi '_j)\;.$$ 

$c_1)$ If we put
$$\mae{\bar{\psi } }{\cbb{\Omega }{F}}{\cbb{\Omega \setminus \Omega _{j_0}}{F}}{x}{x|(\Omega \setminus \Omega _{j_0})}$$
then
$$\oc{\cbb{\Omega _{j_0}}{F}}{\varphi ''_{j_0}}{\cbb{\Omega }{F}}{\bar{\psi } }{\cbb{\Omega \setminus \Omega _{j_0}}{F}}$$ 
is an exact sequence in \frm. Since $\cbb{\Omega \setminus \Omega _{j_0}}{F}$ is K-null, it follows that $K_i(\varphi ''_{j_0})$ is a group isomorphism by the Topological six-term sequence (\pr{24.11} $c_1)$).

$c_2)$ Since $\varphi \circ \varphi' _{j_0}=\varphi ''_{j_0}$,
$$K_i(\varphi )\circ \Phi_i =K_i(\varphi )\circ K_i(\varphi '_{j_0})\circ K_i(\varphi ''_{j_0})^{-1}=K_i(\varphi ''_{j_0})\circ K_i(\varphi ''_{j_0})^{-1}=id_{K_i(\cbb{\Omega }{F})}\;.$$
Since $\psi $ is K-null,
$$\og{K_{i+1}(\cbb{\Gamma }{F})}{\delta _{i+1}}{K_i(\cbb{\Omega \setminus \Gamma }{F})}{K_i(\varphi )}{\Phi_i }{K_i(\cbb{\Omega }{F})}{0}{10}{10}$$
is a split exact sequence and this implies the last assertion.\qed

\begin{p}\label{12.10'b}
If $(\Omega _j)_{j\in J}$, $J\not=\emptyset $, is a finite family of compact spaces belonging to $\Upsilon _1$ then $\pro{j\in J}\Omega _j\in \Upsilon _1$. 
\end{p}

The assertion follows immediately from \pr{10.10'b}.\qed

{\center{\section{Homotopy}}}

\begin{p}\label{24.11a}
Let $\Omega $ be a locally compact space, $\Omega ^*$ its Alexandroff compactification, $(\vartheta _s)_{s\in ]0,1]}$ a family of proper continuous maps $\Omega \rightarrow \Omega $, and for every $s\in ]0,1]$ let $\mac{\vartheta _s^*}{\Omega ^*}{\Omega ^*}$ be the continuous extension of $\vartheta _s$. We assume:
\renewcommand{\labelenumi}{\arabic{enumi})} 
\begin{enumerate}
\item $\mad{\Omega^* \times ]0,1]}{\Omega^* }{(\omega ,s)}{\vartheta^* _s(\omega )}$ is continuous,
\item $\vartheta _1(\omega )=\omega $ for every $\omega \in \Omega $,
\item for every compact set $\Gamma $ of $\Omega $ there is an $\varepsilon \in ]0,1]$ with $\Gamma \cap \vartheta _s(\Omega )=\emptyset $ for all $s\in ]0,\varepsilon [$.
\end{enumerate}
\renewcommand{\labelenumi}{\alph{enumi})} 
\renewcommand{\labelenumii}{\alph{enumi}_{\arabicenumii}))} 
Then $\Omega$ is null-homotopic and $\Omega ^*\in \Upsilon _1$.
\end{p}

We put for every $s\in [0,1]$,
$$\mae{\phi _s}{\cbb{\Omega }{\bc}}{\cbb{\Omega }{\bc}}{x}{\abb{x\circ \vartheta _s}{s\in ]0,1]}{0}{s=0}}\;.$$
Then $(\phi _s)_{s\in [0,1]}$ is a pointwise continuous path in $\cbb{\Omega }{\bc}$ with $\phi _0=0$ and $\phi _1$ the identity map of $\cbb{\Omega }{F}$. Thus $\Omega $ is null-homotopic. By \pr{5.7'a} d), $\Omega$ is $\Upsilon $-null and by Alexandroff's K-theorem, (\h{20.4} c)), $\Omega ^*\in \Upsilon _1$.\qed

\begin{co}\label{22.9'}
Let $J$ be a set and $\Omega :=[0,1]^J$.
Then $\Omega \setminus \z{0}$ is null-homotopic and
$\Omega \in \Upsilon_1$.  
\end{co}

The assertion follows from \pr{24.11a} by using the map
$$\mae{\vartheta }{\Omega \times [0,1]}{\Omega }{(\omega ,s)}{s\omega }\;.\qedd$$

\begin{p}\label{4.9}
Let $\Omega $ be a locally compact space, $\Gamma _0,\Gamma _1$ compact subspaces of $\Omega $, $\mac{\vartheta _0}{\Gamma _0}{\Gamma _1}$ a homeomorphism, and $\mac{\vartheta}{\Gamma _0\times [0,1]}{\Omega }$ a continuous map such that $\vartheta (\omega ,0)=\omega $ and $\vartheta (\omega ,1)=\vartheta _0(\omega )$ for every $\omega \in \Gamma _0$. We put
$$\mae{\psi _j}{\cbb{\Omega }{F}}{\ccb{\Gamma _j}{F}}{x}{x|\Gamma _j}$$
for every $j\in \z{0,1}$ and
$$\mae{\varphi }{\ccb{\Gamma _1}{F}}{\ccb{\Gamma _0}{F}}{x}{x\circ \vartheta _0}\;.$$
\begin{enumerate}
\item $K_i(\varphi )$ is a group isomorphism and $K_i(\psi _0)=K_i(\varphi )\circ K_i(\psi _1)$.
\item For every $j\in \z{0,1}$ let $\mac{\varphi _j}{\cbb{\Omega \setminus \Gamma _j}{F}}{\cbb{\Omega }{F}}$ be the inclusion map and $\oaa{\ccb{\Gamma _j}{F}}{\lambda _j}{\cbb{\Omega }{F}}$ be a morphism in \frm such that $\psi _j\circ \lambda _j=id_{\ccb{\Gamma _j}{F}}$ and $\lambda _1=\lambda _0\circ \varphi $.
\begin{enumerate}
\item For every $j\in \z{0,1}$,
$$\oddg{\kk{i}{\cbb{\Omega \setminus \Gamma _j}{F}}}{\kk{i}{\varphi _j}}{\kk{i}{\cbb{\Omega }{F}}}{\kk{i}{\psi _j}}{\kk{i}{\lambda _j}}{20}{20}{20}$$
$$\odeg{\kk{i}{\psi _j}}{\kk{i}{\lambda _j}}{\kk{i}{\cbb{\Gamma _j}{F}}}{20}{20}$$
is a split exact sequence.
\item $Im\,\kk{i}{\varphi _0}=Im\,\kk{i}{\varphi _1}$.
\item If we put for every $j\in \z{0,1}$
$$\mae{\Psi_{j,i}}{\kk{i}{\cbb{\Omega \setminus \Gamma _j}{F}}}{Im\,\kk{i}{\varphi _j}}{a}{\kk{i}{\varphi _j}a}$$ 
then $\Psi _{j,i}$ and
$$\mac{(\Psi _{1,i})^{-1}\circ \Psi _{0,i}}{\kk{i}{\cbb{\Omega \setminus \Gamma _0}{F}}}{\kk{i}{\cbb{\Omega \setminus \Gamma _1}{F}}}$$
are well-defined group isomorphisms.
\item If $\Omega \setminus \Gamma _0\in \Upsilon $ or $\Omega \setminus \Gamma _1\in \Upsilon $ then
$$\Omega \setminus \Gamma _0,\Omega \setminus \Gamma _1\in \Upsilon \,,\qquad p(\Omega \setminus \Gamma _0)=p(\Omega \setminus \Gamma _1)\,,\qquad q(\Omega \setminus \Gamma _0)=q(\Omega \setminus \Gamma _1)\,,$$
$$\Phi _{i,(\Omega \setminus \Gamma _1),F}=(\Psi _{1,i}^F)^{-1}\circ \Psi _{0,i}^F\circ \Phi _{i,(\Omega \setminus \Gamma ),F}\;.$$
\end{enumerate} 
\item If $\Omega $ is compact and if for every $j\in \z{0,1}$ there is a continuous map $\mac{\vartheta '_j}{\Omega }{\Gamma _j}$ such that $\vartheta '_j(\omega )=\omega $ for every $\omega \in \Gamma _j$ and $\vartheta _0 \circ\vartheta '_0=\vartheta '_1$ then the hypotheses of b) are fulfilled.
\end{enumerate}
\end{p}

a) For every $s\in [0,1]$ put
$$\mae{\nu _s}{\cbb{\Omega }{F}}{\ccb{\Gamma _0}{F}}{x}{x(\vartheta (\cdot\, ,s))}\;.$$
Then $K_i(\nu _0)=K_i(\nu _1)$ by the homotopy axiom (\axi{27.9'b}). $K_i(\varphi )$ is obviously a group isomorphism. For every $x\in \cbb{\Omega }{F}$ and $\omega \in \Gamma _0$,
$$(\nu _0x)(\omega )=x(\vartheta (\omega ,0))=x(\omega )=(\psi _0x)(\omega )\,,$$
$$(\nu _1x)(\omega )=x(\vartheta (\omega ,1))=x(\vartheta _0(\omega ))=(\psi _1x)(\vartheta _0(\omega ))=(\varphi \psi _1x)(\omega )\,,$$
so $\nu _0=\psi _0$, $\nu _1=\varphi \circ \psi _1$,
$$K_i(\psi _0)=K_i(\nu _0)=K_i(\nu _1)=K_i(\varphi )\circ K_i(\psi _1)\;.$$

$b_1$ follows from the split exact axiom (\axi{27.9'a}).

$b_2)$ Let $j\in \z{0,1}$. We want to prove
$$Im\,\kk{i}{\varphi _j}=\me{c-\kk{i}{\lambda _j}\kk{i}{\psi _j}c}{c\in \kk{i}{\ccb{\Omega }{F}}}\;.$$
Let $a\in \kk{i}{\cbb{\Omega \setminus \Gamma _j}{F}}$ and put $c:=\kk{i}{\varphi _j}a$. Then
$$c-\kk{i}{\lambda _j}\kk{i}{\psi _j}c=\kk{i}{\varphi _j}a-\kk{i}{\lambda _j}\kk{i}{\psi _j}\kk{i}{\varphi _j}a=\kk{i}{\varphi _j}a\,,$$
which proves the "$\subset$ "-inclusion. Let $c\in \kk{i}{\ccb{\Omega }{F}}$. Then
$$\kk{i}{\psi _j}(c-\kk{i}{\lambda _j}\kk{i}{\psi _j}c)=$$
$$=\kk{i}{\psi _j}c-\kk{i}{\psi _j}\kk{i}{\lambda _j}\kk{i}{\psi _j}c=\kk{i}{\psi _j}c-\kk{i}{\psi _j}c=0\,,$$
$$c-\kk{i}{\lambda _j}\kk{i}{\psi _j}c\in Ker\,\kk{i}{\psi _j}=Im\,\kk{i}{\varphi _j}\,,$$
which proves the "$\supset $"-inclusion (by $b_1)$).

Since $\lambda _1\circ \psi _1=\lambda _0\circ \varphi \circ \psi _1$, we get by a),
$$\kk{i}{\lambda _1}\circ \kk{i}{\psi _1}=\kk{i}{\lambda _0}\circ \kk{i}{\varphi }\circ \kk{i}{\psi _1}=\kk{i}{\lambda _0}\circ \kk{i}{\psi _0}\;.$$
Thus, by the above, $Im\,\kk{i}{\varphi _0}=Im\,\kk{i}{\varphi _1}$.

$b_3)$ By $b_1)$, $\kk{i}{\varphi _0}$ and $\kk{i}{\varphi _1}$ are injective, the assertion follows from $b_2)$.

$b_4)$ Let $\oaa{F}{\phi }{F'}$ be o morphism in \frm and for every $j\in \z{0,1}$ put 
$$\mae{\mu _j}{\cbb{\Omega \setminus \Gamma _j}{F}}{\cbb{\Omega \setminus \Gamma _j}{F'}}{x}{\phi \circ x}\,,$$
$$\mae{\mu }{\cbb{\Omega }{F}}{\cbb{\Omega}{F'}}{x}{\phi \circ x}\;.$$
We mark by a prime the notation associated to $F$ when applied to $F'$. 
For every $j\in \z{0,1}$ the diagram 
$$\begin{CD}
\cbb{\Omega\setminus \Gamma _j }{F}@>\mu _j >>\cbb{\Omega\setminus \Gamma _j }{F'}\\
@V\varphi _j VV         @VV\varphi '_j   V\\
\cbb{\Omega}{F}@>>\mu  >    \cbb{\Omega }{F'}
\end{CD}$$
is commutative. Thus the diagrams
$$\begin{CD}
\kk{i}{\cbb{\Omega\setminus \Gamma _j }{F}}@>\kk{i}{\mu _j} >>\kk{i}{\cbb{\Omega\setminus \Gamma _j }{F'}}\\
@V\kk{i}{\varphi _j} VV         @VV\kk{i}{\varphi '_j} V\\
\kk{i}{\cbb{\Omega}{F}}@>>\kk{i}{\mu } > \kk{i}{   \cbb{\Omega}{F'}}
\end{CD}$$
$$\begin{CD}
\kk{i}{\cbb{\Omega\setminus \Gamma _j }{F}}@>\kk{i}{\mu _j} >>\kk{i}{\cbb{\Omega\setminus \Gamma _j }{F'}}\\
@V\Psi _{j,i}VV         @VV\Psi '_{j,i}V\\
Im\,\kk{i}{\varphi _j}@>>\Lambda _i> Im\,\kk{i}{\varphi '_j}
\end{CD}$$
are also commutative, where $\Lambda _i$ is the map defined by $\kk{i}{\mu }$.

Assume $\Omega \setminus \Gamma _0\in \Upsilon $ and consider the diagram (by $b_2)$)
$$\begin{CD}
\kk{i}{F}^{p(\Omega \setminus \Gamma _0)}\times \kk{i+1}{F}^{q(\Omega \setminus \Gamma _0)}@>\Delta >>A\\
@V\Phi _{i,(\Omega \setminus \Gamma _0),F}VV@VV\Phi _{i,(\Omega \setminus \Gamma _0),F'}V \\ 
\kk{i}{\cbb{\Omega \setminus \Gamma _0}{F}}@>\kk{i}{\mu _0}>>\kk{i}{\cbb{\Omega \setminus \Gamma _0}{F'}}\\ 
@V\Psi _{0,i}VV@VV\Psi' _{0,i}V\\
Im\,\kk{i}{\varphi _0}@>>\Lambda _{i}>Im\,\kk{i}{\varphi '_0}\\
@A\Psi _{1,i}AA@AA\Psi '_{1,i}A\\
\kk{i}{\cbb{\Omega \setminus \Gamma _1}{F}}@>>\kk{i}{\mu _1}>\kk{i}{\cbb{\Omega \setminus \Gamma _1}{F'}}
\end{CD}$$
where
$$\Delta :=\kk{i}{\phi }^{p(\Omega \setminus \Gamma _0)}\times \kk{i+1}{\phi }^{q(\Omega \setminus \Gamma _0)}\,,$$
$$A:=\kk{i}{F'}^{p(\Omega \setminus \Gamma _0)}\times \kk{i+1}{F'}^{q(\Omega \setminus \Gamma _0)}\;.$$
By the above, this diagram is commutative and the assertion follows from $b_3)$.

c) For every $j\in \z{0,1}$ put
$$\mae{\lambda _j}{\ccb{\Gamma _j}{F}}{\ccb{\Omega }{F}}{x}{x\circ \vartheta '_j}\;.$$
Then $\psi _j\circ \lambda _j=id_{\ccb{\Gamma _j}{F}}$ and for every $x\in \ccb{\Gamma _1}{F}$,
$$\lambda _1x=x\circ \vartheta '_1=x\circ \vartheta _0\circ \vartheta '_0=(\varphi x)\circ \vartheta '_0=\lambda _0(\varphi x)\,,\qquad\qquad \lambda _1=\lambda _0\circ \varphi \;.\qedd$$

\begin{co}\label{5.5'}
Let $\Omega $ be a  compact space and $\omega ,\omega '\in \Omega $ such that there is a continuous path in $\Omega $ from $\omega $ to $\omega '$.
\begin{enumerate}
\item $\kk{i}{\cbb{\Omega \setminus \z{\omega }}{F}}\approx \kk{i}{\cbb{\Omega \setminus \z{\omega' }}{F}}$.
\item If $\Omega \setminus \z{\omega} \in \Upsilon $ then
$$\Omega \setminus \z{\omega'} \in \Upsilon\,,\qquad p(\Omega \setminus \z{\omega '})= p(\Omega \setminus \z{\omega })\,,\qquad q(\Omega \setminus \z{\omega '})= q(\Omega \setminus \z{\omega })\;.$$
\end{enumerate}
\end{co}

a) follows from \pr{4.9} $b_3)$ and c).

b) follows from \pr{4.9} $b_4)$ and c).\qed

 \begin{co}\label{19.5'}
Let $\Omega ,\Omega '$ be compact spaces such that $\Omega '\setminus \z{\omega '}$ is null-homotopic for all $\omega '\in \Omega '$, $\omega \in \Omega $, and $\omega ''\in \Omega \times \Omega '$. Then
$$\kk{i}{\cbb{\Omega \setminus \z{\omega }}{F}}\approx \kk{i}{\cbb{(\Omega \setminus \z{\omega })\times \Omega '}{F}}\approx $$
$$\approx  \kk{i}{\cbb{\Omega \times \Omega '\setminus \z{\omega ''}}{F}}\;.$$
\end{co}

Let $\omega ''=:(\omega _0,\omega '_0)\in \Omega \times \Omega '$. By \cor{28.2'} a),
$$\kk{i}{\cbb{(\Omega \setminus \z{\omega_0 })\times \Omega '}{F}}\approx \kk{i}{\cbb{\Omega \times \Omega '\setminus \z{\omega ''}}{F}}$$
and by \pr{4.9} c),
$$\kk{i}{\cbb{(\Omega \setminus \z{\omega })\times \Omega '}{F}}\approx \kk{i}{\cbb{(\Omega \setminus \z{\omega_0 })\times \Omega '}{F}}\;.$$
By \pr{29.3'} $b_3)$,c),
$$\cbb{(\Omega \setminus \z{\omega })\times (\Omega '\setminus \z{\omega '_0})}{F}\approx \cbb{\Omega '\setminus \z{\omega '_0}}{\bc}\otimes \cbb{\Omega \setminus \z{\omega }}{F}$$
is null-homotopic. Since the sequence in \frm
$$\occ{\cbb{(\Omega \setminus \z{\omega })\times (\Omega '\setminus \z{\omega '_0})}{F}}{}{\cbb{(\Omega \setminus \z{\omega })\times \Omega '}{F}}$$
$$\ocd{\cbb{(\Omega \setminus \z{\omega })\times \Omega '}{F}}{}{\cbb{(\Omega \setminus \z{\omega })\times \z{\omega '_0}}{F}}$$
is exact it follows from the topological six-term sequence (\pr{24.11} $a_1)$),
$$\kk{i}{\cbb{(\Omega \setminus \z{\omega })\times \Omega '}{F}}\approx$$
$$\approx \kk{i}{\cbb{(\Omega \setminus \z{\omega })\times \z{\omega '_0}}{F}}\approx  \kk{i}{\cbb{\Omega \setminus \z{\omega }}{F}}\;.\qedd$$

\begin{co}\label{19.9'}
Let $\Omega $ be a locally compact space and $\omega _1,\omega _2\in \Omega $ and for every $j\in \z{1,2}$ put 
$$\mae{\psi _j}{\cbb{\Omega }{F}}{F}{x}{x(\omega_j )}\;.$$
If there is a continuous path in $\Omega $ from $\omega _1$ to $\omega _2$ then
 $\kk{i}{\psi _1 }=\kk{i}{\psi _2}$.
\end{co}

The assertion follows from \pr{4.9} a).\qed

\begin{co}\label{20.12}
Let $\Omega $ be a locally compact space, $\Gamma $ a finite subset of $\Omega $, $\omega _0\in \Omega $, and
$$\mae{\psi }{\cbb{\Omega }{F}}{\ccb{\Gamma }{F}}{x}{x|\Gamma }\,,$$
$$\mae{\psi_{\omega _0} }{\cbb{\Omega }{F}}{F}{x}{x(\omega _0)}\;.$$
If for every $\omega \in \Gamma $ there is a continuous path in $\Omega $ connecting $\omega _0$ with $\omega $ then
$$\mac{K_i(\psi )}{K_i(\cbb{\Omega }{F})}{K_i(\ccb{\Gamma }{F})\approx K_i(F)^{Card\,\Gamma }}\,,$$
$$a\longmapsto (K_i(\psi _{\omega _0})a)_{\omega \in \Gamma }\;.$$
\end{co}

We put 
$$\mae{\psi _\omega }{\cbb{\Omega }{F}}{F}{x}{x(\omega )}$$
for every $\omega \in \Gamma $. By \cor{19.9'} , $K_i(\psi _\omega )=K_i(\psi _{\omega _0})$ for every $\omega \in \Gamma $ and the assertion follows from the Product Theorem (\pr{14.11} a)).\qed

\begin{p}\label{5.2'}
Let $\Omega $ be a path connected compact space, $\Gamma $ a finite subset of $\Omega $, $\omega _0\in \Gamma $, $\Gamma ':=\Gamma \setminus \z{\omega _0}$,
$$\mac{\varphi }{\cbb{\Omega \setminus \Gamma }{F}}{\ccb{\Omega }{F}}\,,$$
$$\mac{\varphi' }{\cbb{\Omega \setminus \Gamma }{F}}{\cbb{\Omega\setminus \z{\omega _0} }{F}}\,,$$
$$\mac{\varphi ''}{\ccb{\Gamma' }{F}}{\ccb{\Gamma  }{F}}$$
the inclusion maps,
$$\mae{\psi }{\ccb{\Omega }{F}}{\ccb{\Gamma }{F}}{x}{x|\Gamma }\,,$$
$$\mae{\psi' }{\cbb{\Omega\setminus \z{\omega _0} }{F}}{\ccb{\Gamma' }{F}}{x}{x|\Gamma' }\,,$$
$$\mae{\psi_\omega  }{\ccb{\Omega }{F}}{F}{x}{x(\omega )}$$
for every $\omega \in \Gamma $,
and $\delta _i,\delta _i'$ the index maps associated to the exact sequences in \frm
$$\oc{\cbb{\Omega \setminus \Gamma }{F}}{\varphi }{\ccb{\Omega }{F}}{\psi }{\ccb{\Gamma }{F}}\,,$$
$$\oc{\cbb{\Omega \setminus \Gamma }{F}}{\varphi' }{\cbb{\Omega\setminus \z{\omega _0} }{F}}{\psi '}{\ccb{\Gamma '}{F}}\;.$$
\begin{enumerate}
\item $K_i(\ccb{\Omega }{F})\approx K_i(F)\times K_i(\cbb{\Omega \setminus \z{\omega _0}}{F})$.
\item $\psi '$ is K-null.
\item If we use the group isomorphism of a) then
$$\mae{K_i(\psi )}{K_i(\ccb{\Omega }{F})}{\kk{i}{\ccb{\Gamma }{F}}\approx K_i(F)^\Gamma }{(a,b)}{(a)_{\omega \in \Gamma }}\;.$$
\item If we identify $K_i(\ccb{\Gamma }{F})$ with $K_i(F)^\Gamma $ and $K_i(\ccb{\Gamma' }{F})$ with $K_i(F)^{\Gamma'} $ then
$$\mae{\delta _i}{K_i(\ccb{\Gamma }{\cdot })}{K_{i+1}(\cbb{\Omega \setminus \Gamma }{\cdot })}{(a_\omega )_{\omega \in \Gamma }}{(\delta '_i(a_\omega -a_{\omega _0}))_{\omega \in \Gamma '}}\;.$$
\item Assume $\cbb{\Omega \setminus \z{\omega _0}}{F}$ K-null.
\begin{enumerate}
\item $\mac{K_i(\psi _{\omega _0})}{K_i(\cbb{\Omega }{F})}{K_i(F)}$ is a group isomorphism.
\item $\mac{\delta '_i}{K_i(\ccb{\Gamma '}{F})}{K_{i+1}(\cbb{\Omega \setminus \Gamma }{F})}$ is a group isomorphism.
\item If we identify $K_i(\ccb{\Gamma '}{F})$ with $K_i(F)^{\Gamma'} $ and $K_i(\ccb{\Gamma }{F})$ with $K_i(F)^\Gamma $ then for all $(a_\omega )_{\omega \in \Gamma '}$
$$K_i(\varphi '')(a_\omega )_{\omega \in \Gamma '}=(a_\omega )_{\omega \in \Gamma }\,,$$
where $a_{\omega _0}=0$.
\item If we identify $K_{i+1}(\cbb{\Omega \setminus \Gamma }{F})$ with $K_i(\ccb{\Gamma '}{F})$ using $(\delta '_i)^{-1}$ of $e_2)$ then for all $(a_\omega )_{\omega \in \Gamma }\in K_i(\ccb{\Gamma }{F})$,
$$\delta _i(a_\omega )_{\omega \in \Gamma }=(a_\omega -a_{\omega _0})_{\omega \in \Gamma '}\;.$$
\end{enumerate}
\end{enumerate}
\end{p}

a) follows from the Alexandroff K-theorem (\h{20.4} a)).

b) Let $\omega \in \Gamma' $ and let $\mac{\vartheta }{[0,1]}{\Omega }$ be a continuous path in $\Omega $ connecting $\omega $ with $\omega _0$. Then for every $x\in \cbb{\Omega \setminus \z{\omega _0}}{F}$ the map
$$\mad{[0,1]}{\cbb{\Omega \setminus \z{\omega _0}}{F}}{x}{x(\vartheta _s(\omega ))}$$
is continuous. By the homotopy axiom (\axi{27.9'b})  , $K_i(\psi _\omega )=0$ so by the Product Theorem (\pr{14.11} a)), $K_i(\psi' )=0$.

c) follows from a), b), and \cor{20.12}.

d) By the commutativity of the index maps (\axi{27.9'e}), $\delta '_i=\delta _i\circ K_i(\varphi '')$ so by the Product Theorem (\pr{14.11} a)),
$$\delta _i(0,(a_\omega )_{\omega \in \Gamma '})=\delta '_i(a_\omega )_{\omega \in \Gamma '}$$
for all $(a_\omega )_{\omega \in \Gamma '}\in K_i(F)^{\Gamma '}$. For $a\in K_i(F)$, by c) and by the above,
$$0=\delta _i\kk{i}{\psi }a=\delta _i(a)_{\omega \in \Gamma }=\delta _i(a, (a)_{\omega \in \Gamma '})=$$
$$=\delta _i(a,0)+\delta _i(0,(a)_{\omega \in \Gamma '})=\delta _i(a,0)+\delta '_i(a)_{\omega \in \Gamma '}\,,$$
$\delta _i(a,0)=-\delta '_i(a)_{\omega \in \Gamma '}$. It follows for all $(a_\omega )_{\omega \in \Gamma }$,
$$\delta _i(a_\omega )_{\omega \in \Gamma }=\delta _i(a_{\omega _0},0)+\delta _i(0,(a_\omega )_{\omega \in \Gamma '})=$$
$$=-\delta '_i(a_{\omega _0})_{\omega \in \Gamma '}+\delta '_i(a_\omega )_{\omega \in \Gamma '}=\delta '_i(a_{\omega }-a_{\omega _0})_{\omega \in \Gamma }\;.$$

$e_1)$ and $e_2)$ follow from the Topological six-term sequence (\pr{24.11}) $a_1)$ and $b_1)$, respectively.

$e_3)$ follows from the Product Theorem (\pr{14.11} a)).

$e_4)$ follows from d).\qed

\begin{e}\label{20.1'}
Let $\bbn$. We use the notation of \emph{\pr{5.2'}} and put
$$\Omega :=\me{re^{\frac{2\pi ij}{n}}}{r\in [0,1],\,j\in \bnn{n}}\,,\quad \Gamma :=\me{e^{\frac{2\pi ij}{n}}}{j\in \bnn{n}}\,,\quad \omega _0:=1\;.$$
\begin{enumerate}
\item $\Omega \setminus \z{\omega _0}$ is null-homotopic and so K-null. 
\item $\mac{K_i(\psi _{\omega _0})}{K_i(\cbb{\Omega }{F})}{K_i(F)}$ is a group isomorphism.
\item $\mac{\delta '_i}{K_i(\ccb{\Gamma '}{F})\approx K_i(F)^{\Gamma '}}{K_{i+1}(\cbb{\Omega \setminus \Gamma }{F})}$ is a group isomorphism.
\item If we identify $K_i(\ccb{\Gamma '}{F})$ with $K_i(F)^{\Gamma'} $ and $K_i(\ccb{\Gamma }{F})$ with $K_i(F)^\Gamma $ \emph{(using e.g. \lm{14.9'} c))} then for all $(a_\omega )_{\omega \in \Gamma '}$
$$K_i(\varphi '')(a_\omega )_{\omega \in \Gamma '}=(a_\omega )_{\omega \in \Gamma }\,,$$
where $a_{\omega _0}=0$.
\item If we identify $K_{i+1}(\cbb{\Omega \setminus \Gamma }{F})$ with $K_i(F)^{\Gamma '}$ using $(\delta '_i)^{-1}$ of c) then for all $(a_\omega )_{\omega \in \Gamma }$,
$$\delta _i(a_\omega )_{\omega \in \Gamma }=(a_\omega -a_{\omega _0})_{\omega \in \Gamma '}\;.$$
\item $\Omega \in \Upsilon \,,\quad p(\Omega) =1\,,\quad q(\Omega) =0\,,\quad\Phi _{i,\Omega ,F}=\kk{i}{\psi _{\omega _0}}\,,\quad\Omega _\Upsilon =\bc_\Upsilon $. 
\end{enumerate}
\end{e}

a) By \pr{24.11a}, $\Omega \setminus \z{\omega _0}$ is null-homotopic.

b) follows from a) and the Topological six-term sequence (\pr{24.11} a)).

c), d), and e) follow from \pr{5.2'} b), c), and d), respectively.

f) follows from a) and \pr{24.11a}.\qed 

\begin{p}\label{12.1'}
Let $\Omega $ be a locally compact spaces, $\omega \in \Omega $, $\Omega '$ a compact space, and 
$$\mac{\vartheta }{\Omega '\times [0,1]}{\Omega }$$
a continuous map such that $\vartheta (\omega ',0)=\omega $ for all $\omega '\in \Omega '$. Then the map
$$\mad{\cbb{\Omega \setminus \z{\omega }}{F}}{\ccb{\Omega '}{F}}{x}{x\circ \vartheta (\,\cdot \,,1)}$$
is K-null
\end{p}

For every $s\in [0,1]$ put 
$$\mae{\psi _s}{\cbb{\Omega \setminus \z{\omega }}{F}}{\ccb{\Omega '}{F}}{x}{x\circ \vartheta (\,\cdot \,,s)}\;.$$
Then for every $x\in \cbb{\Omega \setminus \z{\omega }}{F}$ the map 
$$\mad{[0,1]}{\ccb{\Omega '}{F}}{s}{\psi _sx}$$
 is continuous and $\psi _0x=0$, so the assertion follows from the homotopy (\axi{27.9'b}).\qed 

\begin{p}\label{26.2'a}
Let $\Omega $ be a locally compact space, $\Delta $ a closed set of $\Omega $, $\Gamma $ a compact set of $\Delta $, $\omega _0\in \Gamma $ such that $\cbb{\Delta \setminus \z{\omega _0}}{F}$ is K-null, and $\mac{\vartheta }{\Gamma \times [0,1]}{\Omega }$ a continuous map such that $\vartheta (\omega ,1)=\omega $ and $\vartheta (\omega ,0)=\omega _0$ for all $\omega \in \Gamma $. Then
$$K_i(\cbb{\Omega \setminus \Gamma }{F})\approx K_i(\cbb{\Omega \setminus \z{\omega _0}}{F})\times K_{i+1}(\cbb{\Gamma \setminus \z{\omega _0}}{F})\;.$$
In particular if $\Gamma $ is finite
$$K_i(\cbb{\Omega \setminus \Gamma }{F})\approx K_i(\cbb{\Omega \setminus \z{\omega _0}}{F})\times K_{i+1}(F)^{Card\,\Gamma -1}\;.$$
\end{p}

We use the notation of the Topological triple (\pr{3.12}) and put
$$\Omega _1:=\Omega \setminus \z{\omega _0},\qquad \Omega _2:=\Omega \setminus \Gamma ,\qquad \Omega _3:=\Omega \setminus \Delta \;.$$
By \pr{12.1'}, $\psi _{1,2}$ is K-null and the first assertion follows from the Topological triple (\pr{3.12} $b_4)$). The last assertion follows from the first one and from the Product Theorem (\pr{14.11} a)).\qed

\begin{center}
\chapter{Some selected locally compact spaces}
\end{center}

\fbox{\parbox{12.05cm}{Throughout this chapter we endow $\{0,1\}$ with a group structure  by identifying it with $\bzz{2}$, $F$ denotes an $E$-C*-algebra,  $i\in \{0,1\}$, and $\bbn$ }}
\vspace{10pt}

{\center{\section{Balls}}}

\begin{de}\label{28.3'b}
We put 
$$\bbb_n:=\me{\alpha \in \br^n}{\n{\alpha} \leq 1}\;.$$
\end{de}

\begin{theo}\label{17.1'b}
Let $\Gamma $ be a closed set of $\bbb_n$, $\omega_0 \in \Gamma $, and $\Gamma ':=\Gamma \setminus \z{\omega_0 }$.  
\begin{enumerate}
\item $\bbb_n\setminus \z{\omega_0 }$ is null-homotopic and so $\Upsilon $-null,
$\bbb_n\in \Upsilon_1$, and every exact sequence in \frm belongs to $(\bbb_n)_\Upsilon $. We use in the sequel the notation of \emph{\pr{12.10'd}} and put there $\Omega :=\bbb_n$.
\item $K_i(\cbb{\bbb_n\setminus \Gamma }{F})\approx K_{i+1}\left(\cbb{\Gamma '}{F}\right)$.
\item  
$$\oddg{K_i(\ccb{\bbb_n}{F})}{K_i(\psi _{1,3})}{K_i(\ccb{\Gamma }{F})}{\delta _{1,3,i}}{\Phi _i}{20}{10}{10}$$
$$\odeg{\delta _{1,3,i}}{\Phi _i}{K_{i+1}(\cbb{\bbb_n\setminus \Gamma }{F})}{10}{10}$$
is a split exact sequence, and the maps
$$K_i(\ccb{\bbb_n}{F})\times K_{i+1}(\cbb{\bbb_n\setminus \Gamma }{F})\longrightarrow K_i(\ccb{\Gamma }{F}),$$
$$(a,b)\longmapsto K_i(\psi _{1,3})a+\Phi _ib\,,$$
$$\mac{\delta _{2,3,i}}{K_{i}\left(\cbb{\Gamma'}{F}\right)}{K_{i+1}(\cbb{\bbb_n\setminus \Gamma }{F})}$$
are group isomorphisms.
\item If $\bbb_n\setminus \Gamma \in \Upsilon $ or $\Gamma '\in \Upsilon $ then with the notation of \emph{\cor{5.9'}}
$$\mac{\delta _i}{\kk{i}{\cbb{\Gamma '}{F}}}{\kk{i+1}{\cbb{\bbb_n\setminus \z{\omega _0}}{F}}}$$
is a group isomorphism and
$$\bbb_n\setminus \Gamma ,\Gamma '\in \Upsilon \,,\qquad p(\bbb_n\setminus \Gamma )=q(\Gamma ')\,,\qquad q(\bbb_n\setminus \Gamma )=p(\Gamma ')\,,$$
$$ \Phi _{i,(\bbb_n\setminus \Gamma ),F}=\delta _{i+1}\circ \Phi _{(i+1),\Gamma ',F}\;.$$
\item Assume $\Gamma $ finite.
\begin{enumerate}
\item 
$$\mac{(\delta _{2,3,i})^{-1}}{K_{i+1}(\cbb{\bbb_n\setminus \Gamma }{F})}{K_{i}(F)^{\Gamma'}}$$ 
is a group isomorphism.
\item
$$\mac{K_i(\psi _{1,3})}{K_i(\ccb{\bbb_n}{F})\approx K_i(F)}{K_i(\ccb{\Gamma }{F})\approx K_i(F)^{\Gamma }}\,,$$
$$a\longmapsto (a)_{\omega \in \Gamma }\,,$$
and, if we identify $K_{i+1}(\cbb{\bbb_n\setminus \Gamma }{F})$ with $K_{i}(F)^{\Gamma'}$ using the above group isomorphism $(\delta _{2,3,i})^{-1}$, then 
$$\mae{\delta _{1,3,i}}{K_i(\ccb{\Gamma }{F})}{K_{i}(F)^{Card\,\Gamma'}}{(a_\omega )_{\omega \in \Gamma }}{(a_\omega -a_{\omega_0} )_{\omega \in \Gamma' }}\;.$$
\item 
$$\bbb_n\setminus \Gamma \in \Upsilon \,,\qquad p(\bbb_n\setminus \Gamma )=0\,,\qquad q(\bbb_n\setminus \Gamma) =Card\,\Gamma '\,,$$
$$\Phi _{i,(\bbb_n\setminus \Gamma ),F}=\delta _{2,3,(i+1)}\circ \Phi _{(i+1),\Gamma ',F}\,,$$
\end{enumerate}
\end{enumerate}
\end{theo}

a) Since $\bbb_n$ is homeomorphic to $[0,1]^n$, it follows from \cor{22.9'} that $\cbb{\Omega \setminus \z{\omega _0}}{\bc}$ is null-homotopic and $\bbb_n\in \Upsilon _1$. By \pr{5.7'a} d), $\bbb_n\setminus \z{\omega _0}$ is $\Upsilon $-null and by \pr{11.10'a}, every exact sequence in \frm belongs to $(\bbb_n)_\Upsilon $.

b), c), d), $e_1)$, and $e_3)$ follow from a) and \pr{12.10'd}.

$e_2)$ follows from a) and \pr{5.2'} $e_3),e_4)$.\qed

{\it Remark.} By b), $K_{i}(\cbb{\bbb_n\setminus \Gamma }{F})$ depends only on $K_{i+1}(\cbb{\Gamma'}{F})$ and not on $n$ or on the embedding of $\Gamma $ in $\bbb_n$. 

\begin{co}\label{9.11}
Let $(\Gamma _j)_{j\in J}$ be a finite family of pairwise disjoint closed sets of $\bbb_n$, $J\not=\emptyset $, and for every $j\in J$ let $\omega _j\in \Gamma _j$ such that $\cbb{\Gamma _j\setminus \z{\omega _j}}{F}$ is K-null. Then
$$K_i\left(\cbb{\bbb_n\setminus \bigcup\limits_{j\in J}\Gamma _j}{F}\right)\approx$$
$$\approx  K_i(\cbb{\bbb_n\setminus \me{\omega _j}{j\in J}}{F})\approx K_{i+1}(F)^{Card\,J-1}\;$$
\end{co}

Put $\Gamma :=\bigcup\limits_{j\in J}(\Gamma _j\setminus \z{\omega _j})$,
$$\mae{\psi }{\cbb{\bbb_n\setminus \me{\omega _j}{j\in J}}{F}}{\cbb{\Gamma }{F}}{x}{x|\Gamma }\,,$$
and denote by $\mac{\varphi }{\cbb{\bbb_n\setminus \bigcup\limits_{j\in J}\Gamma _j}{F}}{\cbb{\bbb_n\setminus \me{\omega _j}{j\in J}}{F}}$ the inclusion map. Then
$$\oc{\cbb{\bbb_n\setminus \bigcup\limits_{j\in J}\Gamma _j}{F}}{\varphi }{\cbb{\bbb_n\setminus \me{\omega _j}{j\in J}}{F}}{\psi }{\cbb{\Gamma }{F}}$$
is an exact sequence in \frm. By the Product Theorem (\pr{14.11} c)), $\cbb{\Gamma }{F}$ is K-null so by the Topological six-term sequence (\pr{24.11} b)) and \h{17.1'b} $e_1)$,
$$K_i\left(\cbb{\bbb_n\setminus \bigcup\limits_{j\in J}\Gamma _j}{F}\right)\approx$$
$$\approx  K_i(\cbb{\bbb_n\setminus \me{\omega _j}{j\in J}}{F})\approx K_{i+1}(F)^{Card\,J-1}\;.\qedd$$

\begin{co}\label{26.11a}
Let $(k_j)_{j\in J}$ be a finite family in $\bn$ and for every $j\in J$ let $\Gamma _j$ be a nonempty finite subset of $\bbb_{k_j}$. If $\Omega $ denotes the Alexandroff compactification of the topological sum of the family $(\bbb_{k_j}\setminus \Gamma _j)_{j\in J}$ then
$$\Omega \in \Upsilon \,,\qquad p(\Omega) =1\,,\qquad q(\Omega) =\si{j\in J}(Card\,\Gamma _j-1)\;.$$
\end{co}

For every $j\in J$ let $\omega _j\in \Gamma _j$. By \h{17.1'b} a), $\bbb_{k_j}\setminus \z{\omega _j}$ is $\Upsilon  $-null and the assertion follows from \cor{26.11} b).\qed

\begin{co}\label{20.5'}
If $\Omega $ is a path connected compact space, $\omega \in \Omega $, and $\omega '\in \bbb_n\times \Omega $ then
$$\kk{i}{\cbb{\Omega \setminus \z{\omega }}{F}}\approx \kk{i}{\cbb{\bbb_n\times \Omega \setminus \z{\omega '}}{F}}\;.$$
\end{co}

B< \h{17.1'b} a), $\cbb{\bbb_n\setminus \z{\omega _0}}{F}$ is K-null for every $\omega _0\in \bbb_n$ and the assertion follows from \cor{19.5'}.\qed

\begin{co}\label{23.5'}
Let $\Gamma $ be a closed set of $\bbb_n$ and $\Omega $ an open set of $\bbb_n$, $\Omega \subset \Gamma $. Then for all $\omega \in \Gamma \setminus \Omega $,
$$\kk{i}{\cbb{(\Gamma \setminus \Omega )\setminus \z{\omega }}{F}}\approx \kk{i}{\cbb{\Gamma \setminus \z{\omega }}{F}}\times \kk{i+1}{\cbb{\Omega }{F}}\,,$$
$$\kk{i}{\cbb{\Gamma \setminus \Omega }{F}}\approx \kk{i}{\ccb{\Gamma }{F}}\times \kk{i+1}{\cbb{\Omega }{F}}\;.$$
\end{co}

By \h{17.1'b} b),
$$\kk{i}{\cbb{\Gamma \setminus \z{\omega }}{F}}\approx \kk{i+1}{\cbb{\bbb_n\setminus \Gamma }{F}}\,,$$
$$\kk{i}{\cbb{(\Gamma \setminus \Omega )\setminus \z{\omega }}{F}}\approx \kk{i+1}{\cbb{\bbb_n\setminus (\Gamma \setminus \Omega )}{F}}$$
and by the Product Theorem (\pr{14.11}a)),
$$\kk{i+1}{\cbb{\bbb_n\setminus (\Gamma\setminus \Omega ) }{F}}\approx \kk{i+1}{\cbb{\bbb_n\setminus \Gamma }{F}}\times \kk{i+1}{\cbb{\Omega }{F}}\,,$$
so
$$\kk{i}{\cbb{(\Gamma \setminus \Omega )\setminus \z{\omega }}{F}}\approx \kk{i}{\cbb{\Gamma \setminus \z{\omega }}{F}}\times \kk{i+1}{\cbb{\Omega }{F}}\;.$$
The last relation follows from the Alexandroff K-theorem (\pr{20.4} a)).\qed

\begin{co}\label{9.6'}
If $\Omega $ is an open set of $\bbb_n$, $\Omega \not=\bbb_n$, and $\Gamma $ a compact set of $\Omega $ then
$$\kk{i}{\cbb{\Omega \setminus \Gamma }{F}}\approx \kk{i}{\cbb{\Omega }{F}}\times \kk{i+1}{\ccb{\Gamma }{F}}\;.$$
\end{co}

Let $\omega \in \bbb_n\setminus \Omega $. By \h{17.1'b} b),
$$\kk{i}{\cbb{\Omega }{F}}\approx \kk{i+1}{\cbb{(\bbb_n\setminus \Omega )\setminus \z{\omega }}{F}}\,,$$
$$\kk{i}{\cbb{\Omega \setminus \Gamma }{F}}\approx \kk{i+1}{\cbb{((\bbb_n\setminus \Omega )\setminus \z{\omega })\cup \Gamma }{F}}\;.$$
By the Product Theorem (\pr{14.11} a)),
$$\kk{i+1}{\cbb{((\bbb_n\setminus \Omega )\setminus \z{\omega })\cup \Gamma }{F}}\approx$$
$$\approx  \kk{i+1}{\cbb{(\bbb_n\setminus \Omega )\setminus \z{\omega }}{F}}\times \kk{i+1}{\ccb{\Gamma }{F}}\,,$$
so
$$\kk{i}{\cbb{\Omega \setminus \Gamma }{F}}\approx \kk{i}{\cbb{\Omega }{F}}\times \kk{i+1}{\ccb{\Gamma }{F}}\;.\qedd$$

{\center{\section{Euclidean spaces and Spheres}}}

\begin{de}\label{28.3'd}
We put 
$$S\hspace{-2mm}S_{n-1}:=\me{\alpha \in \br^n}{\n{\alpha }=1}\,,\qquad\qquad \bt:=\bs_1\;.$$
\end{de}

\begin{theo}\label{17.1'c}
\rule{0mm}{0mm}
\begin{enumerate}
\item 
$$\br^n\in \Upsilon \,,\qquad p(\br^n)=\frac{1+(-1)^n}{2}\,,\qquad q(\br^n)=\frac{1-(-1)^n}{2}\,,$$
$$\br_\Upsilon \subset (\br^n)_\Upsilon \,,\qquad \kk{i}{\cbb{\br^n}{F}}\approx \kk{i+n}{F}\;.$$
\item 
$$\bs_n\in \Upsilon \,,\quad p(\bs_n)=\frac{3+(-1)^n}{2}\,,\quad q(\bs_n)=\frac{1-(-1)^n}{2}\,,\quad 
\br_\Upsilon \subset (\bs_n)_\Upsilon\,, $$
$$K_i(\ccb{S\hspace{-2mm}S_n}{F})\approx $$
$$\approx \ab{K_i(F)^2}{n\quad\emph{is \,\,even}}{K_i(F)\times K_{i+1}(F)}{n\quad\emph{is\,\,odd}}=K_i(F)\times K_{i+n}(F)\,,$$
and the map
$$\mad{K_i(F)\times K_{i+n}(F)}{K_i(\ccb{\bs_n}{F})}{(a,b)}{K_i(\lambda )a+K_{i+n}(\varphi )b}$$
is a group isomorphism, where $\mac{\varphi }{\cbb{\br^n}{F}\approx K_{i+n}(F)}{\ccb{\bs_n}{F}}$ denotes the inclusion map and
$$\mae{\lambda }{F}{\ccb{\bs_n}{F}}{x}{x1_{\ccb{\bs_n}{\bc}}}\;.$$
\item Let $\Gamma $ be a closed set of $\br^n$, $\Gamma \not=\br^n$.
\begin{enumerate}
\item The map
$$\mad{\cbb{\br^n}{F}}{\cbb{\Gamma }{F}}{x}{x|\Gamma }$$
is K-null.
\item If $\Gamma $ is compact then
$$K_i(\cbb{\br^n\setminus \Gamma }{F})\approx K_{i+n}(F)\times K_{i+1}(\ccb{\Gamma }{F})\;.$$
If in addition $\Gamma\in \Upsilon $ then $\br^n\setminus \Gamma \in \Upsilon$, and
$$ p(\br^n\setminus \Gamma )=\ab{q(\Gamma) +1}{n \,\;\emph{is\, even}}{q(\Gamma) }{n\,\;\emph{is odd}}\,,$$
$$ q(\br^n\setminus \Gamma )=\ab{p(\Gamma) }{n \,\;\emph{is\, even}}{p(\Gamma) +1}{n\,\;\emph{is odd}}\;.$$
\end{enumerate}
\item If $\Gamma $ is finite then $\br^n\setminus \Gamma \in \Upsilon$, and
$$ p(\br^n\setminus \Gamma )=\ab{1}{n\; \;\emph{is\, even}}{0}{n\;\;\emph{is odd}}\,,$$
$$ q(\br^n\setminus \Gamma )=\ab{\emph{Card}\,\Gamma}{n \;\;\emph{is\, even}}{{\emph{Card}\,\Gamma }+1 }{n\;\;\emph{is odd}}\;.$$
\item Let $\Gamma $ be a closed set of $\bs_n$, $\Gamma \not=\bs_n$, $\omega \in \Gamma $, and $\Gamma ':=\Gamma \setminus \z{\omega }$.
\begin{enumerate}
\item $K_i(\cbb{\bs_n\setminus \Gamma }{F})\approx K_{i+n}(F)\times K_{i+1}(\cbb{\Gamma \setminus \z{\omega }}{F})$.
\item If $\Gamma '\in \Upsilon $ then $\bs_n\setminus \Gamma \in \Upsilon $, and
$$ p(\bs_n\setminus \Gamma )=\ab{q(\Gamma') +1}{n \,\;\emph{is\, even}}{q(\Gamma') }{n\,\;\emph{is odd}}\,,$$
$$ q(\bs_n\setminus \Gamma )=\ab{p(\Gamma')}{n \,\;\emph{is\, even}}{p(\Gamma')+1 }{n\,\;\emph{is odd}}\;.$$
\item If $\Gamma $ is finite, then $ \bs_n\setminus \Gamma \in \Upsilon $, and
$$ p(\bs_n\setminus \Gamma )=\ab{1}{n \,\;\emph{is\, even}}{0 }{n\,\;\emph{is odd}}\,,$$
$$ q(\bs_n\setminus \Gamma )=\ab{\emph{Card}\,\Gamma'}{n \,\;\emph{is\, even}}{\emph{Card}\,\Gamma }{n\,\;\emph{is odd}}\;.$$
\end{enumerate}
\item If $m\in \bn\,,m<n$, then 
$$K_i(\cbb{\bs_n\setminus \bs_m}{F})\approx K_i(\cbb{\br^n\setminus \br^m}{F})\approx K_i(F)\times K_{i+n-m+1}(F)\;.$$
\item For $m\in \bn\,,m<n$,
$$K_i(\cbb{\bbb_n\setminus \bs_m}{F})\approx K_{i+m+1}(F)\;.$$
\end{enumerate}
\end{theo}

a) Since $\br$ is homeomorphic to $]0,1[=\bbb_1\setminus \z{-1,1}$ we get
$$\br\in \Upsilon \,,\qquad p(\br)=0\,,\qquad q(\br)=1$$
 from \h{17.1'b} $e_3)$  and the assertion follows from \cor{27.3'c}.

b) Since $S\hspace{-2mm}S_n$ is homeomorphic to the Alexandroff compactification of $\br^n$, b) follows from a) and the Alexandroff K-theorem (\h{20.4} a),b)).

$c_1)$ We may assume $0\in \br^n\setminus \Gamma $. Put
$$\mae{\vartheta }{\Gamma \times ]0,1]}{\br^n}{(\omega ,s)}{\frac{1}{s}\,\omega }$$
and for every $s\in [0,1]$
$$\mae{\psi _s}{\cbb{\br^n}{F}}{\cbb{\Gamma }{F}}{x}{\ab{x\circ \vartheta (\cdot\, ,\,s)}{s\not=0}{0}{s=0}}\;.$$
Then for every $x\in \cbb{\br^n}{F}$,
$$\mad{[0,1]}{\cbb{\Gamma }{F}}{s}{\psi _sx}$$
is continuous, $\psi _1x=x|\Gamma $, and $\psi _0x=0$. Thus the assertion follows from the homotopy axiom (\axi{27.9'b}).

$c_2)$  We identify the homeomorphic spaces $\me{\alpha \in \br^n}{\n{\alpha }<1}$ and $\br^n$, put $\omega :=(1,0,\cdots,0)\in \bbb_n$ and
$$\mae{\psi }{\cbb{\bbb_n\setminus \z{\omega }}{F}}{\cbb{(\bs_{n-1}\setminus \z{\omega })\cup \Gamma }{F}}{x}{x|((\bs_{n-1}\setminus \z{\omega})\cup \Gamma) }\,,$$
and denote by $\mac{\varphi }{\cbb{\br^n\setminus \Gamma }{F}}{\cbb{\bbb_n\setminus \z{\omega }}{F}}$ the inclusion map and by $\delta _i$ the index maps associated to the exact sequence in \frm
$$\oc{\cbb{\br^n\setminus \Gamma }{F}}{\varphi }{\cbb{\bbb_n\setminus \z{\omega }}{F}}{\psi }{\cbb{(\bs_{n-1}\setminus \z{\omega })\cup \Gamma }{F}}\;.$$
By \h{17.1'b} a), $\cbb{\bbb_n\setminus \z{\omega }}{F}$ is K-null so by the Topological six-term sequence (\pr{24.11} c)), the map
$$\mac{\delta _{i+1}}{K_{i+1}(\cbb{(\bs_{n-1}\setminus \z{\omega })\cup \Gamma }{F})}{K_i(\cbb{\br^n\setminus \Gamma }{F})}$$
is a group isomorphism. By the Product Theorem (\pr{14.11} a),b)), 
$$K_{i+1}(\cbb{(\bs_{n-1}\setminus \z{\omega })\cup \Gamma }{F})\approx K_{i+1}(\cbb{\bs_{n-1}\setminus \z{\omega }}{F})\times K_{i+1}(\ccb{\Gamma }{F})
\,,$$
and $\Gamma \in \Upsilon $ implies $\br^n\setminus \Gamma \in \Upsilon $. By a), $K_{i+1}(\cbb{\bs_{n-1}\setminus \z{\omega }}{F})\approx K_{i+n}(F)$ so
$$K_i(\cbb{\br^n\setminus \Gamma }{F})\approx K_{i+n}(F)\times K_{i+1}(\ccb{\Gamma }{F})$$
as well as the last assertions.

d) follows from c) and the Product Theorem (\pr{14.11} a),b)).

e) $\bs_n\setminus \Gamma $ is homeomorphic to $\br^n\setminus (\Gamma \setminus \z{\omega })$ and the assertion follows from c) and d).

f) 
\begin{center}
Step 1 $K_i(\cbb{\bs_n\setminus \bs_m}{F})\approx K_i(\cbb{\br^n\setminus \br^m}{F})$
\end{center}

Let $\omega \in \bs_m$. Then $\bs_n\setminus \bs_m=(\bs_n\setminus \z{\omega })\setminus (\bs_m\setminus \z{\omega })$. Since $(\bs_n\setminus \z{\omega })\setminus (\bs_m\setminus \z{\omega })$ is homeomorphic to $\br^n\setminus \br^m$ we get
$$K_i(\cbb{\bs_n\setminus \bs_m}{F})\approx K_i(\cbb{\br^n\setminus \br^m}{F})\;.$$

\begin{center}
Step 2 $K_i(\cbb{\br^n\setminus \br^m}{F})\approx K_i(F)\times K_{i+n-m+1}(F)$
\end{center}

We identify $\br^n\setminus \br^m$ with $\me{\alpha \in \bbb_n}{\n{\alpha }<1\,,\sii{j=m+1}{n}\alpha _j^2\not=0}$, put
$$\mae{\psi }{\bbb_n\setminus \bbb_m}{\bs_{n-1}\setminus \bs_{m-1}}{x}{x|(\bs_{n-1}\setminus \bs_{m-1})}\,,$$
and denote by $\mac{\varphi }{\br^n\setminus \br^m}{\bbb_n\setminus \bbb_m}$ the inclusion map and by $\delta _i$ the index maps associated to the exact sequence in \frm
$$\oc{\cbb{\br^n\setminus \br^m}{F}}{\varphi }{\cbb{\bbb_n\setminus \bbb_m}{F}}{\psi }{\cbb{\bs_{n-1}\setminus \bs_{m-1}}{F}}\;.$$
By \pr{24.11a}, $\cbb{\bbb_n\setminus \bbb_m}{F}$ is K-null so by the Topological six-term sequence (\pr{24.11} c)) and Step 1,
$$K_i(\cbb{\br^n\setminus \br^m}{F})\approx K_{i+1}(\cbb{\bs_{n-1}\setminus \bs_{m-1}}{F})\approx $$
$$\approx K_{i+1}\left(\cbb{\br^{n-1}\setminus \br^{m-1}}{F}\right)\;.$$
For $m=1$, by $e_1)$,
$$K_i(\cbb{\br^n\setminus \br}{F})\approx K_{i+1}(\cbb{\bs_{n-1}\setminus \bs_0}{F})\approx  K_{i+n}(F) \times K_i(F)\;.$$
By induction and by the above,
$$K_i(\cbb{\br^n\setminus \br^m}{F})\approx K_{i+n-m+1}\left(\cbb{\br^{n-m+1}\setminus \br}{F}\right)\approx $$
$$\approx K_{i+n-m+1}(F)\times K_i(F)\;.$$

g) Let $\omega \in \bs_m$. Since $\bs_m\setminus \z{\omega }$ is homeomorph to $\br^m$, by a), 
$$K_i(\cbb{\bs_m\setminus \z{\omega }}{F})\approx K_{i+m}(F)\;.$$
By \h{17.1'b} b),
$$K_i(\cbb{\bbb_n\setminus \bs_m}{F})\approx K_{i+1}(\cbb{\bs_m\setminus \z{\omega }}{F})\approx K_{i+1+m}(F)\;.\qedd$$

\begin{e}\label{26.4'}
Put
$$\Omega _1:=\bs_1\cup \me{re^{\frac{2\pi ij}{n}}}{r\in [0,1],\;j\in \bnn{n}}\,,$$
$$\Omega _2:=\bs_2\cup \me{\alpha \in \bbb_3}{\alpha _3=0}\cup \me{\alpha \in \bbb_3}{\alpha _1=\alpha _2=0}\,,$$
$$\Omega _3:=\bs_{n-1}\cup \left(\bigcup_{j\in \bnn{n}}\me{\alpha \in \bbb_n}{\alpha _j=0} \right)\;.$$
\begin{enumerate}
\item $\kk{i}{\ccb{\Omega _1}{F}}=\kk{i}{F}\times \kk{i+1}{F}^n$.
\item $\kk{i}{\ccb{\Omega _2}{F}}\approx \kk{i}{F}^3\times \kk{i+1}{F}^2$.
\item $\kk{i}{\ccb{\Omega _3}{F}}=\kk{i}{F}\times \kk{i+n+1}{F}^{2^n}$.
\end{enumerate}
\end{e}

a) By \h{17.1'c} b) and the Product Theorem (\pr{14.11} a)), 
$$\kk{i}{\cbb{\bbb_2\setminus \Omega _1}{F}}\approx \kk{i}{F}^n$$
and by \h{17.1'b} a),b),c),
$$\kk{i}{\ccb{\Omega _1}{F}}\approx \kk{i}{\ccb{\bbb_2}{F}}\times \kk{i+1}{\cbb{\bbb_2\setminus \Omega _1}{F}}\approx \kk{i}{F}\times \kk{i+1}{F}^n.$$

b) By \h{17.1'c} a),b),
$$\br^2,\bs_1\in \Upsilon \,,\qquad p(\br^2)=1\,,\qquad q(\br^2)=0\,,\qquad p(\bs_1)=1\,,\qquad q(\bs_1)=1\,,$$
so by \cor{27.3'b} $d_1)$,
$$\kk{i}{\cbb{\br^2\times \bs_1}{F}}\approx \kk{i}{F}\times \kk{i+1}{F}\;.$$
Since $\bbb_3\setminus \Omega _2$ is homeomorphic to the topological sum of two copies of $\br^2\times \bs_1$ we get by the Product Theorem (\pr{14.11} a))
$$\kk{i}{\cbb{\bbb_3\setminus \Omega _2}{F}}\approx \kk{i}{F}^2\times \kk{i+1}{F}^2\;.$$
By \h{17.1'b} a),b),c),
$$\kk{i}{\ccb{\Omega _2}{F}}\approx \kk{i}{\ccb{\bbb_3}{F}}\times \kk{i+1}{\cbb{\bbb_3\setminus \Omega _2}{F}}\approx \kk{i}{F}^3\times \kk{i+1}{F}^2.$$

c) By \h{17.1'c} a), $\kk{i}{\cbb{\br^n}{F}}\approx \kk{i+n}{F}$. Since $\bbb_n\setminus \Omega _3$ is homeomorphic to the topological sum of $2^n$ copies of $\br^n$, we get by the Product Theorem (\pr{14.11} a)) $\kk{i}{\cbb{\bbb_n\setminus \Omega _3}{F}}\approx \kk{i+n}{F}^{2^n}$. By \h{17.1'b} a),b),c),
$$\kk{i}{\ccb{\Omega _3}{F}}\approx \kk{i}{\ccb{\bbb_n}{F}}\times \kk{i+1}{\cbb{\bbb_n\setminus \Omega _3}{F}}\approx $$
$$\approx \kk{i}{F}\times \kk{i+n+1}{F}^{2^n}\;.\qedd$$

{\it Remark.} The above a) and b) will be generalized in \ee{27.4'} b) and c), respectively.

\begin{co}\label{20.1'a}
Let $(k_j)_{j\in J}$ be a finite family in $\bn$ and 
$$p:=Card\,\me{j\in J}{k_j \;\emph{is even}}\,,\qquad q:=Card\,\me{j\in J}{k_j\; \emph{is odd}}\;.$$
\begin{enumerate}
\item If $\Omega $ denotes the Alexandroff compactification of the topological sum of the family $(\br^{k_j})_{j\in J}$ then
$$\Omega \in \Upsilon \,,\qquad \br_\Upsilon \subset \Omega _\Upsilon \,,\qquad p(\Omega) =p+1\,,\qquad q(\Omega) =q\;.$$
\item For every $j\in J$ let $\omega _j\in \bs_{k_j}$ and let $\Omega '$ denote the compact space obtained from the topological sum of the family $(\bs_{k_j})_{j\in J}$ by identifying all the points of the family $(\omega_j)_{j\in J}$. If $J\not=\emptyset $ then
$$\Omega' \in \Upsilon \,,\qquad \br_\Upsilon \subset \Omega' _\Upsilon \,,\qquad p(\Omega') =p+1\,,\qquad q(\Omega') =q\;.$$
In particular if $k_j=1$ for all $j\in J$ then $p(\Omega ')=1$, $q(\Omega ')=Card\,J$.
\end{enumerate}
\end{co}

a) By \h{17.1'c} a), $\br^{k_j}\in \Upsilon $, $\br_\Upsilon \subset (\br^{k_j})_\Upsilon $,
$$p\left(\br^{k_j}\right)=\abb{1}{k_j\; \mbox{is even}}{0}{k_j\;\mbox{is odd}}\,,\qquad\qquad q\left(\br^{k_j}\right)=\abb{0}{k_j\; \mbox{is even}}{1}{k_j\;\mbox{is odd}}$$
for every $j\in J$. The assertion follows now from the Product Theorem (\pr{14.11} b)) and from Alexandroff's K-theorem (\pr{20.4} b)).

b) follows from a) since $\Omega $ and $\Omega '$ are homeomorphic.\qed

\begin{co}\label{27.5'}
Let $(k_j)_{j\in J}$ be a finite family in $\bn$, 
$$p:=Card\,\me{j\in J}{k_j \;\emph{is even}}\,,\qquad q:=Card\,\me{j\in J}{k_j\; \emph{is odd}}\,,$$
$(\Gamma _j)_{j\in J}$ a pairwise disjoint family of closed sets of $\bbb_n$ such that $\Gamma _j$ is homeomorphic to $\bs_{k_j}$ for every $j\in J$, and $\Gamma :=\bigcup_{j\in J}\Gamma _j $.  Then
$$\bbb_n\setminus \Gamma \in \Upsilon \,,\qquad \br_\Upsilon \subset (\bbb_n\setminus \Gamma )_\Upsilon \,,\qquad p(\bbb_n\setminus \Gamma) =q\,,\qquad q(\bbb_n\setminus \Gamma)=2p-1 \;.$$
\end{co}

By \h{17.1'c} a),b), for $j\in J$,
$$\br^{k_j},\bs_{k_j}\in \Upsilon \,,\qquad \br_\Upsilon \subset (\br^{k_j})_\Upsilon \cap (\bs_{k_j})_\Upsilon \,,$$
$$p\left(\br^{k_j}\right)=\frac{1+(-1)^{k_j}}{2}\,,\qquad q\left(\br^{k_j}\right)=\frac{1-(-1)^{k_j}}{2}\,,$$
$$p\left(\bs_{k_j}\right)=\frac{3+(-1)^{k_j}}{2}\,,\qquad q\left(\bs_{k_j}\right)=\frac{1-(-1)^{k_j}}{2}\;.$$
Let $\omega \in \Gamma $ and $\Gamma ':=\Gamma \setminus \z{\omega }$. By the Product Theorem (\pr{14.11} b)),
$$\Gamma '\in \Upsilon \,,\qquad \br_\Upsilon \subset \Gamma '_\Upsilon \,,\qquad p(\Gamma ')=2p-1\,,\qquad q(\Gamma ')=q\,,$$
so by \h{17.1'b} d),
$$\bbb_n\setminus \Gamma \in \Upsilon \,,\qquad \br_\Upsilon \subset (\bbb_n\setminus \Gamma )_\Upsilon \,,\qquad p(\bbb_n\setminus \Gamma) =q\,,\qquad q(\bbb_n\setminus \Gamma)=2p-1 \;.\qedd$$

\begin{co}\label{27.5'a}
If $\Omega $ is a connected closed set of $\bbb_2$ possessing a triangulation with $r_0$ vertices, $r_1$ chords, and $r_2$ triangles then
$$\kk{i}{\ccb{\Omega }{F}}\approx \kk{i}{F}\times \kk{i+1}{F}^{1-r_0+r_1-r_2}\;.$$
\end{co}

{\it Sketch of a proof.} If $\Omega $ has $k$ holes then $r_0-r_1+r_2+k=1$. By \h{17.1'b} c),
$$\kk{i}{\ccb{\Omega }{F}}\approx \kk{i}{F}\times \kk{i+1}{\cbb{\bbb_2\setminus \Omega }{F}}\;.$$
By \h{17.1'c} a) and the Product Theorem (\pr{14.11} a)),
$$\kk{i}{\cbb{\bbb_2\setminus \Omega }{F}}\approx \kk{i}{F}^k$$
so
$$\kk{i}{\ccb{\Omega }{F}}\approx \kk{i}{F}\times \kk{i+1}{F}^{1-r_0+r_1-r_2}\;.\qedd$$

\begin{co}\label{18.1'a}
We identify the homeomorphic spaces $\br^n$ and 
$$\me{\alpha \in \br^n}{\n{\alpha }<1}\;.$$ 
Let $\Gamma $ be a finite subset of $\br^n$, $\Delta $ a subset of $\Gamma $, $\omega \in \Delta $, $\Gamma ':=\Gamma \setminus \z{\omega }$, $\Delta ':=\Delta \setminus \z{\omega }$. We use the notation of the \emph{Topological triple (\pr{3.12})} and put
$$\Omega _1:=\bbb_n\setminus \z{\omega }\,,\qquad \Omega _2:=\br^n\setminus \Delta \,,\qquad \Omega _3:=\br^n\setminus \Gamma \;.$$
\begin{enumerate}
\item $\delta_{1,2,i}$ and $\delta_{1,3,i}$ are group isomorphisms.
\item $\psi_{2,3} $ is K-null.
\item If we put $\Phi _i:=\delta_{1,3,(i+1)}\circ K_{i+1}(\varphi ')\circ (\delta_{1,2,(i+1)})^{-1}$ then 
$$\oddg{K_{i+1}(\ccb{\Gamma \setminus \Delta }{F})}{\delta _{2,3,(i+1)}}{K_i(\cbb{\br^n\setminus \Gamma }{F})}{K_i(\varphi_{2,3} )}{\Phi _i}{20}{20}{20}$$
$$\odeg{K_i(\varphi_{2,3} )}{\Phi _i}{K_i(\cbb{\br^n\setminus \Delta }{F})}{20}{20}$$
is a split exact sequence and  the map
$$K_{i+1}(\ccb{\Gamma \setminus \Delta }{F})\times K_i(\cbb{\br^n\setminus \Delta }{F})\longrightarrow K_i(\cbb{\br^n\setminus \Gamma }{F})\,,$$
$$(a,b)\longmapsto \delta _{2,3,(i+1)}a+\Phi _ib$$
is a group isomorphism.
\end{enumerate}
\end{co}

By \h{17.1'b} a), $\cbb{\Omega _1}{F}$ is K-null and by \pr{12.1'}, $\psi _{2,3}$ is K-null. By the Product Theorem (\pr{14.11} a)), 
$$K_i(\psi \circ \varphi ')=id_{K_i(\cbb{\Omega _1\setminus \Omega _2}{F})}$$
and a) and c) follow from the Topological triple (\pr{3.12}  c)).\qed

\begin{co}\label{31.12}
Let $\omega \in \bs_{n-1}$. We use the notation of the \emph{Topological triple (\pr{3.12})} and put
$$\Omega _1:=\bbb_n\,,\qquad\qquad \Omega _2:=\bbb_n\setminus \z{\omega }\,,\qquad\qquad \Omega_3:=\bbb_n\setminus \bs_{n-1}\;. $$
\begin{enumerate}
\item $\varphi_{1,3} $ is K-null. 
\item $\mac{\delta _{2,3,i}}{K_i(\cbb{\bs_{n-1}\setminus \z{\omega }}{F})}{K_{i+1}(\cbb{\bbb_n\setminus \bs_{n-1}}{F})}$ is a group isomorphism.
\item If we put $\Phi _i:=K_i(\varphi)\circ (\delta _{2,3,i})^{-1}$ then
$$\oddg{K_i(\ccb{\bbb_n}{F})}{K_i(\psi_{1.3} )}{K_i(\ccb{\bs_{n-1}}{F})}{\delta _{1,3,i}}{\Phi _i}{20}{20}{20}$$
$$\odeg{\delta _{1,3,i}}{\Phi _i}{K_{i+1}(\cbb{\bbb_n\setminus \bs_{n-1}}{F})}{20}{20}$$
is a split exact sequence and the map
$$K_i(\ccb{\bbb_n}{F})\times K_{i+1}(\cbb{\bbb_n\setminus \bs_{n-1}}{F})\longrightarrow K_i(\ccb{\bs_{n-1}}{F}),$$
$$(a,b)\longmapsto K_i(\psi_{1,3} )a+\Phi _ib$$
is a group isomorphism.
\item Let $\mac{\phi }{G}{H}$ be a morphism in \frm and put
$$\mae{\phi _{\bbb}}{\ccb{\bbb_n}{G}}{\ccb{\bbb_n}{H}}{x}{\phi\circ x}\,,$$
$$\mae{\phi _{\bs}}{\ccb{\bs_{n-1}}{G}}{\ccb{\bs_{n-1}}{H}}{x}{\phi\circ x}\,,$$
$$\mae{\phi _{\bbb,\bs}}{\cbb{\bbb_n\setminus \bs_{n-1}}{G}}{\cbb{\bbb_n\setminus \bs_{n-1}}{H}}{x}{\phi\circ x}\;.$$
If we identify $K_i(\ccb{\bs_{n-1}}{F})$ with 
$$K_i(\ccb{\bbb_n}{F})\times K_{i+1}(\cbb{\bbb_n\setminus \bs_{n-1}}{F})$$
 for $F\in \z{G,H}$ using the isomorphism of c) then
$$\mac{K_i(\phi _{\bs})}{K_i(\ccb{\bs_{n-1}}{G})}{K_i(\ccb{\bs_{n-1}}{H})},$$
$$(a,b)\longmapsto (K_i(\phi _{\bbb})a,K_{i+1}(\phi _{\bbb,\bs})b)\;.$$
\end{enumerate}
\end{co}

By \h{17.1'b} a), $\cbb{\bbb_n\setminus \z{\omega }}{F}$ is K-null and the assertion follows from the Topological triple (\pr{3.12} a)) and \cor{1.1'} b).\qed

\begin{p}\label{31.1'}
Put
$$\Omega :=\bbb_{n+1}\setminus \me{\alpha \in \bs_n}{\alpha _{n+1}=0}\,,$$
$$\Omega ':=\bs_n\setminus \me{\alpha \in \bs_n}{\alpha _{n+1}=0}\,,$$
$$\mae{\psi }{\cbb{\Omega }{F}}{\cbb{\Omega '}{F}}{x}{x|\Omega '}$$
and denote by
$$\mac{\varphi }{\cbb{\bbb_{n+1}\setminus \bs_n}{F}}{\cbb{\Omega }{F}}$$
the inclusion map and by $\delta _i$ the index maps associated to the exact sequence in \frm
$$\oc{\cbb{\bbb_{n+1}\setminus \bs_n}{F}}{\varphi }{\cbb{\Omega }{F}}{\psi }{\cbb{\Omega '}{F}}\;.$$
\begin{enumerate}
\item 
$$K_i(\cbb{\Omega }{F})\approx K_{i+n}(F)\,,\qquad K_i\left(\cbb{\Omega '}{F}\right)\approx K_{i+n}(F)^2\,,$$
$$K_{i+1}(\cbb{\bbb_{n+1}\setminus \bs_n}{F})\approx K_{i+n}(F)\;.$$
\item If we identify the groups of a) then
$$\mae{\delta _i}{K_i\left(\cbb{\Omega '}{F}\right)}{K_{i+1}(\cbb{\bbb_{n+1}\setminus \bs_n}{F})}{(a,b)}{a+b}\,,$$
$$\oc{K_i(\cbb{\Omega }{\cdot })}{K_i(\psi )}{K_i\left(\cbb{\Omega '}{\cdot }\right)}{\delta _i}{K_{i+1}(\cbb{\bbb_{n+1}\setminus \bs_n}{\cdot })}$$
is an exact sequence, and there is a group automorphism $\mac{\Phi_i }{K_{i+n}(F)}{K_{i+n}(F)}$ such that
$$\mae{K_i(\psi )}{K_i(\cbb{\Omega }{F})}{K_i\left(\cbb{\Omega '}{F}\right)}{a}{(\Phi_i a,-\Phi_i a)}\;.$$
\item If
$$\mac{\lambda '}{\cbb{\Omega }{F}}{\ccb{\bbb_{n+1}}{F}}\,,$$
$$\mac{\lambda ''}{\cbb{\Omega' }{F}}{\ccb{\bs_n}{F}}$$
denote the inclusion maps and if we identify $K_i(\cbb{\Omega '}{F})$ with $K_{i+n}(F)^2$ using a) and $K_i(\ccb{\bs_n}{F})$ with $K_i(F)\times K_{i+n}(F)$ using \emph{\h{17.1'c} b)} then $\lambda '$ is K-null and
$$\mae{K_i\left(\lambda ''\right)}{K_i\left(\ccb{\Omega '}{F}\right)}{K_i(\ccb{\bs_n}{F})}{(a,b)}{(0,a+b)}\;.$$ 
\end{enumerate}
\end{p}

a) By \h{17.1'c} a), $K_i(\cbb{\br^n}{F})\approx K_{i+n}(F)$. Since $\bbb_{n+1}\setminus \bs_n$ is homeomorphic to $\br^{n+1}$, $K_{i+1}(\cbb{\bbb_{n+1}\setminus \bs_n}{F})\approx K_{i+n}(F)$. Since $\Omega '$ is homeomorphic to the topological sum of $\br^n$ and $\br^n$, $K_i(\cbb{\Omega '}{F})\approx K_{i+n}(F)^2$ by the Product Theorem (\pr{14.11} a)). Put 
$$\Gamma :=\me{\alpha \in \Omega }{\alpha _{n+1}=0}$$
 and for every $s\in ]0,1]$
$$\mae{\vartheta _s}{\Omega \setminus \Gamma }{\Omega \setminus \Gamma }{(\alpha _j)_{j\in \bnn{n+1}}}{((\alpha _j)_{j\in \bnn{n}},s\alpha _{n+1})}\;.$$
By \pr{24.11a}, $\cbb{\Omega \setminus \Gamma }{F}$ is K-null, so by the Topological six-term sequence (\pr{24.11} a)), $K_i(\cbb{\Omega }{F})\approx K_i(\cbb{\Gamma }{F})$. Since $\Gamma $ is homeomorphic to $\br^n$,  $K_i(\cbb{\Omega }{F})\approx K_{i+n}(F)$ by the above.

b) Put $\omega :=(1,0,\cdots,0)\in \bbb_{n+1}$,
$$\mae{\psi '}{\cbb{\bbb_{n+1}\setminus \z{\omega }}{F}}{\cbb{\bs_n\setminus \z{\omega }}{F}}{x}{x|(\bs_n\setminus \z{\omega })}\,,$$
and denote by 
$$\mac{\varphi '}{\cbb{\bbb_{n+1}\setminus \bs_n}{F}}{\cbb{\bbb_{n+1}\setminus \z{\omega }}{F}}\,,$$
$$\mac{\varphi ''}{\cbb{\Omega }{F}}{\cbb{\bbb_{n+1}\setminus \z{\omega }}{F}}$$
$$\mac{\varphi '''}{\cbb{\Omega '}{F}}{\cbb{\bs_n\setminus \z{\omega }}{F}}$$
the inclusion maps and by $\delta '_i$ the six-term sequence index maps associated with the exact sequence in \frm
$$\oc{\cbb{\bbb_{n+1}\setminus \bs_n}{F}}{\varphi '}{\cbb{\bbb_{n+1}\setminus \z{\omega }}{F}}{\psi '}{\cbb{\bs_n\setminus \z{\omega }}{F}}\;.$$
By \h{17.1'b} a), $\cbb{\bbb_{n+1}\setminus \z{\omega }}{F}$ is K-null so by the Topological six-term sequence (\pr{24.11} c)),
$$\mac{\delta' _i}{K_i(\cbb{\bs_n\setminus \z{\omega }}{F})}{K_{i+1}(\cbb{\bbb_{n+1}\setminus \bs_n}{F})}$$
is a group isomorphism. By the commutativity of the index maps (\axi{27.9'e}), $\delta _i=\delta '_i\circ K_i(\varphi ''')$. Thus if we identify the above groups using $\delta '_i$ then $\delta _i$ is identified with $K_i(\varphi ''')$. By \cor{30.1'}
$$\mae{K_i(\varphi ''')}{\kk{i}{\cbb{\Omega '}{F}}}{\kk{i}{\cbb{\bs_n\setminus \z{\omega }}{F}}}{(a,b)}{a+b}\;.$$
 Since $\bs_n\setminus \z{\omega }$ is homeomorphic to $\br^n$, we get
$$\mae{\delta _i}{K_i\left(\cbb{\Omega '}{F}\right)}{K_{i+1}(\cbb{\bbb_{n+1}\setminus \bs_n}{F})}{(a,b)}{a+b}\;.$$
 Thus $\delta _i$ is surjective and the other assertions follow from the six-term axiom (\axi{27.9'd}).

c) $\lambda '$ is K-null since it factorizes through null (\h{17.1'b} a)). Put $\omega :=(1,0,\cdots,0)\in \bbb_{n+1}$ and denote by
$$\mac{\lambda '''}{\cbb{\bs_n\setminus \z{\omega }}{F}}{\ccb{\bs_n}{F}}$$
the inclusion map.
By the proof of b), since $\lambda ''=\lambda '''\circ \varphi '''$,
$$\mae{K_i\left(\lambda ''\right)}{K_i\left(\ccb{\Omega '}{F}\right)}{K_i(\ccb{\bs_n}{F})}{(a,b)}{(0,a+b)}$$
by the Alexandroff K-theorem (\h{20.4} a)).\qed

\begin{p}\label{27.11}
Let $\Gamma $ be a closed set of $\br^n$, $\Gamma \not=\br^n$, 
$$\mac{\varphi }{\cbb{\br^n\setminus \Gamma }{F}}{\cbb{\br^n}{F}}$$
the inclusion map,
$$\mae{\psi }{\cbb{\br^n}{F}}{\cbb{\Gamma }{F}}{x}{x|\Gamma }\,,$$
and $\delta _i$ the index maps associated to the exact sequence in \frm
$$\oc{\cbb{\br^n\setminus \Gamma }{F}}{\varphi }{\cbb{\br^n}{F}}{\psi }{\cbb{\Gamma }{F}}\;.$$
\begin{enumerate}
\item $\psi $ is K-null.
\item The sequence
$$\oc{K_{i+1}(\cbb{\Gamma }{F})}{\delta _{i+1}}{K_i(\cbb{\br^n\setminus \Gamma }{F})}{K_i(\varphi )}{\cbb{\Gamma }{F}}$$
is exact.
\item Let $(\Omega _j)_{j\in J}$ be a finite family of pairwise disjoint open sets of $\br^n$ the union of which is $\br^n\setminus \Gamma $. If there is a $j_0\in J$ such that $\cbb{\br^n\setminus \Omega _{j_0}}{F}$ is K-null then for every clopen set $\Gamma '$ of $\Gamma $
$$K_i\left(\cbb{\br^n\setminus \Gamma' }{F}\right)\approx K_{i+1}\left(\cbb{\Gamma' }{F}\right)\times K_{i+n}(F)\;.$$
\end{enumerate}
\end{p}

a) follows from \pr{12.1'}.

b) follows from a) and the six-term axiom (\axi{27.9'd}).

c) We use the notation of \pr{17.11}. For $\Gamma '=\Gamma $ the assertion follows from \pr{17.11} $c_2)$ and \h{17.1'c} b). Let
$$\mac{\tilde{\varphi } }{\cbb{\br^n\setminus \Gamma '}{F}}{\cbb{\br^n}{F}}\,,$$
$$\mac{\tilde{\tilde{\varphi }} }{\cbb{\br^n\setminus \Gamma }{F}}{\cbb{\br^n\setminus \Gamma '}{F}}$$
be the inclusion maps,
$$\mae{\tilde{\psi } }{\cbb{\br^n}{F}}{\cbb{\Gamma '}{F}}{x}{x|\Gamma '}\,,$$
$\tilde{\delta}_i$ the index maps associated to the exact sequence in \frm
$$\oc{\cbb{\br^n\setminus \Gamma '}{F}}{\tilde{\varphi } }{\cbb{\br^n}{F}}{\tilde{\psi } }{\cbb{\Gamma '}{F}}\,,$$
and $\tilde{\Phi_i }:=K_i\left(\tilde{\tilde{\varphi } }\right)\circ \Phi_i   $. Since $\varphi =\tilde{\varphi }\circ \tilde{\tilde{\varphi } }  $,
$$K_i(\tilde{\varphi } )\circ \tilde{\Phi_i} =K_i(\tilde{\varphi } )\circ K_i\left(\tilde{\tilde{\varphi } } \right)\circ \Phi_i =K_i(\varphi) \circ \Phi_i =id_{K_i(\cbb{\br^n}{F})}\;.$$
Thus 
$$\og{K_{i+1}\left(\cbb{\Gamma '}{F}\right)}{\tilde{\delta}_{i+1} }{K_i\left(\cbb{\br^n\setminus \Gamma '}{F}\right)}{K_i(\tilde{\varphi } )}{\tilde{\Phi_i } }{K_i(\cbb{\br^n}{F})}{0}{10}{10}$$
is a split exact sequence and this implies c).\qed

\begin{p}\label{20.5'a}
Let $\Omega ,\Omega '$ be compact spaces and $m\in \bn$. If $\Omega $ is path connected, $\Omega \times \Omega '\subset \bbb_n$, and $\bbb_n\setminus (\Omega \times \Omega ')$ is homeomorphic to the topological sum of $\bbb_n\setminus (\Omega \times \bbb_m)$ and $\Omega \times (\bbb_m\setminus \Omega ')$ then for all $\omega \in \Omega $ and $\omega _0\in \Omega \times \Omega '$
$$\kk{i}{\cbb{(\Omega \times \Omega ')\setminus \z{\omega _0}}{F}}\approx $$
$$\approx \kk{i}{\cbb{\Omega \setminus \z{\omega }}{F}}\times \kk{i+1}{\cbb{\Omega \times (\bbb_m\setminus \Omega ')}{F}}\;.$$
In particular if there is a $p\in \bn$ such that $\bbb_m\setminus \Omega '$ is homeomorphic to $p$ copies of $\br^m$ then
$$\kk{i}{\cbb{(\Omega \times \Omega ')\setminus \z{\omega _0}}{F}}\approx $$
$$\approx \kk{i}{\cbb{\Omega \setminus \z{\omega }}{F}}\times \kk{i+m+1}{\ccb{\Omega }{F}}^p\;.$$
\end{p}

By \h{17.1'b} b) and the Product Theorem (\pr{14.11} a)),
$$\kk{i}{\cbb{(\Omega \times \Omega ')\setminus \z{\omega _0}}{F}}\approx \kk{i+1}{\cbb{\bbb_n\setminus (\Omega \times \Omega ')}{F}}\approx $$
$$\approx \kk{i+1}{\cbb{\bbb_n\setminus (\Omega \times \bbb_m)}{F}}\times \kk{i+1}{\cbb{\Omega \times (\bbb_m\setminus \Omega ')}{F}}\;.$$
By \h{17.1'b} b) and \cor{20.5'},
$$\kk{i+1}{\cbb{\bbb_n\setminus (\Omega \times \bbb_m)}{F}}\approx \kk{i}{\cbb{(\Omega \times \bbb_m)\setminus \z{\omega _0}}{F}}\approx$$
$$\approx  \kk{i}{\cbb{\Omega \setminus \z{\omega }}{F}}$$
and so
$$\kk{i}{\cbb{(\Omega \times \Omega ')\setminus \z{\omega _0}}{F}}\approx $$
$$\approx \kk{i}{\cbb{\Omega \setminus \z{\omega }}{F}}\times \kk{i+1}{\cbb{\Omega \times (\bbb_m\setminus \Omega ')}{F}}\;.$$

We prove now the last assertion. By \h{17.1'b} a),
$$\kk{i+1}{\cbb{\br^m}{\ccb{\Omega }{F}}}\approx \kk{i+m+1}{\ccb{\Omega }{F}}$$
so by the Product Theorem (\pr{14.11} a)),
$$\kk{i+1}{\cbb{\Omega \times (\bbb_m\setminus \Omega ')}{F}}\approx \kk{i+1}{\cbb{\bbb_m\setminus \Omega '}{\ccb{\Omega }{F}}}\approx $$
$$\approx \kk{i+1}{\cbb{\br^m}{\ccb{\Omega }{F}}}^p\approx \kk{i+m+1}{\ccb{\Omega }{F}}^p\,,$$
$$\kk{i}{\cbb{(\Omega \times \Omega ')\setminus \z{\omega _0}}{F}}\approx $$
$$\approx \kk{i}{\cbb{\Omega \setminus \z{\omega }}{F}}\times \kk{i+m+1}{\ccb{\Omega }{F}}^p\;.\qedd$$

\begin{co}\label{13.5'}
Let $\Omega $ be a connected graph contained in $\bbb_2$ and containing $\bs_1$, $r_0$ and $r_1$ the number of vertices and chords of $\Omega $, respectively, and $\Gamma $ a nonempty finite subset of $\bs_n\times \Omega $. Then
$$\kk{i}{\cbb{(\bs_n\times \Omega)\setminus \Gamma  }{F}}\approx$$
$$\approx  \kk{i+n}{F}\times \kk{i+1+n}{F}^{1-r_0+r_1}\times \kk{i+1}{F}^{r_1-r_0+Card\,\Gamma }\;.$$
\end{co}

Assume first $\Gamma =\z{\omega_0 }$ for some $\omega_0 \in \bs_n\times \Omega $. There is an embedding of $\bs_n\times \Omega $ in $\bbb_{n+2}$ such that $\bbb_{n+2}\setminus (\bs_n\times \Omega )$ is homeomorphic to the topological sum of $\bbb_{n+2}\setminus (\bs_n\times \bbb_2)$ and $\bs_n\times (\bbb_2\setminus \Omega )$. Since $\bs_n\times (\bbb_2\setminus \Omega )$ is homeomorphic to $1-r_0+r_1$ copies of $\bs_n\times \br^2$, we get by \pr{20.5'a}, for $\omega \in \bs_n$,
$$\kk{i}{\cbb{(\bs_n\times \Omega) \setminus \z{\omega_0 }}{F}}\approx $$
$$\approx \kk{i}{\cbb{\bs_n\setminus \z{\omega }}{F}}\times \kk{i+1}{\ccb{\bs_n}{F}}^{1-r_0+r_1}\;.$$
By \h{17.1'c} a),b),
$$\kk{i}{\cbb{(\bs_n\times \Omega)\setminus  \z{\omega_0 }}{F}}\approx \kk{i+n}{F}\times \kk{i+1+n}{F}^{1-r_0+r_1}\times \kk{i+1}{F}^{1-r_0+r_1}.$$

By \pr{26.2'a},
$$\kk{i}{\cbb{(\bs_n\times \Omega) \setminus \Gamma }{F}}\approx$$
$$\approx  \kk{i}{\cbb{(\bs_n\times \Omega )\setminus \z{\omega_0 }}{F}}\times \kk{i+1}{F}^{Card\,\Gamma -1}\approx $$
$$\approx \kk{i+n}{F}\times \kk{i+1+n}{F}^{1-r_0+r_1}\times \kk{i+1}{F}^{r_1-r_0+Card\,\Gamma }\;.\qedd$$

\begin{co}\label{21.5'}
If 
$$\Omega :=\bs_{n-1}\cup \left(\bigcup _{j\in \bnn{n}}\me{\alpha \in \bbb_n}{\alpha _j=0}\right)\,,$$
$m\in \bn$, and $\Gamma $ is a finite subset of $\bs_m\times \Omega $ then
$$\kk{i}{\cbb{(\bs_m\times \Omega )\setminus \Gamma }{F}}\approx $$
$$\approx \kk{i+m}{F}\times \kk{i+n+1}{F}^{2^n}\times \kk{i+m+n+1}{F}^{2^n}\times \kk{i+1}{F}^{Card\,\Gamma -1}\;.$$
\end{co}

Assume first $\Gamma =\z{\omega_0 }$ for some $\omega_0 \in \bs_m\times \Omega $. There is an embedding of $\bs_m\times \Omega $ in $\bbb_{m+n+1}$ such that $\bbb_{m+n+1}\setminus (\bs_m\times \Omega )$ is homeomorphic to the topological sum of $\bbb_{m+n+1}\setminus (\bs_m\times \bbb_n)$ and $\bs_m\times (\bbb_n\setminus \Omega )$. Since $\bbb_n\setminus \Omega $ is homeomorphic to the topological sum of $2^n$ copies of $\br^n$, by \pr{20.5'a}, for $\omega \in \bs_m$,
$$\kk{i}{\cbb{(\bs_m\times \Omega) \setminus \z{\omega_0 }}{F}}\approx $$
$$\approx \kk{i}{\cbb{\bs_m\setminus \z{\omega }}{F}}\times \kk{i+n+1}{\ccb{\bs_m}{F}}^{2^n}\;.$$
By \pr{17.1'c} a),b),
$$\kk{i}{\cbb{(\bs_n\times \Omega)\setminus  \z{\omega_0 }}{F}}\approx \kk{i+m}{F}\times \kk{i+1+n}{F}^{2^n}\times \kk{i+1+m+n}{F}^{2^n}.$$

By \pr{26.2'a},
$$\kk{i}{\cbb{(\bs_m\times \Omega) \setminus \Gamma }{F}}\approx$$
$$\approx  \kk{i}{\cbb{(\bs_m\times \Omega) \setminus \z{\omega_0 }}{F}}\times \kk{i+1}{F}^{Card\,\Gamma -1}\approx $$
$$\approx \kk{i+m}{F}\times \kk{i+n+1}{F}^{2^n}\times \kk{i+m+n+1}{F}^{2^n}\times \kk{i+1}{F}^{Card\,\Gamma -1}\;.\qedd$$

\begin{lem}\label{12.5'}
Let $(k_j)_{j\in \bnn{n}}$ be a family in $\bn$, $n\not=1$, and $m:=1+\si{j\in \bnn{n}}k_j$. There is an embedding of $\pro{j\in \bnn{n}}\bs_{k_j}$ in $\bbb_m$ such that \mbox{$\bbb_m\setminus \pro{j\in \bnn{n}}\bs_{k_j}$} has two connected components: one is homeomorphic to $\br^{1+k_n}\times \pro{j\in \bnn{n-1}}\bs_{k_j}$ and the other is homeomorphic to $\bbb_m\setminus \left(\bbb_{1+k_n}\times \pro{j\in \bnn{n-1}}\bs_{k_j}\right)$. 
\end{lem}

We prove the assertion by induction with respect to $n\in \bn\setminus \z{1}$. Assume first $n=2$, put 
$$\Gamma :=\me{\alpha \in \bbb_m}{\n{\alpha}=\frac{1}{2}\,,\;\;\alpha _{2+k_1}=\alpha _{3+k_1}=\cdots=\alpha _m=0}\,,$$
and for every $\alpha \in \bbb_m$ denote by $d(\alpha )$ the distance of $\alpha $ to $\Gamma $. Then 
$$\me{\alpha \in \bbb_m}{d(\alpha )=\frac{1}{4}}$$
is an embedding of $\bs_{k_1}\times \bs_{k_2}$ in $\bbb_m$ with the desired properties.

Let now $n>2$ and assume the assertion holds for $n-1$. Let $\Gamma $ be a closed set of $\bbb_{m-k_n}$ homeomorphic to $\pro{j\in \bnn{n-1}}\bs_{k_j}$. We may assume $\Gamma \subset \bs_{m-k_m}$. We denote for every $\alpha \in \bbb_m$ by $d(\alpha )$ the distance of $\alpha $ to $\frac{1}{2}\Gamma $. Then $\me{\alpha \in \bbb_m}{d(\alpha )=\frac{1}{4}}$ is an embedding with the desired properties.\qed

\begin{p}\label{12.5'a}
Let $(k_j)_{j\in \bnn{n}}$ be a family in $\bn$.
\begin{enumerate}
\item $\qquad\qquad\qquad\qquad\proo{j=1}{n}\bs_{k_j}\in \Upsilon \,,\qquad\br_\Upsilon \subset \left(\proo{j=1}{n}\bs_{k_j}\right)_\Upsilon $,
$$K_i\left(\ccb{\proo{j=1}{n}\bs_{k_j}}{F}\right)\approx $$
$$\approx \ab{K_i(F)^{2^n}}{\emph{all}\;(k_j)_{j\in \bnn{n}}\;\emph{are\,\,even}}{\Big(K_i(F)\times  K_{i+1}(F)\Big)^{2^{n-1}}}{\emph{not\,\,all}\;(k_j)_{j\in \bnn{n}}\;\emph{are\,\,even}}\;.$$
\item If $\Gamma $  is a nonempty finite subset of $\pro{j\in \bnn{n}}\bs_{k_j}$ then
$$\kk{i}{\cbb{\pro{j\in \bnn{n}}\bs_{k_j}\setminus \Gamma }{F}}\approx $$
$$\approx \ab{\kk{i}{F}^{2^n-1}\times \kk{i+1}{F}^{Card\,\Gamma -1}}{\emph{all}\,\, k_j\,\emph{are even}}{\kk{i}{F}^{2^{n-1}-1}\times \kk{i+1}{F}^{2^{n-1}+Card\,\Gamma -2}}{\emph{not all}\,\, k_j\,\emph{are even}}\;.$$
\end{enumerate}
\end{p}

a) By \h{17.1'c} b), $\bs_{k_j}\in \Upsilon ,\,\br_\Upsilon \subset (\bs_{k_j})_\Upsilon $ for every $j\in J$ so by \pr{27.3'b} a),f), $\proo{j=1}{n}\bs_{k_j}\in \Upsilon \,,\br_\Upsilon \subset \left(\proo{j=1}{n}\bs_{k_j}\right)_\Upsilon $. By \h{17.1'c} b), with the notation of \pr{27.3'b} a),f),  
$$p_j=\frac{3+(-1)^{k_j}}{2}\,,\quad q_j=\frac{1-(-1)^{k_j}}{2}\,,\quad p_j+q_j=2\,,\quad p_j-q_j=1+(-1)^{k_j}\,,$$
$$p_J=\frac{1}{2}\left(2^n+\proo{j=1}{n}\left(1+(-1)^{k_j}\right)\right)\,,\quad q_J=\frac{1}{2}\left(2^n-\proo{j=1}{n}\left(1+(-1)^{k_j}\right)\right)\,,$$
and this implies the result.

b) Assume first $\Gamma =\z{\omega_0 }$ for some $\omega_0 \in \pro{j\in \bnn{n}}\bs_{k_j}$. We prove the assertion by induction with respect to $\bbn$. For $n=1$ this follows from \h{17.1'c} $e_1)$. Let $n\not=1$ and assume the assertion holds for $n-1$. By \lm{12.5'}, $\bbb_m\setminus \pro{j\in \bnn{n}}\bs_{k_j}$ is homeomorphic to the topological sum of $\br^{1+k_n}\times \pro{j\in \bnn{n-1}}\bs_{k_j}$ and $\bbb_m\setminus \left(\bbb_{1+k_n}\times \pro{j\in \bnn{n-i}}\bs_{k_j}\right)$. By \pr{20.5'a}, for $\omega \in \pro{j\in \bnn{n-1}}\bs_{k_j}$,
$$\kk{i}{\cbb{\pro{j\in \bnn{n}}\bs_{k_j}\setminus \z{\omega_0 }}{F}}\approx $$
$$\approx \kk{i}{\cbb{\pro{j\in \bnn{n-1}}\bs_{k_j}\setminus \z{\omega }}{F}}\times $$
$$\times \kk{i+1}{\cbb{(\bbb_m\setminus \bs_{k_n})\times \pro{j\in \bnn{n-1}}\bs_{k_j}}{F}}.$$
By a) and \h{17.1'c} g),
$$\kk{i+1}{\cbb{(\bbb_n\setminus \bs_{k_n})\times \pro{j\in \bnn{n-1}}\bs_{k_j}}{F}}\approx$$
$$\approx  \kk{i+k_n}{\ccb{\pro{j\in \bnn{n-1}}\bs_{k_j}}{F}}\approx $$
$$\approx \abb{\kk{i+k_n}{F}^{2^{n-1}}}{\mbox{all}\,\, (k_j)_{j\in \bnn{n-1}}\,\mbox{are even}}{\Big(\kk{i}{F}\times \kk{i+1}{F}\Big)^{2^{n-2} }}{\mbox{not all}\,\, (k_j)_{j\in \bnn{n-1}}\,\mbox{are even}}\;.$$
By the induction hypothesis,
$$\kk{i}{\cbb{\pro{j\in \bnn{n-1}}\bs_{k_j}\setminus \z{\omega }}{F}}\approx $$
$$\approx \abb{\kk{i}{F}^{2^{n-1}-1}}{\mbox{all}\,\, (k_j)_{j\in \bnn{n-1}}\,\mbox{are even}}{\kk{i}{F}^{2^{n-2}-1}\times \kk{i+1}{F}^{2^{n-2} }}{\mbox{not all}\,\, (k_j)_{j\in \bnn{n-1}}\,\mbox{are even}}$$
so
$$\kk{i}{\cbb{\pro{j\in \bnn{n}}\bs_{k_j}\setminus \z{\omega_0 }}{F}}\approx $$
$$\approx \abb{\kk{i}{F}^{2^n-1}}{\mbox{all}\,\, (k_j)_{j\in \bnn{n}}\,\mbox{are even}}{\kk{i}{F}^{2^{n-1}-1}\times \kk{i+1}{F}^{2^{n-1} }}{\mbox{not all}\,\, (k_j)_{j\in \bnn{n}}\,\mbox{are even}}\;.$$
This finishes the inductive proof. 

We prove now the general case and put $\Omega :=\pro{j\in \bnn{n}}\bs_{k_j}$. Since it is possible to find a closed set $\Delta $ of $\Omega $ such that $\Gamma \subset \Delta $ and $\Delta \setminus \z{\omega _0}$ is K-null, the assertion follows from \pr{26.2'a}.\qed

{\center{\section{Some morphisms}}}

\begin{p}\label{3.2'}
We put 
$$\mae{\vartheta }{\bbb_n}{\bbb_n}{(\alpha _j)_{j\in \bnn{n}}}{(\alpha _1,\cdots,\alpha _{n-1},-\alpha _n)}\,,$$
$$\mae{\vartheta' }{\br^n}{\br^n}{(\alpha _j)_{j\in \bnn{n}}}{(\alpha _1,\cdots,\alpha _{n-1},-\alpha _n)}\,,$$
$$\mae{\vartheta'' }{\bs_{n-1}}{\bs_{n-1}}{(\alpha _j)_{j\in \bnn{n}}}{(\alpha _1,\cdots,\alpha _{n-1},-\alpha _n)}\,,$$
$$\mae{\phi }{\ccb{\bbb_n}{F}}{\ccb{\bbb_n}{F}}{x}{x\circ \vartheta }\,,$$
$$\mae{\phi' }{\cbb{\br^n}{F}}{\cbb{\br^n}{F}}{x}{x\circ \vartheta' }\,,$$
$$\mae{\phi'' }{\ccb{\bs_{n-1}}{F}}{\ccb{\bs_{n-1}}{F}}{x}{x\circ \vartheta''}\;.$$
\begin{enumerate}
\item $\mae{K_i(\phi )}{K_i(\ccb{\bbb_n}{F})}{K_i(\ccb{\bbb_n}{F})}{a}{a}$.
\item $\mae{K_i(\phi' )}{K_i(\cbb{\br^n}{F})}{K_i(\cbb{\br^n}{F})}{b}{-b}$.
\item $$\mac{K_i(\phi'' )}{K_i(\ccb{\bs_{n-1}}{F})}{K_i(\ccb{\bs_{n-1}}{F})}\,,$$
$$(a,b)\longmapsto \ab{(b,a)}{n=1}{(a,-b)}{n>1}\,,$$
where we identified $K_i(\ccb{\bs_{n-1}}{F})$ with 
$$K_i(\ccb{\bbb_n}{F})\times K_{i+1}(\cbb{\bbb_{n}\setminus \bs_{n-1}}{F})$$
 using the group isomorphism of \emph{\cor{31.12} d)} if $n>1$.
\end{enumerate}
\end{p}

a) follows from the homotopy axiom (\axi{27.9'b}) since $\phi $ is homotopic to the identity map of $\ccb{\bbb_n}{F}$.

b) We identify $\br^n$ with the homeomorphic space $\bbb_n\setminus \bs_{n-1}$.

Assume first $n=1$. Put
$$\mae{\psi }{\ccb{\bbb_1}{F}}{\ccb{\z{-1,1}}{F}}{x}{x|\z{-1,1}}$$
and denote by $\mac{\varphi }{\cbb{]-1,1[}{F}}{\ccb{\bbb_1}{F}}$ the inclusion map and by $\delta _i$ the index maps associated to the exact sequence in \frm
$$\oc{\cbb{]-1,1[}{F}}{\varphi }{\ccb{\bbb_1}{F}}{\psi }{\ccb{\z{-1,1}}{F}}\;.$$
By \cor{20.12}, $K_i(\psi )a=(a,a)$ for every $a\in K_i(\ccb{\bbb_1}{F})$ so by the six-term axiom (\axi{27.9'd}),
$$\delta _i(a+b,a+b)=0\,,\qquad\qquad \delta _i(a,b)=-\delta _i(b,a)$$
for all $(a,b)\in K_i(\ccb{\z{-1,1}}{F})$. By the commutativity of the index maps (\axi{27.9'e}), $K_{i+1}(\phi ')\circ \delta _i=\delta _i\circ K_i(\phi '')$. For $(a,b)\in K_i(\ccb{\z{-1,1}}{F})$, by the above,
$$K_{i+1}(\phi ')\delta _i(a,b)=\delta _iK_i(\phi '')(a,b)=\delta _i(b,a)=-\delta _i(a,b)\;.$$
Since $\delta _i$ is surjective (because $\varphi $ factorizes through null and is therefore K-null), $K_i(\phi ')b=-b$ for all $b\in \kk{i}{\cbb{]-1,1[}{F}}$.

If $n>1$ then the assertion follows from the case $n=1$, since $\cbb{\br^n}{F}\approx \cbb{\br}{\cbb{\br^{n-1}}{F}}$

c) follows from a), b), and \cor{31.12} c).\qed

\begin{co}\label{4.2'}
If we put 
$$\mae{\vartheta }{\bbb_n}{\bbb_n}{\alpha }{-\alpha }\,,$$
$$\mae{\vartheta' }{\br^n}{\br^n}{\alpha }{-\alpha }\,,$$
$$\mae{\vartheta'' }{\bs_{n-1}}{\bs_{n-1}}{\alpha }{-\alpha }\,,$$
$$\mae{\phi }{\ccb{\bbb_n}{F}}{\ccb{\bbb_n}{F}}{x}{x\circ \vartheta }\,,$$
$$\mae{\phi' }{\cbb{\br^n}{F}}{\cbb{\br^n}{F}}{x}{x\circ \vartheta' }\,,$$
$$\mae{\phi'' }{\ccb{\bs_{n-1}}{F}}{\ccb{\bs_{n-1}}{F}}{x}{x\circ \vartheta''}$$
then
$$\mae{K_i(\phi )}{K_i(\ccb{\bbb_n}{F})}{K_i(\ccb{\bbb_n}{F})}{a}{a}\,,$$
$$\mae{K_i(\phi' )}{K_i(\cbb{\br^n}{F})}{K_i(\cbb{\br^n}{F})}{b}{(-1)^nb}\,,$$
$$\mac{K_i(\phi'' )}{K_i(\ccb{\bs_{n-1}}{F})}{K_i(\ccb{\bs_{n-1}}{F})}\,,$$
$$(a,b)\longmapsto \ab{(b,a)}{n=1}{(a,(-1)^{n+1}b)}{n>1}\,,$$
where we identified $K_i(\ccb{\bs_{n-1}}{F})$ with 
$$K_i(\ccb{\bbb_n}{F})\times K_{i+1}(\cbb{\bbb_{n}\setminus \bs_{n-1}}{F})$$
 using the group isomorphism of \emph{\cor{31.12} c)} if $n>1$.
\end{co}

The assertion for $K_i(\phi )$ follows from the homotopy axiom (\axi{27.9'b}) since $\phi $ is homotopic to the identity map of $\ccb{\bbb_n}{F}$. If $n$ is even then the same holds for $K_i(\phi ')$. Assume now $n$ odd and let us denote by $\bar{\phi }' $ the map denoted by $\phi '$ in \pr{3.2'}. Then $\phi '\circ \bar{\phi }' $ is homotopic to the identity map of $\cbb{\br^n}{F}$ so by \cor{3.2'}, for every $b\in K_i(\cbb{\br^n}{F})$,
$$K_i(\phi ')b=-K_i(\phi ')K_i(\bar{\phi }' )b=-b=(-1)^nb\;.$$
The assertion for $K_i(\phi '')$ follows from the corresponding assertions for $K_i(\phi )$ and $K_i(\phi ')$ and from \cor{31.12} d).\qed

\begin{p}\label{8.10}
Let $\alpha ,\beta \in [0,2\pi [$, $\alpha <\beta $, $\Omega :=\me{e^{i\omega }}{\omega \in ]\alpha ,\beta [}$, $\Gamma :=\bt\setminus \Omega $, $\mac{\varphi }{\cbb{\Omega }{F} }{\ccb{\bt}{F}}$ the inclusion map, and
$$\mae{\psi }{\ccb{\bt}{F}}{\ccb{\Gamma }{F}}{x}{x|\Gamma }\,,$$
$$\mae{\bar{\psi}}{\ccb{\Gamma }{F}}{F}{x}{x(1)}\,,$$
$$\mae{\vartheta }{]0,2\pi [}{]\alpha ,\beta [}{\omega }{\frac{\beta -\alpha }{2\pi }\omega +\alpha }\;.$$
For every $x\in \cbb{\Omega }{F}$ put
$$\mae{\tilde{x} }{\bt}{F}{e^{i\omega }}{\ab{x\left(e^{i\vartheta (\omega )}\right)}{\omega\in ]0,2\pi [}{0}{\omega\in \z{0,2\pi }}}$$
and define
$$\mae{\phi }{\cbb{\Omega }{F}}{\cbb{\bt\setminus \z{1}}{F}}{x}{\tilde{x} }\;.$$
\begin{enumerate}
\item $K_i(\phi )$ and $K_i(\bar{\psi})$ are group isomorphisms and so 
$$K_i(\cbb{\Omega }{F})\approx K_{i+1}(F)\,, \qquad K_i(\ccb{\Gamma }{F})\approx K_i(F)\;.$$
\item If we identify $K_i(\cbb{\Omega }{F})$ with $K_{i+1}(F)$ and $K_i(\ccb{\Gamma }{F})$ with $K_i(F)$ using the isomorphisms from a) and $K_i(\ccb{\bt}{F})$ with $K_i(F)\times K_{i+1}(F)$ using e.g. Alexandroff K-theorem \emph{(\h{20.4} a))} then
$$\mae{K_i(\varphi )}{K_i(\cbb{\Omega }{F})}{K_i(\ccb{\bt}{F})}{b}{(0,b)}\,,$$
$$\mae{K_i(\psi )}{K_i(\ccb{\bt}{F})}{K_i(\ccb{\Gamma }{F})}{(a,b)}{a}\;.$$ 
\end{enumerate}
\end{p} 

a) $\phi$ is an $E$-C*-isomorphism. Put
$$\mae{\tilde{\psi}}{F}{\ccb{\Gamma }{F}}{x}{1_{\ccb{\Gamma }{\bc}}x}\;.$$
Then $\obb{\ccb{\Gamma }{F}}{\bar{\psi}}{F}{\tilde{\psi}}{\ccb{\Gamma }{F}}$ is a homotopy in \frm so $K_i(\phi )$ and $K_i(\bar{\psi})$ are group isomorphisms by the homotopy axiom (\axi{27.9'b}). The
last assertion follows now from \h{17.1'c} a).

b) For every $s\in [0,1]$ put
$$\mae{\vartheta _s}{\bt}{\bt}{e^{i\omega }}{\acc{e^{is\omega }}{\omega \in [0,\alpha ]}{e^{is\omega }e^{\frac{2\pi i(1-s)(\omega -\alpha )}{\beta -\alpha }}}{\omega \in ]\alpha ,\beta [}{e^{is\omega }e^{2\pi i(1-s)}}{\omega \in [\beta ,2\pi ]}}\,,$$
$$\mae{\phi _s}{\ccb{\bt}{F}}{\ccb{\bt}{F}}{x}{x\circ \vartheta _s}\;.$$
Then $(\phi _s)_{s\in [0,1]}$ is a pointwise continuous path in $\ccb{\bt}{F}$ such that $\phi _1$ is the identity map. By the homotopy axiom (\axi{27.9'b}), $K_i(\phi _0)$ is the identity map of $K_i(\ccb{\bt}{F})$. Let
$$\mac{\varphi '}{\cbb{\bt\setminus \z{1}}{F}}{\ccb{\bt}{F}}$$
be the inclusion map and
$$\mae{\psi '}{\ccb{\bt}{F}}{F}{x}{x(1)}\;.$$
Then $\phi _0\circ \varphi =\varphi '\circ \phi $ and $\psi '\circ \phi _0=\bar{\psi }\circ \psi  $ so (by a)) for $a\in K_i(F)$ and $b\in K_{i+1}(F)$,
$$K_i(\varphi )b=K_i(\phi _0)K_i(\varphi )b=K_i(\varphi ')K_i(\phi )b=K_i(\varphi ')b=(0,b)\,,$$
$$K_i(\psi )(a,b)=K_i(\bar{\psi } )K_i(\psi )(a,b)=K_i(\psi ')K_i(\phi _0)(a,b)=K_i(\psi ')(a,b)=a$$
by the Alexandroff K-theorem (\h{20.4} a)).\qed

\begin{p}\label{19.10}
Put $\Gamma :=\me{e^{\frac{2\pi ij}{n}}}{j\in \bnn{n}}$ and  
$$\mae{\psi }{\ccb{\bt}{F}}{\ccb{\Gamma }{F}}{x}{x|\Gamma }\,,$$
and denote by 
$$\mac{\varphi }{\cbb{\bt\setminus \Gamma }{F}}{\ccb{\bt}{F}}$$
the inclusion map and by $\delta _i$ the index maps associated to the exact sequence in \frm
$$\oc{\cbb{\bt\setminus \Gamma }{F}}{\varphi }{\ccb{\bt}{F}}{\psi }{\ccb{\Gamma }{F}}\;.$$
\begin{enumerate}
\item $K_i(\cbb{\bt\setminus \Gamma }{F})\approx K_{i+1}(F)^n,\quad K_i(\ccb{\Gamma }{F})\approx K_i(F)^n$.
\item We identify the isomorphic groups of a) and identify $K_i(\ccb{\bt}{F})$ with $K_i(F)\times K_{i+1}(F)$ \emph{(\h{17.1'c} b))}.
$$\mae{K_i(\varphi )}{K_i(\cbb{\bt\setminus \Gamma }{F})}{K_i(\ccb{\bt}{F})}{(b_j)_{j\in \bnn{n}}}{\left(0,\si{j\in \bnn{n}}b_j\right)},$$
$$\mae{K_i(\psi )}{K_i(\ccb{\bt}{F})}{K_i(\ccb{\Gamma }{F})}{(a,b)}{(a)_{j\in \bnn{n}}}\;.$$
If $n=2$ and $\kk{i}{F}$ is isomorphic to $\bz$ or to $\bz_p$ for some $p\in \bn$ or to the group of rational numbers then there is an automorphism $$\mac{\Phi_i }{\kk{i}{F}}{\kk{i}{F}}$$
 such that
$$\mae{\delta _i}{\kk{i}{\ccb{\Gamma }{F}}}{\kk{i+1}{\cbb{\bt\setminus \Gamma }{F}}}{(a,b)}{(\Phi_i (a-b),\Phi (b-a))}.$$
\item If we put
$$\mae{\vartheta }{\bt\setminus \Gamma }{\bt\setminus \z{1}}{z}{z^n}\,,$$
$$\mae{\vartheta' }{\bt}{\bt}{z}{z^n}\,,$$
$$\mae{\vartheta'' }{\Gamma }{\z{1}}{z}{z^n}\,,$$
$$\mae{\phi }{\cbb{\bt\setminus \z{1}}{F}}{\cbb{\bt\setminus \Gamma }{F}}{x}{x\circ \vartheta }\,,$$
$$\mae{\phi' }{\ccb{\bt}{F}}{\ccb{\bt}{F}}{x}{x\circ \vartheta' }\,,$$
$$\mae{\phi'' }{\ccb{\z{1}}{F}}{\ccb{\Gamma }{F}}{x}{x\circ \vartheta'' }$$
then, with the identifications of a) and b),
$$\mae{K_i(\phi )}{K_i(\cbb{\bt\setminus \z{1}}{F})}{K_i(\cbb{\bt\setminus \Gamma }{F})}{b}{(b)_{j\in \bnn{n}}}\,,$$
$$\mae{K_i(\phi' )}{K_i(\ccb{\bt}{F})}{K_i(\ccb{\bt}{F})}{(a,b)}{(a,nb)}\,,$$
$$\mae{K_i(\phi'' )}{K_i(\ccb{\z{1}}{F})}{K_i(\ccb{\Gamma }{F})}{a}{(a)_{j\in \bnn{n}}}\;.$$
\end{enumerate}
\end{p}

a) Put $\Omega _j:=\me{e^{\frac{2\pi i\omega }{n}}}{\omega \in ]j-1,j[}$ for every $j\in \bnn{n}$. By \pr{8.10} a), for every $j\in \bnn{n}$,
$$K_i(\cbb{\Omega _j}{F})\approx K_{i+1}(F)\;.$$
so
$$K_i(\cbb{\bt\setminus \Gamma }{F})\approx K_{i+1}(F)^n\,,\qquad K_i(\ccb{\Gamma }{F})\approx K_i(F)^n$$
by the Product Theorem (\pr{14.11} a)).

b) By \cor{20.12},
$$\mae{K_i(\psi )}{K_i(\ccb{\bt}{F})}{K_i(\ccb{\Gamma }{F})}{(a,b)}{(a)_{j\in \bnn{n}}}\;.$$
If we denote by
$$\mac{\varphi _j}{\cbb{\Omega _j}{F}}{\ccb{\bt}{F}}$$
the inclusion map then
$$\mae{K_i(\varphi _j)}{K_i(\cbb{\Omega _j}{F})}{K_i(\ccb{\bt}{F})}{b}{(0,b)}$$
by \pr{8.10} b). By \pr{8.10} a) and \cor{30.1'},
$$\mae{K_i(\varphi )}{K_i(\cbb{\bt\setminus \Gamma }{F})}{K_i(\ccb{\bt}{F})}{(b_j)_{j\in \bnn{n}}}
{\left(0,\si{j\in \bnn{n}}b_j\right)}\;.$$

In order to prove the last assertion we define $a',b',a'',b''\in \kk{i}{F}$ by
$$(a',b'):=\delta _i(1,0)\,,\qquad\qquad (a'',b''):=\delta _i(0,1)\;.$$
From
$$0=\delta _i(1,1)=(a',b')+(a'',b'')=(a'+a'',b'+b'')$$
we get $a''=-a'$ and $b''=-b'$. There are $j,k\in \bz$ such that $\delta _i(j,k)=(1,-1)$. Then
$$(1,-1)=\delta _i(j,k)=(ja',jb')-(ka',kb')=((j-k)a',(j-k)b')\,,$$
$$ (j-k)a'=1\,,\qquad\qquad (j-k)b'=-1\;.$$
Thus $a'$ is invertible in the ring $\kk{i}{F}$ and $a'^{-1}=j-k$. It follows $b'=-a'$. If we put
$$\mae{\Phi_i }{\kk{i}{F}}{\kk{i}{F}}{c}{a'c}$$
then $\Phi_i $ is an automorphism and for all $a,b\in \kk{i}{F}$,
$$\delta _i(a,b)=(a'a,-a'a)-(a'b,-a'b)=(a'(a-b),a'(b-a))=(\Phi_i (a-b),\Phi_i (b-a))\;.$$

c) The assertions for $K_i(\phi )$ and $K_i(\phi '')$ follow from the Product Theorem (\pr{14.11} a)). If $\mac{\varphi '}{\cbb{\bt\setminus \z{1}}{F}}{\ccb{\bt}{F}}$ denotes the inclusion map and
$$\mae{\psi '}{\ccb{\bt}{F}}{\ccb{\z{1}}{F}}{x}{x|\z{1}}$$ then the diagram
$$\begin{CD}
K_i(\cbb{\bt\setminus \z{1}}{F})@>K_i(\varphi ')>>K_i(\ccb{\bt}{F})@>K_i(\psi ')>>K_i(\ccb{\z{1}}{F})\\
@VK_i(\phi )VV         @VVK_i(\phi ')V@VVK_i(\phi '')V\\
K_i(\cbb{\bt\setminus\Gamma }{F})@>>K_i(\varphi )>K_i(\ccb{\bt}{F})@>>K_i(\psi)>K_i(\ccb{\Gamma }{F})
\end{CD}$$
is commutative. Let $(a,b)\in K_i(\ccb{\bt}{F})$ and put $(a',b'):=K_i(\phi ')(a,0)$. By b),
$$(a)_{j\in \bnn{n}}=K_i(\phi '')a=K_i(\phi '')K_i(\psi ')(a,0)=$$
$$=K_i(\psi )K_i(\phi ')(a,0)=
K_i(\psi )(a',b')=(a')_{j\in \bnn{n}}\,,$$
$$K_i(\phi ')(0,b)=K_i(\phi ')K_i(\varphi ')b=K_i(\varphi )K_i(\phi )b=K_i(\varphi )(b)_{j\in \bnn{n}}=(0,nb)$$
so $K_i(\phi ')(a,b)=(a,nb)$.\qed 

\begin{co}\label{31.12a}
If we put
$$\mae{\vartheta }{\bbb_2}{\bbb_2}{z}{z^n}\,,$$
$$\mae{\vartheta' }{\bc}{\bc}{z}{z^n}\,,$$
$$\mae{\vartheta'' }{\bs_1}{\bs_1}{z}{z^n}\,,$$
$$\mae{\phi }{\ccb{\bbb_2}{F}}{\ccb{\bbb_2}{F}}{x}{x\circ \vartheta }\,,$$
$$\mae{\phi' }{\cbb{\bc}{F}}{\cbb{\bc}{F}}{x}{x\circ \vartheta '}\,,$$
$$\mae{\phi ''}{\ccb{\bs_1}{F}}{\cbb{\bs_1}{F}}{x}{x\circ \vartheta ''}\;.$$
then $K_i(\phi )$ is the identity map of $K_i(\ccb{\bbb_2}{F})$ and
$$\mae{K_i(\phi ')}{K_i(\cbb{\bc}{F})}{K_i(\cbb{\bc}{F})}{a}{na}\,,$$
$$\mae{K_i(\phi '')}{K_i(\ccb{\bs_1}{F})}{K_i(\ccb{\bs_1}{F})}{(a,b)}{(a,nb)}\,,$$
\end{co}

We identify the homeomorphic spaces $\bc$ and $\bbb_2\setminus \bs_1$. By \cor{31.12} c), 
$$K_i(\ccb{\bs_1}{F})\approx K_i(\ccb{\bbb_2}{F})\times K_{i+1}(\cbb{\bbb_2\setminus \bs_1}{F})$$
and by \pr{19.10} e),
$$\mae{K_i(\phi '')}{K_i(\ccb{\bs_1}{F})}{K_i(\ccb{\bs_1}{F})}{(a,b)}{(a,nb)}\;.$$
By \cor{30.12} b) and \h{17.1'b} a), $K_i(\phi )$ is the identity map of $K_i(\ccb{\bbb_2}{F})$ and
$$\mae{K_i(\phi ')}{K_i(\cbb{\bbb_2\setminus \bs_1}{F})}{K_i(\cbb{\bbb_2\setminus \bs_1}{F})}{a}{na}\;.\qedd$$

\begin{p}\label{4.1'}
Let $m,n\in \bn$ and
$$\mae{\vartheta _1}{\bt}{\bt}{w}{w^m}\,,$$
$$\mae{\vartheta _2}{\bt}{\bt}{z}{z^n}\,,$$
$$\mae{\psi }{\ccb{\bt\times \bt}{F}}{\ccb{\bt\times \bt}{F}}{x}{x\circ (\vartheta _1\times \vartheta _2)}\;.$$
We identify $K_i(\ccb{\bt}{F'})$ with $K_i(F')\times K_{i+1}(F')$ for all $E$-C*-algebras $F'$ by using the group isomorphism of \emph{\h{17.1'c} b)}. Let 
$$a\in K_i(\ccb{\bt\times \bt}{F})\approx K_i(\ccb{\bt}{\ccb{\bt}{F}})\approx K_i(\ccb{\bt}{F})\times K_{i+1}(\ccb{\bt}{F})$$
and put $a_0\in K_i(\ccb{\bt}{F})$, $a_1\in K_{i+1}(\ccb{\bt}{F})$ such that $a=(a_0,a_1)$ and $a_{0,0},\,a_{1,1}\in K_i(F)$ and $a_{0,1},\,a_{1,0}\in K_{i+1}(F)$ such that $a_0=(a_{0,0},\,a_{0,1})$ and $a_1=(a_{1,0},\,a_{1,1})$. Then
$$K_i(\psi )=((a_{0,0},\,mna_{1,1}),\,(na_{0,1},\,ma_{1,0}))\;.$$
\end{p}

We put 
$$\mae{\phi }{\ccb{\bt}{F}}{\ccb{\bt}{F}}{x}{x\circ \vartheta _2}\,,$$
$$\mae{\phi _1}{\ccb{\bt}{\ccb{\bt}{F}}}{\ccb{\bt}{\ccb{\bt}{F}}}{x}{x\circ \vartheta _1}\,,$$
$$\mae{\phi _2}{\ccb{\bt}{\ccb{\bt}{F}}}{\ccb{\bt}{\ccb{\bt}{F}}}{x}{\phi \circ x}\,,$$
By \cor{31.12a},
$$K_i(\phi _1)a=(a_0,\,ma_1),\quad K_i(\phi )a_0=(a_{0,0},\,na_{0,1}),\quad K_{i+1}(\phi )a_1=(a_{1,0},\,na_{1,1})$$
so by \cor{31.12} c), d),
$$K_i(\phi _2)K_i(\phi _1)a=(K_i(\phi )a_0,\,K_{i+1}(\phi )ma_1)=((a_{0,0},\,na_{0,1}),\,(ma_{1,0},\,mna_{1,1}))\;.$$
Since $\psi =\phi _2\circ \phi _1$,
$$K_i(\psi )=((a_{0,0},\,mna_{1,1}),\,(na_{0,1},\,ma_{1,0}))\;.\qedd$$

{\center{\section{Some non-orientable compact spaces}}}

\begin{de}\label{4.2'a}
We denote by $\bp_n$ the $n$-dimensional {\bf{projective space}}, which is obtained from $\bbb_n$ by identifying $\alpha $ with $-\alpha $ for all $\alpha \in \bbb_n$ with $\n{\alpha }=1$.
\end{de}

\begin{p}\label{4.2'b}
Put
$$\Omega :=\bp_{n+1}\setminus \me{\alpha \in \bp_{n+1}}{\n{\alpha }=1,\,\alpha _{n+1}=0}\,,$$
$$\Omega ':=\me{\alpha \in \Omega }{\n{\alpha }=1}\,,$$
$$\mae{\psi }{\cbb{\Omega }{F}}{\cbb{\Omega '}{F}}{x}{x|\Omega '}$$
and denote by
$$\mac{\varphi }{\cbb{\bbb_{n+1}\setminus \bs_n}{F}}{\cbb{\Omega }{F}}$$
the inclusion map and by $\delta _i$ the index maps associated to the exact sequence in \frm
$$\oc{\cbb{\bbb_{n+1}\setminus \bs_n}{F}}{\varphi }{\cbb{\Omega }{F}}{\psi }{\cbb{\Omega '}{F}}\;.$$
\begin{enumerate}
\item $\kk{i}{\cbb{\bbb_{n+1}\setminus \bs_n}{F}}\approx \kk{i+n+1}{F},\, \kk{i}{\cbb{\Omega '}{F}}\approx \kk{i+n}{F}$, and there is an automorphism $\mac{\Phi _i}{\kk{i+n}{F}}{\kk{i+n}{F}}$ such that
$$\mac{\delta _i}{\kk{i}{\cbb{\Omega '}{F}}}{\kk{i+1}{\cbb{\bbb_{n+1}\setminus \bs_n}{F}}},$$
$$ a\longmapsto \Phi _i\left(a-(-1)^na\right)\;.$$
\item If $n$ is even then $\delta _i=0$, $\kk{i}{\varphi }$ is injective, $\kk{i}{\psi }$ is surjective, and
$$\frac{\kk{i}{\cbb{\Omega }{F}}}{\kk{i+1}{F}}\approx \kk{i}{F}\;.$$
\item If $n$ is odd and for a fixed $i\in \z{0,1}$
$$a\in \kk{i}{F},\,2a=0\Longrightarrow a=0$$
then $\kk{i}{\psi }=0$, $\kk{i}{\cbb{\Omega }{F}}\approx \frac{\kk{i}{F}}{2\kk{i}{F}}$,
$$\mac{\kk{i}{\varphi }}{\kk{i}{\cbb{\bbb_{n+1}\setminus \bs_n}{F}}}{\kk{i}{\cbb{\Omega }{F}}}$$
 is the quotient map, and
$$\mae{\delta _i}{\kk{i}{\cbb{\Omega '}{F}}}{\kk{i+1}{\cbb{\bbb_{n+1}\setminus \bs_n}{F}}}{a}{2\Phi _ia}\;.$$
\end{enumerate}
\end{p}

a) By \h{17.1'c} a), $K_i(\cbb{\br^n}{F})\approx K_{i+n}(F)$. Since $\bbb_{n+1}\setminus \bs_n$ is homeomorphic to $\br^{n+1}$, $K_{i}(\cbb{\bbb_{n+1}\setminus \bs_n}{F})\approx K_{i+n+1}(F)$. Since $\Omega '$ is homeomorphic to  $\br^n$, $K_i(\cbb{\Omega '}{F})\approx K_{i+n}(F)$. We use the notation of \pr{31.1'}, which we mark by a bar in order to distinguish it from the present notation. Moreover we denote by $\mac{\vartheta }{\bar{\Omega } }{\Omega }$ and $\mac{\vartheta' }{\bar{\Omega }' }{\Omega' }$ the covering maps and put
$$\mae{\phi }{\cbb{\Omega }{F}}{\cbb{\bar{\Omega } }{F}}{x}{x\circ \vartheta }\,,$$
$$\mae{\phi' }{\cbb{\Omega' }{F}}{\cbb{\bar{\Omega }' }{F}}{x}{x\circ \vartheta' }\;.$$
By the Product Theorem (\pr{14.11} a)), \pr{31.1'} a), and \pr{3.2'} b),
$$\mae{K_i(\phi ')}{K_i\left(\cbb{\Omega '}{F}\right)}{K_i\left(\cbb{\bar{\Omega }' }{F}\right)}{a}{(a,(-1)^na)}\;.$$
By the commutativity of the index maps (\axi{27.9'e}), $\delta _i=\bar{\delta }_i\circ K_i(\phi ') $ so by \pr{31.1'} b),
$$\mae{\delta _i}{K_i\left(\cbb{\Omega '}{F}\right)}{K_{i+1}(\cbb{\bbb_{n+1}\setminus \bs_n}{F})}{a}{\Phi_i (a-(-1)^na)}\;.$$

b) and c) follow from a) and the six-term axiom (\axi{27.9'd}).\qed

\begin{co}\label{21.6'}
We use the notation and the hypothesis of {\emph{\pr{4.2'b}}}, take $n=1$, put $\Gamma :=\me{x\in \bp_2}{\n{x}=1}$,
$$\mae{\psi '}{\ccb{\bp_2}{F}}{\ccb{\Gamma }{F}}{x}{x|\Gamma }\,,$$
and denote by $\mac{\varphi '}{\cbb{\bbb_2\setminus \bs_1}{F}}{\ccb{\bp_2}{F}}$ the inclusion map and by $\delta _i'$ the index maps associated to the exact sequence in \frm
$$\oc{\cbb{\bbb_2\setminus \bs_1}{F}}{\varphi '}{\ccb{\bp_2}{F}}{\psi '}{\ccb{\Gamma }{F}}\;.$$
Then $\kk{i}{\ccb{\bp_2}{F}}\approx \kk{i}{F}\times \frac{\kk{i}{F}}{2\kk{i}{F}}$, $\kk{i}{\ccb{\Gamma }{F}}\approx \kk{i}{F}\times \kk{i+1}{F}$,
$$\mae{\kk{i}{\varphi '}}{\kk{i}{\cbb{\bbb_2\setminus \bs_1}{F}}}{\kk{i}{\kk{i}{\ccb{\bp}{F}}}}{a}{(0,\Phi _ia)}\,,$$
$$\mae{\kk{i}{\psi '}}{\kk{i}{\ccb{\bp_2}{F}}}{\kk{i}{\ccb{\Gamma }{F}}}{(a,c)}{(a,0)}\,,$$
$$\mae{\delta'_i }{\kk{i}{\ccb{\Gamma }{F}}}{\kk{i+1}{\cbb{\bbb_2\setminus \bs_1}{F}}}{(a,b)}{2b}\;.\qedd$$
\end{co}

\begin{p}\label{24.9}
Let 
$$\mae{\vartheta }{[0,1]}{\bt}{\omega }{e^{2\pi i\omega }}\,,$$
$$\mae{\phi }{\ccb{\bt}{F}}{\ccb{[0,1]}{F}}{x}{x\circ \vartheta }\;.$$
If we identify $K_i(\ccb{\bt}{F})$ with $K_i(F)\times K_{i+1}(F)$ \emph{(\h{17.1'c} b))} and $K_i(\ccb{[0,1]}{F})$ with $K_i(F)$ \emph{(\h{17.1'b} a))} then
$$\mae{K_i(\phi )}{K_i(\ccb{\bt}{F})}{K_i(\ccb{[0,1]}{F})}{(a,b)}{a}\;.$$
\end{p}

Put
$$\mae{\vartheta '}{]0,1[}{\bt\setminus \z{1}}{\omega }{e^{2\pi i\omega }}\,,$$
$$\mae{\phi '}{\cbb{\bt\setminus \z{1}}{F}}{\cbb{]0,1[}{F}}{x}{x\circ \vartheta '}$$
and denote by 
$$\mac{\varphi }{\cbb{]0,1[}{F}}{\ccb{[0,1]}{F}}\,, \quad \mac{\varphi '}{\cbb{\bt\setminus \z{1}}{F}}{\ccb{\bt}{F}}$$
the inclusion maps. Then $\phi \circ \varphi '=\varphi\circ \phi ' $, so $K_i(\phi )\circ K_i(\varphi ')=K_i(\varphi )\circ K_i(\phi ')=0$, since $\varphi $ factorizes through $0$. Thus $K_i(\phi )(0,b)=0$ for all $b\in K_{i+1}(F)$.

Put
$$\mae{\psi }{\ccb{[0,1]}{F}}{\ccb{\z{0,1}}{F}\approx F\times F}{x}{x|\z{0,1}}\,,$$
$$\mae{\psi' }{\ccb{\bt}{F}}{F}{x}{x(1)}\,,$$
$$\mae{\mu }{F}{\ccb{\z{0,1}}{F}}{x}{(x,x)}\;.$$
Then $\psi \circ \phi=\mu \circ \psi ' $, so $K_i(\psi )\circ K_i(\phi )=K_i(\mu )\circ K_i(\psi ')$ and we get (by the above)
$$K_i(\psi )K_i(\phi )(a,b)=K_i(\psi )K_i(\phi )(a,0)=K_i(\mu )K_i(\psi ')(a,0)=K_i(\mu )a=(a,a)\,,$$
$$K_i(\phi )(a,b)=a$$
for all $(a,b)\in K_i(F)\times K_{i+1}(F)$.\qed

\begin{de}\label{27.1'}
We denote by $\bm$ the {\bf{M$\ddot{o}$bius band}} obtained from $[0,1]\times [-1,1]$ by identifying the points $(0,\beta )$ and $(1,-\beta )$ for every $\beta \in [-1,1]$. We put for every $j\in \z{-1,0,1}$ 
$$\Gamma^{\bm} _j:=\me{(\alpha ,j)\in \bm}{\alpha \in [0,1]}\;.$$
\end{de}

\begin{p}\label{27.1'a}
For every $j\in \z{-1,0,1}$ put
$$\mae{\psi _j}{\ccb{\bm}{F}}{\ccb{\Gamma^{\bm} _j}{F}}{x}{x|\Gamma^{\bm} _j}\;.$$
\begin{enumerate}
\item $\Gamma^{\bm} _0$ is homeomorphic to $\bt$ and $\Gamma^{\bm} _j$ is homeomorphic to $[0,1]$ for all $j\in \z{-1,1}$.
\item $\cbb{\bm\setminus \Gamma^{\bm} _0}{F}$ is $K$-null and
$$\mac{K_i(\psi _0)}{K_i(\ccb{\bm}{F})}{K_i\left(\ccb{\Gamma^{\bm} _0}{F}\right)\approx K_i(F)\times K_{i+1}(F)}$$
is a group isomorphism.
\item If we identify $K_i(\ccb{\bm}{F})$ with $K_i(F)\times K_{i+1}(F)$ using the group isomorphism $K_i(\psi _0)$ of b)  and $K_i\left(\ccb{\Gamma^{\bm} _1}{F}\right)$ with $K_i(F)$ using a) \emph{(and \h{17.1'b} a))} then
$$\mae{K_i(\psi _1)}{K_i(\ccb{\bm}{F})}
{K_i\left(\ccb{\Gamma^{\bm} _1}{F}\right)}{(a,b)}{a}\;.$$ 
\item If we put $\omega :=(0,0)=(1,0)\in \bm$, $\Gamma :=\me{(\alpha ,0)}{\alpha \in ]0,1[}$, and
$$\mae{\psi }{\cbb{\bm\setminus \z{\omega }}{F}}{\cbb{\Gamma }{F}}{x}{x|\Gamma }$$
then
$$\mac{K_i(\psi )}{K_i(\cbb{\bm\setminus \z{\omega }}{F})}{K_i(\cbb{\Gamma }{F})\approx K_{i+1}(F)}$$
is a group isomorphism.
\item If $\Gamma '$ is a finite subset of $\bm$ then
$$K_i\left(\cbb{\bm\setminus \Gamma '}{F}\right)\approx K_{i+1}(F)^{\Gamma '}$$
\end{enumerate}
\end{p}

a) is easy to see.

b) For every $s\in ]0,1]$ put
$$\mae{\vartheta _s}{\bm\setminus \Gamma^{\bm} _0}{\bm\setminus \Gamma^{\bm} _0}{(\alpha ,\beta )}{(\alpha ,s\beta )}\;.$$
By \pr{24.11a} (replacing there $\Omega $ by $\bm\setminus \Gamma _0^{\bm}$), $\cbb{\bm\setminus \Gamma^{\bm} _0}{F}$ is K-null and the assertion follows from the Topological six-term sequence (\pr{24.11} a)) and a) (and \h{17.1'c} b)). 

c) follows from b) and \pr{24.9}.

d) If $\mac{\varphi }{\cbb{\bm\setminus \Gamma^{\bm} _0}{F}}{\cbb{\bm\setminus \z{\omega }}{F}}$ denotes the inclusion map then
$$\oc{\cbb{\bm\setminus \Gamma^{\bm} _0}{F}}{\varphi }{\cbb{\bm\setminus \z{\omega }}{F}}{\psi }{\cbb{\Gamma }{F}}$$
is an exact sequence in \frm. By b), $\cbb{\bm\setminus \Gamma^{\bm} _0}{F}$ is K-null so by the Topological six-term sequence (\pr{24.11} a)), $K_i(\psi )$ is a group isomorphism. Since $\Gamma $ is homeomorphic to $\br$, $K_i(\cbb{\Gamma }{F})\approx K_{i+1}(F)$ by \h{17.1'c} a).

e) follows from d) and \pr{26.2'a}.\qed

\begin{p}\label{28.1'}
Put 
$$\Gamma':=\Gamma^{\bm} _0\cup \Gamma^{\bm} _1\,,\qquad \Gamma'':=\Gamma^{\bm} _0\cup \Gamma^{\bm} _1\cup \Gamma^{\bm} _{-1}\,,\qquad \Gamma ''':=\Gamma^{\bm} _1\cup \Gamma^{\bm} _{-1}\,,$$ 
$$\bm':=\bm\setminus \Gamma'\,,\qquad \bm'':=\bm\setminus \Gamma''\qquad \bm''':=\bm\setminus \Gamma '''\;.$$ 
Let 
$$\mac{\varphi'}{\cbb{\bm'}{F}}{\ccb{\bm}{F}}\,,$$
$$\mac{\varphi''}{\cbb{\bm''}{F}}{\ccb{\bm}{F}}\,,$$
$$\mac{\bar{\varphi }'}{\cbb{\bm'}{F}}{\cbb{\bm\setminus \Gamma^{\bm} _0}{F}}\,,$$
$$\mac{\bar{\varphi }''}{\cbb{\bm''}{F}}{\cbb{\bm\setminus \Gamma^{\bm} _0}{F}}\,,$$
$$\mac{\varphi'''}{\cbb{\bm''}{F}}{\cbb{\bm'''}{F}}\,,$$
$$\mac{\lambda '}{\ccb{\Gamma^{\bm} _1}{F}}{\ccb{\Gamma '}{F}}\,,$$
$$\mac{\lambda ''}{\ccb{\Gamma'''}{F}}{\ccb{\Gamma ''}{F}}$$
$$\mac{\lambda'''}{\ccb{\Gamma^{\bm} _0}{F}}{\ccb{\Gamma ''}{F}}\,,$$
be the inclusion maps,
$$\mae{\psi'}{\ccb{\bm}{F}}{\ccb{\Gamma'}{F}}{x}{x|\Gamma '}\,,$$
$$\mae{\psi''}{\ccb{\bm}{F}}{\ccb{\Gamma''}{F}}{x}{x|\Gamma''}\,,$$
$$\mae{\bar{\psi}'}{\cbb{\bm\setminus \Gamma^{\bm} _0}{F}}{\ccb{\Gamma^{\bm} _1}{F}}{x}{x|\Gamma^{\bm} _1}\,,$$
$$\mae{\bar{\psi}''}{\cbb{\bm\setminus \Gamma^{\bm} _0}{F}}{\ccb{\Gamma'''}{F}}{x}{x|\Gamma '''}\,,$$
$$\mae{\psi'''}{\cbb{\bm'''}{F}}{\ccb{\Gamma^{\bm}_0}{F}}{x}{x|\Gamma^{\bm}_0}\,,$$
and $\delta'_i\,,\delta''_i\,,\bar{\delta}'_i\,,\bar{\delta }''_i\,,\delta'''_i $ the index maps associated to the exact sequences in \frm
$$\oc{\cbb{\bm'}{F}}{\varphi'}{\ccb{\bm}{F}}{\psi'}{\ccb{\Gamma'}{F}}\,,$$
$$\oc{\cbb{\bm''}{F}}{\varphi''}{\ccb{\bm}{F}}{\psi''}{\ccb{\Gamma''}{F}}\,,$$
$$\oc{\cbb{\bm'}{F}}{\bar{\varphi}'}{\ccb{\bm\setminus \Gamma^{\bm} _0}{F}}{\bar{\psi}'}{\ccb{\Gamma^{\bm} _1}{F}}\,,$$
$$\oc{\cbb{\bm''}{F}}{\bar{\varphi}''}{\ccb{\bm\setminus \Gamma^{\bm} _0}{F}}{\bar{\psi}''}{\ccb{\Gamma'''}{F}}\,,$$
$$\oc{\cbb{\bm''}{F}}{\varphi'''}{\cbb{\bm'''}{F}}{\psi'''}{\ccb{\Gamma^{\bm}_0}{F}}\,,$$
respectively.
\begin{enumerate}
\item $\Gamma'''$ is homeomorphic to $\bt$.
\item The maps
$$\mac{\bar{\delta }'_i}{K_i\left(\ccb{\Gamma^{\bm} _1}{F}\right)\approx K_i(F)}{K_{i+1}\left(\cbb{\bm'}{F}\right)}\,,$$
$$\mac{\bar{\delta }''_i}{K_i\left(\ccb{\Gamma'''}{F}\right)\approx K_i(F)\times K_{i+1}(F)}{K_{i+1}\left(\cbb{\bm''}{F}\right)}$$
are group isomorphisms.
\item If we put $\Phi '_i:=K_i(\lambda ')\circ (\bar{\delta }'_i)^{-1}$, $\Phi ''_i:=K_i(\lambda '')\circ (\bar{\delta }''_i)^{-1}$ (using b)) then the sequences
$$\og{K_i(\ccb{\bm}{F})}{K_i(\psi ')}{K_i\left(\ccb{\Gamma '}{F}\right)}{\delta '_i}{\Phi '_i}{K_{i+1}\left(\cbb{\bm'}{F}\right)}{10}{0}{0}\,,$$
$$\og{K_i(\ccb{\bm}{F})}{K_i(\psi '')}{K_i\left(\ccb{\Gamma ''}{F}\right)}{\delta ''_i}{\Phi ''_i}{K_{i+1}\left(\cbb{\bm''}{F}\right)}{10}{0}{0}$$
are split exact and the maps
$$K_i(\ccb{\bm}{F})\times K_{i+1}\left(\cbb{\bm'}{F}\right)\longrightarrow K_i\left(\ccb{\Gamma '}{F}\right),$$
$$(a,b)\longmapsto K_i(\psi ')a+\Phi '_ib\,,$$
$$K_i(\ccb{\bm}{F})\times K_{i+1}\left(\cbb{\bm''}{F}\right)\longrightarrow K_i\left(\ccb{\Gamma ''}{F}\right),$$
$$(a,b)\longmapsto K_i(\psi '')a+\Phi ''_ib$$
are group isomorphisms.
\item $\delta'''_i=0$ and the sequence
$$\occ{K_i\left(\cbb{\bm''}{F}\right)}{K_i(\varphi''' )}{K_i\left(\cbb{\bm'''}{F}\right)}$$
$$\ocd{K_i\left(\cbb{\bm'''}{F}\right)}{K_i(\psi''')}{K_i\left(\ccb{\Gamma^{\bm} _0}{F}\right)}$$
is exact.
\end{enumerate}
\end{p}

a) is easy to see.

b) By \pr{27.1'a} b), $\cbb{\bm\setminus \Gamma^{\bm} _0}{F}$ is K-null and the assertion follows from a), the Topological six-term sequence (\pr{24.11} b)), and \pr{27.1'a} a) (and \h{17.1'b} a), \h{17.1'c} b)).

c) If we put $\Omega _1:=\bm$, $\Omega _2:=\bm\setminus \Gamma^{\bm} _0$, and $\Omega _3:=\bm'$ (respectively $\Omega _3:=\bm''$) then the assertion follows from the Topological triple (\pr{3.12} a)).

d) By the commutativity of the index maps (\axi{27.9'e}), $\delta_i'''=\delta ''_i\circ K_i(\lambda''')$. By c), $Im\,(\Phi ''_i\circ \delta ''_i)\subset Im\,K_i(\lambda '')$. Since $Im\,K_i(\lambda''')=K_i\left(\ccb{\Gamma^{\bm} _0}{F}\right)$ we get
$$\Phi ''_i\circ \delta'''_i=\Phi ''_i\circ \delta ''_i\circ K_i(\lambda''')=0\;.$$
Thus $\delta'''_i=\delta'' _i\circ \Phi ''_i\circ \delta'''_i=0$
and the assertion follows from the six-term axiom (\axi{27.9'd}).\qed

\begin{de}\label{28.1'a}
We denote by $\bk$ the {\bf{Klein bottle}} obtained from the M$\ddot{o}$bius band $\bm$ by identifying the points $(\alpha ,-1)$ and $(\alpha ,1)$ for all $\alpha \in [0,1]$ and put for every $j\in \z{0,1}$
$$\Gamma _j^{\bk}:=\me{(\alpha ,j)\in \bk}{\alpha \in [0,1]}\;.$$
\end{de}

\begin{p}\label{28.1'b}
We put $\bk':=\bk\setminus \Gamma _0^{\bk}$, $\bk'':=\bk\setminus (\Gamma _0^{\bk}\cup \Gamma _1^{\bk})$,
$$\mae{\psi }{\cbb{\bk'}{F}}{\ccb{\Gamma _1^{\bk}}{F}}{x}{x|\Gamma _1^{\bk}}$$
and denote by $\mac{\varphi }{\cbb{\bk''}{F}}{\cbb{\bk'}{F}}$ the inclusion map and by $\delta _i$ the index maps associated to the exact sequence in \frm
$$\oc{\cbb{\bk''}{F}}{\varphi }{\cbb{\bk'}{F}}{\psi }{\ccb{\Gamma _1^{\bk}}{F}}\;.$$
We use the notation of \emph{\pr{28.1'}} (so $\Gamma _0^{\bk}=\Gamma _0^{\bm}$ and $\bk''=\bm''$).
\begin{enumerate}
\item $\Gamma _0^{\bk}$ and $\Gamma _1^{\bk}$ are homeomorphic to $\bt$.
\item The map
$$\mac{(\bar{\delta }_{i+1}'' )^{-1}}{K_i\left(\cbb{\bm''}{F}\right)}{K_{i+1}\left(\ccb{\Gamma '''}{F}\right)\approx K_i(F)\times K_{i+1}(F)}$$
is a group isomorphism.
 \item If we identify $K_i(\ccb{\Gamma^{\bk} _1}{F})$ with $K_i(F)\times K_{i+1}(F)$ using a) and \emph{\h{17.1'c} b)} and $K_{i+1}(\cbb{\bk''}{F})$ with $K_i(F)\times K_{i+1}(F)$  using b) then
$$\mae{\delta _i}{K_i\left(\ccb{\Gamma^{\bk} _1}{F}\right)}{K_{i+1}\left(\cbb{\bk''}{F}\right)}{(a,b)}{(a,2b)}\;.$$
\item If $\delta _i$ is injective then $\psi $ is K-null and
 $K_i(\cbb{\bk'}{F})\approx\frac{K_i(F)}{2K_i(F)} $ and if we denote by 
$$\mac{\Phi _i}{K_i(F)}{\frac{K_i(F)}{2K_i(F)}}$$
the quotient map then
$$\mae{K_i\left(\varphi \right)}{K_i\left(\cbb{\bk''}{F}\right)}{K_i\left(\cbb{\bk'}{F}\right)}{(a,b)}{\Phi _ib}\;.$$
\end{enumerate}
\end{p}

a) is easy to see.

b) follows from \pr{28.1'} b).

c) We denote by 
$$\mac{\vartheta }{\bm\setminus \Gamma^{\bm} _0}{\bk'}$$
the covering map, by 
$$\mac{\vartheta '}{\Gamma '''}{\Gamma^{\bk} _1}$$
 the map defined by $\vartheta $, and put 
$$\mae{\phi }{\cbb{\bk'}{F}}{\cbb{\bm\setminus \Gamma^{\bm} _0}{F}}{x}{x\circ \vartheta }\,,$$
$$\mae{\phi' }{\cbb{\Gamma^{\bk} _1}{F}}{\cbb{\Gamma'''}{F}}{x}{x\circ \vartheta' }\;.$$
With the identifications of $\Gamma'''$ and $\Gamma^{\bk} _1$  with $\bt$ (by a) and \pr{28.1'} a)),
$$\mae{\vartheta '}{\Gamma'''}{\Gamma^{\bk} _1}{z}{z^2}\;.$$
By the commutativity of the index maps (\axi{27.9'e}) the diagrams
$$\begin{CD}
K_i(\cbb{\bm''}{F})@>K_i(\varphi)>>K_i(\cbb{\bk'}{F})@>K_i(\psi )>>K_i\left(\ccb{\Gamma^{\bk} _1}{F}\right)\\
@VV=V         @VVK_i(\phi )V    @VVK_i(\phi ')V \\
K_i(\cbb{\bm''}{F})@>>K_i(\bar{\varphi}'')>K_i(\cbb{\bm''')}{F})@>>K_i(\bar{\psi}'')>K_i(\ccb{\Gamma'''}{F})
\end{CD}$$
$$\begin{CD}
K_i\left(\ccb{\Gamma^{\bk} _1}{F}\right)@>\delta _i>>K_{i+1}(\cbb{\bm''}{F})\\
@VVK_i(\phi ')V      @VV=V\\
K_i(\ccb{\Gamma'''}{F})@>>\bar{\delta }''_i>K_{i+1}(\cbb{\bm''}{F})
\end{CD}$$
are commutative. By \pr{19.10} c),
$$\mae{K_i(\phi ')}{K_i\left(\ccb{\Gamma^{\bk}}{F}\right)}{K_i\left(\ccb{\Gamma'''}{F}\right)}{(a,b)}{(a,2b)}\;.$$
By b),
$$\mae{\delta _i}{K_i\left(\ccb{\Gamma^{\bk} _1}{F}\right)}{K_{i+1}\left(\cbb{\bk''}{F}\right)}{(a,b)}{(a,.2b)}\;.$$

d) By  the six-term axiom (\axi{27.9'd}), $\psi $ is K-null. The other assertions follow from c) and the six-term axiom (\axi{27.9'd}).\qed

{\center{\section{Pasting locally compact spaces}}}

\begin{p}\label{19.11}
Let $\Omega _1,\Omega _2$ be locally compact spaces, $\Gamma _1$ and $\Gamma _2$ closed sets of $\Omega _1$ and $\Omega _2$, respectively, $\mac{\vartheta }{\Gamma _1}{\Gamma _2}$ a homeomorphism, $\Omega '$ the topological sum of $\Omega _1\setminus \Gamma _1$ and $\Omega _2\setminus \Gamma _2$, $\Omega $ the locally compact space obtained from the topological sum of $\Omega _1$ and $\Omega _2$ by identifying the points $\omega $ and $\vartheta (\omega )$ for all $\omega \in \Gamma _1$, $\Gamma $ the closed set of $\Omega $ corresponding to the identified $\Gamma _1$ and $\Gamma _2$ (so $\Omega \setminus \Gamma =\Omega '$), $\mac{\varphi }{\cbb{\Omega \setminus \Gamma }{F}}{\cbb{\Omega }{F}}$ the inclusion map,
$$\mae{\psi }{\cbb{\Omega }{F}}{\cbb{\Gamma }{F}}{x}{x|\Gamma }\,,$$
and $\delta _i$ the index maps associated to the exact sequence in \frm
$$\oc{\cbb{\Omega \setminus \Gamma }{F}}{\varphi }{\cbb{\Omega }{F}}{\psi }{\cbb{\Gamma }{F}}\;.$$  
Let $J:=\z{1,2}$ and for every $j\in J$ let
$$\mac{\varphi _j}{\cbb{\Omega _j\setminus \Gamma _j}{F}}{\cbb{\Omega _j}{F}}\,,$$
$$\mac{\varphi' _j}{\cbb{\Omega _j\setminus \Gamma _j}{F}}{\cbb{\Omega'}{F}}\,,$$
$$\mac{\varphi'' _j}{\cbb{\Omega _j\setminus \Gamma _j}{F}}{\cbb{\Omega}{F}}$$
be the inclusion maps,
$$\mae{\psi _j}{\cbb{\Omega _j}{F}}{\cbb{\Gamma _j}{F}}{x}{x|\Gamma _j}\,,$$
$$\mae{\psi' _j}{\cbb{\Omega'}{F}}{\cbb{\Omega _j\setminus \Gamma _j}{F}}{x}{x|(\Omega _j\setminus \Gamma _j)}\,,$$
and $\delta _{j,i}$ the index maps associated to the exact sequence in \frm
$$\oc{\cbb{\Omega _j\setminus \Gamma _j}{F}}{\varphi _j}{\cbb{\Omega _j}{F}}{\psi _j}{\cbb{\Gamma _j}{F}}\;.$$
\begin{enumerate}
\item $\delta _{j,i}=K_{i+1}(\psi '_j)\circ \delta _i$ for every $j\in J$ and
$$\delta _i=K_{i+1}(\varphi '_1)\circ \delta _{1,i}+K_{i+1}(\varphi '_2)\circ \delta _{2,i}\;.$$
\item Assume $\cbb{\Omega _1}{F}$ K-null.
\begin{enumerate}
\item $\mac{\delta _{1,i}}{K_{i}(\cbb{\Gamma _1}{F})}{K_{i11}(\cbb{\Omega _1\setminus \Gamma _1}{F})}$ is a group isomorphism.
\item $\delta _i$ is injective.
\item $\psi $ is K-null.
\item $\mac{K_i(\varphi ''_2)}{K_i(\cbb{\Omega _2\setminus \Gamma _2}{F})}{K_i(\cbb{\Omega }{F})}$ is a group isomorphism.
\item If we put
$$\mac{\Phi_i :=K_i(\varphi '_2)\circ K_i(\varphi ''_2)^{-1}}{K_i(\cbb{\Omega }{F})}{K_i\left(\cbb{\Omega '}{F}\right)}$$
then the map
$$K_{i+1}(\cbb{\Gamma }{F})\times K_i(\cbb{\Omega }{F})\longrightarrow K_i\left(\cbb{\Omega '}{F}\right)\,,$$
$$(a,b)\longmapsto {\delta _{i+1}a+\Phi_i b}$$
is a group isomorphism.
\item If also $\cbb{\Omega _2}{F}$ is K-null then
$$K_i(\cbb{\Omega }{F})\approx K_{i+1}(\cbb{\Gamma }{F})\,,$$
$$K_i\left(\cbb{\Omega '}{F}\right)\approx K_{i+1}(\cbb{\Gamma }{F})^2\;.$$
\end{enumerate}
\end{enumerate}
\end{p}

a) follows from \pr{17.11} a), since $\psi ''_j$ of this Proposition is the identity map in the present case.

$b_1)$ follows from the Topological six-term sequence (\pr{24.11} a)).

$b_2)$ Let $a\in K_i(\cbb{\Gamma }{F})$ such that $\delta _ia=0$. By a), $\delta _{1,i}a=K_{i+1}(\psi '_1)\delta _ia=0$ and by $b_1)$, $a=0$.

$b_3)$ follows from $b_2)$ and the six-term axiom (\axi{27.9'd}).

$b_4)$ and $b_5)$ follow from $b_3)$ and \pr{17.11} $c_1)$, $c_2)$.

$b_6)$ follows from $b_1)$, $b_4)$, and the Product Theorem (\pr{14.11} a)).\qed

\begin{co}\label{19.11a}
Let $\Gamma $ be a locally compact space, $(\Omega _j)_{j\in J}$ a nonempty finite family of locally compact spaces such that $\cbb{\Omega _j}{F}$ is K-null for every $j\in J$, and for every $j\in J$ let $\Gamma _j$ be a closed set of $\Omega _j$ and $\mac{\vartheta _j}{\Gamma }{\Gamma _j}$ a homeomorphism. Let $\Omega '$ the topological sum of the family $(\Omega _j\setminus \Gamma _j)_{j\in J}$, and $\Omega $ the locally compact space obtained from the topological sum of the family $(\Omega _j)_{j\in J}$ by identifying for every $\omega \in \Gamma $ all the points $\vartheta_j(\omega )$ $(j\in J)$. Then
$$K_i(\cbb{\Omega }{F})\approx K_{i+1}(\cbb{\Gamma }{F})^{Card\,J-1}\,,$$
$$K_i\left(\cbb{\Omega '}{F}\right)\approx K_i(\cbb{\Omega }{F})\times K_{i+1}(\cbb{\Gamma }{F})\approx K_{i+1}(\cbb{\Gamma }{F})^{Card\,J}\;.$$ 
\end{co}

We prove the Corollary by induction with respect to $Card\,J$. For $Card\,J\in \z{1,2}$ the assertion follows from \pr{19.11} $b_1)$,$b_5),b_6)$. Let $k\in J$, assume the assertion holds for $J':=J\setminus \z{k}$, and denote by $\Omega ''$ the topological sum of the family $(\Omega _j\setminus \Gamma _j)_{j\in J'}$. By \pr{19.11} $b_4),b_5)$ and the induction hypothesis,
$$K_i(\cbb{\Omega }{F})\approx K_i\left(\cbb{\Omega ''}{F}\right)\approx K_{i+1}(\cbb{\Gamma }{F})^{Card\,J-1}\,,$$
$$K_i\left(\cbb{\Omega '}{F}\right)\approx K_{i+1}(\cbb{\Gamma }{F})\times K_i\left(\cbb{\Omega ''}{F}\right)\approx K_{i+1}(\cbb{\Gamma }{F})^J\;.\qedd$$

\begin{co}\label{19.11b}
Let $m,n\in \bn$,
$$\Gamma _+:=\me{\alpha \in \bbb_n}{\n{\alpha}=1,\,\alpha _n>0 },\, \Gamma _-:=\me{\alpha \in \bbb_n}{\n{\alpha}=1,\,\alpha _n\leq 0 }\,,$$
and $\Omega $ the locally compact space obtained from the topological sum of the family $(\bbb_n\setminus \Gamma _-)_{j\in \bnn{m}}$ by identifying all the $\Gamma _+$. Then
$$K_i(\cbb{\Omega }{F})\approx K_{i+n}(F)^{m-1}\;.$$
\end{co}

By \pr{24.11a}, $\cbb{\bbb_n\setminus \Gamma _-}{F}$ is null-homotopic and so K-null. For $n>1$, $\Gamma _+$ is homeomorphic to $\br^{n-1}$ so by \h{17.1'c} a), 
$$K_i(\cbb{\Gamma _+}{F})\approx K_{i+n-1}(F)$$
and this relation obviously holds also for $n=1$. Then by \cor{19.11a}, 
$$K_i(\cbb{\Omega }{F})\approx K_{i+n}(F)^{m-1}\;.\qedd$$

{\it Remark.} The above result can be deduced also from \ee{20.1'} by using \pr{27.3'b} d).

\begin{co}\label{3.9}
Let $\Omega ',\Omega ''$ be locally compact spaces, $\omega '\in \Omega '$, $\omega ''\in \Omega ''$, and $\Omega $ the locally compact space obtained from the topological sum of $\Omega '$ and $\Omega ''$ by identifying $\omega '$ and $\omega ''$. If $\cbb{\Omega ''}{F}$ is K-null then
$$\kk{i}{\cbb{\Omega }{F}}\approx \kk{i}{\cbb{\Omega'\setminus \z{\omega '} }{F}}\;.$$
\end{co}

The assertion follows from \pr{19.11} $b_4)$.\qed

\begin{p}\label{15.9b}
Let $\Omega ',\Omega ''$ be compact spaces, $\omega '\in \Omega '$, $\omega ''\in \Omega ''$, and $\Omega $ the compact space obtained by identifying the points $\omega '$ and $\omega ''$ in the topological sum of $\Omega '$ and $\Omega ''$. Then
$$K_i(\ccb{\Omega }{F})\approx K_i\left(\cbb{\Omega \setminus \Omega '}{F}\right)\times K_i\left(\ccb{\Omega '}{F}\right)\;.$$
\end{p}

Let $\mac{\varphi }{\cbb{\Omega \setminus \Omega '}{F}}{\ccb{\Omega }{F}}$ be the inclusion map and
$$\mae{\psi }{\ccb{\Omega }{F}}{\ccb{\Omega '}{F}}{x}{x|\Omega '}\;.$$
We put for every $x\in \ccb{\Omega '}{F}$,
$$\mae{\lambda x}{\Omega }{F}{\omega }{\abb{x(\omega )}{\omega \in \Omega '}{x(\omega _0)}{\omega \in \Omega ''}}\,,$$
where $\omega _0\in \Omega $ denotes the point corresponding to the identified points $\omega '$ and $\omega ''$. Then
$$\od{\cbb{\Omega \setminus \Omega '}{F}}{\varphi }{\ccb{\Omega }{F}}{\psi }{\lambda }{\ccb{\Omega '}{F}}$$
is a split exact sequence in \frm and the assertion follows from the split exact axiom (\axi{27.9'a}).\qed

\begin{p}\label{16.5}
Let $(\Omega _j)_{j\in \bnn{n}}$ be a family of compact spaces and for every $j\in \bnn{n}$ let $\omega _j,\omega '_j$ be distinct points of $\Omega _j$. If $\Omega $ denotes the compact space obtained from the topological sum of the family $(\Omega _j)_{j\in \bnn{n}}$ by identifying $\omega '_j$ with $\omega _{j+1}$ for all $j\in \bnn{n-1}$ then 
$$K_i(\ccb{\Omega }{F})\approx K_i(F)\times \proo{j=1}{n}K_i(\cbb{\Omega _j\setminus \z{\omega _j}}{F})\;.$$
If $(k_j)_{j\in \bnn{n}}$ is a family in $\bn$, $\Omega _j=\bs_{k_j}$ for every $j\in \bnn{n}$, and 
$$p:=Card\,\me{j\in \bnn{n}}{k_j \;\emph{is\,\,even}}\,,\qquad q:=Card\,\me{j\in \bnn{n}}{k_j \;\emph{is\,\,odd}}$$
then
$$K_i(\ccb{\Omega }{F})\approx K_i(F)^{p+1}\times K_{i+1}(F)^q\;.$$
\end{p}

We put $\bar{\Omega }_n:=\Omega  $ and prove the assertion by induction with respect to $n\in \bn$. For $n=1$ the assertion follows from the Alexandroff K-theorem (\h{20.4} a)). Assume the assertion hods for an $n\in \bn$. By \pr{15.9b} and the induction hypothesis,
$$K_i\left(\ccb{\bar{\Omega }_{n+1} }{F}\right)\approx K_i\left(\cbb{\bar{\Omega }_{n+1}\setminus \bar{\Omega }_n  }{F}\right)\times K_i\left(\ccb{\bar{\Omega }_n }{F}\right)\approx $$
$$\approx K_i(\cbb{\Omega_{n+1}\setminus \z{\omega _{n+1}}  }{F})\times K_i\left(\ccb{\bar{\Omega }_n }{F}\right)\approx $$
$$\approx K_i(\cbb{\Omega _{n+1}\setminus \z{\omega _{n+1}}}{F})\times K_i(F)\times \proo{j=1}{n}K_i(\cbb{\Omega _j\setminus \z{\omega _j}}{F})\approx $$
$$\approx K_i(F)\times \proo{j=1}{n+1}K_i(\cbb{\Omega _j\setminus \z{\omega _j}}{F})\,,$$
which finishes the inductive proof. The last assertion follows now from \h{17.1'c} a), since $\bs_{k_j}\setminus \z{\omega _j}$ is homeomorphic to $\br^{k_j}$.\qed

\begin{p}\label{31.10}
Let $\Omega _1,\Omega _2$ be locally compact spaces such that the $\eo$ $\cbb{\Omega _2}{F}$ is K-null, $\Gamma $ a compact set of $\Omega _1$, and $\mac{\vartheta }{\Gamma }{\Omega _2}$ a continuous map. We denote by $\Omega $ the locally compact space obtained from the topological sum of $\Omega _1$ and $\Omega _2$ by identifying the points $\omega $ and $\vartheta (\omega )$ for all $\omega \in \Gamma $. 
\begin{enumerate}
\item If 
$$\mac{\varphi }{\cbb{\Omega _1\setminus \Gamma }{F}}{\cbb{\Omega }{F}}$$
denotes the inclusion map then
$$\mac{K_i(\varphi )}{K_i(\cbb{\Omega_1 \setminus \Gamma }{F})}{K_i(\cbb{\Omega }{F})}$$
is a group isomorphism. If in addition $\Omega \in \Upsilon $ or $\Omega _1\setminus \Gamma \in \Upsilon $ then
$$\Omega,\Omega _1\setminus \Gamma  \in \Upsilon \,,\; p(\Omega) =p(\Omega _1\setminus \Gamma )\,,\; q(\Omega) =q(\Omega _1\setminus \Gamma) \,,\; \Omega _\Upsilon =(\Omega _1\setminus \Gamma )_\Upsilon \;.$$
\item If $\Omega ^*$ denotes the Alexandroff compactification of $\Omega $ then
$$K_i(\ccb{\Omega^*}{F})\approx K_i(F)\times K_i(\cbb{\Omega_1 \setminus \Gamma }{F})\;.$$
\end{enumerate}.
\end{p}

a) If we put
$$\mae{\psi }{\cbb{\Omega }{F}}{\cbb{\Omega _2}{F}}{x}{x|\Omega _2}$$
then
$$\oc{\cbb{\Omega _1\setminus \Gamma }{F}}{\varphi }{\cbb{\Omega }{F}}{\psi }{\cbb{\Omega _2}{F}}$$
is an exact sequence in \frm. Since $\cbb{\Omega _2}{F}$ is K-null, the assertion follows from the Topological six-term sequence (\pr{24.11} c)).

b) follows from a) and Alexandroff's K-theorem (\h{20.4} a)).\qed

\begin{co}\label{31.10a}
Let $(\Omega _j)_{j\in J}$ be a finite family of locally compact spaces, $\omega _j\in \Omega _j$ for all $j\in J$, and $\Omega $ the locally compact space obtained from the topological sum of the family $(\Omega _j)_{j\in J}$ by identifying the points $\omega _j$ for all $j\in J$.
\begin{enumerate}
\item If there is a $j_0\in J$ such that $\cbb{\Omega _{j_0}}{F}$ is K-null then
$$K_i(\cbb{\Omega }{F})\approx \pro{j\in J\setminus \z{j_0}}K_i(\cbb{\Omega _j\setminus \z{\omega _j}}{F})\;.$$
\item If $\Omega _j:=[0,1[$ for all $j\in J$ and $n:=Card\,J$ then 
$$K_i(\cbb{\Omega }{F})\approx K_{i+1}(F)^{n-1}\;.$$
\item Let $j_0\in J$ and $\Omega _{j_0}:=[0,1[$. If $(k_j)_{j\in J\setminus \z{j_0}}$ is a family in $\bn$,
$$p:=Card\,\me{j\in J\setminus \z{j_0}}{k_j\, \emph{is even}}\,,$$
$$q:=Card\,\me{j\in J\setminus \z{j_0}}{k_j \,\emph{is odd}}\,,$$
and $\Omega _j:=\bs_{k_j}$ for every $j\in J\setminus {j_0}$ then
$$K_i(\cbb{\Omega }{F})\approx K_i(F)^p\times K_{i+1}(F)^q\;.$$
\end{enumerate}
\end{co}

a) Let $\Omega '$ be the locally compact space obtained from the topological sum of the family $(\Omega _j)_{J\setminus \z{j_0}}$ by identifying the points $\omega _j$ for all $j\in J\setminus \z{j_0}$ and let $\bar{\omega } $ denote the point obtained by this identification. If we replace in \pr{31.10} $\Omega _1$ by $\Omega '$, $\Omega _2$ by $\Omega _{j_0}$, $\Gamma $ by $\bar{\omega } $, and take $\vartheta(\bar{\omega } ):=\omega _{j_0} $ 
then we get
$$K_i(\cbb{\Omega }{F})\approx K_i\left(\cbb{\Omega '\setminus \z{\bar{\omega } }}{F}\right)\;.$$
$\Omega '\setminus \z{\bar{\omega } }$ is the topological sum of the family $(\Omega _j\setminus \z{\omega _j})_{j\in J\setminus \z{j_0}}$ so by the Product Theorem (\pr{14.11} a)),
$$K_i\left(\cbb{\Omega '\setminus \z{\bar{\omega } }}{F}\right)\approx \pro{j\in J\setminus \z{j_0}}K_i(\cbb{\Omega _j\setminus \z{\omega _j}}{F})\;.$$

b) follows immediately from a) since $\cbb{[0,1[}{F}$ is K-null and
$$K_i(\cbb{[0,1[\setminus \z{\omega }}{F})\approx K_{i+1}(F)$$
for all $\omega \in [0,1[$.

c) For $j\in J\setminus \z{j_0}$, $\bs_{k_j}\setminus \z{\omega _j}$ is homeomorphic to $\br^{k_j}$ and so by \h{17.1'c} a), 
$K_i\left(\cbb{\bs_{k_j}\setminus \z{\omega _j}}{F}\right)\approx K_{i+k_j}(F)$. Since $\cbb{[0,1[}{F}$ is K-null, we get from a),
$$K_i(\cbb{\Omega }{F})\approx K_i(F)^p\times K_{i+1}(F)^q\;.\qedd$$

\begin{co}\label{3.11}
Let $J_1,J_2,J_3$ be pairwise disjoint finite sets and let $\Omega $ be the locally compact space (the {\bf{graph}}) obtained from the topological sum of $[0,1]\times J_1$, $[0,1[\times J_2$, and $]0,1[\times J_3$ by identifying some of the points of the set
$$\me{(0,j)}{j\in J_1\cup J_2}\cup \me{(1,j)}{j\in J_1}\;.$$
If $s$ denotes the number of compact connected components of $\Omega $ and $r_0$ and $r_1$ denote the number of vertices and chords of the graph $\Omega $, respectively, then
$$K_i(\cbb{\Omega }{F})\approx K_i(F)^s\times K_{i+1}(F)^{s+r_1-r_0}\;.$$
\end{co}

By the Product Theorem (\pr{14.11} a)), we may assume $\Omega $ connected.

Assume first there is a $j\in J_3$ such that $\Omega $ contains $]0,1[\times \z{j}$. Since $\Omega $ is connected, $\Omega =]0,1[\times \z{j}$. Thus $\Omega $ is homeomorphic to $\br$, $r_1-r_0=1$, and the assertion follows from \h{17.1'c} a).

Assume now there is a $j\in J_2$ such that $\Omega $ contains $[0,1[\times \z{j}$. By \pr{31.10} a),
$$K_i(\cbb{\Omega }{F})\approx K_i(\cbb{\Omega \setminus ([0,1[\times \z{j})}{F})\;.$$
$\Omega $ and $\Omega \setminus ([0,1[\times \z{j})$ have the same $r_1-r_0$, so we may replace $\Omega $ by $\Omega \setminus ([0,1[\times \z{j})$. Repeating the operation, we obtain finally a locally compact space, which is the topological sum of a finite family $(]0,1[)_{j\in J}$, and in this case the assertion follows from the Product Theorem (\pr{14.11} a)) and \h{17.1'c} a).

Finally assume $\Omega $ compact. Then there is a $j\in J_1$ such that $\Omega $ contains $[0,1]\times \z{j}$. By the above and by Alexandroff's K-theorem
$$K_i(\ccb{\Omega }{F})\approx K_i(F)\times K_i(\cbb{\Omega \setminus \z{(1,j)}}{F})\;.$$
If $s',r'_0,r'_1$ denote the corresponding numbers associated to $\Omega \setminus \z{(1,j)}$ then $s'=0$, $r'_0=r_0-1$, and $r'_1=r_1$. All the connected components of $\Omega \setminus \z{(1,j)}$ satisfy the condition of the above paragraphs, so
$$K_i(\cbb{\Omega \setminus \z{(1,j)}}{F})\approx K_{i+1}(F)^{r'_1-r'_0}\approx K_{i+1}(F)^{1+r_1-r_0}\,,$$
$$K_i(\ccb{\Omega }{F})\approx K_i(F)\times K_{i+1}(F)^{1+r_1-r_0}\;.\qedd$$

\begin{co}\label{22.5'}
If $\Omega $ is a compact graph contained in $\bbb_n$ then
$$\kk{i}{\cbb{\bbb_n\setminus \Omega }{F}}\approx \kk{i}{F}^{s-r_0+r_1}\times \kk{i+1}{F}^{s-1}\,,$$
where $s$ denotes the number of connected components of $\Omega $ and $r_0$ and $r_1$ the munber of vertices and chords of $\Omega $, respectively.
\end{co}

Let $\omega $ be a vertex of $\Omega $. By \cor{3.11} and \cor{5.5'} a),
$$\kk{i}{\cbb{\Omega \setminus \z{\omega }}{F}}\approx \kk{i}{F}^{s-1}\times \kk{i+1}{F}^{s-r_0+r_1}$$
and by \h{17.1'b} b),
$$\kk{i}{\cbb{\bbb_n\setminus \Omega }{F}}\approx \kk{i+1}{\cbb{\Omega \setminus \z{\omega }}{F}}\approx  \kk{i}{F}^{s-r_0+r_1}\times \kk{i+1}{F}^{s-1}\;.\qedd$$

\begin{e}\label{27.4'}
Let $\bbn$, $\Gamma $ a closed set of $\bs_n$, $\emptyset \not=\Gamma \not=\bs_n$, $\omega \in \Gamma $, $\Gamma '$ the compact space obtained from $\Gamma \times [0,1]$ by identifying the points of $\Gamma \times 0$, and $\Omega $ the compact space obtained from the topological sum of $\bs_n$ and $\Gamma '$ by identifying the points of $\Gamma \subset \bs_n$ with the points of $\Gamma \times \z{1}\subset \Gamma '$.
\begin{enumerate}
\item $\kk{i}{\ccb{\Omega }{F}}\approx \kk{i}{F}\times \kk{i+n}{F}\times \kk{i+1}{\cbb{\Gamma \setminus \z{\omega }}{F}}$.
\item If $\Gamma $ is finite then
$$\kk{i}{\ccb{\Omega }{F}}\approx \kk{i}{F}\times \kk{i+n}{F}\times \kk{i+1}{F}^{Card\,\Gamma -1}\;.$$
\item If $\Gamma $ is a graph then
$$\kk{i}{\ccb{\Omega }{F}}\approx \kk{i}{F}^{1+s+r_1-r_0}\times \kk{i+n}{F}\times \kk{i+1}{F}^{s-1}\,,$$
where $s$ denotes the number of connected components of $\Omega $ and $r_0$ and $r_1$ denote the number of vertices and chords of the graph $\Gamma $, respectively,
\end{enumerate}
\end{e}

a) By \h{17.1'c} $e_1)$,
$$\kk{i}{\cbb{\bs_n\setminus \Gamma }{F}}\approx \kk{i+n}{F}\times \kk{i+1}{\cbb{\Gamma \setminus \z{\omega }}{F}}\;.$$
By \pr{24.11a}, $\cbb{\Gamma '\setminus \z{0}}{F}$ is K-null, where $0$ is the point obtained from the identification of the points of $\Gamma \times \z{0}$. By \pr{31.10} a),
$$\kk{i}{\cbb{\Omega \setminus \z{0}}{F}}\approx \kk{i}{\cbb{\bs_n\setminus \Gamma }{F}}\,,$$
so by Alexandroff's K-theorem (\h{20.4} a)),
$$\kk{i}{\ccb{\Omega }{F}}\approx \kk{i}{F}\times \kk{i+n}{F}\times \kk{i+1}{\cbb{\Gamma \setminus \z{\omega }}{F}}\;.$$

b) follows from a) and the Product Theorem (\pr{14.11} a)).

c) By \cor{3.11} and Alexandroff's K-theorem (\h{20.4} a)),
$$\kk{i}{\ccb{\Gamma }{F}}\approx \kk{i}{F}^s\times \kk{i+1}{F}^{s+r_1-r_0}\,,$$
$$\kk{i}{\cbb{\Gamma \setminus \z{\omega }}{F}}\approx \kk{i}{F}^{s-1}\times \kk{i+1}{F}^{s+r_1-r_0}\,,$$
so by a),
$$\kk{i}{\ccb{\Omega }{F}}\approx \kk{i}{F}^{1+s+r_1-r_0}\times \kk{i+n}{F}\times \kk{i+1}{F}^{s-1}\;.\qedd$$

\begin{p}\label{15.4'}
Let $(p_j)_{j\in J}$ be a finite family in $\bn$, $(J\not=\emptyset )$, and for every $j\in J$ put $\Omega _j:=\bs_{p_j}$. Let $\Omega '$ be the topological sum of the family $(\Omega _j)_{j\in J}$, $(\Gamma _k)_{k\in K}$ a finite family of pairwise disjoint nonempty finite subsets of $\Omega '$, $\Gamma :=\bigcup\limits _{k\in K}\Gamma _k$, and $\Omega $ the compact space obtained from $\Omega '$ by identifying for every $k\in K$ the points of $\Gamma _k$. If $\Omega $ is connected then
$$\kk{i}{\ccb{\Omega }{F}}\approx \kk{i}{F}\times \kk{i+1}{F}^{Card\,\Gamma -Card\,J-Card\,K+1}\times \pro{j\in J}\kk{i+p_j}{F}\;.$$
\end{p}

If $K=\emptyset $, since $\Omega $ is connected, $J$ is a one-point set and the assertion holds by \h{17.1'c} b). Thus we may assume $K=\bnn{n}$ for some $\bbn$. Take $k_1\in K$ and put $J_1:=\me{j\in J}{\Omega _j\cap \Gamma _{k_1}\not=\emptyset }$. We define recursively an injective family $(k_m)_{m\in \bnn{n}}$ in $K$ and an increasing family $(J_m)_{m\in \bnn{n}}$ of subsets of $J$ in the following way. Let $m\in \bnn{n}$, $m>1$, and assume the families were defined up to $m-1$. Since $\Omega $ is connected there is a $k_m\in K\setminus \me{k_q}{q\in \bnn{m-1}}$ such that $\Gamma _{k_m}\cap J_{m-1}\not=\emptyset $. We put 
$$J_m:=\me{j\in J}{\Omega _j\cap \left(\bigcup_{q=1}^m\Gamma _{k_q} \right)\not=\emptyset }\;.$$
It is easy to prove by induction with respect to $m\in \bnn{n}$ that
$$Card\,\left(\bigcup_{q=1}^m \Gamma _{k_q}\right)-Card\,J_m-m+1\geq 0$$
for every $m\in \bnn{n}$. In particular,
$$Card\,\Gamma -Card\,J-Card\,K+1\geq 0\;.$$

For every $j\in J$, by \pr{26.2'a} and \h{17.1'c} a),
$$\kk{i}{\cbb{\Omega _j\setminus \Gamma }{F}}\approx \kk{i+1}{F}^{Card\,(\Gamma \cap \Omega _j)-1}\times \kk{i+p_j}{F}$$
so that by the Product Theorem (\pr{14.11} a)),
$$\kk{i}{\cbb{\Omega '\setminus \Gamma }{F}}\approx \kk{i+1}{F}^{Card\,\Gamma -Card\,J}\times \pro{j\in J}\kk{i+p_j}{F}\;.$$
For every $k\in K$ let $\omega _k$ be the point of $\Omega $ corresponding to the unified points of $\Gamma _k$ and put $\Delta :=\me{\omega _k}{k\in K}$. Then by \pr{26.2'a},
$$\kk{i}{\cbb{\Omega \setminus \Delta }{F}}\approx \kk{i}{\cbb{\Omega \setminus \z{\omega _{k_0}}}{F}}\times \kk{i+1}{F}^{Card\,K-1}\,,$$
where $k_0\in K$. By the above and by Alexandroff's K-theorem, since $\Omega \setminus \Delta =\Omega '\setminus \Gamma $,
$$\kk{i}{\ccb{\Omega }{F}}\times \kk{i+1}{F}^{Card\,K-1}\approx $$
$$\approx \kk{i}{F}\times \kk{i}{\cbb{\Omega \setminus \z{\omega _{k_0}}}{F}}\times \kk{i+1}{F}^{Card\,K-1}\approx $$
$$\approx \kk{i}{F}\times \kk{i}{\cbb{\Omega \setminus \Delta }{F}}\approx \kk{i}{F}\times \kk{i}{\cbb{\Omega '\setminus \Gamma }{F}}\approx $$
$$\approx \kk{i}{F}\times \kk{i+1}{F}^{Card\,\Gamma -Card\,J-Card\,K+1}\times \kk{i+1}{F}^{Card\,K-1}\times \pro{j\in J}\kk{i+p_j}{F}\,,$$
$$\kk{i}{\ccb{\Omega }{F}}\approx \kk{i}{F}\times \kk{i+1}{F}^{Card\,\Gamma -Card\,J-Card\,K+1}\times \pro{j\in J}\kk{i+p_j}{F}\;.\qedd$$

\begin{co}\label{21.8}
Let $(p_j)_{j\in \bnn{n}}$ be a family in $\bn$ and for every $j\in \bnn{n}$ put $\Omega _j:=\bs_{p_j}$. For every $j\in \bnn{n}$ let $\Gamma _j$ and $\Gamma '_j$ be disjoint nonempty finite subsets of $\Omega _j$ such that $k_j:=Card\,\Gamma '_j=Card\,\Gamma _{j+1}$ for every $j\in \bnn{n-1}$. We denote by $\Omega $ the compact space obtained from the topological sum of the family $(\Omega _j)_{j\in \bnn{n}}$ by identifying in a bijective way $\Gamma '_j$ with $\Gamma _{j+1}$ for all $j\in \bnn{n-1}$. Then
$$K_i(\ccb{\Omega }{F})\approx K_i(F)\times  K_{i+1}(F)^{\sii{j=1}{n-1}(k_j-1)}\times \proo{j=1}{n}K_{i+p_j}(F)\;.\qedd$$
\end{co}

\begin{p}\label{11.6'}
Let $\Omega _1,\,\Omega _2$ be locally compact spaces and for every $j\in \z{1,2}$ let $\Gamma _j$ be a compact set of $\Omega _j$ and $\mac{\vartheta _j}{\bbb_n}{\Gamma _j}$ a homeomorphism such that $\Delta _j:=\vartheta _j(\bbb_n\setminus \bs_{n-1})$ is an open set of $\Omega _j$. We denote by $\Omega $ the locally compact space obtained from the topological sum of $\Omega _1\setminus \Delta _1$ and $\Omega _2\setminus \Delta _2$ by identifying the points $\vartheta _1(\omega )$ and $\vartheta _2(\omega )$ for all $\omega \in \bs_{n-1}$. Then for every $\omega \in \bs_{n-1}$,
$$\kk{i}{\cbb{\Omega \setminus \z{\vartheta _1(\omega )}}{F}}\approx$$
$$\approx  \kk{i}{\cbb{\Omega _1\setminus \Gamma _1}{F}}\times\kk{i}{\cbb{\Omega _2\setminus \Gamma _2}{F}}\times \kk{i+n-1}{F}\;. $$
\end{p}

We use the notation of the topological triple (\pr{3.12}), which we mark with a prime in order to distinguish them from the present notation. We put $\Omega '_2:=\Omega \setminus \z{\vartheta _1(\omega )}$ and take as $\Omega '_3$ the topological sum of $\Omega _1\setminus \Gamma _1$ and $\Omega _2\setminus \Gamma _2$ and as $\Omega '_1$ the locally compact space obtained from $\Omega $ by completing first $\vartheta _1(\bs_{n-1})$ to $\vartheta _1(\bbb_n)$ and deleting then $\omega $. By the Product Theorem (\pr{14.11} a)),
$$\kk{i}{\cbb{\Omega '_3}{F}}\approx \kk{i}{\cbb{\Omega _1\setminus \Gamma _1}{F}}\times\kk{i}{\cbb{\Omega _2\setminus \Gamma _2}{F}}\;.$$
Since $\Omega '_2\setminus \Omega '_3$ is homeomorphic to $\bs_{n-1}\setminus \z{\omega }$, we get by \h{17.1'c} $e_1)$,
$$\kk{i}{\cbb{\Omega '_2\setminus \Omega '_3}{F}}\approx \kk{i+n-1}{F}\;.$$
Thus by the topological triple (\pr{3.12} $b_3)$) (and \h{17.1'b} b)),
$$\kk{i}{\cbb{\Omega \setminus \z{\vartheta _1(\omega )}}{F}}\approx\kk{i}{\cbb{\Omega '_2}{F}}\approx$$
$$\approx  \kk{i}{\cbb{\Omega '_3}{F}}\times \kk{i}{\cbb{\Omega '_2\setminus \Omega '_3}{F}}\approx $$
$$\approx \kk{i}{\cbb{\Omega _1\setminus \Gamma _1}{F}}\times\kk{i}{\cbb{\Omega _2\setminus \Gamma _2}{F}}\times \kk{i+n-1}{F}\;.\qedd$$

\begin{co}\label{10.6'}
If $S_g$ is an orientable compact connected surface of genus $g\in \bn$ and $\Gamma $ is a nonempty finite subset of $S_g$ then
$$\kk{i}{\ccb{S_g}{F}}\approx \kk{i}{F}^{g+1}\times \kk{i+1}{F}^{3g-1}\,,$$
$$\kk{i}{\ccb{S_g\setminus \Gamma }{F}}\approx \kk{i}{F}^{g}\times \kk{i+1}{F}^{3g-2+Card\,\Gamma }\;.$$
\end{co}

Assume first $\Gamma $ is a one-point set $\z{\omega }$. We prove the second assertion in this case by induction with respect to $g\in\bn$. By \pr{12.5'a} b), the assertion holds for $g=1$. Assume now the assertion holds for $g\in \bn$. Let $\Delta _1$ be a closed disc of $S_1$, $\Delta _g$ a closed disc of $S_g$, $\omega \in \Delta _1$, and $\omega \in \Delta _g$. $S_{g+1}\setminus \z{\omega }$ can be obtained from the topological sum of $S_1\setminus \Delta _1$, $S_g\setminus \Delta _2$, and $\bs_1\setminus \z{\omega }$ by pasting $\bs_1\setminus \z{\omega }$ in the the boundaries of $\Delta _1\setminus \z{\omega }$ and $\Delta _g\setminus \z{\omega }$. By the induction hypothesis, since $S_g\setminus \Delta _g$ is homeomorphic to $S_g\setminus \z{\omega }$,
$$\kk{i}{\cbb{S_g\setminus \Delta _g}{F}}\approx \kk{i}{F}^g\times \kk{i+1}{F}^{3g-1}\;.$$
By \pr{11.6'}, 
$$\kk{i}{\cbb{S_{g+1}\setminus \z{\omega }}{F}} \approx \kk{i}{F}^{g+1}\times \kk{i+1}{F}^{3g+2}\,,$$ 
which finishes the inductive proof.

The first assertion follows now from Alexandroff's K-theorem (\pr{20.4} a)) and the second one from \pr{26.2'a}.\qed 

The following Example shows a way to generalize \cor{10.6'}.

\begin{e}\label{10.6'a}
Let $\Omega $ be the compact space obtained from the topological sum of $\bs_1\times \bs_2\setminus \Delta $, $\bs_1\times \bs_1\times \bs_1\setminus \Delta '$, and $\bs_2$, where $\Delta $ and $\Delta '$ denote balls homeomorphic to $\bbb_3$ by pasting $\bs_2$ in the boundaries of $\Delta $ and $\Delta '$. Then for every nonempty finite subset $\Gamma $ of $\Omega $,
$$\kk{i}{\ccb{\Omega }{F}}\approx \kk{i}{F}^5\times \kk{i+1}{F}^6\,,$$
$$\kk{i}{\cbb{\Omega \setminus \Gamma }{F}}\approx \kk{i}{F}^4\times \kk{i+1}{F}^{5+Card\,\Gamma }\;.\qedd$$
\end{e}

{\it Remark.} Let 
$$\oc{F_1}{\varphi }{F_2}{\psi }{F_3}\,,$$
$$\oc{G_1}{\varphi '}{G_2}{\psi '}{G_3}$$
be exact sequences in \frm and $\mac{\lambda }{F_3}{G_3}$ and isomorphism in \frm. Then
$$H:=\me{(x,y)\in F_2\times G_2}{\psi 'y=\lambda \psi x}$$
is a C*-subalgebra of $F_2\times G_2$ containing the ideal $F_1\times G_1$ of $F_2\times G_2$. $H$ corresponds to the operation of pasting $F_2$ and $G_2$ in \frm.

\begin{center}
\chapter{Some supplementary results}
\end{center}

\begin{center}
\fbox{\parbox{8.8cm}{Throughout this chapter $F$ denotes an $E$-C*-algebra}}

\section{Full $E$-C*-algebras}
\end{center}

\begin{de}\label{10.3'r}
A {\bf{\en}} is a unital C*-algebra $F$ for which $E$ is a \emph{canonical} unital C*-subalgebra such that $\alpha x=x\alpha $ for all $(\alpha ,x)\in E\times F$. Every \en is canonically an $E$-C*-algebra, the exterior multiplication being the restriction of the interior multiplication. We denote by \frc the category of full $E$-C*-algebras for which the morphisms are the unital $E$-linear C*-homomorphisms. In particular $\fr{C}_{\bc}$ is the category of all unital C*-algebras with unital C*-homomorphisms. A {\bf full $E$-C*-subalgebra of $F$} is a C*-subalgebra of $F$ containing $E$. An isomorphism of full $E$-C*-algebras is also called {\bf $E$-C*-isomorphism}. 

If $\pro{j\in J}F_j$ is a finite family of \en $\!\!$s, $J\not=\emptyset $,  then $\pro{j\in J}F_j$ is a \en, the canonical embedding $E\rightarrow \pro{j\in J}F_j$ being given by
$$\mad{E}{\pro{j\in J}F_j}{\alpha }{(\alpha)_{j\in J}}\;.$$
\end{de}

If $F$ is a  full $E$-C*-algebra and $G$ a unital C*-algebra then the map
$$\mad{E}{F\otimes G}{\alpha }{\alpha \otimes 1_G}$$
is an injective C*-homomorphism. In particular, the $E$-C*-algebra $F\otimes G$ has a canonical structure of a full $E$-C*-algebra.

\begin{p}\label{10.3'a}
Let $F$ be an \eo. We denote by $\check{F} $ the vector space $E\times F$ endowed with the bilinear map
$$\mad{(E\times F)\times (E\times F)}{E\times F}{((\alpha ,x),(\beta ,y))}{(\alpha \beta ,\alpha y+\beta x+xy)}$$ and with the involution 
$$\mad{E\times F}{E\times F}{(\alpha ,x)}{(\alpha ^*,x^*)}\;.$$
\begin{enumerate}
\item $\check{F} $ is an involutive unital algebra with $(1_E,0)$ as unit and $\me{(\alpha ,0)}{\alpha \in E}$ is a unital involutive subalgebra of $\check{F} $ isomorphic to $E$. 
\item If $E$ and $F$ are C*-subalgebras of a C*-algebra $G$ then the map 
$$\mae{\varphi }{\check{F} }{E\times G}{(\alpha ,x)}{(\alpha ,\alpha +x)}$$
is an injective involutive algebra homomorphism with closed image 
$$\me{(\alpha ,y)\in E\times G}{\alpha -y\in F}\;.$$
In particular $\varphi (\check{F} )$ is a C*-subalgebra of $E\times G$ and there is a norm on $\check{F} $ with respect to which $\check{F} $ is a C*-algebra.
\item There is a unique C*-norm on $\check{F} $ making it a C*-algebra. Moreover $\check{F} $ is a \en and $F$ may be identified with the closed ideal 
$$\me{(0,x)}{x\in F}$$
of $\check{F} $. We shall always consider $\check{F} $ endowed with the structure of a \en.
\item If $F$ is a \en then the map
$$\mad{\check{F} }{E\times F}{(\alpha ,x)}{(\alpha ,\alpha +x)}$$
is an isomorphism of $E$-C*-algebras with inverse
$$\mad{E\times F}{\check{F} }{(\alpha ,x)}{(\alpha ,x-\alpha )}\;.$$
\item If $E=\bc$ then $\check{F} $ is the unitization $\tilde{F} $ of $F$.
\end{enumerate}
\end{p} 

a) is easy to verify.

b) Only the assertion that the image of $\varphi $ is closed needs a proof. Let $(\alpha ,x)\in \overline{\varphi (\check F)}$. There are sequences $(\alpha _n)_{n\in \bn}$ and $(x_n)_{n\in \bn}$ in $E$ and $F$, respectively, such that
$$\lim_{n\rightarrow \infty }(\alpha _n,\alpha _n+x_n)=(\alpha ,x)\;.$$
It follows
$$\alpha  =\lim_{n\rightarrow \infty }\alpha _n\in E\,,\quad x-\alpha =\lim_{n\rightarrow \infty }x_n\in F\,,\quad (\alpha ,x)=\varphi (\alpha ,x-\alpha )\in \varphi (\check F)\;.$$
Thus $\varphi (\check F)$ is closed.

c) Let $\Omega $ be the spectrum of $E$ and $\tilde{F} $ the unitization of $F$. Then $E$ and $F$ are C*-subalgebras of the C*-algebra $\ccb{\Omega }{\tilde{F} }$ and the assertion follows from b).

d) follows from c) and b).

e) is obvious.\qed

\begin{e}\label{22.4'}
Let $F$ be a commutative $E$-C*-algebra.
\begin{enumerate}
\item $\check{F} $ is commutative. We denote by $\Omega _E$, $\Omega _F$, and $\Omega _{\check{F} }$ the spectra of $E$, $F$, and $\check{F} $, respectively.
\item $\Omega _F$ is homeomorphic to an open set $\Omega '$ of $\Omega _{\check{F} }$ such that $F\approx \cbb{\Omega '}{\bc}$.
\item There is a unique surjective continuous map $\mac{\vartheta }{\Omega _{\check{F} }}{\Omega _E}$ such that if we put
$$\mae{\phi }{E\approx \ccb{\Omega _E}{\bc}}{\check{F}\approx \ccb{\Omega _{\check{F} }}{\bc} }{\alpha }{\alpha \circ \vartheta }$$
then $\phi $ is an injective continuous C*-homomorphism (so we may identify $E$ with $\phi (E)$).
\item The restriction of $\vartheta $ to $\Omega _{\check{F} }\setminus \Omega '$ is a homeomorphism.
\item If $F$ is unital then $\Omega _{\check{F} }$ is homeomorphic to the topological sum of $\Omega _E$ and $\Omega _F$.
\end{enumerate}
\end{e}

a) is easy to see.

b) follows from the fact that $F$ may be identified with a closed ideal of $\check{F} $ (\pr{10.3'a} c)).

c) is proved in [C1] Proposition 4.1.2.15.

d) Let $\omega \in \Omega _E$ and put
$$\mae{\omega '}{\check{F} }{\bc}{(\alpha ,x)}{\alpha (\omega )}\;.$$
Then $\omega '\in \Omega _{\check{F} }\setminus \Omega '$ and $\vartheta (\omega ')=\omega $, so $\vartheta |(\Omega _{\check{F}}\setminus \Omega ' )$ is surjective.

Let $\omega _1,\omega _2\in \Omega _{\check{F}}\setminus \Omega ' $, $\omega _1\not=\omega _2$. There is an $(\alpha ,x)\in \check{F} $ with 
$$\langle(\alpha ,x),\omega _1\rangle\not=\langle(\alpha ,x),\omega _2\rangle\;.$$
Since $\langle(\alpha ,x),\omega _j\rangle=\langle\alpha ,\omega_j\rangle$ for every $j\in \z{1,2}$, $\vartheta |(\Omega _{\check{F} }\setminus \Omega ')$ is injective.

e) follows from d) since in this case $\Omega '$ is clopen.\qed

{\it Remark.} The above d) may be seen as a kind of generalization of Alexandroff's compactification.

\begin{de}\label{924}
We put for every $E$-C*-algebra $F$ 
$$\mae{\iota ^F}{F}{\check{F} }{x}{(0,x)}\,,$$
$$\mae{\pi^F}{\check{F} }{E}{(\alpha ,x)}{\alpha }\,,$$
$$\mae{\lambda ^F}{E}{\check{F} }{\alpha }{(\alpha ,0)}\,,$$
$$\sigma ^F:=\lambda ^F\circ \pi ^F\;.$$
\end{de}

If $E=\bc$ then 
$$\check{F}=\tilde{F}\,,\qquad \iota _F=\iota ^F\,,\qquad \pi _F=\pi ^F\,,\qquad \lambda _F=\lambda ^F\;.$$

All these maps are $E$-linear C*-homomorphisms,
$$\pi ^F\circ \iota ^F=0,\qquad \pi ^F\circ \lambda ^F=id_E,\qquad \pi ^F\circ \sigma ^F=\pi ^F\,,$$
$\iota ^F$ and $\lambda ^F$ are injective, $\pi ^F,\;\lambda ^F$, and $\sigma ^F$ are unital, and
$$\od{F}{\iota ^F}{\check{F} }{\pi^F}{\lambda ^F}{E}$$
is a split exact sequence in $\frm$.

\begin{p}\label{34}
\rule{1ex}{0em}
\begin{enumerate}
\item If $\oaa{F}{\varphi }{F'}$ is a morphism in \frm then the map
$$\mae{\check \varphi }{\check F}{\check F'}{(\alpha ,x)}{(\alpha ,\varphi x)}$$
is an involutive unital algebra homomorphism, injective or surjective if $\varphi $ is so. If $F=F'$ and if $\varphi $ is the identity map then $\check \varphi $ is also the identity map.
\item Let $F_1,F_2,F_3$ be $E$-C*-algebras and let $\varphi :F_1\rightarrow F_2$ and $\psi :F_2\rightarrow F_3$ be $E$-linear C*-homomorphisms. Then $\check {\overbrace{\psi \circ \varphi }} =\check \psi \circ \check \varphi $.\qed 
\end{enumerate}
\end{p}

{\it Remark.} If $E=\bc$ then $\check{\varphi }=\tilde{\varphi }  $.

\begin{e}\label{18.3'}
Let $F$ be a \en and $F'$ a closed ideal of $F$.
\begin{enumerate}
\item $F'$ endowed with the exterior multiplication
$$\mad{E\times F'}{F'}{(\alpha ,x)}{\alpha x}$$
is an \eo.
\item The map
$$\mad{\check{F'} }{E\times F}{(\alpha ,x)}{(\alpha ,\alpha +x)}$$
is an injective $E$-linear C*-homomorphism with image 
$$\me{(\alpha ,x)\in E\times F}{\alpha -x\in F'}\;.$$
\item \frc is a full subcategory of \frm.\qed
\end{enumerate}
\end{e}

\begin{p}\label{12.3'}
Let $F$ be a \en and $J$ a finite set.
\begin{enumerate}
\item $F^J=F\otimes l^2(J)$ endowed with the maps
$$\mad{F\times F^J}{F^J}{(x,\xi )}{(x\xi _j)_{j\in J}}\,,$$
$$\mad{F^J\times F}{F^J}{(\xi,x )}{(\xi _jx)_{j\in J}}\,,$$
$$\mad{F^J\times F^J}{F}{(\xi ,\eta )}{\si{j\in J}\eta _j^*\xi _j}$$
is a unital Hilbert $F$-module \emph{([C1] Proposition 5.6.4.2 c))}.
\item Let $\ccc{L}(F^J)$ be the Banach space of operators on $F^J$. The set $\ccc{L}_F(F^J)$ of adjointable operators on $F^J$ is a Banach subspace of $\ccc{L}(F^J)$. $\ccc{L}_F(F^J)$ endowed with the restriction of the norm of $\ccc{L}(F^J)$ it is a \en \emph{([C1] Theorem 5.6.1.11 d), [C1] Proposition 5.6.1.8 g),h))}.\qed
\end{enumerate}
\end{p}

\begin{p}\label{28.9'g}
For every $E$-C*-algebra $F$ the sequence
$$\og{\kk{i}{F}}{\kk{i}{\iota ^F}}{\kk{i}{\check{F} }}{\kk{i}{\pi ^F}}{\kk{i}{\lambda ^F}}{\kk{i}{E}}{20}{20}{20}$$
is split exact and the map
$$\mad{\kk{i}{F}\times \kk{i}{E}}{\kk{i}{\check{F} }}{(a,b)}{\kk{i}{\iota ^F}a+\kk{i}{\lambda ^F}b}$$
is a group isomorphism.
\end{p}

Since the sequence in \frm
$$\od{F}{\iota ^F}{\check{F} }{\pi ^F}{\lambda ^F}{E}$$
is split exact the assertion follows from the split exact axiom (\axi{27.9'a}).\qed

\begin{co}\label{23.4'a}
Let $G$ be a C*-algebra. 
\begin{enumerate}
\item The sequence in \frm
$$\og{F\otimes G}{\iota ^F\otimes id_G}{\check{F}\otimes G}{\pi ^F\otimes id_G}{\lambda ^F\otimes id_G}{E\otimes G}{20}{20}{20}$$
is split exact.
\item The sequence
$$\og{K_i(F\otimes G)}{K_i\left(\iota ^F\otimes id_G\right)}{K_i\left(\check{F}\otimes G \right) }{K_i\left(\pi ^F\otimes id_G\right)}{K_i\left(\lambda ^F\otimes id_G\right)}{K_i(E\otimes G)}{30}{30}{30}$$
is split exact and the map
$$K_i(E\otimes G)\times K_i(F\otimes G)\longrightarrow K_i\left(\check{F}\otimes G \right)\,,$$
$$(a,b)\longmapsto K_i\left(\lambda ^F\otimes id_G\right)a+K_i\left(\iota ^F\otimes id_G\right)b$$
is a group isomorphism.
\item Let $\oaa{F}{\varphi }{F'}$ be a morphism in \frm and $\oaa{G}{\psi }{G'}$ a morphism in $\fr{M}_{\bc}$. If we identify the isomorphic groups of b) then
$$\mac{\kk{i}{\check{\varphi }\otimes \psi  }}{\kk{i}{\check{F}\otimes G }}{\kk{i}{\check{F'}\otimes G' }},$$
$$(a,b)\longmapsto (\kk{i}{id_E\otimes \psi }a,\kk{i}{\varphi \otimes \psi }b)$$
is a group isomorphism.
\end{enumerate}
\end{co}

a) follows from \pr{23.4'} a).

b) follows from a) and the split exact axiom (\axi{27.9'a}).

c) follows from b) and the commutativity of the following diagram:
$$\begin{CD}
F\otimes G@>\iota ^F\otimes id_G>>\check{F}\otimes G @<\lambda ^F\otimes id_G<<E\otimes G\\
@V\varphi \otimes \psi VV@V\check{\varphi }\otimes \psi  VV         @VVid_E\otimes \psi V\\
F'\otimes G'@>>\iota ^{F'}\otimes id_{G'}>    \check{F'}\otimes G'  @<<\lambda ^{F'}\otimes id_{G'}<E\otimes G'\\
\end{CD}\hspace{1cm}.\qedd$$

\begin{co}\label{3.11'}
Let $\oaa{F}{\phi _1}{F'}$ and $\oaa{F}{\phi _2}{F'}$ be morphisms in \frm. If $F$ is K-null then $\kk{i}{\check{\phi _1} }=\kk{i}{\check{\phi _2} }$.
\end{co}

By \pr{28.9'g}, the map
$$\mad{\kk{i}{F}\times \kk{i}{E}}{\kk{i}{\check{F} }}{(a,b)}{\kk{i}{\iota ^F}a+\kk{i}{\lambda ^F}b}$$
is a group isomorphism. Since $F$ is K-null, $\kk{i}{\lambda ^F}$ is a group isomorphism. We get from $\check{\phi _1}\circ \lambda ^F=\check{\phi _2}\circ \lambda ^F $,
$$\kk{i}{\check{ \phi} _1}\circ \kk{i}{\lambda ^F}=\kk{i}{\check{ \phi}  _2}\circ \kk{i}{\lambda ^F}\,,\qquad \kk{i}{\check{\phi _1} }=\kk{i}{\check{\phi _2} }\;.\qedd$$

{\center{\section{Continuity and stability}}}

\begin{ax}[Continuity axiom]\label{5.10'}
If $\{(F_j)_{j\in J},\,(\varphi _{j,k})_{j,k\in I}\}$ is an inductive system in \frm such that $\varphi _{j,k}$ are injective for all $j,k\in J$, $k<j$, and if $\{F,\,(\varphi _j)_{j\in J}\}$ denotes its inductive limit in \frm then
$\{\kk{i}{F},\,(\kk{i}{\varphi _j})_{j\in J}\}$ is the inductive limit of the inductive system $\{(\kk{i}{F_j})_{j\in J},\,(\kk{i}{\varphi _{j,k}})_{j,k\in J}\}$.
\end{ax}

\begin{p}\label{19.2}
If $\Omega $ is a totally disconnected compact space then 
$$K_i(\ccc{C}(\Omega ,F))\approx \me{a\in K_i(F)^\Omega }{a(\Omega)\; {\emph{is finite}}}\;.$$ 
\end{p}

Let $\Xi $ be the set of clopen partitions of $\Omega $ ordered by fineness and for every $\Theta :=(\Omega _j)_{j\in J}\in \Xi $ and $x\in F^\Theta $ put
$$\mae{\tilde{x} }{\Omega }{F}{\omega }{x(j)} \quad \mbox{for}\quad \omega \in \Omega _j\;.$$
Then the map
$$\mad{F^\Theta }{\ccc{C}(\Omega ,F)}{x}{\tilde{x} }$$
is an injective $E$-C*-homomorphism for every $\Theta \in \Xi $ and $\ccc{C}(\Omega ,F)$ is isomorphic to the corresponding inductive limit in \frm of $(F^\Theta  )_{\Theta \in \Xi }$. By \lm{14.9'} c), $K_i(F^\Theta   )\approx K_i(F)^\Theta $ for every $\Theta \in \Xi $ and the assertion follows from the continuity axiom (\axi{5.10'}).\qed

\begin{p}\label{5.9}
Let $\xi $ be an ordinal number, $(\Omega _\eta )_{\eta <\xi }$ a family of path connected, non-compact, locally compact spaces, and $\omega _\eta \in \Omega _\eta $ for every $\eta <\xi $. We denote by $\Omega^ \xi $ the locally compact space obtained by endowing the disjoint union of the family of sets $(\Omega _\eta )_{\eta <\xi }$ with the topology for which a subset $U$ of $\Omega^ \xi $ is open if it has the following properties:
\renewcommand{\labelenumi}{\arabic{enumi})}
\begin{enumerate}
\item $\Omega _\eta \cap U$ is open for every $\eta <\xi $.
\item If $\omega _\eta \in U$ for some $\eta <\xi $ and if there is a $\zeta <\eta $ with $\eta =\zeta +1$ then $\Omega _\zeta \setminus U$ is compact.
\item If $\omega _\eta \in U$ for some limit ordinal number $\eta <\xi $ then there is a $\zeta <\eta $ such that $\bigcup_{\zeta <\zeta '<\eta }\Omega _{\zeta '}\subset U $.
\end{enumerate}
\renewcommand{\labelenumi}{\alph{enumi})} 
\renewcommand{\labelenumii}{\alph{enumi}_{\arabicenumii}))}
If $K_i(\cbb{\Omega _\eta }{F})=0$ for all $\eta <\xi $ then $K_i\left(\cbb{\Omega^ \xi }{F}\right)=0$.
\end{p}

The assertion is trivial for $\xi =0$. We prove the general case by transfinite induction. If $\xi =\eta +1$ for some $\eta <\xi $ for which the assertion holds then by \cor{3.9}, the assertion holds also for $\xi $. If $\xi $ is a limit ordinal number and the assertion holds for every $\eta <\xi $ then by the continuity axiom (\axi{5.10'}) the assertion holds also for $\xi $ since $\cbb{\Omega^ \xi }{F}$ is the inductive limit of the inductive system $\me{\cbb{\Omega ^\eta }{F}}{\eta <\xi }$.\qed

{\it Remark.} If $\Omega _\eta =[0,1[$ for every $\eta <\xi $ then $\Omega ^\xi $ is "one-dimensional".

\begin{lem}\label{7.10'}
Let $\{(F_j)_{j\in J},\,(\varphi _{j,k})_{j,k\in J}\}$ be an inductive system in \frm, $\{F,\,(\varphi _j)_{j\in J}\}$ its inductive limit in \frm, $G$ an $E$-C*-algebra, and for every $j\in J$ an injective morphism $\mac{\psi _j}{F_j}{G}$ in \frm such that $\psi _j=\psi _k\circ \varphi _{k,j}$ for all $j,k\in J$, $j<k$. Then the morphism $\mac{\psi }{F}{G}$ in \frm such that $\psi _j=\psi \circ \varphi _j$ for all $j,k\in J$, $j<k$, \emph{([W] Theorem L.2.1)} is injective.
\end{lem}

For $j\in J$ and $x\in F_j$,
$$\n{\varphi _jx}\leq \n{x}=\n{\psi _jx}=\n{\psi \varphi _jx}\leq \n{\varphi _jx}\,,$$
so $\psi $ preserves the norms on $\varphi _j(F_j)$. Since $\bigcup\limits_{j\in J}\varphi _j(F_j)$ is dense in $F$, $\psi $ preserves the norms, i.e. it is injective.\qed 

\begin{p}\label{7.10'a}
Let $\{(G_j)_{j\in J},\,(\varphi _{j,k})_{j,k\in J}\}$ be an inductive system in \frcc such that $\varphi _{j,k}$ are injective for all $j,k\in J$, $k<j$, and let $\{G,\,(\varphi _j)_{j\in J}\}$ be its inductive limit in \frcc. If $\{F',\,(\psi _j)_{j\in J}\}$ denotes the inductive limit in \frm of the inductive system $\{(F\otimes G_j)_{j\in J},\,(id_F\otimes \varphi _{j,k})_{j,k\in J}\}$ in \frm and $\mac{\psi }{F'}{F\otimes G}$ denotes the morphism in \frm such that $\psi \circ \psi _j=id_F\otimes \varphi _j$ for all $j\in J$ \emph{([W] Theorem L.2.1)} then $\psi $ is an isomorphism.
\end{p}

By [W] Corollary T.5.19, $id_F\otimes \varphi _j$ are injective for all $j\in J$. By \lm{7.10'}, $\psi $ is injective. Since
$$F\otimes\left( \bigcup_{j\in J}G_j\right)\subset Im\,\psi\,,$$ 
$\psi $ is  surjective and so it is an isomorphism.\qed

\begin{co}\label{7.10'b}
If $\{(G_j)_{j\in J},\,(\varphi _{j,k})_{j,k\in J}\}$ is an inductive system in \frcc $\;$ such that $\varphi _{j,k}$ are injective for all $j,k\in J$, $k<j$, and if $\{G,\,(\varphi _j)_{j\in J}\}$ is its inductive limit in \frcc then $\{\kk{i}{F\otimes G},\,(\kk{i}{id_F\otimes \varphi _j})_{j\in J}\}$ is the inductive limit of the inductive system $\{(\kk{i}{F\otimes G_j})_{j\in J},\,(\kk{i}{id_F\otimes \varphi _{j,k}})_{j,k\in J}\}$. In particular if $G_j$ is $\Upsilon $-null for every $j\in J$ then $G$ is also $\Upsilon $-null.
\end{co}

By [W] Corollary T.5.19, $id_F\otimes \varphi _{j,k}$ are injective for all $j,k\in J$, $k<j$. By \pr{7.10'a}, $\{F\otimes G,\,(id_F\otimes \varphi_ j)_{j\in J}\}$ may be identified with the inductive limit in \frm of the inductive system $\{(F\otimes G_j)_{j\in J},\,(id_F\otimes \varphi _{j,k})_{j,k\in J}\}$ in \frm and the assertion follows from the continuity axiom (\axi{5.10'}).\qed

\begin{co}\label{11.10'}
Let $(G_j)_{j\in J}$ be an infinite family in $\Upsilon _1$, $\fr{J}$ the set of nonempty finite subsets of $J$ ordered by inclusion, and for all $K,L\in \fr{J}$, $K\subset L$, put $G_K:=\bigotimes\limits _{j\in K}G_j$ and
$$\mae{\varphi (L,K)}{G_K}{G_L}{\bigotimes _{j\in K}x_j}{\bigotimes _{j\in L}y_j}\,,$$
where 
$$y_j:=\ab{x_j}{j\in K}{1_{G_j}}{j\in L\setminus K}\;.$$
Then $\{(G_K)_{K\in \fr{J}},\,(\varphi (L,K))_{K,L\in \fr{J}}\}$ is an inductive system in \frcc and its limit belongs to $\Upsilon _1$.
\end{co}

We denote by $\{G,\,(\varphi (K))_{K\in \fr{J}}\}$ the above inductive limit. By \pr{10.10'b}, $G_K\in \Upsilon _1$ for all $K\in \fr{J}$ so by \cor{7.10'b}, $p(G)=1$, $q(G)=0$. Let $\oaa{F}{\phi }{F'}$ be a morphism in \frm and let $K\in \fr{J}$. Then the diagram
$$\begin{CD}
F@>\phi _{G_K,F}>>F\otimes G_K@>id_F\otimes \varphi(K)>>F\otimes G\\
@V\phi VV @VV\phi \otimes id_{G_K}V @VV\phi \otimes id_G V\\
F'@>>\phi _{G_K,F'}>F'\otimes G_K@>>id_{F'}\otimes \varphi (K)>F'\otimes G
\end{CD}$$
is commutative. Since
$$\phi _{G,F}=(id_F\otimes \varphi (K))\circ \phi _{G_K,F}\,,\qquad \phi _{G,F'}=(id_{F'}\otimes \varphi (K))\circ \phi _{G_K,F'}\,,$$
the diagrams

\parbox{2cm}{
$$\begin{CD}
F@>\phi_{G,F} >>F\otimes G\\
@V\phi VV @VV\phi \otimes id_GV\\
F'@>>\phi_{G,F'} >F'\otimes G  
\end{CD}$$}
\hspace{3cm}
\parbox{3cm}{
$$\begin{CD} 
\kk{i}{F}@>\kk{i}{\phi_{G,F}} >>\kk{i}{F\otimes G}\\
@V\kk{i}{\phi} VV @VV\kk{i}{\phi \otimes id_G}V\\
\kk{i}{F'}@>>\kk{i}{\phi_{G,F'}} >\kk{i}{F'\otimes G}  
\end{CD}$$}

\hspace{-0.6cm} are commutative and so $G\in \Upsilon _1$.\qed

\begin{co}\label{9.11'}
Let $\{(G_j)_{j\in J},\,(\varphi _{j,k})_{j,k\in J}\}$ be an inductive system in \frcc $\;$ such that $\varphi _{k,j}$ are injective for all $j,k\in J$, $j<k$, and let $\{G,\,(\varphi _j)_{j\in J}\}$ be its inductive limit. We assume that for all $j,k\in J$, $j<k$,
$$G_j,G_k\in \Upsilon \,,\qquad\qquad \Phi _{i,G_k,F}=\kk{i}{id_F\otimes \varphi _{k,j}}\circ \Phi _{i,G_j,F}\;.$$
Then
$$G\in \Upsilon \,,\qquad\qquad \Phi _{i,G,F}=\kk{i}{id_F\otimes \varphi _j}\circ \Phi _{i,G_j,F}$$
for all $j\in J$.
\end{co}

By \cor{7.10'b}, $\{\kk{i}{F\otimes G},\,(\kk{i}{id_F\otimes \varphi _j})_{j\in J}\}$ is the inductive limit of the inductive system $\{(\kk{i}{F\otimes G_j})_{j\in J},\,(\kk{i}{id_F\otimes \varphi _{j,k}})_{j,k\in J}\}$. By the hypothesis of the Corollary,
$$\mac{\kk{i}{id_F\otimes \varphi _{k,j}}}{\kk{i}{F\otimes G_j}}{\kk{i}{F\otimes G_k}}$$
is a group isomorphism for all $j,k\in J$, $j<k$, so
$$\mac{\kk{i}{id_F\otimes \varphi _j}}{\kk{i}{F\otimes G_j}}{\kk{i}{F\otimes G}}$$
is also a group isomorphism for all $j\in J$. Let $\oaa{F}{\phi }{F'}$ be a morphism in \frm. The assertion follows from the commutativity of the diagram
$$\begin{CD}
\kk{i}{F}^{p(G_j)}\times \kk{i+1}{F}^{q(G_j)}@>\kk{i}{\phi }^{p(G_j)}\times \kk{i+1}{\phi }^{q(G_j)}>>A\\
@V\Phi _{i,G_j,F}VV@V\Phi _{i,G_j,F'}VV \\ 
\kk{i}{F\otimes G_j}@>\kk{i}{\phi \otimes G_j}>>\kk{i}{F'\otimes G_j}\\ 
@V\kk{i}{id_F\otimes \varphi _j}VV@V\kk{i}{id_{F'}\otimes \varphi _j}VV\\
\kk{i}{F\otimes G}@>>\kk{i}{\phi \otimes id_G}>\kk{i}{F'\otimes G}
\end{CD}$$
where $A:=\kk{i}{F'}^{p(G_j)}\times \kk{i+1}{F'}^{q(G_j)}$.\qed

\begin{de}\label{18.5}
We denote for every family $(\ccc{G}_j)_{j\in J}$ of additive groups by $\si{j\in J}\ccc{G}_j$ its direct sum i.e.
$$\si{j\in J}\ccc{G}_j:=\me{a\in \pro{j\in J}\ccc{G}_j}{\me{j\in J}{a_j\not=0}\;\emph{is finite}}\;.$$
\end{de}

\begin{p}\label{5.5}
If $(F_j)_{j\in J}$ is a family of $E$-C*-algebras and $F$ is its C*-direct sum \emph{([C1] Example 4.1.1.6)} then 
$$\kk{i}{F}\approx \si{j\in J}\kk{i}{F_j}\;.$$
In particular, the C*-direct sum of a family of K-null E-C*-algebras is K-null.
\end{p}

If $J$ is finite then the assertion follows from \pr{28.9'a}. The general case follows now from the continuity (\axi{5.10'}).\qed

\begin{co}\label{5.10'a}
If $(\Omega _j)_{j\in J}$ is a family of locally compact spaces and $\Omega $ is its topological sum then 
$$K_i(\cbb{\Omega }{F})\approx\si{j\in J}\kk{i}{\cbb{\Omega _j}{F}}\approx $$
$$\approx  \me{a\in \pro{j\in J}K_i(\cbb{\Omega _j}{F})}{\me{j\in J}{a_j\not=0} \emph{is finite}}\;.$$
\end{co}

By \pr{7.10'a}, $\cbb{\Omega }{F}$ is the direct sum of the family $(\cbb{\Omega _j}{F})_{j\in J}$ and the assertion follows from \pr{5.5}.\qed

\begin{p}\label{17.1'a}
If $\xi $ is an ordinal number endowed with its usual topology then $K_i(\cbb{\xi }{F})\approx \si{\eta <\xi }K_i(F)$.
\end{p}  

We prove the assertion by transfinite induction. If $\xi $ is not a limit ordinal number then the assertion follows from \cor{5.5a} a). Assume $\xi $ is a limit ordinal number and for all $\eta <\zeta <\xi $ let $\mac{\varphi _{\zeta,\eta  }}{\cbb{\eta }{F}}{\cbb{\zeta }{F}}$ be the inclusion map.
By \pr{7.10'a}, $\cbb{\xi }{F}$ may be identified with the inductive limit in \frm of the inductive system $\z{(\cbb{\eta }{F})_{\eta <\xi },\;(\varphi _{\zeta,\eta  })_{\eta <\zeta <\xi }}$
  in \frm. Thus the assertion follows from the continuity axiom (\axi{5.10'}) and the induction hypothesis.\qed
  
\begin{de}\label{9.11'a}
We denote for every $n\in \bn$ by $M(n)$ the C*-algebra of $n\times n$-matrices with entries in $\bc$.
\end{de}

\begin{ax}[Stability axiom]\label{5.10'b}
There is an $h\in \bn$, $h\not=1$, such that 
$$M(h)\in \Upsilon \,,\qquad p(M(h))=1\,,\qquad q(M(h))=0\,,$$
$$\Phi _{i,M(h),F}=\kk{i}{id_F\otimes \varphi }\circ \Phi _{i,\bc,F}\,,$$
where
$$\mae{\varphi }{\bc}{M(h)}{\alpha }{\left(\begin{array}{cccc}
\alpha &0&\cdots&0\\
0&0&\cdots&0\\
\vdots&\vdots&&\vdots\\
0&0&\cdots&0
\end{array}\right)}\;.$$
\end{ax}

\begin{p}\label{5.11'}
We put for all $j,k\in \bn^*$, $j<k$,
$$\mae{\varphi _{k,j}}{M(h^j)}{M(h^k)}{x}{\left(\begin{array}{cccc}
x&0&\cdots&0\\
0&0&\cdots&0\\
\vdots&\vdots&&\vdots\\
0&0&\cdots&0
\end{array}\right)}\;.$$
\begin{enumerate}
\item For all $j\in \bn$,
$$M(h^j)\in \Upsilon \,,\qquad p(M(h^j))=1\,,\qquad q(M(h^j))=0\,,$$
$$\Phi _{i,M(h^j),F}=\kk{i}{id_F\otimes \varphi _{j,0}}\circ \Phi _{i,\bc,F}\;.$$
\item For all $j,k\in \bn^*$, $j<k$,
$$\Phi _{i,M(h^k),F}=\kk{i}{id_F\otimes \varphi _{k,j}}\circ \Phi _{i,M(h^j),F}$$
and $\kk{i}{id_F\otimes \varphi _{k,j}}$ is a group isomorphism.
\end{enumerate}
\end{p}

a) We prove the assertion by induction with respect to $j\in \bn$. For $j=1$ the assertion is exactly the Stability axiom (\axi{5.10'b}). Let $j>1$ and assume the assertion holds for $j-1$. With the notation of \pr{5.7'a} b),
$$\left(id_{F\otimes M(h)}\otimes \varphi _{(j-1),0}\right)\circ \phi _{\bc,F\otimes M(h)}\circ (id_F\otimes \varphi _{1,0})=id_F\otimes \varphi _{j,0}\,,$$
so by the above and by the induction hypothesis,
$$\kk{i}{id_F\otimes \varphi _{j,0}}\circ \Phi _{i,\bc,F}=$$
$$=\kk{i}{id_{F\otimes M(h)}\otimes \varphi _{(j-1),0}}\circ \Phi _{i,\bc,F\otimes M(h)}\circ \kk{i}{id_F\otimes \varphi _{1,0}}\circ \Phi _{i,\bc,F}=$$
$$=\Phi _{i,M(h^{j-1}),F\otimes M(h)}\circ \Phi _{i,M(h),F}\;.$$
Thus
$$\mac{\kk{i}{id_F\otimes \varphi _{j,0}}\circ \Phi _{i,\bc,F}}{\kk{i}{F}}{\kk{i}{F\otimes M(h^j)}}$$
is a group isomorphism. Let $\oaa{F}{\phi }{F'}$ be a morphism in \frm. Since the diagram
$$\begin{CD}
\kk{i}{F}@>\Phi _{i,\bc,F}>>\kk{i}{F\otimes M(1)}@>\kk{i}{id_F\otimes \varphi _{j,0}}>>\kk{i}{F\otimes M(h^j)}\\
@V\kk{i}{\phi} VV @VV\kk{i}{\phi \otimes id_{M(1)}}V @VV\kk{i}{\phi \otimes id_{M(h^j)}} V\\
\kk{i}{F'}@>>\Phi _{i,\bc,F'}>\kk{i}{F'\otimes M(1)}@>>\kk{i}{id_{F'}\otimes \varphi (j,0)}>\kk{i}{F'\otimes M(h^j)}
\end{CD}$$
is commutative, we may take
$$\Phi _{i,M(h^j),F}=\kk{i}{id_F\otimes \varphi _{j,0}}\circ \Phi _{i,\bc,F}\;.$$

b) By a),
$$\kk{i}{id_F\otimes \varphi _{k,j}}\circ \Phi _{i,M(h^j),F}=\kk{i}{id_F\otimes \varphi _{k,j}}\circ \kk{i}{id_F\otimes \varphi _{j,0}}\circ \Phi _{i,\bc,F}=$$
$$=\kk{i}{id_F\otimes \varphi _{k,0}}=\Phi _{i,M(h^k),F}\;.\qedd$$

\begin{theo}\label{5.10'c}
Let $H$ be an infinite-dimensional Hilbert space and $\ccc{K}(H)$ the C*-algebra of compact operators on $H$. Then 
$$\ccc{K}(H)\in \Upsilon \,,\qquad p(\ccc{K}(H))=1\,,\qquad q(\ccc{K}(H))=0\,,$$
$$\Phi _{i,\ccc{K}(H),F}=\kk{i}{id_F\otimes \varphi }\circ \Phi _{i,\bc,F}\,,$$
where $\mac{\varphi  }{\bc}{\ccc{K}(H)}$ is an inclusion map. 
\end{theo}

Let $\Xi $ be the set of subspaces of $H$ of dimension $h^j$ for some $j\in \bn^*$ ordered by inclusion and for every $K\in \Xi $ let $\pi _K$ be the orthogonal projection of $H$ on $K$ and $G_K:=\pi _K\ccc{K}(H)\pi _K$. We denote for all $K,L\in \Xi $, $K\subset L$, by 
$$\mac{\varphi _{L,K}}{G_K}{G_L}\,,\qquad\qquad \mac{\varphi _K}{G_K}{\ccc{K}(H)}$$
the inclusion maps. Then $\{(G_K)_{K\in\, \Xi }\,,\;(\varphi _{L,K})_{L,K\in\, \Xi }\}$ is an inductive system in \frcc and $\{\ccc{K}(H)
\,,(\varphi _K)_{K\in\, \Xi }\}$ is its inductive limit. By \pr{5.11'}, for $K,L\in \Xi $, $K\subset L$,
$$G_K,G_L\in \Upsilon \,,\qquad p(G_K)=p(G_L)=1\,,\qquad q(G_K)=q(G_L)=0\,,$$
$$\Phi _{i,G_L,F}=\kk{i}{id_F\otimes \varphi _{L,K}}\circ \Phi _{i,G_K,F}\,,$$
and $\kk{i}{id_F\otimes \varphi _{L,K}}$ is a group isomorphism. By \cor{9.11'}, for $K\in \Xi $,
$$\ccc{K}(H)\in \Upsilon \,,\qquad \Phi _{i,\ccc{K}(H),F}=\kk{i}{id_F\otimes \varphi _K}\circ \Phi _{i,G_K,F}\,,$$
so $p(\ccc{K}(H))=1\,,\;q(\ccc{K}(H))=0$.\qed
  
\begin{center}
\part{Projective K-theory}
\end{center}

\begin{center}
\fbox{\parbox{10.5cm}{Throughout this part we use the following notation: $T$ is a group, 1 is its neutral element, $K$ is the complex Hilbert space $l^2(T)$, $(T_n)_{n\in \bn}$ is an increasing sequence of finite subgroups of $T$ the union of which is $T$, $T_0:=\{1\}$, $E$ is a unital commutative C*-algebra, and $f$ is a Schur $E$-function for $T$ (\dd{703})}}
\end{center}

\thispagestyle{empty}

In the usual K-theory the orthogonal projections (used for $K_0$) and the unitaries (used for $K_1$) are identified with elements of the square matrices, which is not a  very elegant procedure from the mathematical point of view, but is justified as a very efficient pragmatic solution. It seems to us that in the present more complicated construction the danger of confusion produced by these identifications is greater and we decided to separate these three domains. Unfortunately this separation complicates the presentation and the notation. Moreover, we also do identifications! In general the stability does not hold.  We present in \h{949} (as an example) some strong conditions under which stability holds for $K_0$.

For projective representations of groups we use [C2] (but the groups will be finite here) and for the K-theory we use [R], the construction of which we follow step by step. In the sequel we give a list of notation used in this Part.

\renewcommand{\labelenumi}{\arabic{enumi})} 
\begin{enumerate}
\item We put for every involutive algebra $F$,
$$Pr\,F:=\me{P\in F}{P=P^*=P^2}$$
and for every $A\subset F$,
$$A^c:=\me{x\in F}{y\in A\Longrightarrow xy=yx}\;.$$
\item We denote for every unital involutive algebra $F$ by $1_F$ its unit and set
$$\unn{F}:=\me{U\in F}{UU^*=U^*U=1_F}\;.$$
\item If $F$ is a unital C*-algebra and $U,V\in \unn{F}$ then we denote by $U\sim _hV$ the assertion $U$ and $V$ are homotopic in $\unn{F}$ and put
$$\unm{F}:=\me{U\in \unn{F}}{U\sim _h1_F}\;.$$
Moreover $GL(F)$ denotes the group of invertible elements of $F$ and $GL_0(F)$ the elements of $GL(F)$ which are homotopic to $1_F$ in $GL(F)$.
\item If $F$ is a unital C*-algebra and $G$ is a unital C*-subalgebra of $F$ then we denote by  $Un_G\,F$ the set of elements of $\unn{F}$ which are homotopic to an element of $\unn{G}$ in $\unn{F}$ and by $GL_G(F)$ the set of elements of $GL(F)$ which are homotopic to an element of $GL(G)$ in $GL(F)$.
\item If $\Omega $ is a topological space, $F$ a C*-algebra, and $A\subset F$ then we put
$$\ccb{\Omega }{A}:=\me{X\in \ccb{\Omega }{F}}{\omega \in \Omega \Longrightarrow X(\omega )\in A}\;.$$
\item Hilbert \eo ([C1] Definition 5.6.1.4).
\item $\ccc{L}_E(H)$ ([C1] Definition 5.6.1.7).
\end{enumerate}

\renewcommand{\labelenumi}{\alph{enumi})}

\begin{center}
\chapter{Some notation and the axiom}
\end{center}

{\center{\section{Some notation and the axiom}}}

\begin{de}\label{703}
Let $S$ be a group and let $1$ be its neutral element. A {\bf Schur $E$-function for $S$ } is a map
$$\mac{f}{S\times S}{\unn{E}}$$
such that $f(1,1)=1_E$ and
$$f(r,s)f(rs,t)=f(r,st)f(s,t)$$
for all $r,s,t\in T$. We denote by $\f{S}{E}$ the set of Schur $E$-functions for $S$.
\end{de}

Schur functions are also called normalized factor set or multiplier or two-co-cycle (for $S$ with values in $\unn{E}$) in the literature. 

\begin{de}\label{a}Let $F$ be an \en and $n\in \bn^*$. We put
 for every $t\in T_n$, $\xi \in F^{T_n}=F\otimes l^2(T_n)$, and $x\in F$,
$$\mae{V_t\xi :=V_t^F\xi }{T_n}{F}{s}{f(t,t^{-1}s)\xi (t^{-1}s)}\,,$$
$$\mae{x\otimes id_K}{F^{T_n}}{F^{T_n}}{\xi }{(x\xi _s)_{s\in T_n}}\,,$$
so we have
$$\mae{(x\otimes id_K)V_t\xi }{T_n}{F}{s}{f(t,t^{-1}s)x\xi (t^{-1}s)}\;.$$
We define
$$F_n:=\me{\si{t\in T_n}(X_t\otimes id_K)V_t}{(X_t)_{t\in T_n}\in F^{T_n}}\;.$$
If $\oa{F}{\varphi }{G}$ is a morphism in \frc then we put
$$\mae{\varphi _n}{F_n}{G_n}{X}{\si{t\in T_n}((\varphi X_t)\otimes id_{K_n})V_t}\;.$$
\end{de}

$F_n$ is a full $E$-C*-subalgebra of $\ccc{L}_F(F^{T_n})$ (\pr{12.3'} b), [C2] Theorem 2.1.9 h), k)), so $1_{F_n}=1_E$, and $\varphi _n$ is an $E$-C*-homomorphism, injective or surjective if $\varphi $ is so ([C2] Corollary 2.2.5). Moreover
 $F_m$ is canonically a full $E$-C*-subalgebra of $F_n$ for every $m\in \bn^*$, $m<n$ ([C2] Proposition 2.1.2). For every $n\in \bn$, $F_n\times G_n\approx (F\times G)_n$.

\begin{ax}\label{b} We fix in Part II a sequence $(C_n)_{n\in \bn}\in \pro{n\in \bn}E_n$, put
$$A_n:=C_n^*C_n\,,\qquad\qquad B_n:=C_nC_n^*\,,$$
and assume $A_n,\,B_n\in Pr\,E_n$, $A_n+B_n=1_E=1_{E_n}$, and $C_n\in (E_{n-1})^c$ for every $n\in \bn$ (where we used the inclusion $E_{n-1}\subset E_n$ in the last relation).
\end{ax}

From
$$A_n=A_n(A_n+B_n)=A_n^2+A_nB_n=A_n+A_nB_n\,,$$
$$C_n=C_n(A_n+B_n)=C_nA_n+C_nB_n=C_n+C_n^2C_n^*$$
we get $A_nB_n=C_n^2=0$ for every $n\in \bn$.

We have $C_n\in (F_{n-1})^c$ for every $n\in \bn$ and for every full $E$-C*-algebra $F$ (where we used the inclusion $F_{n-1}\subset F_n$). 

\begin{e}\label{13.4}
Let $(S_m)_{m\in \bn}$ be a sequence of finite groups and $(k_n)_{\bbn}$ a strictly increasing sequence in $\bn$ such that $T_n=\proo{m=1}{k_n}S_m$ for all \bbn. We identify $S_m$ with a subgroup of $T$ for every $m\in \bn$. Assume that for every $m\in \bn$ there is a  $g_m\in \f{S_m}{E}$ such  that 
$$f(s,t)=\pro{m\in \bn}g_m(s_m,t_m)$$
for all $s,t\in T$. For every $n\in \bn$ let $m\in \bn$, $k_{n-1}<m\leq k_n$, let $\mac{\chi }{\bzz{2}\times \bzz{2}}{S_m}$ be an injective group homomorphism, and $\beta _1,\beta _2\in \unn{E}$. We put
$$a:=\chi (1,0)\,,\qquad b:=\chi (0,1)\,,\qquad \alpha _1:=f(a,a)\,,\qquad \alpha _2:=f(b,b)\,,$$
$$C_n:=\frac{1}{2}((\beta _1\otimes id_K)V_a^f+(\beta _2\otimes id_K)V_b^f)\;.$$
If $f(a,b)=-f(b,a)=1_E$ and
 $\alpha _1\beta _1^2+\alpha _2\beta _2^2=0$ then $(C_n)_{n\in \bn}$ fulfills the conditions of \emph{\axi{b}}.
\end{e}

The assertion follows from [C2] Theorem 2.2.18 a), b).\qed

{\it Remark 1.} If $E=\bc$, $S_m=\bzz{2}\times \bzz{2}$, and $k_m=m$ for every $m\in \bn$  then (by [C2] Proposition 3.2.1 c) and [C2] Corollary 3.2.2 d)) we may choose $(C_n)_{\bbn}$ in such a way that the corresponding K-theory coincides with the classical one.

{\it Remark 2.} Denote by $T_n$ the set of permutations $p$ of $\bn$ such that  $\me{j\in \bn}{p(j)\not=j}\subset \bnn{4n}$ so $T$ is the set of permutations $p$ of $\bn$ such that $\me{j\in \bn}{p(j)\not=j}$ is finite. This example shows that the given conditions for $T_n$ in \ee{13.4} are not automatically fulfilled.

\begin{center}
\chapter{The functor $K_0$}

\section{$K_0$ for $\frc$}
\end{center}

\begin{center}
\fbox{Throughout this section $F$ denotes a full $E$-C*-algebra}
\end{center}

\vspace{5ex}
\begin{p}\label{911}
Let $n\in \bn$.
\begin{enumerate}
\item $A_n,\,B_n\in (F_{n-1})^c$ (where we used the inclusion $F_{n-1}\subset F_n$).
\item $A_nF_nA_n$ is a unital C*-algebra with $A_n$ as unit.
\item The map
$$\mae{\bar{\rho }_n^F }{F_{n-1}}{F_n}{X}{A_nX=XA_n=A_nXA_n=C_n^*XC_n}$$
(where we used the inclusion $F_{n-1}\subset F_n$) is an $E$-linear injective C*-homomorphism.
\end{enumerate}
\end{p}
 
Only the injectivity of $\bar{\rho }_n^F $ needs a proof. Let $X\in F_{n-1}$ with $\bar{\rho }_n^FX=0 $. Then
$$C_n^*C_nX=0\,,\qquad XC_n=C_nX=0\,,$$
$$XB_n=XC_nC_n^*=0\,,\qquad X=X(A_n+B_n)=0\;.\qedd$$

{\it Remark.} $ \bar{\rho }_n^F$ is not unital since $\bar{\rho }_n^F1_E=A_n $. 

\begin{de}\label{912}
We put for all $m,n\in \bn,\,m<n$,
$$\rho _{n,m}^F:=\bar{\rho }_n^F\circ \bar{\rho }_{n-1}^F\circ \cdots \circ \bar{\rho }_{m+1}^F:F_m\longrightarrow F_n\;.   $$
Then $\{(F_n)_{n\in \bn},\,(\rho _{n,m}^F)_{n,m\in \bn}\}$ is an inductive system of full $E$-C*-algebras with injective $E$-linear (but not unital) maps. We denote by $\{F_\rightarrow ,\,(\rho _n^F)_{n\in \bn}\}$ its algebraic inductive limit. $F_\rightarrow $ is an involutive (but not unital) algebra endowed with the structure of an algebraic $E$-C*-algebra, $\rho _n^F$ is injective and $E$-linear for every $n\in \bn$, and $(Im\,\rho _n^F)_{n\in \bn}$ is an increasing sequence of involutive subalgebras and algebraic $E$-C*-subalgebras of $F_\rightarrow $ the union of which is $F_\rightarrow $. We put for every $X\in F_n$,
$$X_\rightarrow :=X_{\rightarrow \,n}:=X_{\rightarrow \,n}^F:=\rho _n^FX\,,$$
and
$$1_{\rightarrow \,n}:=1_{\rightarrow \,n}^F:=\rho _n^F1_{F_n}=\rho _n^F1_E \,,$$
$$F_{\rightarrow \,n}:=Im\,\rho _n^F\;.$$
In particular
$$(A_n)_\rightarrow =\rho _n^FA_n=1_{\rightarrow ,n-1},\qquad (B_n)_\rightarrow =\rho _n^FB_n,\qquad (C_n)_\rightarrow =\rho _n^FC_n\;.$$

We put
$$Pr\,F_\rightarrow :=\me{P\in F_\rightarrow }{P=P^*=P^2}=\bigcup _{n\in \bn}(Pr\,F_{\rightarrow \,n})\;.$$
For $P,Q\in Pr\,F_\rightarrow $ we put $P\sim _0Q$ if there is an $X\in F_\rightarrow $ with $X^*X=P$, $XX^*=Q$ (in this case there is an $n\in \bn$ such that $P,Q,X\in F_{\rightarrow \,n}$); $\sim _0$ is the Murray - von Neumann equivalence relation, which we shall use also in the case of C*-algebras. For every $P\in Pr\,F_\rightarrow $ we denote by $\dot{P} $ its equivalence class in $Pr\,F/\!\sim _0$.
\end{de}

Often we shall identify $F_n$ with $F_{\rightarrow \,n}$ by using $\rho _n^F$. By this identification $F_{\rightarrow \,n}$  is a full $E$-C*-algebra with $1_{\rightarrow \,n}$ as unit.

$F_\rightarrow $ is also endowed with a C*-norm and its completion in this norm is the C*-inductive limit of the above inductive system, but we shall not use this supplementary structure in the sequel. 

\begin{p}\label{914}
If $n\in \bn$ and $P\in Pr\,F_{\rightarrow ,\,n-1}$ then
$$P=(A_n)_\rightarrow P\sim _0(B_n)_\rightarrow P=(C_n)_\rightarrow P(C_n)_\rightarrow ^*\;.$$
\end{p}

We have
$$((C_n)_\rightarrow P)^*((C_n)_\rightarrow P)=P(C_n)_\rightarrow ^*(C_n)_\rightarrow P=(A_n)_\rightarrow P\,,$$
$$((C_n)_\rightarrow P)((C_n)_\rightarrow P)^*=P(C_n)_\rightarrow (C_n)_\rightarrow ^*P=(B_n)_\rightarrow P\,,$$
so $(A_n)_\rightarrow P\sim _0(B_n)_\rightarrow P$.\qed

\begin{p}\label{915}
For every finite family $(P_i)_{i\in I}$ in $Pr\,F_\rightarrow $ there is a family $(Q_i)_{i\in I}$ in $Pr\,F_\rightarrow $ such that $P_i\sim _0Q_i$ for every $i\in I$ and $Q_iQ_j=0$ for all distinct $i,j\in I$.
\end{p}

We prove the assertion by complete induction with respect to $Card\, I$. Let $i_0\in I$ and put $J:=I\setminus \{i_0\}$. We may assume, by the induction hypothesis, that there is an $n\in \bn$ with $P_i\in Pr\,F_{\rightarrow  ,\,n-1}$ for all $i\in I$ and $P_iP_j=0$ for all distinct $i,j\in J$. By \pr{914}, 
$$P_{i_0}=(A_n)_\rightarrow P_{i_0}\sim _0(C_n)_\rightarrow P_{i_0}(C_n)_\rightarrow ^*=:Q_{i_0}\,,$$
and
$$Q_{i_0}P_j=(C_n)_\rightarrow P_{i_0}(C_n)_\rightarrow ^*(A_n)_\rightarrow P_j=(C_n)_\rightarrow P_{i_0}(C_n^*A_n)_\rightarrow P_j=0$$
for all $j\in J$.\qed

\begin{p}\label{916}
Let $P,Q\in Pr\,F_\rightarrow $.
\begin{enumerate}
\item If $P',P'',Q',Q''\in Pr\,F_\rightarrow $ such that
$$P\sim _0P'\sim _0P'',\qquad Q\sim _0Q'\sim _0Q'',\qquad P'Q'=P''Q''=0$$
then
$$P'+Q'\sim _0P''+Q''\;. $$
We put 
$$\dot{P}\oplus \dot{Q}:=\dot{\overbrace{P'+Q'}}\;.   $$
\item $Pr\,F_\rightarrow /\!\sim _0$ endowed with the above composition law $\oplus $ is an additive semi-group with $\dot{0} $ as neutral element. We denote by $K_0(F)$ its associated Grothendieck group and by 
$$\mac{[\;\cdot \;]_0}{Pr\,F_\rightarrow }{K_0(F)}$$
the Grothendieck map \emph{([R] 3.1.1)}.
\item $K_0(F)=\me{[P]_0-[Q]_0}{P,Q\in Pr\,F_\rightarrow }$.
\item For every $a\in K_0(F)$ there are $P,Q\in Pr\,F_\rightarrow $ and $n\in \bn$ such that
$$P=P(A_n)_\rightarrow,\qquad Q=Q(B_n)_\rightarrow ,\qquad a=[P]_0-[Q]_0\;. $$
\end{enumerate}
\end{p}

a) Let $X,Y\in F_\rightarrow $ with
$$X^*X=P',\qquad XX^*=P'',\qquad Y^*Y=Q',\qquad YY^*=Q''\;.$$
Then
$$0=P'Q'=X^*XY^*Y,\qquad 0=P''Q''=XX^*YY^*$$
so 
$$XY^*=X^*Y=0,\qquad (X+Y)^*(X+Y)=X^*X+Y^*Y=P'+Q'\,,$$
$$(X+Y)(X+Y)^*=XX^*+YY^*=P''+Q''\,,\qquad P'+Q'\sim _0P''+Q''\;.$$

b) and c) follow from a) and \pr{915}.

d) follows from c) and \pr{914}.\qed

\begin{co}\label{918}
The following are equivalent for all $n\in \bn$ and $P,Q\in Pr\,F_{\rightarrow\, n}$.
\begin{enumerate}
\item $[P]_0=[Q]_0$.
\item There is an $R\in Pr\,F_\rightarrow $ such that
$$PR=QR=0\,,\qquad\qquad P+R\sim _0Q+R\;.$$
\item There is an $m\in \bn$, $m>n+1$, such that
$$P+(B_m)_\rightarrow \sim _0Q+(B_m)_\rightarrow$$
or (by identifying $F_m$ with $F_{\rightarrow \,m}$)
$$\left(\proo{i=n+1}{m}A_i\right)P+\left(1_E-\proo{i=n+1}{m}A_i\right)\sim _0
\left(\proo{i=n+1}{m}A_i\right)Q+\left(1_E-\proo{i=n+1}{m}A_i\right)\;.$$
\end{enumerate}
\end{co} 

$a\Rightarrow b$ follows from \pr{915} (and from the definition of the Grothendieck group).

$b\Rightarrow c$. We may assume $R\in F_{\rightarrow ,\,m-1}$ for some $m>n+1$. By \pr{914},
$$P+(B_m)_\rightarrow R\sim _0P+R\sim _0 Q+R\sim _0Q+(B_m)_\rightarrow R\,,$$
so
$$P+(B_m)_\rightarrow =P+(B_m)_\rightarrow R+((B_m)_\rightarrow -(B_m)_\rightarrow R)\sim _0$$
$$\sim _0Q+(B_m)_\rightarrow R+((B_m)_\rightarrow -(B_m)_\rightarrow R)=Q+(B_m)_\rightarrow \;.$$
It follows
$$\left(\proo{i=n+1}{m}A_i\right)P+\left(1_E-\proo{i=n+1}{m}A_i\right)=\rho _{m,n}^FP+B_m+\left(A_m-\proo{i=n+1}{m}A_i\right)\sim _0$$
$$\sim _0\rho _{m,n}^FQ+B_m+\left(A_m-\proo{i=n+1}{m}A_i\right)=\left(\proo{i=n+1}{m}A_i\right)Q+\left(1_E-\proo{i=n+1}{m}A_i\right)\;.$$

c$\Rightarrow $a is trivial.\qed

\begin{co}\label{30.3'}
If for every $\bbn$ and $P\in \pp{F_{\rightarrow \,n}}$ there is an $m\in \bn$, $m>n+1$, such that $P+(B_m)_\rightarrow \sim _01_E$ then $K_0(F)=\z{0}$.
\end{co}

Let $P,Q\in \pp{F_\rightarrow }$. By our hypothesis there is an $m\in \bn$ such that $P+(B_m)_\rightarrow \sim _0Q+(B_m)_\rightarrow $. By \cor{918} $c\Rightarrow a$, $[P]_0=[Q]_0$. Thus by \pr{916} c), $K_0(F)=\z{0}$.\qed

\begin{co}\label{21.3'}
$K_0(E)\not=\z{0}$.
\end{co}

Assume $K_0(E)=\z{0}$. Then $[1_E]_0=[0]_0$, so by \cor{918} $a\Rightarrow c$, there is an $\bbn$ such that 
$$1_E\sim _01_E-\proo{i=1}{n}A_i\;.$$
Let $\omega $ be a point of the spectrum of $E$. Since $E_n(\omega )$ is a product of square matrices the above relation leads to a contradiction by using the trace function.\qed

\begin{p}\label{920}
Let $\ccc{G}$ be an additive group and $\nu :Pr\,F_\rightarrow \rightarrow \ccc{G}$ a map such that
\renewcommand{\labelenumi}{\arabic{enumi})}
\begin{enumerate}
\item $P,Q\in Pr\,F_\rightarrow ,\, PQ=0\;\Longrightarrow  \;\nu (P+Q)=\nu (P)+\nu (Q)$. 
\item $P,Q\in Pr\,F_\rightarrow ,\,P\sim _0Q\;\Longrightarrow \;\nu (P)=\nu (Q)$.
\renewcommand{\labelenumi}{\alph{enumi})} 
\renewcommand{\labelenumii}{\alph{enumi}_{\arabicenumii}))}
\end{enumerate}
Then there is a unique group homomorphism $\mu :K_0(F)\rightarrow \ccc{G}$ such that $\mu [P]_0=\nu (P)$ for every $P\in Pr\,F_\rightarrow $.
\end{p}

By 2), $\nu $ is well-defined on $Pr\,F_\rightarrow /\!\sim _0$ and by 1) and \pr{916} a),b), $\nu $ is an additive map on $Pr\,F_\rightarrow /\!\sim _0$. By 2) and \cor{918} a$\Rightarrow $b, $\nu $ is well-defined on $K_0(F)$. The existence and uniqueness of $\mu $ with the given properties follows now from \pr{916} c).\qed

\begin{p}\label{921}
Let $\oa{F}{\varphi }{G}$ be a morphism in $\frc$.
\begin{enumerate}
\item For $m,n\in \bn$, $m<n$, the diagram
\[ \begin{CD} 
F_m @>\rho _{n,m}^F>> F_n\\
@V\varphi _mVV        @VV\varphi _nV\\
G_m @>>\rho _{n,m}^G> G_n
\end{CD} \]
is commutative. Thus there is a unique $E$-linear involutive algebra homomorphism $\mac{\varphi _\rightarrow }{F_\rightarrow }{G_\rightarrow }$ with
$$\varphi _\rightarrow \circ \rho _n^F=\rho _n^G\circ \varphi _n$$
for every $n\in \bn$.

\item $\varphi _\rightarrow $ is injective or surjective if $\varphi $ is so.
\item There is a unique group homomorphism $\mac{K_0(\varphi )}{K_0(F)}{K_0(G)}$ such that
$$K_0(\varphi )[P]_0=[\varphi _\rightarrow P]_0$$
for every $P\in Pr\,F_\rightarrow $.
\item If $\varphi $ is the identity map then $K_0(\varphi )$ is also the identity map.
\item If $\varphi =0$ then $K_0(\varphi )=0$.
\end{enumerate}
\end{p}

a) It is sufficient to prove the assertion for $n=m+1$. For $X\in F_m$,
$$\varphi _n\bar{\rho }_n^FX=\varphi _n(A_nX)=A_n\varphi _nX=\bar{\rho }_n^G\varphi _nX  $$
(where we used the inclusion $F_m\subset F_n$).

b) follows from the fact that for every $n\in \bn$, $\varphi _n$ is injective or surjective if $\varphi $ is so ([C2] Theorem 2.1.9 a))).

c) By a) and \pr{914}, the map
$$\mad{Pr\,F_\rightarrow }{K_0(G)}{P}{[\varphi _\rightarrow P]_0}$$
possesses the properties from \pr{920}.

d) and e) are obvious.\qed

\begin{co}\label{923}
If $\ob{F}{\varphi }{G}{\psi }{H}$ are morphisms in $\frc$ then
$$(\psi \circ \varphi )_\rightarrow =\psi _\rightarrow \circ \varphi _\rightarrow ,\qquad K_0(\psi \circ \varphi )=K_0(\psi )\circ K_0(\varphi )\;.\qedd$$
\end{co}

\begin{p}\label{925}
\rule{0mm}{0mm}
\begin{enumerate}
\item The maps
$$\mae{\mu }{\check{F} }{F}{(\alpha ,x)}{\alpha +x}\,,$$
$$\mae{\lambda '}{E}{\check{F} }{\alpha }{(\alpha ,-\alpha )}$$
are $E$-C*-homomorphisms.
\item $$\mu \circ \iota ^F=id_F,\qquad\qquad \iota ^F\circ \mu +\lambda '\circ \pi ^F=id_{\check{F} }\,,$$
$$K_0(\iota ^F)\circ K_0(\mu )+K_0(\lambda ')\circ K_0(\pi ^F)=id_{K_0(\check{F} )}\;.$$
\item $$\og{K_0(F)}{K_0(\iota ^F)}{K_0(\check{F} )}{K_0(\pi ^F)}{K_0(\lambda ^F)}{K_0(E)}{20}{20}{20}$$
is a split exact sequence.
\end{enumerate}
\end{p}

a) is easy to see.

b) For $(\alpha ,x),(\beta ,y)\in \check{F} $,
$$\iota ^F\mu (\alpha ,x)=(0,\alpha +x),\qquad\quad \lambda '\pi ^F(\alpha ,x)=(\alpha ,-\alpha )\,,$$
$$(\iota ^F\mu (\alpha ,x))(\lambda '\pi ^F(\beta ,y))=(0,\alpha +x)(\beta ,-\beta )=(0,0)\,,$$
$$(\iota ^F\mu +\lambda '\pi ^F)(\alpha ,x)=(\alpha ,x)$$
so $\iota ^F\circ \mu +\lambda '\circ \pi ^F$ is a full $E$-C*-homomorphism and
$$\iota ^F\circ \mu +\lambda '\circ \pi ^F=id_{\check{F} }\;.$$
By a) and \cor{923},
$$\iota _\rightarrow ^F\circ \mu _\rightarrow +\lambda '_\rightarrow \circ \pi _\rightarrow ^F=id_{\check{F}_\rightarrow  }\;.$$
By \pr{921} c),d) and \cor{923}, for $P\in Pr\,\check{F} _\rightarrow $,
$$(K_0(\iota ^F)\circ K_0(\mu )+K_0(\lambda ')\circ K_0(\pi ^F))[P]_0=K_0(\iota ^F\circ \mu )[P]_0+K_0(\lambda '\circ \pi ^F)[P]_0=$$
$$=[\iota _\rightarrow ^F\mu _\rightarrow P]_0+[\lambda '_\rightarrow \pi _\rightarrow ^FP]_0=[(\iota ^F\circ \mu +\lambda '\circ \pi ^F)_\rightarrow P]_0=[P]_0$$
so by \pr{916} c),
$$K_0(\iota ^F)\circ K_0(\mu )+K_0(\lambda ')\circ K_0(\pi ^F)=id_{K_0(\check{F} )}\;.$$

c) By b), \pr{921} d),e), and \cor{923},
$$K_0(\pi ^F)\circ K_0(\iota ^F)=K_0(\pi ^F\circ \iota ^F)=0\,,$$
$$K_0(\pi ^F)\circ K_0(\lambda ^F)=K_0(\pi ^F\circ \lambda ^F)=id_ {K_0(E)}\,,$$
$$K_0(\mu )\circ K_0(\iota ^F)=K_0(\mu \circ \iota ^F)=id_{K_0(F)}$$
and so $K_0(\iota ^F)$ is injective. By b), for $a\in K_0(\check{F} )$,
$$a=K_0(\iota ^F)K_0(\mu )a+K_0(\lambda ')K_0(\pi ^F)a\;.$$
Thus if $a\in Ker\,K_0(\pi ^F)$ then $a=K_0(\iota ^F)K_0(\mu )a\in Im\,K_0(\iota ^F)$, and so $Ker\,K_0(\pi ^F)=Im\,K_0(\iota ^F)$.

\begin{center}
\section{$K_0$ for \frm}
\end{center}

\begin{de}\label{926}
Let $F$ be an $E$-C*-algebra and consider the split exact sequence
$$\od{F}{\iota ^F}{\check{F} }{\pi ^F}{\lambda ^F}{E}$$
introduced in \emph{\dd{924}}. We put
$$K_0(F):=Ker\,K_0(\pi ^F)\;.$$
\end{de}

By \pr{925} c), this definition does not contradict the definition given in \pr{916} b) for the case that $F$ is an \en.

$K_0(\{0\})=\{0\}$ since $\pi ^{\{0\}}$ is bijective.

\begin{p}\label{927}
Let \oa{F}{\varphi }{G} be a morphism in \frm.
\begin{enumerate}
\item The diagram
\[ \begin{CD}
F @>\iota ^F>> \check{F} @>\pi ^F>> E\\
@V\varphi VV   @VV\check{\varphi }V \parallel \\
G @>>\iota ^G> \check{G} @>>\pi^G>  E           
\end{CD} \]
is commutative.
\item The diagram
\[ \begin{CD}
K_0(F) @>\subset >> K_0(\check{F} ) @>K_0(\pi ^F)>> K_0(E)\\
@VK_0(\varphi )VV   @VVK_0(\check{\varphi } )V \parallel \\
K_0(G) @>>\subset > K_0(\check{G} ) @>>K_0(\pi ^G)> K_0(E)
\end{CD} \]
is commutative, where $K_0(\varphi )$ is defined by $K_0(\check{\varphi } )$.
\item If $P\in Pr\,F_\rightarrow $ then
$$K_0(\varphi )[P]_0=[\varphi _\rightarrow P]_0\;.$$
\item $K_0(id_F)=id_{K_0(F)}$.
\item If $\varphi =0$ then $K_0(\varphi )=0$.
\end{enumerate}
\end{p}

a) is obvious.

b) By a) and \cor{923}, the right part of the diagram is commutative. This implies the existence (and uniqueness) of $K_0(\varphi )$.

c) By a),\,b), \pr{921} a),c), and \cor{923},
$$K_0(\varphi )[P]_0=K_0(\check{\varphi } )[\iota _\rightarrow ^FP]_0=[\check{\varphi} _\rightarrow \iota _\rightarrow ^F P]_0=[\iota _\rightarrow ^G\varphi _\rightarrow P]_0=[\varphi _\rightarrow P]_0\;.$$

d) and e) follow from c) and \pr{916} c). \qed
 
\begin{co}\label{928}
Let $\ob{F}{\varphi }{G}{\psi }{H}$ be morphisms in \frm.
\begin{enumerate}
\item $K_0(\psi )\circ K_0(\varphi )=K_0(\psi \circ \varphi )$.
\item If $\varphi $ is an isomorphism then $K_0(\varphi )$ is also an isomorphism and 
$$K_0(\varphi )^{-1}=K_0(\varphi ^{-1})\;.$$
\end{enumerate}
\end{co}

a) follows from \pr{34} b), \cor{923}, and \pr{927} b).

b) follows from a) and \pr{927} d).\qed

\begin{p}\label{929}
For every \eo \,$F$,
$$K_0(F)=\me{[P]_0-[\sigma _\rightarrow ^FP]_0}{P\in Pr\,\check{F}_\rightarrow  }\;.$$
\end{p}

For $P\in Pr\,\check{F}_\rightarrow  $, by \pr{927} c) and \cor{923} (since $\pi ^F=\pi ^F\circ\, \sigma ^F$),
$$K_0(\pi ^F)[\sigma _\rightarrow ^FP]_0=[\pi ^F_\rightarrow \sigma _\rightarrow ^FP]_0=[\pi _\rightarrow ^FP]_0=K_0(\pi ^F)[P]_0$$
so
$$[P]_0-[\sigma _\rightarrow ^FP]_0\in Ker\,K_0(\pi ^F)=K_0(F)\;.$$

Let $a\in K_0(F)$. By \pr{916} d), there are $Q,R\in Pr\,\check{F}_\rightarrow  $ and $n\in \bn$ such that
$$Q=Q(A_n)_\rightarrow \,,\qquad R=R(B_n)_\rightarrow\, ,\qquad a=[Q]_0-[R]_0\;.$$
Then
$$a=[Q(A_n)_\rightarrow ]_0+[(B_n)_\rightarrow -R(B_n)_\rightarrow ]_0-([R(B_n)_\rightarrow ]_0-[(B_n)_\rightarrow -R(B_n)_\rightarrow ]_0)=$$
$$=[Q(A_n)_\rightarrow +((B_n)_\rightarrow -R(B_n)_\rightarrow )]_0-[(B_n)_\rightarrow ]_0\;.$$
If we put
$$P:=Q(A_n)_\rightarrow +((B_n)_\rightarrow -R(B_n)_\rightarrow )$$
then
$$a=[P]_0-[(B_n)_\rightarrow ]_0\;.$$
By \pr{927} c) and \cor{923} (and \dd{924})
$$0=K_0(\pi^F)a=K_0(\pi ^F)[P]_0-K_0(\pi ^F)[(B_n)_\rightarrow ]_0=[\pi _\rightarrow ^FP]_0-[\pi _\rightarrow ^F(B_n)_\rightarrow ]_0\,,$$
$$[\sigma _\rightarrow ^FP]_0=[\lambda _\rightarrow ^F\pi _\rightarrow ^FP]_0=K_0(\lambda ^F)[\pi _\rightarrow ^FP]_0=K_0(\lambda ^F)[\pi _\rightarrow ^F(B_n)_\rightarrow ]_0=$$
$$=[\lambda _\rightarrow ^F\pi _\rightarrow ^F(B_n)_\rightarrow ]_0=[\sigma _\rightarrow ^F(B_n)_\rightarrow ]_0=[(B_n)_\rightarrow ]_0\,,$$
$$a=[P]_0-[\sigma _\rightarrow ^FP]_0\;.\qedd$$

\begin{p}\label{931}
Let $F$ be an \en \,and $n\in \bn$. 
\begin{enumerate}
\item $C_n+C_n^*\in \unm{E_n}$. 
\item For $X,Y\in F_{n-1}$,
$$(C_n+C_n^*)(A_nX+B_nY)(C_n+C_n^*)=B_nX+A_nY\;.$$
\item If $U,V\in \unn{F_{n-1}}$ then $A_nU+B_nV\in \unn{F_n}$.
\item If $U\in \unn{F_{n-1}}$ then $A_nU+B_n\in \unn{F_n}$ and $A_nU+B_nU^*\in \unm{F_n}$. 
\end{enumerate}
\end{p}

a) From
$$(C_n+C_n^*)(C_n+C_n^*)=B_n+A_n=1_E$$
it follows that $C_n+C_n^*$ is unitary. Being selfadjoint, its spectrum is contained in $\{-1,+1\}$ and so it belongs to \unm{E_n} ([R] Lemma 2.1.3 (ii)). 

b) We have
$$(C_n+C_n^*)(A_nX+B_nY)(C_n+C_n^*)=(C_nX+C_n^*Y)(C_n+C_n^*)=B_nX+A_nY\;.$$

c) We have
$$(A_nU+B_nV)(A_nU+B_nV)^*=A_n+B_n=1_E\,,$$
$$(A_nU+B_nV)^*(A_nU+B_nV)=A_n+B_n=1_E\;.$$

d) By c), $A_nU+B_n\in \unn{F_n}$. By b),
$$(C_n+C_n^*)(A_nU^*+B_n)(C_n+C_n^*)=B_nU^*+A_n\,,$$
so it follows from a), that $A_nU^*+B_n$ is homotopic  to $B_nU^*+A_n$ in $\unn{F_n}$ and so
$$A_nU+B_nU^*=(A_nU+B_n)(A_n+B_nU^*)$$
is homotopic in $\unn{F_n}$ to
$$(A_nU+B_n)(A_nU^*+B_n)=A_n+B_n=1_E\,,$$
i.e. $A_nU+B_nU^*\in \unm{F_n}$.\qed

\begin{p}\label{930}
Let $F$ be a full $E$-C*-algebra, \bbn, $P,Q\in \pp{F_n}$, and $X\in F_n$ with $X^*X=P$, $XX^*=Q$. Then there is a $U\in \unm{F_{n+2}}$ with 
$$U(A_{n+2}A_{n+1}P)U^*=A_{n+2}A_{n+1}Q\,, \qquad \emph{i.e.}\qquad U_\rightarrow P_\rightarrow U^*_\rightarrow =Q_\rightarrow \;.$$
\end{p}

We have $X(1_E-P)=(1_E-Q)X=0$. Put
$$V:=A_{n+1}X+C_{n+1}(1_E-P)+C_{n+1}^*(1_E-Q)+B_{n+1}X^*\qquad (\in F_{n+1})\;.$$
Then
$$V^*=A_{n+1}X^*+C_{n+1}^*(1_E-P)+C_{n+1}(1_E-Q)+B_{n+1}X\,,$$
$$VV^*=A_{n+1}Q+B_{n+1}(1_E-P)+A_{n+1}(1_E-Q)+B_{n+1}P=A_{n+1}+B_{n+1}=1_E\,,$$
$$V^*V=A_{n+1}P+A_{n+1}(1_E-P)+B_{n+1}(1_E-Q)+B_{n+1}Q=A_{n+1}+B_{n+1}=1_E$$
so $V\in \unn{F_{n+1}}$. Moreover
$$VA_{n+1}P=A_{n+1}X\,,\qquad\qquad A_{n+1}XV^*=A_{n+1}Q\;.$$
Put
$$U:=A_{n+2}V+B_{n+2}V^*\;.$$
By \pr{931} d), $U\in \unm{F_{n+2}}$. We have
$$U(A_{n+2}A_{n+1}P)U^*=(A_{n+2}V+B_{n+2}V^*)A_{n+2}A_{n+1}P(A_{n+2}V^*+B_{n+2}V)=$$
$$=A_{n+2}A_{n+1}X(A_{n+2}V^*+B_{n+2}V)=A_{n+2}A_{n+1}Q\;.\qedd$$

\begin{p}\label{932}
Let\, $\oa{F}{\varphi }{G}$\, be a morphism in\, \frm\, and\, $a\in Ker\,K_0(\varphi )$.
\begin{enumerate}
\item There are $\bbn,\,P\in Pr\,\check{F}_{\rightarrow \,n }$, and $U\in \unm{\check{G}_{\rightarrow ,\,n+2} }$ such that
$$a=[P]_0-[\sigma _\rightarrow ^FP]_0\qquad\qquad U(\check{\varphi }_\rightarrow P )U^*=\sigma _\rightarrow ^G\check{\varphi }_\rightarrow P\;. $$
\item If $\varphi $ is surjective then there is a $P\in Pr\,\check{F}_\rightarrow  $ such that
$$a=[P]_0-[\sigma _\rightarrow ^FP]_0,\qquad\qquad \check{\varphi }_\rightarrow P=\sigma _\rightarrow ^G\check{\varphi }_\rightarrow P\;.  $$
\end{enumerate}
\end{p}

a) By \pr{929}, there are $m\in \bn$ and $Q\in Pr\,\check{F}_{\rightarrow ,\,m-1} $ such that
$$a=[Q]_0-[\sigma _\rightarrow ^FQ]_0\;.$$
Since $\check{\varphi }\circ \sigma ^F=\sigma ^G\circ \check{\varphi }  $, by \pr{921} c) and \cor{923},
$$0=K_0(\varphi )a=[\check{\varphi }_\rightarrow Q ]_0-[\check{\varphi }_\rightarrow \sigma _\rightarrow ^FQ ]_0=[\check{\varphi }_\rightarrow Q ]_0-[\sigma _\rightarrow ^G\check{\varphi }_\rightarrow Q ]_0\;.$$
By \cor{918} a$\Rightarrow $c, there is an \bbn, $n>m$, such that
$$\check{\varphi }_\rightarrow Q+(B_n)_\rightarrow \sim _0\sigma _\rightarrow ^G\check{\varphi }_\rightarrow Q+(B_n)_\rightarrow =\sigma _\rightarrow ^G(\check{\varphi }_\rightarrow Q+(B_n)_\rightarrow  )\;.  $$
Put
$$P:=Q+(B_n)_\rightarrow \in Pr\,\check{F}_{\rightarrow \,n}\;. $$
Then
$$[P]_0-[\sigma _\rightarrow ^FP]_0=[Q]_0+[(B_n)_\rightarrow ]_0-[\sigma _\rightarrow ^FQ]_0-[(B_n)_\rightarrow ]_0=a\,,$$
$$[\check{\varphi }_\rightarrow P ]_0-[\sigma _\rightarrow ^G\check{\varphi }_\rightarrow P ]_0=[\check{\varphi }_\rightarrow Q ]_0+[(B_n)_\rightarrow ]_0-[\sigma _\rightarrow ^G\check{\varphi }_\rightarrow Q ]_0-[(B_n)_\rightarrow ]_0=0\;.$$
By \cor{918} a$\Rightarrow $b and \pr{930}, there is a $U\in \unm{\check{G}_{\rightarrow ,\,n+2} }$ with
$$U(\check{\varphi }_\rightarrow P )U^*=\sigma _\rightarrow ^G\check{\varphi }P\;. $$

b) By a), there are \bbn, $n>2$, $Q\in Pr\,\check{F}_{\rightarrow, \,n-2} $, and $U\in \unm{\check{G}_{\rightarrow \,n} }$ such that
$$a=[Q]_0-[\sigma _\rightarrow ^FQ]_0,\qquad\qquad U(\check{\varphi }_\rightarrow Q )U^*=\sigma _\rightarrow ^G\check{\varphi }_\rightarrow Q\;. $$
Since $\mac{\varphi _n }{\check{F}_n }{\check{G}_n }$ is surjective, by [R] Lemma 2.1.7 (i), there is a $V\in \unn{\check{F}_{\rightarrow \,n} }$ with $\check{\varphi }_nV=U $. We put
$$P:=VQV^*\sim _0Q$$
so 
$$a=[P]_0-[\sigma _\rightarrow ^FP]_0$$
and
$$\check{\varphi }_\rightarrow P=(\check{\varphi }_\rightarrow V )(\check{\varphi }_\rightarrow Q )(\check{\varphi }_\rightarrow V^* )=U(\check{\varphi }_\rightarrow Q )U^*=\sigma _\rightarrow ^G\check{\varphi }_\rightarrow Q\,, $$
$$\sigma _\rightarrow ^G\check{\varphi }_\rightarrow P=\sigma _\rightarrow ^G\check{\varphi }_\rightarrow Q=\check{\varphi }_\rightarrow P\;.\qedd    $$

\renewcommand{\labelenumi}{\alph{enumi})}
\renewcommand{\labelenumii}{$\alph{enumi}_{\arabic{enumii}}$)}

\begin{p}\label{934}
Let
$$\oc{F}{\varphi }{G}{\psi }{H}$$
be an exact sequence in \frm.
\begin{enumerate}
\item $\check{\varphi }_\rightarrow  $ is injective.
\item The following are equivalent for all $X\in \check{G}_\rightarrow  $:
\begin{enumerate}
\item $X\in Im\,\check{\varphi }_\rightarrow $.
\item $\check{\psi }_\rightarrow X=\sigma _\rightarrow ^H\check{\psi }_\rightarrow X  $. 
\end{enumerate}
\item $\ob{K_0(F)}{K_0(\varphi )}{K_0(G)}{K_0(\psi )}{K_0(H)}$ is exact.
\end{enumerate}
\end{p}

a) $\check{\varphi } $ is injective (\pr{34} a)) and the assertion follows from \pr{921} b).

$b_1\Rightarrow b_2$ follows from $\psi \circ \varphi =0$.

$b_2\Rightarrow b_1$. Let \bbn\, such that $X\in \check{G}_{\rightarrow \,n} $, which we identify with $\check{G}_n$. Then $X$ has the form
$$X=\si{t\in T_n}((\alpha _t,Y_t)\otimes id_K)V_t^{\check{G} }\,,$$
where $(\alpha _t,Y_t)\in \check{G} $ for every $t\in T_n$, and so by $b_2$),
$$\si{t\in T_n}((\alpha _t,\psi Y_t)\otimes id_K)V_t^{\check{H} }=\check{\psi }_nX=\sigma _n^H\check{\psi }_nX=\si{t\in T_n}((\alpha _t,0)\otimes id_K)V_t^{\check{H} }\;.  $$
It follows $\psi Y_t=0$ for every $t\in T_n$ ([C2] Theorem 2.1.9 a)). Thus for every $t\in T_n$ there is a $Z_t\in F$ with $\varphi Z_t=Y_t$ and we get
$$X=\si{t\in T_n}((\alpha _t,\varphi Z_t)\otimes id_K)V_t^{\check{G} }=$$
$$=\check{\varphi }_n\left(\si{t\in T_n}((\alpha _t,Z_t)\otimes id_K)V_t^{\check{F} }\right)\in Im\,\check{\varphi }_n\subset Im\,\check{\varphi }_\rightarrow \;.   $$ 

c) By \cor{928} a) and \pr{927} e), 
$$K_0(\psi )\circ K_0(\varphi )=K_0(\psi \circ \varphi )=0$$
so $Im\,K_0(\varphi )\subset Ker\,K_0(\psi )$. Let $a\in Ker\,K_0(\psi )$. By \pr{932} b), there is a $P\in Pr\,\check{G}_\rightarrow  $ such that
$$a=[P]_0-[\sigma _\rightarrow ^GP]_0,\qquad\qquad \check{\psi }_\rightarrow P=\sigma _\rightarrow ^H\check{\psi }_\rightarrow P\;.  $$
Then $P$ has the form
$$P=\si{t\in T_n}((\alpha _t,X_t)\otimes id_K)V_t^{\check{G} }$$
for some \bbn \,with $(\alpha _t,X_t)\in E\times G$ for every $t\in T_n$, where we identified $\check{G}_n $ with $\check{G}_{\rightarrow \,n} $. We get
$$\si{t\in T_n}((\alpha _t,\psi X_t)\otimes id_K)V_t^{\check{H} }=\check{\psi }_\rightarrow P=\sigma _\rightarrow ^H\check{\psi }_\rightarrow P=\si{t\in T_n}((\alpha _t,0)\otimes id_K)V_t^{\check{H} }\;.  $$
Thus $\psi X_t=0$ ([C2] Theorem 2.1.9 a)) and there is an $Y_t\in F$ with $\varphi Y_t=X_t$ for every $t\in T_n$. We put
$$Q:=\si{t\in T_n}((\alpha _t,Y_t)\otimes id_K)V_t^{\check{F} }\in Pr\,\check{F}_\rightarrow  $$
with the usual identification ($\check{\varphi }$ is an embedding !). Then
$$\check{\varphi }_\rightarrow Q=\si{t\in T_n}((\alpha _t,\varphi Y_t)\otimes id_K)V_t^{\check{G}}=\si{t\in T_n}((\alpha _t,X_t)\otimes id_K)V_t^{\check{G} }=P $$
and by \pr{927} c) (since $\check{\varphi} \circ \sigma ^F=\sigma ^G\circ \check{\varphi }  $),
$$K_0(\varphi )([Q]_0-[\sigma _\rightarrow^FQ]_0)=[\check{\varphi }_\rightarrow Q ]_0-[\check{\varphi} _\rightarrow \sigma _\rightarrow ^FQ ]_0=$$
$$=[\check{\varphi} _\rightarrow Q ]_0-[\sigma _\rightarrow ^G\check{\varphi }_\rightarrow Q ]_0=[P]_0-[\sigma _\rightarrow ^GP]_0=a\;.$$
Thus $Ker\,K_0(\psi )\subset Im\,K_0(\varphi )$, $Ker\,K_0(\psi )= Im\,K_0(\varphi )$.\qed

\begin{p}[Split Exact Theorem for $K_0$]\label{936}
If
$$\od{F}{\varphi }{G}{\psi }{\lambda }{H}$$
is a split exact sequence in \frm\, then
$$\og{K_0(F)}{K_0(\varphi )}{K_0(G)}{K_0(\psi )}{K_0(\lambda )}{K_0(H)}{15}{15}{15}$$
is also split exact. In particular the map
$$\mad{K_0(F)\times K_0(H)}{K_0(G)}{(a,b)}{K_0(\varphi )a+K_0(\lambda )b}$$
is a group isomorphism and $K_0(\check{F} )\approx K_0(E)\times K_0(F)$ for every \eo $F$.
\end{p}

By \pr{934} c), the second sequence is exact at $K_0(G)$. From
$$K_0(\psi )\circ K_0(\lambda )=K_0(\psi \circ \lambda )=K_0(id_H)=id_{K_0(H)}$$
(\cor{928} a) and \pr{927} d)) it follows that this sequence is (split) exact at $K_0(H)$.

Let $a\in Ker\,K_0(\varphi )$. By \pr{932} a), there are \bbn, $P\in Pr\,\check{F} _{\rightarrow \,n}$, and $U\in \unm{\check{G}_{\rightarrow, \,n+2}  }$ such that
$$a=[P]_0-[\sigma _\rightarrow ^FP]_0,\qquad\qquad U(\check{\varphi }_\rightarrow P )U^*=\sigma _\rightarrow ^G\check{\varphi }_\rightarrow P\;. $$
Put
$$V:=(\check{\lambda} _\rightarrow \check{\psi} _\rightarrow U^*  )U\in \unn{\check{G} _{\rightarrow, \,n+2}}\;.$$
Then
$$\check{\psi }_\rightarrow V=(\check{\psi }_\rightarrow U^* )(\check{\psi }_\rightarrow U )=1_{\rightarrow, \,n+2},\qquad \sigma _\rightarrow ^H\check{\psi }_\rightarrow V=\check{\psi }_\rightarrow V\;.   $$
By \pr{934} $b_2\Rightarrow b_1$, there is a $W\in \unn{\check{F} _{\rightarrow, \,n+2}}$ with $\check{\varphi }_\rightarrow W=V $ ($\check{\varphi } $ is an embedding). We have
$$\check{\varphi }_\rightarrow (WPW^*)=V(\check{\varphi }_\rightarrow P )V^*=(\check{\lambda }_\rightarrow \check{\psi}_\rightarrow  U^*)U(\check{\varphi }_\rightarrow P )U^*(\check{\lambda }_\rightarrow \check{\psi }_\rightarrow U  )= $$
$$=(\check{\lambda }_\rightarrow \check{\psi }_\rightarrow U^*  )(\sigma _\rightarrow ^G\check{\varphi}_\rightarrow  P )(\check{\lambda }_\rightarrow \check{\psi }_\rightarrow U  )=\check{\lambda }_\rightarrow \check{\psi }_\rightarrow (U^*(\sigma _\rightarrow ^G\check{\varphi }_\rightarrow P )U)=  $$
$$=\check{\lambda }_\rightarrow \check{\psi }_\rightarrow \check{\varphi }_\rightarrow P=\sigma _\rightarrow ^G\check{\varphi }_\rightarrow P=\check{\varphi }_\rightarrow \sigma _\rightarrow ^FP\;.     $$
Since $\check{\varphi }_\rightarrow  $ is injective (\pr{934} a)),
$$P\sim _0WPW^*=\sigma _\rightarrow ^FP,\qquad\qquad a=0$$
and $K_0(\varphi )$ is injective.

The last assertion follows since
$$\od{F}{\iota ^F}{\check{F} }{\pi ^F}{\lambda ^F}{E}$$
is a split exact sequence.\qed

\begin{co}\label{938}
Let $F,G$ be \eo $\!$s. 
\begin{enumerate}
\item If we put
$$\mae{\iota _1}{F}{F\times G}{x}{(x,0)},\qquad\mae{\pi _1}{F\times G}{F}{(x,y)}{x}\,,$$
$$\mae{\iota _2}{G}{F\times G}{y}{(0,y)},\qquad\mae{\pi _2}{F\times G}{F}{(x,y)}{y}\,,$$
then the sequences
$$\og{K_0(F)}{K_0(\iota _1)}{K_0(F\times G)}{K_0(\pi _2)}{K_0(\iota _2)}{K_0(G)}{20}{20}{20}\,,$$
$$\og{K_0(G)}{K_0(\iota _2)}{K_0(F\times G)}{K_0(\pi _1)}{K_0(\iota _1)}{K_0(F)}{20}{20}{20}$$
are split exact.
\item The map
$$\mad{K_0(F)\times K_0(G)}{K_0(F\times G)}{(a,b)}{K_0(\iota _1)a+K_0(\iota _2)b}$$
is a group isomorphism {\bf (Product Theorem for $K_0$)}.
\end{enumerate}
\end{co}

a) is easy to see.

b) follows from a) and \pr{936}.\qed

\begin{theo}\label{940}{\bf (Homotopy invariance of $K_0$)}
\begin{enumerate}
\item If $\mac{\varphi ,\psi }{F}{G}$ are homotopic morphisms in \frm, then $K_0(\varphi )=K_0(\psi )$.
\item If $\oa{F}{\varphi }{G},\,\oa{G}{\psi }{F}$ is a homotopy in \frm\, then
$$K_0(\varphi )\circ K_0(\psi )=id_{K_0(G)},\qquad\qquad K_0(\psi )\circ K_0(\varphi )=id_{K_0(F)}\;.$$
\item If $F$ and $G$ are homotopic \eo $\!$s then $K_0(F)$ and $K_0(G)$ are isomorphic.
\item If $F$ is an $E$-C*-algebra such that $id_F$ is homotopic to
$$\mae{0_F}{F}{F}{x}{0}$$
then $F$ is homotopic to $\{0\}$.
\item If the $E$-C*-algebra $F$ is homotopic to $\{0\}$ then $K_0(F)=\{0\}$.
\end{enumerate}
\end{theo}

a) Let
$$\mac{\phi _s}{F}{G},\qquad\qquad s\in [0,1]$$
be a pointwise continuous path of morphisms in \frm\, such that $\phi _0=\varphi ,\,\phi _1=\psi $. Then
$$\mac{\check{\phi _s} }{\check{F} }{\check{G} },\qquad\qquad s\in [0,1]$$
is a pointwise continuous path of morphisms in \frc\, with $\check{\phi} _0=\check{\varphi  },\,\check{\phi} _1=\check{\psi }   $ and for every \bbn,
$$\mac{(\check{\phi _s})_{\rightarrow \,n} }{(\check{F} )_{\rightarrow \,n}}{(\check{G} )_{\rightarrow \,n}},\qquad\qquad s\in [0,1]$$
is a pointwise continuous path in \frc\, with $(\check{\phi _0})_{\rightarrow \,n}=(\check{\varphi } )_{\rightarrow \,n} $ and $(\check{\phi _1})_{\rightarrow \,n}=(\check{\psi  } )_{\rightarrow \,n} $. For every $P\in Pr\,\check{F}_{\rightarrow \,n} $,
$$\mad{[0,1]}{Pr\,(\check{G} )_{\rightarrow \,n}}{s}{(\check{\phi }_s )_{\rightarrow \,n}P}$$
is continuous so (by [R] Proposition 2.2.7)
$$K_0(\varphi )[P]_0=[\varphi _\rightarrow P]_0=[\psi _\rightarrow P]_0=K_0(\psi )[P]_0$$
(\pr{927} c)). By \pr{929}, $K_0(\varphi )=K_0(\psi )$.

b) follows from a), \cor{928} a), and \pr{927} d).

c) follows from b).

d) If we put $\mac{\varphi }{F}{\{0\}}$ and $\mac{\psi }{\{0\}}{F}$ then $\psi \circ \varphi =0_F$ is homotopic to $id_F$ and $\varphi \circ \psi $ is homotopic to $id_{\{0\}}$, so $F$ is homotopic to  $\{0\}$.

e) follows from c).\qed

We show now that $K_0$ is continuous with respect to inductive limits.

\begin{theo}[Continuity of $K_0$]\label{941}
Let $\{(F_i)_{i\in I},\,(\varphi _{ij})_{i,j\in I}\}$ be an inductive system in \frm\, and let $\{F,\,(\varphi _i)_{i\in I}\}$ be its inductive limit in \frm\,. By \emph{\cor{928} a)}, 
$$\{(K_0(F_i))_{i\in I},\,(K_0(\varphi _{ij}))_{i,j\in I}\}$$
 is an inductive system in the category of additive groups. Let $\{\ccc{G},\,(\psi _i)_{i\in I}\}$ be its limit in this category and let $\mac{\psi }{\ccc{G}}{K_0(F)}$ be the group homomorphism such that $\psi \circ \psi _i=K_0(\varphi _i)$ for every $i\in I$. Then $\psi $ is a group isomorphism.
\end{theo}

$\{(\check{F_i})_{i\in I},\,(\check{\varphi _{ij}} )_{i,j\in I} \}$ is an inductive system in \frc\, and by [C2] Proposition 1.2.9 b), $\{\check{F},\,(\check{\varphi _i} )_{i\in I} \}$ may be identified with its inductive limit in \frc. By[C2] Proposition 2.3.5, for every \bbn, $\{((\check{F_i} )_{\rightarrow \,n})_{i\in I},\,((\check{\varphi _{ij}} )_{\rightarrow \,n})_{i,j\in I}\}$ is an inductive system in \frc\, and $\{(\check{F}_{\rightarrow \,n},\,((\check{\varphi _i} )_{\rightarrow \,n} )_{i\in I}\}$ may be identified with its inductive limit in \frc.

\begin{center}
Step 1 $\psi $ is surjective
\end{center}

Let $Q\in Pr\,(\check{F} )_{\rightarrow \,n}$. By [W] L.2.2, there are $i\in I$ and  $P\in Pr\,(\check{F_i} )_{\rightarrow \,n}$ such that $\n{(\check{\varphi }_i )_{\rightarrow \,n}P-Q}<1$, so by [R] Proposition 2.2.4, $(\check{\varphi }_i )_{\rightarrow \,n}P\sim _0Q$. By \pr{927} b),c)
$$\psi \psi _i[P]_0=K_0(\varphi _i)[P]_0=K_0(\check{\varphi }_i )[P]_0=[(\check{\varphi _i})_{\rightarrow \,n}P ]_0=[Q]_0\;.$$
Since
$$Pr\,\check{F}_\rightarrow =\bigcup_{\bbn}Pr\,(\check{F} )_{\rightarrow \,n} \,, $$
$\psi $ is surjective.

\begin{center}
Step 2 $\psi $ is injective
\end{center}

Let $a\in \ccc{G}$ with $\psi a=0$. Since $\ccc{G}=\bigcup_{i\in I}Im\;\psi_i 
 $, there is an $i\in I$ and an $a_i\in K_0(F_i)$ with $a=\psi _i a_i$. There are \bbn\, and $P,Q\in Pr\,(\check{F_i})_{\rightarrow \,n} $ such that
$$a_i=[P]_0-[Q]_0$$
(by \pr{916} c)). By \pr{927} c),
$$0=\psi a=\psi \psi _ia=K_0(\varphi _i)a=K_0(\varphi _i)[P]_0-K_0(\varphi _i)[Q]_0=$$
$$=[(\check{\varphi _i})_{\rightarrow \,n} P]_0-[(\check{\varphi _i})_{\rightarrow \,n}Q ]_0\;.$$
By \cor{918} a$\Rightarrow $b, there is an $R\in Pr\,(\check{F_i} )_{\rightarrow }$ such that
$$PR=QR=0,\qquad\qquad P+R\sim _0Q+R$$
and we get
$$a=[P]_0+[R]_0-[Q]_0-[R]_0=[P+R]_0-[Q+R]_0=0\;.\qedd$$

\begin{center}
\section{Stability of $K_0$}
\end{center}

The stability of $K_0$ holds only under strong supplementary hypotheses. We present below such possible hypotheses, which we fix for this section. We shell give only a sketch of the proof.

Let $S$ be a finite group, $\mac{\chi }{\bzz{2}\times \bzz{2}}{S}$ an injective group homomorphism,
$$a:=\omega (1,0),\qquad b:=\omega (0,1),\qquad c:=\omega (1,1)\,,$$
and $g$ a Schur $E$-function for $S$ such that
$$g(a,b)=g(a,c)=g(b,c)=-g(b,a)=1_E\;.$$
We put for every \bbn,
$$T_n:=S^n=\me{t\in S^{\bn}}{m\in \bn,\,m>n\Rightarrow t_m=1}\,,$$
$$T:=\bigcup_{\bbn}T_n=\me{t\in S^{\bn}}{\{\bbn,\,t_n\not=1\}\,\mbox{is finite}}\,, $$
$$\mae{f}{T\times T}{E}{(s,t)}{\pro{\bbn}g(s_n,t_n)}\,,$$
$$\mae{\stackrel{\underline{n}}{s}}{\bn}{S}{m}{\abb{s}{m=n}{1}{m\not=n}}\,,$$
for every $s\in S$, and
$$C_n:=\frac{1}{2}(V_{\stackrel{\underline{n}}{a}}^f+V_{\stackrel{\underline{n}}{b}}^f),\qquad A_n:=C_n^*C_n,\qquad B_n:=C_nC_n^*\;.$$

Then $f$ is a Schur $E$-function for $T$ and the following hold for all $s,t\in S$ and \bbn:
$$f(\stackrel{\underline{n}}{s},\stackrel{\underline{n}}{t})=g(s,t)\,,$$
$$\stackrel{\underline{n}}{t}=1\;\Longrightarrow\; V_{\stackrel{\underline{n}}{s}}^fV_{\stackrel{\underline{n}}{t}}^f=V_{\stackrel{\underline{n}}{t}}^fV_{\stackrel{\underline{n}}{s}}^f\,,$$
$$s\in T_{n-1}\;\Longrightarrow \;V_s^fV_{\stackrel{\underline{n}}{t}}=V_{\stackrel{\underline{n}}{t}}V_s^f\,,$$
$$A_n=\frac{1}{2}(V_1^f+V_{\stackrel{\underline{n}}{c}}^f)\in Pr\,E_n,\qquad\qquad B_n=\frac{1}{2}(V_1^f-V_{\stackrel{\underline{n}}{c}}^f)\in Pr\,E_n\,,$$
$$A_n+B_n=V_1^f=1_E\,,$$
so the assumptions of \axi{b} are fulfilled.

{\it Remark.} If $\chi $ is bijective and $E=\bc$ then the corresponding projective K-theory coincides with the usual K-theory.

\begin{p}\label{945}
Let $F$ be an \en\, and $m,n\in \bn$. We define
$$\mac{\alpha :=\alpha _{m,n}^F}{(F_m)_n}{F_{m+n}}\,,$$ 
$$\mac{\beta  :=\beta _{m,n}^F}{F_{m+n}}{(F_m)_n}\,,$$  
by
$$(\alpha X)_{(s,t)}:=(X_t)_s,\qquad\qquad ((\beta Y)_t)_s:=Y_{(s,t)}$$
for every $X\in (F_m)_n,\, Y\in F_{m+n}$, and $(s,t)\in S^m\times S^n=S^{m+n}$, where the identification is given by the bijective map
$$\mad{S^m\times S^n}{S^{m+n}}{(s,t)}{(s_1,\cdots,s_m,t_1,\cdots,t_n)}\;.$$
\begin{enumerate}
\item $\alpha $ and $\beta $ are $E$-C*-isomorphisms and $\alpha =\beta ^{-1}$.
\item $\alpha A_n=A_{m+n}$.
\item The diagram
\[ \begin{CD}
(F_m)_{n-1} @>\alpha _{m,n-1}^F>> F_{m+n-1}\\
@V\bar{\rho }_n^{F_m}VV           @VV\bar{\rho }_{m+n}^FV\\
(F_m)_n @>>\alpha _{m,n}^F>       F_{m+n}  
\end{CD} \]
is commutative.
\end{enumerate}
\end{p}

It is obvious that $\alpha $ and $\beta $ are $E$-linear and $\alpha \circ \beta =id_{F_{m+n}}$, $\beta \circ \alpha =id_{(F_m)_n}$. Thus $\alpha $ and $\beta $ are bijective and $\alpha =\beta ^{-1}$.

For $X,Y\in (F_m)_n$ and $(s,t)\in S^m\times S^n$, by [C2] Theorem 2.1.9 c),g),
$$(\alpha X^*)_{(s,t)}=((X^*)_t)_s=(\tilde{f}(t)(X_{t^{-1}})^* )_s=\tilde{f}(s)\tilde{f}(t)((X_{t^{-1}})_{s^{-1}})^*  =$$
$$=\tilde{f}((s,t))(\alpha X_{(s,t)^{-1}})^*=((\alpha X)^*)_{(s,t)}\,, $$
$$((\alpha X)(\alpha Y))_{(s,t)}=$$
$$=\si{(u,v)\in S^m\times S^n}f((u,v),(u^{-1}s,v^{-1}t))(\alpha X)_{(u,v)}(\alpha Y)_{(u^{-1}s,v^{-1}t)}=$$
$$=\si{(u,v)\in S^m\times S^n}f(u,u^{-1}s)f(v,v^{-1}t)(X_v)_u(Y_{v^{-1}t})_{u^{-1}s}=$$
$$=\si{v\in S^n}f(v,v^{-1}t)(X_vY_{v^{-1}t})_s=$$
$$=\left(\si{v\in S^n}f(v,v^{-1}t)X_vY_{v^{-1}t}\right)\!_s=((XY)_t)_s=(\alpha (XY))_{(s,t)}$$
so $\alpha $ is a C*-homomorphism and the assertion follows.

b) follows from the definition of $A_n$ and $A_{m+n}$.

c) follows from b).\qed

\begin{p}\label{948}
Let \oa{F}{\varphi }{G} be a morphism in \frc\, and $m,n\in \bn$. With the notation of \emph{\pr{945}} the diagram
\[ \begin{CD}
(F_m)_n @>\alpha _{m,n}^F>> F_{m+n}\\
@V(\varphi _m)_nVV          @VV\varphi _{m+n}V\\
(G_m)_n @>>\alpha _{m,n}^G>  G_{m+n}
\end{CD} \]
is commutative.
\end{p}

For $X\in (F_m)_n$ and $(s,t)\in S^m\times S^n=S^{m+n}$,
$$(\varphi _{m+n}\alpha _{m,n}^FX)_{(s,t)}=\varphi (\alpha _{m,n}^FX)_{(s,t)}=\varphi (X_t)_s=$$
$$=(\varphi _mX_t)_s=(((\varphi _m)_nX)_t)_s=(\alpha _{m,n}^G(\varphi _m)_nX)_{(s,t)}$$
so
$$\varphi _{m+n}\circ \alpha _{m,n}^F=\alpha _{m,n}^G\circ (\varphi _m)_n\;.\qedd$$

\begin{theo}\label{949}{\bf (Stability for $K_0$)}
If $\oa{F}{\varphi }{G}$ is a morphism in \frm\, and \bbn\, then
$$K_0(F_n)\approx K_0(F),\qquad K_0(G_n)\approx K_0(G),\qquad K_0(\varphi _n)\approx K_0(\varphi )\;.\qedd$$
\end{theo}

{\it Remark.} If $(F_\infty ,(\rho _n^F)_{\bbn})$ and $(G_\infty ,(\rho _n^G)_{\bbn})$ denote the inductive limits in \frm of the corresponding inductive systems $((F_n)_{\bbn},(\rho _{n,m}^F)_{n,m\in \bn})$ and  $((G_n)_{\bbn},(\rho _{n,m}^G)_{n,m\in \bn})$ then, with obvious notation,
$$K_0(F_\infty )\approx K_0(F),\qquad K_0(G_\infty )\approx K_0(G),\qquad K_0(\varphi _\infty )\approx K_0(\varphi )\;.$$

\begin{center}
\chapter{The functor $K_1$}

\section{Definition of $K_1$}

\end{center}

\begin{p}\label{950}
If $F$ is a full $E$-C*-algebra and $n\in \bn$ then
$$\mae{\bar{\tau }_n^F }{\unn{F_{n-1}}}{\unn{F_n}}{U}{A_nU+B_n} $$
is an injective group homomorphism with
$$\bar{\tau }_n^F(\unhh{F_{n-1}}{n-1})\subset \unhh{F_n}{n}\;. $$
For $U,V\in \unn{F_n}$ we put $U\sim _1V$ if $UV^*,U^*V\in \unn{E_n}$. $\sim _1$ is an equivalence relation and $\sim _h$ implies $\sim _1$.
\end{p}

For $U,V\in \unn{F_{n-1}}$,
$$\bar{\tau }_n^FU^*=A_nU^*+B_n=(\bar{\tau }_n^FU )^*\,, $$
$$(\bar{\tau }_n^FU )(\bar{\tau }_n^FV )=(A_nU+B_n)(A_nV+B_n)=A_nUV+B_n=\bar{\tau }_n^F(UV)\,, $$
$$(\bar{\tau }_n^FU )(\bar{\tau }_n^FU )^*=(\bar{\tau }_n^FU )^*(\bar{\tau }_n^FU )=A_n+B_n=1_{F_n}\,,$$
i.e. $\bar{\tau }_n^F $ is well-defined and it is a group homomorphism. If $\bar{\tau }_n^FU=1_{F_n} $ then
$$A_nU+B_n=\bar{\tau }_n^FU=1_{F_n}=1_E=A_n+B_n,\qquad\qquad A_nU=A_n\,, $$
so by \pr{911} c), $U=1_{F_{n-1}}=1_E$ and $\bar{\tau }_n^F $ is injective.

The other assertions are obvious.\qed

\begin{de}\label{951}
Let $F$ be a full $E$-C*-algebra. We put for all $m,n\in \bn$, $m<n$,
$$\mac{\tau _{n,m}^F:=\bar{\tau }_n^F\circ \bar{\tau }_{n-1}^F\circ \cdots\circ \bar{\tau} _{m+1}^F   }{\unn{F_m}}{\unn{F_n}}\;.$$
Then $\{(\unn{F_n})_{n\in \bn},\;(\tau _{n,m})_{m,n\in \bn}\}$ is an inductive system of groups with injective maps. We denote by $\{un\,F,\;(\tau _n^F)_{n\in \bn}\}$ its inductive limit.
$\tau _n^F$ is injective for every $n\in \bn$, so $(\tau _n^F(\unn{F_n}))_{n\in \bn}$ is an increasing sequence of subgroups of $un\,F$, the union of which is $un\,F$. We put for every $n\in \bn$ and $U\in \unn{F_n}$,
$$\unn{F_{\leftarrow \,n}}:=\tau _n^F(\unn{F_n}),\qquad\qquad U_\leftarrow :=U_{\leftarrow \,n}:=U_{\leftarrow \,n}^F:=\tau_n^FU\,,$$
$$1_{\leftarrow \,n}:=1_{\leftarrow \,n}^F:=\tau _n^F1_{F_n}\,(=\tau _n^F1_E)\;.$$
$(\tau _n^F(\unhh{F_n}{n}))_{n\in \bn}$ is an increasing sequence of subgroups of $un\,F$; we denote by $\zh{F}$ their union.  
\end{de} 

We often identify $\unn{F_n}$ with $\unn{F_{\leftarrow \,n}}$.

\begin{p}\label{951a}
For $m,n\in \bn$, $m<n$, and $U\in \unn{F_m}$,
$$\tau _{n,m}^FU=\left(\proo{i=m+1}{n}A_i\right)U+\left(1_E-\proo{i=m+1}{n}A_i\right)\;.$$
\end{p}

We prove this identity by induction with respect to $n$. The identity holds for $n:=m+1$. Assume it holds for $n-1\geq m$. Then
$$\tau _{n,m}^FU=\bar{\tau }_n^F\tau _{n-1,m}^FU=A_n\tau _{n-1,m}^FU+B_n= $$
$$=A_n\left(\left(\proo{i=m+1}{n-1}A_i\right)U+\left(1_E-\proo{i=m+1}{n-1}A_i \right)\right)+B_n=$$
$$=\left(\proo{i=m+1}{n}A_i\right)U+\left(1_E-\proo{i=m+1}{n}A_i\right)\;.\qedd$$

\begin{p}\label{952}
Let $F$ be a full $E$-C*algebra.
\begin{enumerate}
\item If $U,V\in \unn{F_{n-1}}$ for some $n\in \bn$ then 
$$\bar{\tau }_n^F(UV)\sim _h\bar{\tau }_n^F(VU)\,,\qquad\qquad \bar{\tau }_n^F(UVU^*)\sim _h\bar{\tau }_n^F(V)\;.$$
\item $\zh{F}$ is a normal subgroup of $\zn{F}$ and $\zn{F}/\zh{F}$ is commutative.
\item For all $U,V\in \zn{F}$,
$$UV^*\in \zh{F}\Longleftrightarrow U^*V\in \zh{F}\;.$$
We put $U\sim _1V$ if $UV^*\in \zh{F}$. $\sim _1$ is an equivalence relation.
\end{enumerate}
\end{p}

a) By \pr{931} a),b),
$$\bar{\tau }_n^F(UV)=A_nUV+B_n=(A_nU+B_n)(A_nV+B_n)\sim _h$$
$$\sim _h(A_nU+B_n)(A_n+B_nV)=A_nU+B_nV\sim _hA_nV+B_nU\sim _h\bar{\tau }_n^F(VU)\;.  $$
It follows
$$\bar{\tau }_n^F(UVU^*)\sim _h\bar{\tau }_n^F(U^*UV)=\bar{\tau }_n^F(V)\;.$$

b) $\zh{F}$ is obviously a subgroup of $\zn{F}$. The other assertions follow from a).

c) Let $q:\zn{F}\rightarrow \zn{F}/\zh{F}$ be the quotient map. If $UV^*\in \zh{F}$ then by b),
$$q(UV^*)=q(U)q(V^*)=q(V^*)q(U)=q(V^*U)\,,$$
$$V^*U\in \zh{F}\,,\qquad\qquad U^*V=(V^*U)^*\in \zh{F}\;.\qedd$$

\begin{de}\label{953}
We denote for every \eo $F$ by $K_1(F)$ the additive group obtained from the commutative group $un\,\check{F}/\zh{\check{F}}  $ \emph{(\pr{952} b))} by replacing the multiplication with the addition $\oplus $; by this the neutral element (which corresponds to $1_E$) is denoted by $0$. For every $U\in un\,\check{F} $ we denote by $[U]_1$ its equivalence class in $K_1(F)$.
\end{de}

{\it Remark.} Let $F$ be a \en. By \pr{10.3'a} d), $\check{F} $ is isomorphic to $E\times F$, so in this case we may define $K_1$ using $F$ instead of $\check{F} $ (as we did for $K_0$).

\begin{p}\label{954}
Let $\oa{F}{\varphi }{G}$ be a morphism in \frm.
\begin{enumerate}
\item For $m,n\in \bn$, $m<n$, the diagram
$$\begin{CD}
\unn{\check{F}_m }@>\tau _{n,m}^{\check{F} }>>\unn{\check{F}_n }\\
@V\check{\varphi }_m VV                     @VV\check{\varphi }_nV\\
\unn{\check{G}_m }@>>\tau _{n,m}^{\check{G} }>\unn{\check{G}_n }  
\end{CD}$$
is commutative. Thus there is a unique group homomorphism
$$\mac{\check{\varphi }_\leftarrow  }{un\,\check{F} }{un\,\check{G} }$$
such that
$$\check{\varphi }_\leftarrow \circ \tau _n^{\check{F} }=\tau _n^{\check{G} }\circ \check{\varphi }_n  $$
for every $n\in \bn$.
\item $\varphi _\leftarrow (\zh{\check{F}} )\subset \zh{\check{G}} $; if $\varphi $ is surjective then $\varphi _\leftarrow (\zh{\check{F}} )=\zh{\check{G}} $.
\item There is a unique group homomorphism
$$\mac{K_1(\varphi )}{K_1(F)}{K_1(G)}$$
such that
$$K_1(\varphi )[U]_1=[\check{\varphi }_\leftarrow U ]_1$$
for every $U\in un\,\check{F} $.
\item $K_1(id_F)=id_{K_1(F)}$.
\item $K_1(\{0\})=\{0\}$.
\end{enumerate}
\end{p}

a) It is sufficient to prove the assertion for $n=m+1$. For $U\in \unn{\check{F}_m }$,
$$\tau _{n,m}^{\check{G} }\check{\varphi }_mU=A_n(\check{\varphi }_mU)+B_n=\check{\varphi }_n(A_nU+B_n) =\check{\varphi }_n\tau _{n,m}^{\check{F} } U\;.$$

b) Since $\check{\varphi }_n(\unhh{\check{F}_n}{n})\subset \unhh{\check{G}_n}{n}$ for every $n\in \bn$, it follows $\varphi _\leftarrow (\zh{\check{F}} )\subset \zh{\check{G}} $. If $\varphi $ is surjective then by [R] Lemma 2.1.7 (iii), we may replace the above inclusion relation by $=$.

c) follows from a) and b).

d) is obvious.

e) follows from $\zn{E}=\zh{E}$.\qed

\begin{de}\label{9.7}
An \eo $F$ is called {\bf{K-null}} if 
$$K_0(F)=K_1(F)=0\;.$$
 Let $\oaa{F}{\varphi }{G}$ be a morphism in \frm. We say that $\varphi$ is {\bf{ K-null}} if 
$$K_0(\varphi )=K_1(\varphi )=0\;.$$ 
We say that {\bf{$\varphi $ factorizes through null}} if there  are morphisms \obb{F}{\varphi '}{H}{\varphi ''}{G} in \frm such that $\varphi =\varphi ''\circ \varphi '$ and $H$is K-null. 
\end{de}

\begin{p}\label{956}
\rule{0mm}{0mm}
\begin{enumerate}
\item If $\ob{F}{\varphi }{G}{\psi }{H}$ are morphisms in \frm \,then
$$\check{\psi} _\leftarrow \circ \check{ \varphi} _\leftarrow =(\check{ \psi} \circ \check{ \varphi} )_\leftarrow=\left(\check{ \overbrace{\psi \circ \varphi }}\right)\!_\leftarrow \;\; ,\qquad K_1(\psi )\circ K_1(\varphi )=K_1(\psi \circ \varphi )\;.$$
\item If $\varphi =0$ then $K_1(\varphi )=0$.
\item {\bf{(Homotopy invariance of $K_1$)}} If $\mac{\varphi ,\psi }{F}{G}$ are homotopic morphisms in \frm\, then 
$$K_1(\varphi )=K_1(\psi )\;.$$ 
\item {\bf{(Homotopy invariance of $K_1$)}} If $\ob{F}{\varphi }{G}{\psi }{F}$ is a homotopy in \frm \,then
$$\mac{K_1(\varphi )}{K_1(F)}{K_1(G)},\qquad\qquad \mac{K_1(\psi  )}{K_1(G)}{K_1(F)}$$
are isomorphisms and $K_1(\psi )=K_1(\varphi )^{-1}$.
\item If the $E$-C*-algebra $F$ is homotopic to $\{0\}$ then $F$ is K-null.
\item If a morphism in \frm factorizes through null then it is K-null.
\end{enumerate}
\end{p}

a) Since 
$$\check{\psi }_n\circ \check{\varphi }_n=(\check{\psi }\circ \check{\varphi }  )_n=\left(\check{\overbrace{\psi \circ \varphi }} \right)\!_n  $$
 for every $n\in \bn$ we get 
$$\check{\psi }_\leftarrow \circ \check{\varphi }_\leftarrow =(\check{\psi }\circ \check{\varphi }  )_\leftarrow =\left(\check{\overbrace{\psi \circ \varphi }} \right)\!_\leftarrow   \;.$$ 
For $U\in un\,\check{F} $, by \pr{954} c),
$$K_1(\psi )K_1(\varphi )[U]_1=K_1(\psi )[\check{\varphi }_\leftarrow U ]_1=[\check{\psi }_\leftarrow \check{\varphi }_\leftarrow U  ]_1=$$
$$=[(\check{\psi }\circ \check{\varphi }  )_\leftarrow U]_1=\left[\left(\overbrace{\psi \circ \varphi }^{\check{\rule{0mm}{0mm}} }\right)\!_\leftarrow  U \right]_1=K_1(\psi \circ \varphi )[U]_1\,,$$
so $K_1(\psi )\circ K_1(\varphi )=K_1(\psi \circ \varphi )$.

b) If we put $\mac{\vartheta }{F}{\{0\}}$, $\mac{\iota }{\{0\}}{G}$ then $\varphi =\iota \circ \vartheta $ and by a) and \pr{954} e), $K_1(\varphi )=0$.

c) Let
$$\mac{\phi _s}{F}{G},\qquad\qquad s\in [0,1]$$
be a pointwise continuous path of morphisms in \frm with $\phi _0=\varphi $ and $\phi _1=\psi $. Let \bbn .\,Then
$$\mac{(\check{\phi }_s)_n }{\check{F}_n }{\check{G}_n },\qquad\qquad s\in [0,1]$$
is a pointwise continuous path of $E$-C*-homomorphisms with $(\check{\phi }_0)_n=\check{\varphi }_n  $ and $(\check{\phi }_1 )_n=\check{\psi }_n $. For every $U\in \unn{\check{F}_n }$, the map
$$\mae{\vartheta }{[0,1]}{\unn{\check{G}_n }}{s}{(\check{\phi }_s )_nU}$$
is continuous and $\vartheta (0)=\check{\varphi }_nU $, $\vartheta (1)=\check{\psi }_nU $, i.e. $\check{\varphi }_nU$ and $\check{\psi }_nU  $ are homotopic in $\unn{\check{G}_n }$. It follows
$$K_1(\varphi )[\tau _n^{\check{F} }U]_1=K_1(\psi )[\tau _n^{\check{F} }U]_1\,,$$
which implies $K_1(\varphi )=K_1(\psi )$.

d) follows from c) and \pr{954} d). 

e) By d) and \pr{954} e), $K_1(F)=\z{0}$. By the Homotopy invariance of $K_0$ (\h{940} e)), $F$ is K-null.

f) follows immediately from a), e), and \cor{928} a).\qed

\begin{p}\label{958}
If 
$$\oc{F}{\varphi }{G}{\psi }{H}$$
is an exact sequence in \frm\, then
$$\ob{K_1(F)}{K_1(\varphi )}{K_1(G)}{K_1(\psi )}{K_1(H)}$$
is also exact.
\end{p}

Let $a\in Ker\,K_1(\psi )$ and let $U\in un\,\check{G} $ with $a=[U]_1$. By \pr{954} c),
$$0=K_1(\psi )a=[\check{\psi }_\leftarrow U ]_1\,,\qquad\qquad \check{\psi }_\leftarrow U\in \zh{\check{H}} \;. $$
By \pr{954} b), there is a $V\in \zh{\check{G}}$ with $\check{\psi }_\leftarrow V=\check{\psi }_\leftarrow U  $. We put $W:=UV^*$. By \pr{952} c), $[W]_1=a$ and so
$$\check{\psi }_\leftarrow W=(\check{\psi }_\leftarrow U )(\check{\psi }_\leftarrow V )^*=1_E\;. $$
$W$ has the form
$$W=\si{t\in T_n}((\alpha _t,X_t)\otimes id_K)V_t^{\check{G} }$$
for some $\bbn$, where $(\alpha _t,X_t)\in E\times G$ for every $t\in T_n$. We get
$$1_E=\check{\psi }_nW=\si{t\in T_n}((\alpha _t,\psi X_t)\otimes id_K)V_t^{\check{H} } $$
and so by [C2] Theorem 2.1.9 a), $\psi X_t=0$ for every $t\in T_n$. For every $t\in T_n$, let $Y_t\in F$ with $\varphi Y_t=X_t$ and put
$$W':=\si{t\in T_n}((\alpha _t,Y_t)\otimes id_K)V_t^{\check{F} }\;.$$
Since $\mac{\check{\varphi } }{\check{F} }{\check{G} }$ is an embedding, $W'\in \unn{\check{F}_{\leftarrow \,n} }$ and by \pr{954} c),
$$K_1(\varphi )[W']_1=[\check{\varphi }_nW' ]_1=[W]_1=a\;.$$
Thus $Ker\,K_1(\psi )\subset Im\,K_1(\varphi )$.

Let now $U\in un\,\check{F}_\leftarrow  $. By \pr{956} a),b),
$$K_1(\psi )K_1(\varphi )[U]_1=K_1(\psi \circ \varphi )[U]_1=K_1(0)[U]_1=0$$
so $Im\,K_1(\varphi )\subset Ker\,K_1(\psi )$.\qed

\begin{p}\label{9.6}
The following are equivalent for every \en $F$.
\begin{enumerate}
\item $K_1(F)=\{0\}$.
\item For every \bbn\, and $U\in \unn{F_n}$ there is an $m\in \bn$, $m>n$, with $\tau _{m,n}^FU\sim _h1_E$ in $\unn{F_m}$.
\end{enumerate}
\end{p} 

$a\Rightarrow b$ Since
$$(1_E,U)\in \unn{E_n}\times \unn{F_n}=\unn{(E_n\times F_n)}=\unn{(E\times F)_n}\,,$$
it follows from \pr{10.3'a} d), $(1_E,U-1_E)\in \unn{\check{F_n} }$. By a), there is an $m\in \bn$, $m>n$, with
$$U_0:=(1_E,\tau _{m,n}^FU-1_E)=\tau _{m,n}^{\check{F} }(1_E,U-1_E)\in Un_{E_m}\,\check{F}_m\;. $$
Thus there is a continuous map
$$\mad{[0,1]}{\unn{\check{F}_m }}{s}{U_s}$$
with $U_1\in \unn{E_m}\,(\subset \unn{\check{F}_m })$. We put
$$U'_s:=U_s(\sigma _m^FU_s)^*\,(\in \unn{\check{F}_m })$$
for every $s\in [0,1]$. Then the map
$$\mad{[0,1]}{\unn{\check{F}_m }}{s}{U'_s}$$
is continuous and $U'_0=U_0$, $U'_1=1_E$. Let 
$$\mae{\varphi }{\check{F} }{E\times F}{(\alpha ,x)}{(\alpha ,x+\alpha )}$$
be the $E$-C*-isomorphism of \pr{10.3'a} d). Then
$$\mae{U''}{[0,1]}{\unn{E_n}\times \unn{F_n}}{s}{\varphi _mU'_s}$$
is continuous and
$$U''_0=\varphi _mU'_0=(1_E,\tau _{m,n}^FU)\,,\qquad U''_1=\varphi _mU'_1=(1_E,1_E)\;.$$
Thus $\tau _{m,n}^FU\sim _h1_E$ in $\unn{F_m}$.

$b\Rightarrow a$ Let $a\in K_1(F)$. There are $n\in \bn$ and $U\in \unn{\check{F}_n }$ with $a=[U]_1$. Since $U(\sigma _n^FU)^*\sim _1U$, we may assume $U=U(\sigma _n^FU)^*$, i.e. $\sigma _n^FU=1_E$. Thus there is a unique $X\in F_n$ with $\iota _n^FX=U-1_E$. Then
$$U':=X+1_E\in \unn{F_n}\;.$$
By b), there is an $m\in \bn$, $m>n$, with $\tau _{m,n}^FU'\sim _h1_E$. By \pr{10.3'a} d),
$$U=(1_E,X)=(1_E,U'-1_E)\,,\qquad \tau _{m,n}^{\check{F} }U=(1_E,\tau _{m,n}^FU'-1_E)\sim _h(1_E,0)\,,$$
i.e. $a=[U]_1=0$.\qed

\begin{co}\label{10.6}
If $F$ is a finite-dimensional \en then $K_1(F)=\{0\}$.
\end{co}

For every \bbn, $F_n$ is finite-dimensional and so there is a finite family $(k_i)_{i\in I}$ in \bn\, such that $F_n\approx \pro{i\in I}\bc_{k_i,k_i}$. Thus every $U\in\unn{F_n} $ is homotopic to $1_E$ in $\unn{F_n}$. By \pr{9.6} $b\Rightarrow a$, $K_1(F)=\{0\}$.\qed

\begin{co}\label{o}
If the spectrum of $E$ is totally disconnected (this happens e.g. if $E$ is a W*-algebra \emph{([C1] Corollary 4.4.1.10})) then $\unn{E_n}=\unm{E_n}$ for every $\bbn$ and so $K_1(E)=\z{0}$.
\end{co}

Let $\Omega $ be the spectrum of $E$ and let $U\in \unn{E_n}$. $U$ has the form
$$U=\si{t\in T_n}(U_t\otimes id_K)V_t\,,$$
with $U_t\in E$ for every $t\in T_n$. We put
$$U(\omega):=\si{t\in T_n}(U_t(\omega)\otimes id_K)V_t$$
for every $\omega \in \Omega $ and denote by $\sigma (U(\omega ))$ its spectrum, which is finite. Let $\omega _0\in \Omega $ and let $\theta _0\in [0,2\pi [$ such that 
 $e^{i\theta _0}\not\in \sigma (U(\omega _0))$. By [C1] Corollary 2.2.5.2, there is o clopen neighborhood $\Omega _0$ of $\omega _0$ such that $e^{i\theta _0}$ does not belong to the spectrum of $U(\omega )$ for all $\omega \in \Omega _0$. Assume for a moment $\Omega_0 =\Omega $ and put for every $s\in [0,1]$,
$$\mae{h_s}{\bt\setminus \{\alpha \}}{\bt}{e^{i\vartheta }}{e^{i\vartheta s}},\qquad\qquad W_s:=h_s(U) \,,$$
where $\vartheta \in ]\vartheta _0-2\pi ,\vartheta _0[$. Then
$$\mad{[0,1]}{\unn{E_n}}{s}{W_s}$$
is a continuous path in $\unn{E_n}$ ([C1] Corollaries 4.1.2.13 and 4.1.3.5) with $W_1=U$ and $W_0=1_E$. Thus $U\in \unm{E_n}$.

Since $\Omega $ is the union of a finite family of pairwise disjoint clopen sets of the above form $\Omega _0$, $U\in \unm{E_n}$.

By \pr{9.6} $b\Rightarrow a$, $K_1(E)=\z{0}$.\qed

\begin{center}
\section{The index map}
\end{center}

\begin{center}
\fbox{\parbox{8cm}{Throughout this section 
$$\oc{F}{\varphi }{G}{\psi }{H}$$
 denotes an exact sequence in \frm\, and $n\in \bn$.}}
 \end{center}
 
\begin{p}\label{960}
Let $U\in \unn{\check{H}_{n-1} }$.
\begin{enumerate}
\item There are $V\in \unn{\check{G}_n }$ and $P\in Pr\,\check{F}_n $ such that
$$\check{\psi }_nV=A_nU+B_nU^*,\qquad\qquad \check{\varphi }_nP=VA_nV^*\;.  $$
\item If $W\in \unn{\check{G}_n }$ and $Q\in Pr\,\check{F}_n $ such that
$$\check{\psi }_nW=A_nU+B_nU^*,\qquad\qquad \check{\varphi }_nQ=WA_nW^*  $$
then $\sigma _n^FQ=A_n$ and $P\sim _0Q$.
\item Let $U_0\in \unn{\check{H}_{n-1} }$, $V_0\in \unn{\check{G}_n }$, and $P_0\in Pr\,\check{F}_n $ with
$$U_0\sim _1U,\qquad\qquad\check{\psi }_nV_0=A_nU_0+B_nU_0^*,\qquad\qquad \check{\varphi}_n P_0=V_0A_nV_0^*\;.  $$
Then $P_0\sim _0P$.
\item If $U\in \unhh{\check{H}_{n-1} }{n-1}$ then $P\sim _0A_n$.
\end{enumerate}
\end{p}

a) By \pr{931} d), $A_nU+B_nU^*\in \unm{\check{H}_n }$ so by [R] Lemma 2.1.7 (i) (and [C2] Theorem 2.1.9 a)), there is a $V\in \unm{\check{G}_n }$ with $\check{\psi }_nV=A_nU+B_nU^* $. We have
$$\check{\psi }_n(VA_nV^*)=(A_nU+B_nU^*)A_n(A_nU^*+B_nU)=A_n\,, $$
$$\sigma_n^H\check{\psi }_n(VA_nV^*)=\sigma_n^HA_n=A_n=\check{\psi }_n(VA_nV^*)\,,  $$
so by \pr{934} $b_2\Rightarrow b_1$, there is a $P\in \pp{\check{F}_n }$ with $\check{\varphi }_nP=VA_nV^* $. 

b) Since $\pi ^F=\pi ^H\circ \check{\psi }\circ \check{\varphi }  $, we have
$$\pi _n^FQ=\pi _n^H\check{\psi }_n\check{\varphi }_nQ=\pi _n^H\check{\psi }_n(WA_nW^*)=$$
$$=\pi_n^H((A_nU+B_nU^*)A_n(A_nU^*+B_nU))=\pi _n^HA_n=A_n\,,$$
$\sigma _n^FQ=A_n$. Since
$$\check{\psi }_n(WV^*)=(A_nU+B_nU^*)(A_nU^*+B_nU)=A_n+B_n=1_E=\sigma _n^H\check{\psi }_n(WV^*)\,,  $$
by \pr{934} $b_2\Rightarrow b_1$, there is a $Z\in \unn{\check{F}_n }$ with $\check{\varphi }_nZ=WV^* $. Then
$$\check{\varphi }_n(ZPZ^*)=(WV^*)(VA_nV^*)(VW^*)=WA_nW^*=\check{\varphi }_nQ\,,$$
$$ ZPZ^*=Q,\qquad\qquad P\sim _0Q\;.  $$

c) By \pr{952} c), $U^*U_0,UU_0^*\in \unhh{\check{H}_{n-1} }{n-1}$ so by [R] Lemma 2.1.7 (iii), there are $X,Y\in \unn{\check{G}_{n-1} }$ such that 
$$\check{\psi }_{n-1}X=U^*U_0,\qquad\qquad \check{\psi }_{n-1}Y=UU_0^*\;.  $$
We put
$$Z:=V(A_nX+B_nY)\;.$$
By \pr{931} c), $Z\in \unn{\check{G}_n }$. We have
$$\check{\psi }_nZ=(A_nU+B_nU^*)(A_nU^*U_0+B_nUU_0^*)=A_nU_0+B_nU_0^*\,, $$
$$\check{\psi }_n(ZA_nZ^*)=(A_nU_0+B_nU_0^*)A_n(A_nU_0^*+B_nU_0)=A_n=\sigma _n^H\check{\psi }_n(ZA_nZ^*)\;. $$
By \pr{934} $b_2\Rightarrow b_1$, there is a $Q\in Pr\,\check{F}_n $ with $\check{\varphi }_nQ=ZA_nZ^* $. By b), $Q\sim _0P_0$. From
$$\check{\varphi }_nQ=ZA_nZ^*=V(A_nX+B_nY )A_n(A_nX^*+B_nY^*)V^*=VA_nV^*=\check{\varphi }_nP $$
it follows $P_0\sim _0Q=P$ (by [C2] Theorem 2.1.9 a)).

d) By c), we may take $U=1_E$. Further we may take $W=1_E$ and $Q=A_n$ in b), so $P\sim A_n$.\qed

\begin{p}\label{962}
For every $i\in \{1,2\}$ let $U_i\in \unn{\check{H}_{n-1} }$, $V_i\in \unn{\check{G}_n} $, and $P_i\in Pr\,\check{F}_n $ such that
$$\check{\psi }_nV_i=A_nU_i+B_nU_i^*,\qquad\qquad \check{\varphi }_nP_i=V_iA_nV_i^*\;.  $$
Put
$$X:=A_{n+1}A_n+C_{n+1}^*C_n+C_{n+1}C_n^*+B_{n+1}B_n,\qquad U:=A_nU_1+B_nU_2\,,$$
$$V:=X(A_{n+1}V_1+B_{n+1}V_2)X,\qquad\qquad P:=X(A_{n+1}P_1+B_{n+1}P_2)X\;,$$
\begin{enumerate}
\item $X\in \unm{E_{n+1}},\;U\in \unn{}\check{H}_n,\;V\in \unn{\check{G}_{n+1} },\;P\in Pr\,\check{F}_{n+1}  $.
\item $\check{\psi }_{n+1}V=A_{n+1}U+B_{n+1}U^*,\;\check{\varphi }_{n+1}P=VA_{n+1}V^*  $.
\end{enumerate}
\end{p}

a) We have
$$X^2=A_{n+1}A_n+A_{n+1}B_n+B_{n+1}A_n+B_{n+1}B_n=1_E\;.$$
Since $X$ is selfadjoint it follows $X\in \unm{E_{n+1}}$ ([R] Lemma 2.1.3 (ii)) and so $P\in Pr\,\check{F}_{n+1} $. By \pr{931} c), $U\in \unn{\check{H}_n }$ and $V\in \unn{\check{G}_{n+1} }$.

b) We have
$$XA_{n+1}X=(A_{n+1}A_n+C_{n+1}C_n^*)X=A_{n+1}A_n+B_{n+1}A_n=A_n\,,$$
$$XB_{n+1}X=(C_{n+1}^*C_n+B_{n+1}B_n)X=A_{n+1}B_n+B_{n+1}B_n=B_n\,,$$
$$XA_nX=A_{n+1},\qquad\qquad XB_nX=B_{n+1}\,, $$
$$XA_{n+1}A_nX=A_{n+1}A_n,\qquad XA_{n+1}B_nX=B_{n+1}A_n\,,$$
$$XB_{n+1}A_nX=A_{n+1}B_n,\qquad XB_{n+1}B_nX=B_{n+1}B_n\,,$$
$$\check{\psi }_{n+1}V=X(A_{n+1}(A_nU_1+B_nU_1^*)+B_{n+1}(A_nU_2+B_nU_2^*))X=$$
$$=A_{n+1}A_nU_1+B_{n+1}A_nU_1^*+A_{n+1}B_nU_2+B_{n+1}B_nU_2^*=A_{n+1}U+B_{n+1}U^*\,, $$
$$VA_{n+1}V^*=X(A_{n+1}V_1+B_{n+1}V_2)XA_{n+1}X((A_{n+1}V_1^*+B_{n+1}V_2^*)X=$$
$$=X(A_{n+1}V_1+B_{n+1}V_2)A_n(A_{n+1}V_1^*+B_{n+1}V_2^*)X=$$
$$=X(A_{n+1}V_1A_nA_{n+1}V_1^*+B_{n+1}V_2A_nB_{n+1}V_2^*)X=$$
$$=X(A_{n+1}V_1A_nV_1^*+B_{n+1}V_2A_nV_2^*)X=$$
$$=X(A_{n+1}\check{\varphi }_nP_1+B_{n+1}\check{\varphi }_nP_2  )X=$$
$$=\check{\varphi }_{n+1}(X(A_{n+1}P_1+B_{n+1}P_2)X)= \check{\varphi }_{n+1}P\;.\qedd $$

\begin{co}\label{964}
There is a unique group homomorphism, called {\bf{the index map}},
$$\mac{\delta _1}{K_1(H)}{K_0(F)}$$
such that
$$\delta _1[U]_1=[P]_0-[\sigma _\rightarrow ^FP]_0$$
for every $U\in un\,{\check{H} }$, where $P$ satisfies the conditions of \emph{\pr{960} a)}. 
\end{co}

By \pr{960} a),b), the map
$$\mae{\nu _n}{\unn{\check{H}_{n-1} }}{K_0(F)}{U}{[P]_0-[\sigma_n^FP]_0}$$
is well-defined for every \bbn, where $P$ is associated to $U$ as in \pr{960} a). By \pr{960} c), $\nu _nU=\nu _nU_0$ for all $U,U_0\in \unn{\check{H}_{n-1} }$ with $U\sim _1U_0$. With the notation of \pr{962},
$$\nu _{n+1}(A_nU_1+B_nU_2)=\nu _{n+1}U=[P]_0-[\sigma_{n+1}^FP]_0=$$
$$=[A_{n+1}P_1+B_{n+1}P_2]_0-[\sigma_{n+1}^F(A_{n+1}P_1+B_{n+1}P_2)]_0=$$
$$=[P_1]_0+[P_2]_0-[\sigma_n^FP_1]_0-[\sigma_n^FP_2]_0=\nu _nU_1+\nu _nU_2\;.$$
Thus by \pr{960} d) (and \pr{962}), for $U\in \unn{\check{H}_{n-1} }$,
$$\nu _{n+1}(\bar{\tau }_n^{\check{H} }U )=\nu _{n+1}(A_nU+B_n)=\nu _nU+\nu _n1_E=\nu _nU\;.$$
Hence the map
$$\mae{\nu }{un\,\check{H} }{K_0(F)}{U}{\nu _nU}$$
is well-defined, where $U\in \unn{\check{H}_{n-1} }$ for some \bbn. By \pr{960} d), again, $\nu $ induces a map $\mac{\delta _1}{K_1(H)}{K_0(F)}$, which is additive by the above considerations. The uniqueness follows from the fact that the map $\mac{[\cdot ]_1}{un\,\check{H} }{K_1(H)}$ is surjective.\qed

\begin{p}\label{966}
Let
$$\oc{F'}{\varphi '}{G'}{\psi '}{H'}$$
be an exact sequence in \frm \, and $\delta '_1$ its associated index map. If the diagram in \frm
$$\begin{CD}
0@>>>F@>\varphi >>G@>\psi >>H@>>>0\\
@.      @V\gamma VV  @V\alpha VV @VV\beta V@.\\
0@>>>F'@>>\varphi' >G'@>>\psi' >H'@>>>0
\end{CD}$$
is commutative then the diagram
$$\begin{CD}
K_1(H)@>\delta _1>>K_0(F)\\
@VK_1(\beta )VV    @VVK_0(\gamma )V\\
K_1(H')@>>\delta '_1>K_0(F')
\end{CD}$$
is also commutative.
\end{p}

Let $U\in \unn{\check{H}_{n-1} }$, $V\in \unn{\check{G}_n }$, and $P\in Pr\,\check{F}_n $ with
$$\check{\psi }_nV=A_nU+B_nU^*,\qquad\qquad \check{\varphi }_nP=VA_nV^*\;.  $$
Put
$$V':=\check{\alpha }_nV\in \unn{\check{G'}_n },\qquad\qquad P':=\check{\gamma }_nP\in Pr\,\check{F'}_n\;.   $$
Then
$$\check{\psi '}_nV'=\check{\psi '}_n\check{\alpha }_nV=\check{\beta }_n\check{\psi }_nV=A_n\check{\beta }_{n-1}U+B_n\check{\beta }_{n-1}U^*\,,     $$
$$\check{\varphi '}_nP'=\check{\varphi '}_n\check{\gamma }_nP=\check{\alpha }_n\check{\varphi }_nP=\check{\alpha }_n(VA_nV^*)=V'A_nV'^*\;.$$
By \cor{964} for $\delta '_1$, \pr{954} c), and \pr{927} c),
$$\delta '_1K_1(\beta )[U]_1=\delta '_1[\check{\beta }_{n-1}U ]_1=[P']_0-[\sigma _n^{F'}P']_0=[\check{\gamma }_nP ]_0-[\sigma _n^{F'}\check{\gamma }_nP ]_0=$$
$$=[\check{\gamma }_nP]_0-[\check{\gamma }_n\sigma _n^FP ]_0=K_0(\gamma )([P]_0-[\sigma _n^FP]_0)=K_0(\gamma )\delta _1[U]_1\;.\qedd$$

\begin{p}\label{965}
\rule{0mm}{0mm}
\begin{enumerate}
\item $\delta _1\circ K_1(\psi )=0$.
\item $K_0(\varphi )\circ \delta _1=0$.
\end{enumerate}
\end{p}

a) Let $U\in \unn{\check{G}_{n-1} }$ and put
$$V:=\bar{\tau }_n^{\check{G} }U =A_nU+B_n\in \unn{\check{G}_n }\;.$$
Then
$$\check{\psi }_nV=A_n(\check{\psi }_{n-1}U )+B_n\,, $$
$$(\check{\psi }_nV )A_n(\check{\psi }_nV )^*
=(A_n(\check{\psi }_{n-1}U )+B_n)A_n(A_n(\check{\psi }_{n-1}U )^*+B_n)=A_n\,,$$
so (by \pr{954} c))
$$\delta _1K_1(\psi )[U]_1=\delta _1[\check{\psi }_{n-1}U ]_1=[A_n]_0-[\sigma _n^FA_n]_0=0\;.$$

b) Let $U\in \unn{\check{H}_{n-1} }$, $V\in \unn{\check{G}_n }$, and $P\in Pr\,\check{F}_n $ with
$$\check{\psi }_nV=A_nU+B_nU^*,\qquad\qquad \check{\varphi }_nP=VA_nV^*\;.  $$
By \pr{927} c) (since $\check{\varphi} \circ \sigma ^F=\sigma ^G\circ \check{\varphi }  $),
$$K_0(\varphi )\delta _1[U]_1=K_0(\varphi )([P]_0-[\sigma _n^FP]_0)=$$
$$=[\check{\varphi }_nP]_0-[\check{\varphi }_n\sigma _n^FP ]_0=[\check{\varphi }_nP ]_0-[\sigma _n^G\check{\varphi }_nP ]_0=$$
$$=  [VA_nV^*]_0-[(\sigma _n^GV)A_n(\sigma _n^GV)^*]_0=[A_n]_0-[A_n]_0=0\;.\qedd$$

\begin{p}\label{967}
Let $U\in \unn{\check{H}_{n-1} }$. There are $V\in \check{G}_n $ and $P,Q\in \pp{\check{F}_n }$ such that
$$V^*V\in \pp{\check{G}_n },\qquad\qquad \check{\psi }_nV=A_nU\,,$$
$$ \check{\varphi }_nP=1_E-V^*V,\qquad  \check{\varphi }_nQ=1_E-VV^*,\qquad \delta _1[U]_1=[P]_0-[Q]_0\;.$$
\end{p}

By \pr{931} d), $A_nU+B_nU^*\in \unm{\check{H}_n }$. Since $\check{\psi }_n $ is surjective, by [R] Lemma 2.1.7 (i), there is a $V_0\in \unn{\check{G}_n }$ with $\check{\psi }_nV_0=A_nU+B_nU^* $. Put $V:=V_0A_n\in \check{G}_n $. Then
$$V^*V=A_nV_0^*V_0A_n=A_n\in \pp{\check{G}_n } $$
and
$$\check{\psi }_nV=(\check{\psi }_nV_0 )A_n=(A_nU+B_nU^*)A_n=A_nU\;. $$

We have
$$\check{\psi }_n(1_E-V^*V)=1_E-A_n=B_n=\check{\psi }_n(1_E-VV^*)\;.  $$
By \pr{934} $b_2\Rightarrow b_1$, there are $P,Q\in \pp{\check{F}_n }$ with
$$\check{\varphi }_nP=1_E-V^*V,\qquad\qquad \check{\varphi }_nQ=1_E-VV^*\;.  $$
Put
$$W:=A_{n+1}V+C_{n+1}(1_E-V^*V)+C_{n+1}^*(1_E-VV^*)+B_{n+1}V^*\in \check{G}_{n+1}\,, $$
$$Z:=A_n+(C_{n+1}+C_{n+1}^*)B_n\in E_{n+1}\;.$$
Since $VV^*V=V$, $V^*VV^*=V^*$, and
$$W^*=A_{n+1}V^*+C_{n+1}^*(1_E-V^*V)+C_{n+1}(1_E-VV^*)+B_{n+1}V\,,$$
we get
$$WW^*=A_{n+1}VV^*+B_{n+1}(1_E-V^*V)+A_{n+1}(1_E-VV^*)+B_{n+1}V^*V=$$
$$=A_{n+1}+B_{n+1}=1_E\,,$$
$$W^*W=A_{n+1}V^*V+A_{n+1}(1_E-V^*V)+B_{n+1}(1_E-VV^*)+B_{n+1}VV^*=$$
$$=A_{n+1}+B_{n+1}=1_E\;.$$
By \pr{931} a),
$$Z^2=A_n+B_n=1_E$$
so $W\in \unn{\check{G}_{n+1} }$, $Z\in \unn{E_{n+1}}$, and $ZW\in \unn{\check{G}_{n+1} }$. By the above and \pr{931} a),
$$\check{\psi }_{n+1}W=A_{n+1}A_nU+(C_{n+1}+C_{n+1}^*)B_n+B_{n+1}A_nU^*\,, $$
$$\check{\psi }_{n+1}(ZW)=Z\check{\psi }_{n+1}W=$$
$$=(A_n+(C_{n+1}+C_{n+1}^*)B_n)(A_{n+1}A_nU+(C_{n+1}+C_{n+1}^*)B_n+B_{n+1}A_nU^*)=  $$
$$=A_{n+1}A_nU+B_{n+1}A_nU^*+B_n=A_{n+1}A_nU+B_{n+1}A_nU^*+(A_{n+1}+B_{n+1})B_n=$$
$$=A_{n+1}(A_nU+B_n)+B_{n+1}(A_nU^*+B_n)\;.$$
We put
$$R:=A_{n+1}(1_E-Q)+B_{n+1}P\in \pp{\check{F}_{n+1} }\;.$$
Using again $VV^*V=V$ and $V^*VV^*=V^*$,
$$\check{\varphi }_{n+1}R=A_{n+1}VV^*+B_{n+1}(1_E-V^*V)\,, $$
$$WA_{n+1}=A_{n+1}V+C_{n+1}(1_E-V^*V)\,,$$
$$WA_{n+1}W^*=A_{n+1}VV^*+B_{n+1}(1_E-V^*V)=\check{\varphi }_{n+1}R\,, $$
$$ZWA_{n+1}W^*Z=Z(\check{\varphi }_{n+1}R )Z=\check{\varphi }_{n+1}(ZRZ)\;. $$
Since $ZRZ\sim _0R$ and $U\sim _1A_nU+B_n$, by the definition of $\delta _1$,
$$\delta _1[U]_1=\delta _1[A_nU+B_n]_1=[R]_0-[\sigma _{n+1}^FR]_0\;.$$
Since $\pi ^H\circ \check{\psi }\circ \check{\varphi }=\pi ^F  $, by the above,
$$\pi _n^FP=\pi _n^H\check{\psi }_n\check{\varphi }_nP=\pi _n^H\check{\psi }_n(1_E-V^*V)=\pi _n^HB_n=B_n=\pi _n^FQ\;.   $$
Thus by \pr{914} (and \pr{960} b)),
$$\sigma _{n+1}^FR=A_{n+1}(1_E-B_n)+B_{n+1}B_n\sim _0A_{n+1}B_n+A_{n+1}A_n=$$
$$=A_{n+1}=\bar{\rho }_{n+1}^{\check{F} }1_E\sim _01_E $$
and we get
$$[R]_0=[1_E-Q]_0+[P]_0=[1_E]_0+[P]_0-[Q]_0\,,$$
$$\delta _1[U]_1=[1_E]_0+[P]_0-[Q]_0-[1_E]_0=[P]_0-[Q]_0\;.\qedd$$

\begin{p}\label{969}
$Ker\;\delta _1\subset Im\;K_1(\psi )$.
\end{p}

Let $a\in Ker\;\delta _1$ and let $U\in \unn{\check{H}_{n-1} }$ with $a=[U]_1$. By \pr{967}, there are $V\in \check{G}_n $ and $P,Q\in \pp{\check{F}_n }$ such that $V^*V\in \pp{\check{G}_n }$, $\check{\psi }_nV=A_nU $,  
$$\check{\varphi }_nP=1_E-V^*V,\qquad \check{\varphi }_nQ=1_E-VV^*,\qquad \delta _1[U]_1=[P]_0-[Q]_0\;.  $$
Then $[P]_0=[Q]_0$. By \cor{918} a$\Rightarrow $c, there is an $m\in \bn$, $m>n+1$, and an $X\in \check{F}_m $ such that
$$X^*X=\left(\proo{i=n+1}{m}A_i\right)P+\left(1_E-\proo{i=n+1}{m}A_i\right)\,,$$
$$XX^*=\left(\proo{i=n+1}{m}A_i\right)Q+\left(1_E-\proo{i=n+1}{m}A_i\right)\;.$$
Put $W:=\check{\varphi }_mX $. Then 
$$W^*W=\check{\varphi }_m(X^*X)=\left(\proo{i=n+1}{m}A_i\right)(1_E-V^*V)+\left(1_E-\proo{i=n+1}{m}A_i\right)=$$
$$=1_E-\left(\proo{i=n+1}{m}A_i \right)V^*V\,,$$
$$WW^*=1_E-\left(\proo{i=n+1}{m}A_i\right)VV^*\,,$$
$$\left(\proo{i=n+1}{m}A_i\right)VV^*WW^*=\left(\proo{i=n+1}{m}A_i\right)V^*VW^*W=0\,,$$
$$\left(\proo{i=n+1}{m}A_i\right)V^*W=\left(\proo{i=n+1}{m}A_i\right)VW^*=0\,,$$
$$\left(\left(\proo{i=n+1}{m}A_i\right)V+W\right)^*\left(\left(\proo{i=n+1}{m}A_i\right)V+W\right)=$$
$$=\left(\proo{i=n+1}{m}A_i\right)V^*V+W^*W=1_E\,,$$
$$\left(\left(\proo{i=n+1}{m}A_i\right)V+W\right)\left(\left(\proo{i=n+1}{m}A_i\right)V+W\right)^*=$$
$$=\left(\proo{i=n+1}{m}A_i\right)VV^*+WW^*=1_E\,,$$
$$\left(\proo{i=n+1}{m}A_i\right)V+W\in \unn{\check{G}_m }\;.$$
From
$$\check{\psi }_m(W^*W)=1_E-\left(\proo{i=n+1}{m}A_i\right)\check{\psi }_m(V^*V)=$$
$$=1_E-\left(\proo{i=n+1}{m}A_i\right)A_n=\check{\psi }_m(WW^*)\,,  $$
since $\check{\psi }_mW=\check{\psi }_m\check{\varphi }_mX\in E_m$, it follows 
$$\check{\psi }_mW+\left(\proo{i=n}{m}A_i\right)\in \unn{E_m} \;.$$
By the above,
$$\left(\proo{i=n}{m}A_i\right)U\check{\psi }_mW^*=\left(\proo{i=n+1}{m}A_i\right)(\check{\psi }_mV )(\check{\psi }_mW^* )=$$
$$=\check{\psi }_m\left(\left(\proo{i=n+1}{m}A_i\right)VW^*\right)= 0\,,$$
$$(\check{\psi }_mW )^*(\check{\psi }_mW )\left(\proo{i=n}{m}A_i\right)=0\,,\qquad\qquad (\check{\psi }_mW )\left(\proo{i=n}{m}A_i\right)=0\,, $$
$$\check{\psi }_m\left(\left(\proo{i=n+1}{m}A_i\right)V+W\right)=\left(\proo{i=n}{m}A_i\right)U+\check{\psi }_mW\sim _1  $$
$$\sim _1\left(\left(\proo{i=n}{m}A_i\right)U+\check{\psi }_mW \right)\left(\left(\proo{i=n}{m}A_i\right)+\check{\psi }_mW^* \right)=$$
$$=\left(\left(\proo{i=n}{m}A_i\right)U+\left(1_E-\proo{i=n}{m}A_i\right)\right)\;.$$
By \pr{951a} and \pr{954} c),
$$a=[U]_1=\left[\left(\proo{i=n}{m}A_i\right)U+\left(1_E-\proo{i=n}{m}A_i\right)\right]_1=$$
$$=\left[\check{\psi }_m\left(\left(\proo{i=n+1}{m}A_i\right)V+W\right) \right]_1=$$
$$=K_1(\psi )\left[\left(\proo{i=n+1}{m}A_i\right)V+W\right]_1\in Im\;K_1(\psi )\;.\qedd$$

\begin{p}\label{971}
$Ker\;K_0(\varphi )\subset Im\;\delta _1$.
\end{p}

Let $a\in Ker\;K_0(\varphi )$. By \pr{929}, there is a $P\in \pp{\check{F}_\rightarrow  }$ with
$$a=[P]_0-[\sigma _\rightarrow ^FP]_0\;.$$
By \pr{927} c),
$$0=K_0(\varphi )a=[\check{\varphi }_\rightarrow P ]_0-[\check{\varphi }_\rightarrow \sigma _\rightarrow ^FP ]_0\;.$$
Let \bbn\, such that $P\in \pp{\check{F}_{\rightarrow \,n} }$. Then $[\check{\varphi }_{\rightarrow \,n}P ]_0=[\check{\varphi }_{\rightarrow \,n}\sigma _{\rightarrow \,n}^FP ]_0$. By \cor{918} a$\Rightarrow $c, there is an $m\in \bn$, $m>n+1$, such that
$$\check{\varphi }_{\rightarrow \,n}P+(B_m)_\rightarrow \sim _0\check{\varphi }_{\rightarrow \,n}\sigma _{\rightarrow \,n}^FP+(B_m)_\rightarrow \;.  $$
Put
$$Q:=P+(B_m)_\rightarrow \in \pp{\check{F}_{\rightarrow \,m} }\;.$$
Then
$$a=[Q]_0-[\sigma _\rightarrow ^FQ]_0,\qquad\qquad \check{\varphi }_{\rightarrow \,m}Q\sim _0\check{\varphi }_{\rightarrow \,m}\sigma _{\rightarrow \,m}^FQ=\sigma _{\rightarrow\,m}^FQ\;.  $$
By \pr{930}, there are $k\in \bn$, $k\geq m+2$, and $W\in \unn{\check{G}_{\rightarrow \,k} }$ with
$$W(\check{\varphi }_{\rightarrow \,m}Q )W^*=\sigma _{\rightarrow \,m}^FQ\;. $$
It follows
$$(\sigma ^F_{\rightarrow m}Q)W=W(\check{\varphi }_{\rightarrow m}Q )W^*W=W(\check{\varphi }_{\rightarrow m}Q )\,,$$
$$(\check{\psi }_{\rightarrow \,k}W )(\sigma _{\rightarrow \,k}^FQ)=(\check{\psi}_{\rightarrow \,k}W)(\check{\psi }_{\rightarrow \,k}\check{\varphi }_{\rightarrow \,k}Q  )=\check{\psi }_{\rightarrow \,k}(W\check{\varphi }_{\rightarrow \,k}Q ) =$$
$$=\check{\psi }_{\rightarrow \,k}((\sigma _{\rightarrow \,k}^FQ)W)= 
(\sigma _{\rightarrow \,k}^FQ )(\check{\psi }_{\rightarrow \,k}W)\;.$$
Put
$$U:=(\check{\psi }_{\rightarrow \,k}W )(1_E-\sigma _{\rightarrow \,k}^FQ)+\sigma _{\rightarrow \,k}^FQ\in \check{H}_{\rightarrow \,k}\;. $$
Then
$$UU^*=U^*U=1_E,\qquad\qquad\qquad U\in \unn{\check{H}_{\rightarrow \,k} }\;.$$
Put
$$V_1:=(A_{k+1})_\rightarrow (1_E-\sigma _{\rightarrow \,k}^FQ)W+(B_{k+1})_\rightarrow \sigma _{\rightarrow \,k}^FQ\in \check{G}_{k+1}\;. $$
Then
$$V_1^*=(A_{k+1})_\rightarrow W^*(1_E-\sigma _{\rightarrow \,k}^FQ)+(B_{k+1})_\rightarrow \sigma _{\rightarrow \,k}^FQ\,,$$
$$V_1V_1^*=(A_{k+1})_\rightarrow (1_E-\sigma _{\rightarrow \,k}^FQ)+(B_{k+1})_\rightarrow \sigma _{\rightarrow \,k}^FQ\in \pp{E_{k+1}}\,,$$
$$V_1^*V_1=(A_{k+1})_\rightarrow W^*(1_E-\sigma _{\rightarrow \,k}^FQ)W+(B_{k+1})_\rightarrow \sigma _{\rightarrow \,k}^FQ=$$
$$=(A_{k+1})_\rightarrow (1_E-W^*(\sigma _{\rightarrow \,k}^FQ)W)+(B_{k+1})_\rightarrow \sigma _{\rightarrow \,k}^FQ\;.$$
Put
$$Z:=(1_E-\sigma _{\rightarrow \,k}^FQ)+((C_{k+1})_\rightarrow +(C_{k+1}^*)_\rightarrow )\sigma _{\rightarrow \,k}^FQ\in E_{k+1}\;.$$
By \pr{931} a),
$$Z^2=(1_E-\sigma _{\rightarrow \,k}^FQ)+\sigma _{\rightarrow \,k}^FQ=1_E,\qquad\qquad Z\in \unn{E_{k+1}}\,,$$
$$ZV_1=(A_{k+1})_\rightarrow (1_E-\sigma _{\rightarrow \,k}^FQ)W+(C_{k+1}^*)_\rightarrow \sigma _{\rightarrow \,k}^FQ\,,$$
$$V:=ZV_1Z=(A_{k+1})_\rightarrow (1_E-\sigma _{\rightarrow \,k}^FQ)W(1_E-\sigma _{\rightarrow \,k}^FQ)+$$
$$+(C^*_{k+1})_\rightarrow (1_E-\sigma _   {\rightarrow k}^FQ)W\sigma _{\rightarrow k}^FQ+(A_{k+1})_\rightarrow \sigma _{\rightarrow \,k}^FQ\in \check{G}_{\rightarrow \,k+1}\,, $$
$$\check{\psi }_\rightarrow V=(A_{k+1})_\rightarrow (1_E-\sigma _{\rightarrow \,k}^FQ)\check{\psi }_{\rightarrow \,k}W+(A_{k+1})_\rightarrow \sigma _{\rightarrow \,k}^FQ=(A_{k+1})_\rightarrow U\,, $$
$$VV^*=ZV_1V_1^*Z\in \pp{E_{k+1}},\qquad\qquad V^*V=ZV_1^*V_1Z\,,$$
$$1_E-VV^*=Z(1_E-V_1V_1^*)Z=$$
$$=Z((A_{k+1})_\rightarrow \sigma _{\rightarrow \,k}^FQ+(B_{k+1})_\rightarrow (1_E-\sigma _{\rightarrow \,k}^FQ))Z\,,$$
$$1_E-V^*V=Z(1_E-V_1^*V_1)Z=$$
$$=Z((A_{k+1})_\rightarrow W^*(\sigma _{\rightarrow \,k}^FQ)W+(B_{k+1})_\rightarrow (1_E-\sigma _{\rightarrow \,k}^FQ))Z=$$
$$=Z((A_{k+1})_\rightarrow \check{\varphi }_{\rightarrow \,k}Q+(B_{k+1})_\rightarrow (1_E-\sigma _{\rightarrow \,k}^FQ) )Z\,,$$
$$\check{\varphi }_{\rightarrow ,k+1}(Z((A_{k+1})_\rightarrow Q+(B_{k+1})_\rightarrow (1_E-\sigma _{\rightarrow \,k}^FQ))Z)= $$
$$=Z((A_{k+1})_\rightarrow \check{\varphi }_kQ+(B_{k+1})_\rightarrow (1_E-\sigma _{\rightarrow \,k}^FQ))Z=1_E-V^*V\,, $$
$$\check{\varphi }_{\rightarrow ,k+1}(Z((A_{k+1})_\rightarrow \sigma _{\rightarrow\,k} ^FQ+(B_{k+1})_\rightarrow (1_E-\sigma _{\rightarrow\,k} ^FQ))Z))=1_E-VV^*\;. $$
By \pr{967},
$$\delta _1[U]_1=
[Z((A_{k+1})_\rightarrow Q+(B_{k+1})_\rightarrow (1_E-\sigma _{\rightarrow \,k}^FQ))Z]_0-$$
$$-[Z((A_{k+1})_\rightarrow \sigma _{\rightarrow\,k} ^FQ+(B_{k+1})_\rightarrow (1_E-\sigma _{\rightarrow \,k}^F)Q)Z]_0=
[Q]_0-[\sigma _\rightarrow ^FQ]_0=a\;.$$
Thus $a\in Im\;\delta _1$.\qed

\begin{theo}\label{974}
The sequence
$$K_1(F)\stackrel{\scriptscriptstyle K_1(\varphi )}{\longrightarrow }K_1(G)\stackrel{\scriptscriptstyle K_1(\psi )}{\longrightarrow }K_1(H)\stackrel{\delta _1}{\longrightarrow }K_0(F)\stackrel{\scriptscriptstyle K_0(\varphi )}{\longrightarrow }K_0(G)\stackrel{\scriptscriptstyle K_0(\psi )}{\longrightarrow }K_0(H)$$
is exact.
\end{theo}

The exactness was proved: for $K_1(G)$ in \pr{958}, for $K_1(H)$ in \pr{969} and \pr{965} a), for $K_0(F)$ in \pr{971} and \pr{965} b), and for $K_0(G)$ in \pr{934} c).\qed 

\begin{center}
\section{$K_1(F)\approx K_0(S F)$}
\end{center} 

\begin{de}\label{975}
Let $F$ be an \eo. We denote by $CF$ the \eo\, of continuous maps $\mac{x}{[0,1]}{F}$ with $x(0)=0$ and by $SF$ its $E$-C*-subalgebra $\me{x\in CF}{x(1)=0}$ \emph{(\dd{28.3'a} or [C2] Corollary 1.2.5 a),d))}. Moreover we denote by $\mac{\theta _F}{K_1(F)}{K_0(SF)}$ the index map associated to the exact sequence
$$\oc{SF}{i_F}{CF}{j_F}{F}\,,$$
in \frm, where $i_F$ is the inclusion map and
$$\mae{j_F}{CF}{F}{x}{x(1)}\;.$$
If $\oa{F}{\varphi }{G}$ is a morphism in \frm\, then we put
$$\mae{S\varphi }{SF}{SG}{x}{\varphi \circ x}\,,$$
$$\mae{C\varphi }{CF}{CG}{x}{\varphi \circ x}\;.$$
\end{de}

If $\ob{F}{\varphi }{G}{\psi }{H}$ are morphisms in \frm then $S(\psi )\circ S(\varphi )=S(\psi \circ \varphi )$.

\begin{theo}\label{976}
$\theta _F$ is a group isomorphism for every \eo\, F.
\end{theo}

$CF$ is null-homotopic ([R] Example 4.1.5 or \pr{24.11a}), so by the Homotopy invariance (\h{940} e), \pr{956} e)), it is $K$-null.
By \h{974}, the sequence
$$K_1(CF)\stackrel{\scriptscriptstyle K_1(j_F)}{\longrightarrow }K_1(F)\stackrel{\theta _F}{\longrightarrow }K_0(SF)\stackrel{\scriptscriptstyle K_0(i_F)}{\longrightarrow }K_0(CF)$$
is exact, so $\theta _F$ is a group isomorphism.\qed

\begin{p}\label{975a}
Let $F$ and $G$ be \eo $\!$s.
\begin{enumerate}
\item For all $(x,y)\in (SF)\times (SG)$ put
$$\mae{\overbrace{(x,y)}}{[0,1]}{F\times G}{s}{(x(s),y(s))}\;.$$
Then the map
$$\mad{(SF)\times (SG)}{S(F\times G)}{(x,y)}{\overbrace{(x,y)} }$$
is an isomorphism in \frm \emph{(\dd{8.6})}.
\item $K_1(F)\times K_1(G)\approx K_1(F\times G)$ {\bf (Product Theorem)}.
\end{enumerate}
\end{p}

a) is easy to see.

b) By \h{976}, the maps
$$\oa{K_1(F)\times K_1(G)}{\theta _F\times \theta _G}{K_0(SF)\times K_0(SG)},\qquad \oa{K_1(F\times G)}{\theta _{F\times G}}{K_0(S(F\times G))}$$
are group isomorphisms. By a), $K_0((SF)\times (SG))\approx  K_0(S(F\times G))$  and by \cor{938} b),  $K_0((SF)\times (SG))\approx K_0(SF)\times K_0(SG)$. Thus 
$$K_1(F)\times K_1(G)\approx K_1(F\times G)\;.\qedd$$

\begin{co}\label{13.11}
Let \oaa{F}{\varphi }{F'}, $\oaa{G}{\psi }{G'}$ be morphisms in \frm and
$$\mae{\varphi \times \psi }{F\times G}{F'\times G'}{(x,y)}{(\varphi x,\psi y)}\;.$$
Then $\varphi \times \psi $ is a morphism in \frm and
$$K_i(\varphi \times \psi )=K_i(\varphi) \times K_i(\psi )$$
for all $i\in \z{0,1}$.
\end{co}

The assertion follows easily from \cor{938} b) and \pr{975a} b).\qed

\begin{p}[Product Theorem]\label{28.3'}
Let $(F_j)_{j\in J}$ be a finite family of $E$-C*-algebras, $F:=\pro{j\in J}F_j$ \emph{(\dd{8.6})}, and for every $j\in J$ let $\mac{\varphi _j}{F_j}{F}$ be the canonical inclusion and $\mac{\psi _j}{F}{F_j}$ the projection. Then for every $i\in \z{0,1}$,
$$\mae{\Phi }{\pro{j\in J}K_i(F_j)}{K_i(F)}{(a_j)_{j\in J}}{\si{j\in J}}K_i(\varphi _j)a_j$$
is a group isomorphism and
$$\mae{\Psi }{K_i(F)}{\pro{j\in J}K_i(F_j)}{a}{(K_i(\psi _j)a)_{j\in J}}$$
is its inverse.
\end{p}

$\Phi $ and $\Psi $ are obviously group homomorphisms. For $j,k\in J$, $\psi _j\circ \varphi _k=0$ if $j\not=k$ and $\psi _j\circ \varphi _j=id_{F_j}$. Thus for $(a_j)_{j\in J}\in \pro{j\in J}K_i(F_j)$ and $k\in J$,
$$(\Psi \Phi (a_j)_{j\in J})_k=K_i(\psi _k)\si{j\in J}K_i(\varphi _j)a_j=a_k$$
i.e. $\Psi \circ \Phi $ is the identity map of $\pro{j\in J}K_i(F_j)$. Since $\si{j\in J}\varphi _j\circ \psi _j=id_F$, for $a\in K_i(F)$,
$$\Phi \Psi a=\Phi (K_i(\psi _j)a)_{j\in J}=\si{j\in J}K_i(\varphi _j)K_i(\psi _j)a=K_i\left(\si{j\in J}\varphi _j\circ \psi _j\right)a=a$$
i.e. $\Phi \circ \Psi=id_{K_i(F)} $.\qed

\begin{theo}[Continuity of $K_1$]\label{17.2}
Let $\{(F_i)_{i\in I},\,(\varphi _{ij})_{i,j\in I}\}$ be an inductive system in \frm\, and let $\{F,\,(\varphi _i)_{i\in I}\}$ be its limit in \frm\,. By \emph{\pr{956} a)}, 
$$\{(K_1(F_i))_{i\in I},\,(K_1(\varphi _{ij}))_{i,j\in I}\}$$
 is an inductive system in the category of additive groups. Let $\{\ccc{G},\,(\psi _i)_{i\in I}\}$ be its limit in this category and let $\mac{\psi }{\ccc{G}}{K_1(F)}$ be the group homomorphism such that $\psi \circ \psi _i=K_1(\varphi _i)$ for every $i\in I$. Then $\psi $ is a group isomorphism.
\end{theo}

By [R] Exercise 10.2, $\{SF,(S\varphi _i)_{i\in I}\}$ is the limit in \frm of the inductive system $\{(SF_i)_{i\in I},(S\varphi _{ij})_{i,j\in I}\}$. By \h{941}, $\{K_0(SF),(K_0(S\varphi _i))_{i\in I}\}$ may be identified with the inductive limit in the category of additive groups of the inductive system $\{K_0(SF_i)_{i\in I},(K_0(S\varphi _{ij}))_{i,j\in I}\}$ and the assertion follows from \h{976}.\qed

\begin{p}\label{977}
Let $F$ be an \eo, $n\in \bn$, $U\in \unn{\check{F}_{n-1} }$, $V\in \unn{(\check{\overbrace{CF}})_n }$, and $P\in \pp{({\check{\overbrace{SF}}})_n}$ such that
$$\check{j}_FV=A_nU+B_nU^*,\qquad\qquad \check{i}_FP=VA_nV^*\;.  $$
Then
$$\theta _F[U]_1=[P]_0-[\sigma _n^{SF}P]_0.$$
\end{p}

The assertion follows from \cor{964} and \dd{975}.\qed

\begin{p}\label{978}
If \oa{F}{\varphi }{G} is a morphism in \frm\, then the diagram
$$\begin{CD}
K_1(F)@>K_1(\varphi )>>K_1(G)\\
@V\theta _FVV         @VV\theta _GV\\
K_0(SF)@>>K_0(S\varphi )>    K_0(SG)
\end{CD}$$
is commutative.
\end{p}

The diagram
$$\begin{CD}
0@>>>SF@>i_F>>CF@>j_F>>F@>>>0\\
@.@VS\varphi VV  @VC\varphi VV  @VV\varphi V@.\\
0@>>>SG@>>i_G>CG@>>j_G>G@>>>0
\end{CD}$$
is commutative and the assertion follows from \pr{966}.\qed

{\it Remark.} By \h{976} and \pr{978}, the functor $K_1$ is determined by the functor $K_0$.

\begin{co}[Split Exact Theorem]\label{987}
If 
$$\od{F}{\varphi }{G}{\psi }{\gamma }{H}$$
is a split exact sequence in \frm then
$$\og{K_1(F)}{K_1(\varphi )}{K_1(G)}{K_1(\psi )}{K_1(\gamma )}{K_1(H)}{20}{20}{20}$$
is also split exact. In particular the map
$$\mad{K_1(F)\times K_1(H)}{K_1(G)}{(a,b)}{K_1(\varphi )a+K_1(\lambda )b}$$
is a group isomorphism and $\kk{1}{\check{F} }\approx \kk{1}{E}\times \kk{1}{F}$.  
\end{co}

By \h{974}, the sequence
$$K_1(F)\stackrel{\scriptscriptstyle K_1(\varphi )}{\longrightarrow }K_1(G)\stackrel{\scriptscriptstyle K_1(\psi )}{\longrightarrow }K_1(H)\stackrel{\delta _1}{\longrightarrow }K_0(F)\stackrel{\scriptscriptstyle K_0(\varphi )}{\longrightarrow }K_0(G)\stackrel{\scriptscriptstyle K_0(\psi )}{\longrightarrow }K_0(H)$$
is exact and by \pr{956} a) and \pr{954} d),
$$K_1(\psi )\circ K_1(\gamma )=K_1(\psi \circ \gamma )=K_1(id_H)=id_{K_1(H)}\;.$$
It remains only to prove that $K_1(\varphi )$ is injective.

It is easy to see that
$$\od{SF}{S\varphi }{SG}{S\psi }{S\gamma }{SH}$$
is split exact. By \pr{936}, $K_0(S\varphi )$ is injective and by \pr{978}, the diagram
$$\begin{CD}
K_1(F)@>K_1(\varphi )>>K_1(G)\\
@V\theta _FVV   @VV\theta _GV\\
K_0(SF)@>>K_0(S\varphi )>K_0(SG)
\end{CD}$$
is commutative. Since $\theta _F$ is injective (\h{976}), $K_1(\varphi )$ is also injective.

The last assertion follows from the fact that
$$\od{F}{\iota ^F}{\check{F} }{\pi ^F}{\lambda ^F}{E}$$
is split exact.\qed

\begin{co}\label{24.10}
Let 
$$\od{F}{\varphi }{G}{\psi }{\gamma }{H}\,,\qquad\qquad 
\od{F'}{\varphi' }{G'}{\psi' }{\gamma' }{H'}$$
be split exact sequences in $\frm$ and 
$$\oa{F}{\lambda }{F'},\qquad \oa{G}{\mu }{G'},\qquad \oa{H}{\nu }{H'}$$
morphisms in $\frm$ such that the corresponding diagram is commutative and let $i\in \z{0,1}$.
\begin{enumerate}
\item If we denote by
$$\mae{\phi }{K_i(F)\times K_i(H)}{K_i(G)}{(a,b)}{K_i(\varphi )a+K_i(\gamma )b}\,,$$
$$\mae{\phi' }{K_i(F')\times K_i(H')}{K_i(G')}{(a',b')}{K_i(\varphi' )a'+K_i(\gamma ')b'}$$
the group isomorphisms \emph{(\pr{936}, \cor{987})} then
$$K_i(\mu )\circ K_i(\phi )=K_i(\phi ')\circ (K_i(\lambda )\times K_i(\nu ))\;.$$
\item If we identify $K_i(G)$ with $K_i(F)\times K_i(H)$ using $\phi $ and $K_i(G')$ with $K_i(F')\times K_i(H')$ using $\phi' $ then
$$\mae{K_i(\mu )}{K_i(G)}{K_i(G')}{(a,b)}{(K_i(\lambda )a,K_i(\nu )b)}\;.$$
\end{enumerate}
\end{co}

a) For $(a,b)\in K_i(F)\times K_i(H)$,
$$K_i(\mu )K_i(\phi )(a,b)=K_i(\mu )(K_i(\varphi )a+K_i(\gamma )b)=$$
$$=K_i(\varphi ')K_i(\lambda )a+K_i(\gamma ')K_i(\nu )b=
K_i(\phi ')(K_i(\lambda )\times K_i(\nu ))(a,b)\;.$$

b) follows from a).\qed

\begin{center}
\chapter{Bott periodicity}

\section{The Bott map}
\end{center}

\begin{lem}\label{979}
Let $F$ be a full $E$-C*-algebra and $\bbn$. We identify $SF$ with $\cbb{\bt\setminus \z{1}}{F}$ in an obvious way. 
\begin{enumerate}
\item $F_{\btt}:=\me{X\in \ccc{C}(\bt,F)}{X(1)\in E}$ is a full $E$-C*-subalgebra of $\ccc{C}(\bt,F)$.
\item If we put for every $(\alpha ,x)\in \check{\overbrace{SF}}$
$$\mae{\overbrace{(\alpha ,x)}}{\bt}{F}{z}{\alpha +x(z)}$$
then the map
$$\mae{\psi }{\check{\overbrace{SF}}}{F_{\btt}}{(\alpha ,x)}{\overbrace{(\alpha ,x)}}$$
is an $E$-C*-isomorphism. Thus the map
$$\mac{\psi _n}{\left(\check{\overbrace{SF}}\right)\!_n}{(F_{\btt})_n}$$
is also an $E$-C*-isomorphism.
\item For every $Y\in (F_{\btt})_n$ put
$$\mae{\ddot{Y} }{\bt}{F_n}{z}{\si{t\in T_n}(Y_t(z)\otimes id_K)V_t}\;.$$
Then $\ddot{Y}\in \me{X\in \ccc{C}(\bt,F_n)}{X(1)\in E_n} $ for every $Y\in (F_{\btt})_n$ and the map
$$\mae{\phi ^n}{(F_{\btt})_n}{\me{X\in \ccc{C}(\bt,F_n)}{X(1)\in E_n}}{Y}{\ddot{Y}  }$$
is an $E$-C*-isomorphism.
\item The map
$$\mac{\phi ^n\circ \psi _n}{\left(\check{\overbrace{SF}}\right)\!_n}{\me{X\in \ccc{C}(\bt,F_n)}{X(1)\in E_n}}$$
is an $E$-C*-isomorphism. 
We identify these two full $E$-C*-algebras by using this isomorphism.The map
$$\unn{\left(\check{\overbrace{SF}}\right)\!_n}\longrightarrow \me{X\in \ccc{C}(\bt,\unn{F_n})}{X(1)\in \unn{E_n}}$$
defined by $\phi ^n\circ \psi _n$ is a homeomorphism.
\item  For every 
$$X:=\si{t\in T_n}((\alpha _t,X_t)\otimes id_K)V_t\in {\left(\check{\overbrace{SF}}\right)\!_n}$$
and $z\in \bt$,
$$(\phi ^n\psi _nX)(z)=\si{t\in T_n}((\alpha _t+X_t(z))\otimes id_K)V_t\in F_n\,,$$
$$(\phi ^n\psi _nX)(1)=\si{t\in T_n}(\alpha _t\otimes id_K)V_t\in E_n\;.$$
\item Consider the split exact sequence in \frm \emph{(\dd{924})}
$$\og{SF}{\iota ^{SF}}{\check{\overbrace{SF}}}{\pi ^{SF}}{\lambda ^{SF}}{E}{0}{10}{15}\;.$$
Then
$$\left(\pi ^{SF}\right)\!_nX=\left(\phi ^n\psi _nX\right)(1)$$
for every $X\in \left(\check{\overbrace{SF}}\right)\!_n$.
\item If $\oaa{F}{\varphi }{G}$ is a morphism in \frc then, by the identification of d), for every $X\in \ccc{C}(\bt,F_n)$ with $X(1)\in E_n$ and for every $z\in \bt$,
$$\left(\left(\overbrace{S\varphi }^{\check{} }\right)\!_nX\right)(z)=\varphi _nX(z)\;.$$ 
\end{enumerate}
\end{lem}

a) is obvious.

b) For $(\alpha ,x),(\beta ,y)\in \check{\overbrace{SF}}$, $\gamma \in E$, and $z\in \bt$,
$$(\overbrace{(\alpha ,x)})^*(z)=\alpha ^*+x(z)^*=\overbrace{(\alpha ,x)^*}(z)\,,$$
$$(\overbrace{(\alpha ,x)}(z))(\overbrace{(\beta ,y)}(z))=(\alpha +x(z))(\beta +y(z))=\alpha \beta +\alpha y(z)+x(z)\beta +x(z)y(z)=$$
$$=\overbrace{(\alpha \beta ,\alpha y+\beta x+xy)}(z)=\overbrace{(\alpha ,x)(\beta ,y)}(z)\,,$$
$$\overbrace{(\gamma ,0)}(z)=\gamma\,, $$
so $\psi $ is an $E$-C*-homomorphism. If $\overbrace{(\alpha ,x)}=0$ then for all $z\in \bt$
$$\alpha =\alpha +x(1)=0\,,\qquad x(z)=\alpha +x(z)=0\,,\qquad x=0\,,$$
so $\psi $ is injective.

Let $X\in F_{\btt}$ and put $\alpha :=X(1)\in E$ and
$$\mae{x}{\bt}{F}{z}{X(z)-X(1)}\;.$$
Then $(\alpha ,x)\in \check{\overbrace{SF}}$ and for $z\in \bt$,
$$\overbrace{(\alpha ,x)}(z)=\alpha +x(z)=X(1)+X(z)-X(1)=X(z)\;.$$
Thus $\overbrace{(\alpha ,x)}=X$ and $\psi $ is surjective.

By [C2] Corollary 2.2.5 and [C2] Theorem 2.1.9 a), $\psi _n$ is an isomorphism.

c) follows from [C2] Proposition 2.3.7 and [C2] Theorem 2.1.9 a).

d) follows from b) and c).

e) We have
$$\psi _nX=\si{t\in T_n}(\overbrace{(\alpha _t,X_t)}\otimes id_K)V_t\,,$$
$$(\phi ^n\psi _nX)(z)=\si{t\in T_n}((\alpha _t+X_t(z))\otimes id_K)V_t\in F_n\,,$$
$$(\phi ^n\psi _nX)(1)=\si{t\in T_n}(\alpha _t\otimes id_K)V_t\in E_n\;.$$

f) and g) follow from e).\qed

\begin{de}\label{982}
We put for every full $E$-C*-algebra $F$, $\bbn$, and $P\in F_n$,
$$\mae{\widetilde{P} }{\bt}{F_n}{z}{zP+(1_E-P)}\;.$$
\end{de}

By the identification of \lm{979} d), 
$$\widetilde{P}\in \me{X\in \ccb{\bt}{\unn{F_n}}}{X(1)\in E_n}=\unn{\left(\check{\overbrace{SF}}\right)\!_n}  $$
 for every $P\in \pp{F_n}$. Obviously, $\widetilde{0}=1_E $ and $\widetilde{1_E}=z1_E $. 

\begin{p}\label{983}
If $F$ is a full $E$-C*-algebra, $\bbn$, and $P\in \pp{F_{n-1}}$ then
$$\bar{\tau }_n^{\check{\overbrace{SF}} }\widetilde{P}=\widetilde{\bar{\rho }_n^FP }\,,$$
(with the identification of \emph{\lm{979} d)}). Thus we get a well-defined map
$$\mac{\nu _F}{\pp{F_\rightarrow }}{un\check{\overbrace{SF}}}$$
with $\nu _FP=\widetilde{P} $ for every $P\in \pp{F_\rightarrow }=\bigcup\limits_{n\in \bn}Pr\,F_{\rightarrow n} $.
\end{p}

For $z\in \bt$, 
$$(\bar{\tau }_n^{\check{\overbrace{SF}}}\widetilde{P} )(z)=(A_n\widetilde{P}+B_n)(z)=A_n(zP+(1_E-P))+B_n=$$
$$=zA_nP+(1_E-A_nP)=\widetilde{\bar{\rho }_n^FP }(z)\;.\qedd  $$

\begin{p}\label{984}
For every full $E$-C*-algebra $F$ there is a unique group homomorphism
$$\mac{\beta _F}{K_0(F)}{K_1(SF)}\qquad\qquad \mbox{\bf{(the Bott map)}}$$
such that for every $P\in \pp{F_\rightarrow }$,
$$\beta _F[P]_0=(\nu _FP)/\sim _1=\left[\tilde{P} \right]_1\;.$$
\end{p}

Let $P,Q\in \pp{F_\rightarrow }$ with $P\sim _0Q$. By \pr{930}, there are $m,n\in \bn$, $m\geq n+2$, and $U\in \unm{F_m}$ with $P,Q\in \pp{F_n}$ and $UPU^*=Q$ and so
$$(U\widetilde{P}U^* )(z)=U\widetilde{P}(z)U^*=zUPU^*+(1_E-UPU^*)=\widetilde{Q}(z)  $$
for every $z\in \bt$. Thus $U\widetilde{P}U^*=\widetilde{Q}  $, $\widetilde{P}\sim _h\widetilde{Q}  $, and $\widetilde{P}\sim _1\widetilde{Q}  $. 

Let $P,Q\in \pp{{F_\rightarrow }}$ with $PQ=0$. We may assume $P,Q\in \pp{F_{n-1}}$ with $P=PA_n$ and $Q=QB_n$ for some $\bbn$ (\pr{914}). For every $z\in \bt$,
$$\widetilde{P}(z)=zPA_n+(1_E-PA_n),\qquad\qquad \widetilde{Q}(z)=zQB_n+(1_E-QB_n)\,,  $$
$$(\widetilde{P}\widetilde{Q})(z)=\widetilde{P}(z)\widetilde{Q}(z)=zPA_n+zQB_n+1_E-QB_n-PA_n=$$
$$=z(P+Q)+(1_E-(P+Q))=
\widetilde{(P+Q)}(z),\qquad\qquad 
 \widetilde{P}\widetilde{Q}=\widetilde{P+Q}\;.$$

By \pr{920}, there is a unique group homomorphism 
$$\mac{\beta _F}{K_0(F)}{K_1(SF)}$$
 with the required property.\qed

\begin{p}\label{985}
Let $F$ be an \eo.
\begin{enumerate}
\item There is a unique map $\mac{\beta _F}{K_0(F)}{K_1(SF)}$ (called {\bf{the Bott map}}) such that the diagram 
$$\begin{CD}
K_0(F)@>K_0(\iota ^F)>>K_0(\check{F} )\\
@V\beta _FVV       @VV\beta _{\check{F} }V\\
K_1(SF)@>>K_1(S\iota ^F)>K_1(S\check{F} )
\end{CD}$$
is commutative. $\beta _F$ is a group homomorphism.
\item If $F$ is a full $E$-C*-algebra then the above map $\beta _F$ coincides with the map $\beta _F$ defined in \emph{\pr{984}}.
\item If $\oa{F}{\varphi }{G}$ is a morphism in \frm then the diagram
$$\begin{CD}
K_0(F)@>K_0(\varphi )>>K_0(G)\\
@V\beta _FVV         @VV\beta _GV\\
K_1(SF)@>>K_1(S\varphi )>K_1(SG)
\end{CD}$$
is commutative. 
\end{enumerate}
\end{p}

c) for \frc with $\oa{F}{\varphi }{G}$ unital. For $\bbn$, $P\in \pp{F_n}$, and $z\in \bt$, by \lm{979} g),
$$\left(\left(\check{\overbrace{S\varphi }}\right)\!_n\widetilde{P}\right)(z)=z\varphi _nP+(1_E-\varphi _nP)=\left(\widetilde{\varphi _nP}\right)(z)\,, $$
$$\left(\check{\overbrace{S\varphi }}\right)\!_n\widetilde{P}=\widetilde{\varphi _nP}\;. $$
By \pr{921} c), \pr{984}, and \pr{954} c),
$$K_1(S\varphi )\beta _F[P]_0=K_1(S\varphi )\left[\widetilde{P} \right]_1=$$
$$=\left[\left(\check{\overbrace{S\varphi }}\right)\!_n\widetilde{P} \right]_1=
\left[\widetilde{\varphi _nP}\right]_1=\beta _G[\varphi _nP]_0=\beta _GK_0(\varphi )[P]_0\,,$$
$$K_1(S\varphi )\circ \beta _F=\beta _G\circ K_0(\varphi )\;.$$

a) By c) for \frc, the diagram
$$\begin{CD}
K_0(\check{F} )@>K_0(\pi ^F)>>K_0(E)\\
@V\beta _{\check{F} }VV     @VV\beta _EV\\
K_1(S\check{F} )@>>K_1(S\pi ^F)>K_1(SE)
\end{CD}$$
is commutative. By \pr{925} c) and \cor{987} the sequences
$$\of{K_0(F)}{\scriptscriptstyle K_0(\iota ^F)}{K_0(\check{F} )}{\scriptscriptstyle K_0(\pi ^F)}{K_0(E)}{20}{20}\,,$$
$$\of{K_1(SF)}{\scriptscriptstyle K_1(S\iota ^F)}{K_1(S\check{F} )}{\scriptscriptstyle K_1(S\pi ^F)}{K_1(SE)}{20}{20}$$
are exact, since the sequence
$$\og{SF}{S\iota ^F}{S\check{F} }{S\pi ^F}{S\lambda^F}{SE}{0}{10}{10}$$
is split exact. By the above c) for \frc, \cor{928} a), and \pr{927} e),
$$K_1(S\pi ^F)\circ \beta _{\check{F} }\circ K_0(\iota ^F)=\beta _E\circ K_0(\pi ^F)\circ K_0(\iota ^F)=$$
$$=\beta _E\circ K_0(\pi ^F\circ \iota ^F)=\beta _E\circ K_0(0)=0\;.$$
Thus
$$Im\,(\beta _{\check{F} }\circ K_0(\iota ^F))\subset Ker\,K_1(S\pi ^F)=Im\,K_1(S\iota ^F)\;.$$
The assertion follows now from the fact that $K_1(S\iota ^F)$ is injective.

b) By c) for \frc, the diagram
$$\begin{CD}
K_0(F)@>K_0(\iota ^F)>>K_0(\check{F} )\\
@V\beta _FVV           @VV\beta _{\check{F} }V\\
K_1(SF)@>>K_1(S\iota ^F)>K_1(S\check{F} )
\end{CD}$$
is commutative, with $\beta _F$ defined in \pr{984}. By a), this $\beta _F$ coincides with $\beta _F$ defined in a).

c) The following diagrams

\parbox{2cm}{
$$\begin{CD}
F@>\varphi >>G\\
@V\iota ^FVV @VV\iota ^GV\\
\check{F}@>>\check{\varphi } >\check{G}  
\end{CD}$$}
\hspace{0.7cm}
\parbox{3cm}{
$$\begin{CD}
SF@>S\varphi >>SG\\
@VS\iota ^FVV  @VVS\iota ^GV\\
S\check{F}@>>S\check{\varphi } >S\check{G}  
\end{CD}$$}
\hspace{0.9cm}
\parbox{3.5cm}{
$$\begin{CD}
K_1(SF)@>K_1(S\varphi )>>K_1(SG)\\
@VK_1(S\iota ^F)VV  @VVK_1(S\iota ^G)V\\
K_1(S\check{F} )@>>K_1(S\check{\varphi } )>K_1(S\hat{G} )
\end{CD}$$}

\hspace{-5.5mm}are obviously commutative (\pr{956} a)). So by a) and c) for \frc (and \cor{928} a), \pr{956} a)),
$$K_1(S\iota ^G)\circ \beta _G\circ K_0(\varphi )=\beta _{\check{G} }\circ K_0(\iota ^G)\circ K_0(\varphi )=
\beta _{\check{G} }\circ K_0(\check{\varphi } )\circ K_0(\iota ^F)=$$
$$=K_1(S\check{\varphi } )\circ \beta _{\check{F} }\circ K_0(\iota ^F)=
K_1(S\check{\varphi } )\circ K_1(S\iota ^F)\circ \beta _F=K_1(S\iota ^G)\circ K_1(S\varphi )\circ \beta _F\;.$$
The assertion follows now from the fact that $K_1(S\iota ^G)$ is injective.\qed

\begin{center}
\section{Higman's linearization trick}
\end{center}

\begin{center}
\fbox{\parbox{7cm}{Throughout this section $F$ denotes a full $E$-C*-algebra, $m,n\in \bn$, and $l:=2^m-1$}}
\end{center}

\begin{de}\label{76}
We shall use the following notation \emph{([R] 11.2)}:
$$Trig(n):=\me{X\in \ccc{C}(\bt,GL_{E_n}(F_n))}{X(z)=\sii{p=-m}{m}a_pz^p,\,a_p\in F_n}\,,$$
$$Pol(n,m):=\me{X\in \ccc{C}(\bt,GL_{E_n}(F_n))}{X(z)=\sii{p=0}{m}a_pz^p,\,a_p\in F_n}\,,$$
$$Pol(n):=\bigcup_{m\in \bn}Pol(n,m),\qquad 
Lin(n):=Pol(n,1)\,,$$
$$Proj(n):=\me{\widetilde{P} }{P\in \pp{F_n}}\;.$$
\end{de}

\begin{lem}\label{78}
\rule{0mm}{0mm}
\begin{enumerate}
\item If $X\in \ccc{C}(\bt,GL_{E_n}(F_n))$ then there are $k\in \bn$ and $Y\in Pol(n)$ such that $z^kX$ is homotopic to $Y$ in $\ccc{C}(\bt,GL_{E_n}(F_n))$.
\item If $P,Q\in \pp{F_n}$ such that $\widetilde{P} $ and $\widetilde{Q} $ are homotopic in $\ccc{C}(\bt,GL_{E_n}(F_n))$ then there are $k,m\in \bn$ such that $z^k\widetilde{P}$ is homotopic to $z^k\widetilde{Q} $ in $Pol(n,l)$.
\end{enumerate}
\end{lem}

a) It is possible to adapt [R] Lemma 11.2.3 to the present situation in order to find a $Z\in Trig(n)$ such that
$$\n{X-Z}<\n{X^{-1}}^{-1}\;.$$
By [R] Proposition 2.1.11, $X$ and $Z$ are homotopic in $\ccc{C}(\bt,GL_{E_n}(F_n))$. There is a $k\in \bn$ such that $Y:=z^kZ\in Pol(n)$. Then $z^kX$ and $Y$ are homotopic in $\ccc{C}(\bt,GL_{E_n}(F_n))$.

b) The proof of [R] Lemma 11.2.4 (ii) works in this case too.\qed 

\begin{de}\label{79}
The map
$$\mad{\{0,1\}^m}{\bnn{l}\cup \{0\}}{j}{\sii{i=1}{m}j_i\,2^{i-1}}$$
is bijective. We denote by
$$\mad{\bnn{l}\cup \{0\}}{\{0,1\}^m}{p}{|p|}$$
its inverse. For every $i\in \bnn{m}$ and $p,q\in \bnn{l}\cup \{0\}$ we put
$$(p,q)_i:=\ad{A_{n+i}}{|p|_i=|q|_i=0}{C_{n+i}^*}{|p|_i=0,\,|q|_i=1}{C_{n+i}}{|p|_i=1,\,|q|_i=0}{B_{n+i}}{|p|_i=|q|_i=1}\;.$$
\end{de}

\begin{lem}\label{80}
\rule{0mm}{0mm}
\begin{enumerate}
\item For $p,q,r,s\in \bnn{l}\cup \{0\}$ and $i\in \bnn{m}$,
$$(p,q)_i(r,s)_i=\ab{0}{|q|_i\not=|r|_i}{(p,s)_i}{|q|_i=|r|_i}\;.$$
In particular
$$\proo{i=1}{m}((p,q)_i(r,s)_i)=\ab{0}{q\not=r}{\proo{i=1}{m}(p,s)_i}{q=r}\;.$$
\item For $p,q\in \bnn{l}\cup \{0\}$ and $i\in \bnn{m}$,
$$A_{n+i}(p,q)_i=\ab{(p,q)_i}{|p|_i=0}{0}{|p|_i=1}\,,$$
$$(p,q)_iA_{n+i}=\ab{(p,q)_i}{|q|_i=0}{0}{|q|_i=1}\;.$$
In particular
$$p\not=0\Longrightarrow \proo{i=1}{m}(A_{n+i}(p,q)_i)=0\,,$$
$$q\not=0\Longrightarrow \proo{i=1}{m}((p,q)_iA_{n+i})=0\,,$$
$$\sii{r=q}{l}\proo{i=1}{m}(A_{n+i}(r,r-q)_i)=\ab{0}{q\not=0}{\proo{i=1}{m}A_{n+i}}{q=0}\;.$$
\item $\sii{p=0}{l}\proo{i=1}{m}(p,p)_i=1_E$.
\end{enumerate}
\end{lem}

a) and b) is a long verification.

c) For every $p\in \bnn{l}\cup \{0\}$ put
$$J_p:=\me{i\in \bnn{m}}{|p|_i=0},\qquad\qquad K_p:=\me{i\in \bnn{m}}{|p|_i=1}\;.$$
Then
$$1_E=\proo{i=1}{m}(A_{n+i}+B_{n+i})=\sii{p=0}{l}\left(\pro{i\in J_p}A_{n+i}\right)\left(\pro{i\in K_p}B_{n+i}\right)=\sii{p=0}{l}\proo{i=1}{m}(p,p)_i\;.\qedd$$

\begin{lem}\label{82}
Let $a\in (F_n)^l$ and
$$X:=\sii{p=1}{l}a_p\sii{q=p}{l}\proo{i=1}{m}(q,q-p)_i\qquad\qquad (X\in F_{m+n})\;.$$
\begin{enumerate}
\item $X^{2^m}=0$.
\item $1_E-X$ is invertible.
\end{enumerate}
\end{lem}

a) We put $D:=\bnn{l}$ and for every $k\in \bn$ and $p\in D^k$,
$$p^{(k)}:=\sii{j=1}{k}p_j,\qquad\qquad a_p^{(k)}:=\proo{j=1}{k}a_{p_j}\;.$$
We want to prove by induction that for every $k\in \bn$,
$$X^k=\si{p\in D^k}a_p^{(k)}\sii{q=p^{(k)}}{l}\proo{i=1}{m}(q,q-p^{(k)})_i\;.$$
The assertion holds for $k=1$. Assume the assertion holds for $k\in \bn$. Then
$$X^{k+1}
=\si{p\in D^k}\si{p'\in D}a_p^{(k)}a_{p'}\sii{q=p^{(k)}}{l}\sii{q'=p'}{l}\proo{i=1}{m}((q,q-p^{(k)})_i(q',q'-p')_i)\;.$$
By \lm{80} a),
$$X^{k+1}=\si{p\in D^k}\si{p'\in D}a_p^{(k)}a_{p'}\sii{q=p^{(k)}+p'}{l}\;\proo{i=1}{m}(q,q-p^{(k)}-p')_i=$$
$$=\si{p\in D^{k+1}}a_p^{(k+1)}\sii{q=p^{(k+1)}}{l}\proo{i=1}{m}(q,q-p^{(k+1)})_i\,,$$
which finishes the inductive proof. Since $p^{(k)}\geq k$ for every $k\in \bn$ we get $X^{2^m}=0$.

b) By a), $1_E+\sii{k=1}{l}X^k $ is the inverse of $1_E-X$.\qed

\begin{p}\label{83}{\bf{(Higman's linearization trick)}}
There is a continuous map
$$\mac{\mu }{Pol(n,l)}{Lin(n+m)}$$
such that $\mu X$ is homotopic to $X\left(\proo{i=1}{m}A_{n+i}\right)+\left(1_E-\proo{i=1}{m}A_{n+i}\right)$ in $Pol(n+m,2l+1)$ for every $X\in Pol(n,l)$. If $X\in Proj(n)$ then the above homotopy takes place in $Lin(n+1)$.
\end{p}

Assume $X\in Pol(n,l)$ is given by 
$$X=\sii{p=0}{l}a_pz^p\,,$$
where $a_p\in F_n$ for every $p\in \bn_l\cup \{0\}$. Put
$$X_p:=\sii{q=p}{l}a_qz^{q-p}\qquad\qquad (\in \ccc{C}(\bt,F_n))$$
for all $p\in \bn_l\cup \{0\}$ and for all $s\in [0,1]$,
$$Y_s:=1_E-s\sii{p=1}{l}X_p\proo{i=1}{m}(0,p)_i\qquad\qquad (\in \ccc{C}(\bt,F_{n+m}))\,,$$
$$Z_s:=1_E+s\sii{q=1}{l}z^q\sii{r=q}{l}\proo{i=1}{m}(r,r-q)_i\qquad\qquad (\in \ccc{C}(\bt,F_{n+m}))\;.$$
By \lm{80} a),
$$Y_s(1_E+s\sii{p=1}{l}X_p\proo{i=1}{m}(0,p)_i)=(1_E+s\sii{p=1}{l}X_p\proo{i=1}{m}(0,p)_i)Y_s=$$
$$=1_E+s^2\sii{p,q=1}{l}X_pX_q\proo{i=1}{m}((0,p)_i(0,q)_i)=1_E\,,$$
so $Y_s$ is invertible. By \lm{82} b), $Z_s$ is also invertible. Thus for every $s\in [0,1]$, $Y_s$ and $Z_s$ are homotopic to $1_E$ in $\ccc{C}(\bt,GL(F_{n+m}))$ and belong therefore to $Pol(n+m,l)$. By \lm{80} c),
$$Z_1=\sii{q=0}{l}z^q\sii{r=q}{l}\proo{i=1}{m}(r,r-q)_i\;.$$

Put
$$\mu X:=1_E-\proo{i=1}{m}A_{n+i}+\sii{p=0}{l}a_p\proo{i=1}{m}(0,p)_i-z\sii{p=1}{l}\proo{i=1}{m}(p,p-1)_i\; \;(\in \ccc{C}(\bt,F_{n+m}))\;.$$
For $z\in \bt$,
$$((\mu X)Z_1)(z)=\sii{p=0}{l}z^p\sii{q=p}{l}\proo{i=1}{m}(q,q-p)_i-\sii{p=0}{l}z^p\sii{q=p}{l}\proo{i=1}{m}(A_{n+i}(q,q-p)_i)+$$
$$+\sii{p,q=0}{l}a_pz^q\sii{r=q}{l}\proo{i=1}{m}((0,p)_i(r,r-q)_i)-\sii{q=0}{l}z^{q+1}\sii{p=1}{l}\sii{r=q}{l}\proo{i=1}{m}((p,p-1)_i(r,r-q)_i)\;.$$
By \lm{80} b),
$$\sii{p=0}{l}z^p\sii{q=p}{l}\proo{i=1}{m}(A_{n+i}(q,q-p)_i)=\proo{i=1}{m}A_{n+i}$$
and by \lm{80} a),
$$\sii{p,q=0}{l}a_pz^q\sii{r=q}{l}\proo{i=1}{m}((0,p)_i(r,r-q)_i)=\sii{q=0}{l}z
^q\sii{p=q}{l}a_p\proo{i=1}{m}(0,p-q)_i=$$
$$=\sii{q=0}{l}z^q\sii{r=0}{l-q}a_{q+r}\proo{i=1}{m}(0,r)_i=\sii{r=0}{l}\sii{q=0}{l-r}z^qa_{q+r}\proo{i=1}{m}(0,r)_i=$$
$$=\sii{r=0}{l}\sii{s=r}{l}z^{s-r}a_s\proo{i=1}{m}(0,r)_i=\sii{r=0}{l}X_r\proo{i=1}{m}(0,r)_i\,,$$
$$\sii{q=0}{l}z^{q+1}\sii{p=1}{l}\sii{r=q}{l}\proo{i=1}{m}((p,p-1)_i(r,r-q)_i)=$$
$$=\sii{q=0}{l}z^{q+1}\sii{p=q+1}{l}\proo{i=1}{m}(p,p-q-1)_i=\sii{q=1}{l}z^q\sii{p=q}{l}\proo{i=1}{m}(p,p-q)_i\;.$$
Thus by \lm{80} c),
$$((\mu X)Z_1)(z)=
\sii{q=0}{l}z^q\sii{p=q}{l}\proo{i=1}{m}(p,p-q)_i-\proo{i=1}{m}A_{n+i}+$$
$$+\sii{r=0}{l}X_r\proo{i=1}{m}(0,r)_i-\sii{q=1}{l}z^q\sii{p=q}{l}\proo{i=1}{m}(p,p-q)_i=$$
$$=\sii{p=0}{l}\proo{i=1}{m}(p,p)_i-\proo{i=1}{m}A_{n+i}+\sii{r=0}{l}X_r\proo{i=1}{m}(0,r)_i=1_E-\proo{i=1}{m}A_{n+i}+\sii{p=0}{l}X_p\proo{i=1}{m}(0,p)_i\;.$$
By \lm{80} a),b), for $z\in \bt$,
$$(Y_1(\mu X)Z_1)(z)=1_E-\proo{i=1}{m}A_{n+i}+\sii{p=0}{l}X_p\proo{i=1}{m}(0,p)_i-\sii{p=1}{l}X_p\proo{i=1}{m}(0,p)_i+$$
$$+\sii{p=1}{l}X_p\proo{i=1}{m}((0,p)_iA_{n+i})-\sii{p=1}{l}\sii{q=0}{l}X_pX_q\proo{i=1}{m}((0,p)_i(0,q)_i)=$$
$$=1_E-\proo{i=1}{m}A_{n+i}+X_0\proo{i=1}{m}(0,0)_i=1_E-\proo{i=1}{m}A_{n+i}+X\proo{i=1}{m}A_{n+i}\;.$$
Since $1_E-\proo{i=1}{m}A_{n+i}+X^{-1}\proo{i=1}{m}A_{n+i}$ is the inverse of $Y_1(\mu X)Z_1$ it follows that $Y_1(\mu X)Z_1$ and $\mu X$ are invertible, i.e. they belong to $\ccc{C}(\bt,GL(F_{n+m}))$. Thus for every $s\in [0,1]$, $Y_s(\mu X)Z_s\in \ccc{C}(\bt,GL(F_{n+m}))$. Let $z\in \bt$ and let
$$\mad{[0,1]}{GL(F_n)}{s}{x_s}$$
be a continuous map with $x_0=X(z)$ and $x_1=1_E$. Since $1_E-\proo{i=1}{m}A_{n+i}+x_s^{-1}\proo{i=1}{m}A_{n+i}$ is the inverse of $1_E-\proo{i=1}{m}A_{n+i}+x_s\proo{i=1}{m}A_{n+i}$ for every $s\in [0,1]$ it follows that the map
$$\mad{[0,1]}{GL(F_{n+m})}{s}{1_E-\proo{i=1}{m}A_{n+i}+x_s\proo{i=1}{m}A_{n+i}}$$
is well-defined and it is a homotopy from $(Y_1(\mu X)Z_1)(z)$ to $1_E$ i.e. $Y_1(\mu X)Z_1\in \ccc{C}(\bt,GL_0(F_{n+m}))$ and $Y_1(\mu X)Z_1\in Pol(n+m,l)$. By the above, for every $s\in [0,1]$, $Y_s(\mu X)Z_s\in \ccc{C}(\bt,GL_0(F_{n+m}))$, so $Y_s(\mu X)Z_s\in Pol(n+m,2l+1)$. Hence $\mu X$ is homotopic to $X\left(\proo{i=1}{m}A_{n+i}\right)+\left(1_E-\proo{i=1}{m}A_{n+i}\right)$ in $Pol(n+m,2l+1)$ and $\mu X\in Lin(n+m)$.

In order to prove the last assertion remark that there is a $P\in \pp{F_n}$ with 
$X=\widetilde{P}=(1_E-P)+zP$. Then $m=l=1$, $a_0=1_E-P$, $a_1=P$,  $X_1=a_0=P$,
$$\mu X=1_E-PA_{n+1}+PC_{n+1}^*-zC_{n+1}\,,$$
and for every $s\in [0,1]$,
$$Y_s=1_E-sPC_{n+1}^*\,,\qquad
 Z_s:=1_E+szC_{n+1}\,,\qquad Y_s(\mu X)Z_s\in Lin(n+1)\;.$$
Thus $\mu X$ is homotopic to $Y_1(\mu X)Z_1$ in $Lin(n+1)$.\qed

\begin{center}
\section{The periodicity}
\end{center}

\begin{center}
\fbox{\parbox{7cm}{Throughout this section $F$ denotes a full $E$-C*-algebra, $m,n\in \bn$, and $l:=2^m-1$}}
\end{center}

\begin{lem}\label{77}
If $X\in \ccc{C}(\bt,GL(F_n))$ and $X(1)\in GL_{E_n}(F_n)$ then
 $$X\in \ccc{C}(\bt,GL_{E_n}(F_n))\;.$$
\end{lem}

Let $\theta \in [0,2\pi [$ and for every $s\in [0,1]$ put
$$\mae{Y_s}{\bt}{GL(F_n)}{z}{X(e^{-is}z)}\;.$$
Then $Y_0(e^{i\theta })=X(e^{i\theta })$ and $Y_\theta (e^{i\theta })=X(1)$ so $X(e^{i\theta })$ is homotopic to $X(1)$ in $GL(F_n)$. Thus $X(e^{i\theta })\in GL_{E_n}(F_n)$ and $X\in \ccc{C}(\bt,GL_{E_n}(F_n))$.\qed

\begin{p}\label{86}
The following are equivalent for every $X\in F_n$.
\begin{enumerate}
\item $\widetilde{X}\in Lin(n) $.
\item $z\in \bt\setminus \{1\}\Longrightarrow \widetilde{X}(z)\in GL\,(F_n) $.
\item $\widetilde{X} $ is a generalized idempotent of $F_n$ \emph{([R] Definition 11.2.8)}.
\end{enumerate}
\end{p}

$a\Rightarrow b$ is trivial.

$b\Rightarrow a$. By \lm{77}, since $\widetilde{X}(1)=1_E $, $\tilde{X}\in \ccc{C}(\bt,GL_{E_n}(F_n)) $  so $\tilde{X}\in Lin(n) $.

$b\Leftrightarrow c$. For $z\in \bt\setminus \{1\}$,
$$\widetilde{X}(z)=(z-1)X+1_E=(z-1)\left(X-\frac{1}{1-z}1_E\right)\;. $$
Since
$$\me{\frac{1}{1-z}}{z\in \bt\setminus \{1\}}=\me{\alpha \in \bc}{real(\alpha) =\frac{1}{2}}\,,$$
b) holds iff $X-\alpha 1_E$ is invertible for every $\alpha \in \bc$ with $real(\alpha )=\frac{1}{2}$, which is equivalent to c).\qed

\begin{lem}\label{87}
For $z\in \bt$,
$$zA_n+B_n\sim _hA_n+zB_n\quad\emph{in}\quad\unn{E_n}\;.$$
\end{lem}

We have 
$$(C_n+C_n^*)(zA_n+B_n)(C_n+C_n^*)=(zC_n+C_n^*)(C_n+C_n^*)=zB_n+A_n$$
and the assertion follows from \pr{931} a).\qed

\begin{lem}\label{88}
For $z\in \bt$,
$$z^l\proo{i=1}{m}A_{n+i}+\sii{p=1}{l}\proo{i=1}{m}(p,p)_i\sim _h\proo{i=1}{m}A_{n+i}+z\sii{p=1}{l}\proo{i=1}{m}(p,p)_i\quad \emph{in} \quad\unn{E_{n+m}}\;.$$
\end{lem}

Let $k\in \bn_l$ and let $j\in \bn_m$ with $|k|_j=1$.  By \lm{87},

$$z^{l-k+1}\prom A_{n+i}+z\sii{p=1}{k-1}\prom (p,p)_i+\sii{p=k}{l}\prom (p,p)_i=$$
$$=\left(z^{l-k}\prom A_{n+i}+\prom (k,k)_i\right)(zA_{n+j}+(k,k)_j)+$$
$$+z\sii{p=1}{k-1}\prom (p,p)_i+\sii{p=k+1}{l}\prom (p,p)_i\sim _h$$
$$\sim _h\left(z^{l-k}\prom A_{n+i}+\prom (k,k)_i\right)(A_{n+j}+
z(k,k)_j)+$$
$$+z\sii{p=1}{k-1}\prom (p,p)_i+\sii{p=k+1}{l}\prom (p,p)_i=$$
$$=z^{l-k}\prom A_{n+i}+z\sii{p=1}{k}\prom (p,p)_i+\sii{p=k+1}{l}\prom (p,p)_i$$
in $\unn{E_{n+m}}$. The assertion follows now by induction on $k\in \bn_l$.\qed

\begin{lem}\label{d}
Let $P,Q\in \pp{F_n}$.
\begin{enumerate}
\item For every $z\in \bt$,
$$\overbrace{P\left(\prom A_{n+i}\right)+\left(1_E-\prom A_{n+i}\right)}^{\widetilde{} }\;(z)=$$
$$=\widetilde{P}(z)\left(\prom A_{n+i}\right)+z\left(1_E-\prom A_{n+i}\right)\;.$$
\item If (with the identification of \emph{\lm{979} d)})
$$\widetilde{P}\left(\prom A_{n+i}\right)+\left(1_E-\prom A_{n+i}\right)\sim _h$$
$$\sim _h\widetilde{Q}\left(\prom A_{n+i}\right)+\left(1_E-\prom A_{n+i}\right)\quad\emph{in}\quad \unn{\left(\check{\overbrace{SF}}\right)\!_{n+m}}\,,$$
then
$$\overbrace{P\left(\prom A_{n+i}\right)+\left(1_E-\prom A_{n+i}\right)}^{\widetilde{} }\sim _h$$
$$\sim _h\overbrace{Q\left(\prom A_{n+i}\right)+\left(1_E-\prom A_{n+i}\right)}^{\widetilde{} }\quad\emph{in}\quad \unn{\left(\check{\overbrace{SF}}\right)\!_{n+m}}\;.$$
\end{enumerate}
\end{lem}

a) We have
$$\overbrace{P\left(\prom A_{n+i}\right)+\left(1_E-\prom A_{n+i}\right)}^{\widetilde{} }\;(z)=$$
$$=zP\left(\prom A_{n+i}\right)+z\left(1_E-\prom A_{n+i}\right)+\prom A_{n+i}+$$
$$+\left(1_E-\prom A_{n+i}\right)-
P\left(\prom A_{n+i}\right)-\left(1_E-\prom A_{n+i}\right)=$$
$$=\widetilde{P}(z)\left(\prom A_{n+i}\right)+z\left(1_E-\prom A_{n+i}\right)\;.$$

b) Let 
$$\mad{[0,1]}{\unn{\left(\check{\overbrace{SF}}\right)\!_{n+m}}}{s}{U_s}$$
be a continuous map with
$$U_0=\widetilde{P}\left(\prom A_{n+i}\right)+\left(1_E-\prom A_{n+i}\right)\,,$$
$$U_1=\widetilde{Q}\left(\prom A_{n+i}\right)+\left(1_E-\prom A_{n+i}\right)\;.$$
Put $U'_s:=U_s\left(\prom A_{n+i}\right)+z\left(1_E-\prom A_{n+i}\right)$ for every $s\in [0,1]$. Then $s\mapsto U'_s$ is a continuous path in $\unn{\left(\check{\overbrace{SF}}\right)\!_{n+m}}$ and by a),
$$U'_0=U_0\left(\prom A_{n+i}\right)+z\left(1_E-\prom A_{n+i}\right)=$$
$$=\widetilde{P}\left(\prom A_{n+i}\right)+z\left(1_E-\prom A_{n+i}\right)=$$
$$=
\overbrace{P\left(\prom A_{n+i}\right)+\left(1_E-\prom A_{n+i}\right)}^{\widetilde{} }\;(z)\,,$$
$$U'_1=\overbrace{Q\left(\prom A_{n+i}\right)+\left(1_E-\prom A_{n+i}\right)}^{\widetilde{} }\;(z)\;.                      \qedd$$

\begin{p}\label{89}
\rule{0mm}{0mm}
\begin{enumerate}
\item If $U\in \unn{\left(\check{\overbrace{SF}}\right)\!_n}$ then there are $k,m\in \bn$ and $P\in \pp{F_{n+m}}$ such that (with the identification of \emph{\lm{979} d)}) 
$$(z^kU)\left(\prom A_{n+i}\right)+\left(1_E-\prom A_{n+i}\right)\sim _h\widetilde{P}\quad \emph{in}\quad  \unn{\left(\check{\overbrace{SF}}\right)\!_{n+m}}\;.$$
\item Let $P,Q\in \pp{F_n}$ with $\widetilde{P}\sim _h\widetilde{Q}  $ in $\unn{\left(\check{\overbrace{SF}}\right)\!_n}$. Then there is an $m\in \bn$ such that
$$P\left(\prom A_{n+i}\right)+\left(1_E-\prom A_{n+i}\right)\sim _h$$
$$\sim _hQ\left(\prom A_{n+i}\right)+\left(1_E-\prom A_{n+i}\right)\quad\emph{in}\quad \pp{F_{n+m}}\;.$$
\end{enumerate}
\end{p}

a) By \pr{78} a), there are $k,m\in \bn$, $k<2^m$, and $X\in Pol(n,l)$ such that $z^kU$ is homotopic to $X$ in $\ccc{C}(\bt,GL_E(F_n))$. By \pr{83}, there is a $Y\in Lin(n+m)$ with
$$X\left(\prom A_{n+i}\right)+\left(1_E-\prom A_{n+i}\right)\sim _h Y\quad \mbox{in}\quad Pol(n+m,2l+1)\;.$$
By [R] Lemma 11.2.12 (i), there is a $P\in \pp{F_{n+m}}$ with $Y\sim _h\widetilde{P} $ in $Lin(n+m)$. Thus
$$(z^kU)\left(\prom A_{n+i}\right)+\left(1_E-\prom A_{n+i}\right)\sim _h$$
$$\sim _hX\left(\prom A_{n+i}\right)+\left(1_E-\prom A_{n+i}\right)\sim _hY\sim _h\widetilde{P} $$
in $\ccc{C}(\bt,GL_E(F_{n+m}))$. By [R] Proposition 2.1.8 (iii) and the identification of \lm{979} d),
$$(z^kU)\left(\prom A_{n+i}\right)+\left(1_E-\prom A_{n+i}\right)\sim _h\widetilde{P}\quad\mbox{in}\quad  \unn{\left(\check{\overbrace{SF}}\right)\!_{n+m}}\;.$$

b) By \pr{78} b), there are $k,m\in \bn$, $k<2^m$, such that $z^k\widetilde{P}\sim _hz^k\widetilde{Q}  $ in $Pol(n,l)$. By \lm{88} and \lm{80} c),
$$z^l\left(\prom A_{n+i}\right)+\left(1_E-\prom A_{n+i}\right)\sim _h\left(\prom A_{n+i}\right)+z\left(1_E-\prom A_{n+i}\right)$$
in $\unn{E_{n+m}}$. By \lm{d} a),
$$\overbrace{P\left(\prom A_{n+i}\right)+\left(1_E-\prom A_{n+i}\right)}^{\widetilde{} }\;(z)=$$
$$=\left(\left(\prom A_{n+i}\right)+z\left(1_E-\prom A_{n+i}\right)\right)\times $$
$$\times \left(\widetilde{P}(z)\left(\prom A_{n+i}\right)+\left(1_E-\prom A_{n+i}\right)\right)\sim _h$$
$$\sim _h\left(z^l\left(\prom A_{n+i}\right)+\left(1_E-\prom A_{n+i}\right)\right)\left(\widetilde{P}(z)+\left(1_E-\prom A_{n+i}\right) \right)=$$
$$=z^l\widetilde{P}(z)\left(\prom A_{n+i}\right)+ \left(1_E-\prom A_{n+i}\right)\sim _h$$
$$\sim _hz^l\widetilde{Q}(z)\left(\prom A_{n+i}\right)+ \left(1_E-\prom A_{n+i}\right)\sim _h$$
$$\sim _h\overbrace{Q\left(\prom A_{n+i}\right)+\left(1_E-\prom A_{n+i}\right)}^{\widetilde{} }\;(z)$$
in $Pol(n+m,l)$. By \pr{83},
$$\widetilde{P}\left(\prom A_{n+i}\right)+\left(1_E-\prom A_{n+i}\right)=$$
$$=\overbrace{P\left(\prom A_{n+i}\right)+\left(1_E-\prom A_{n+i}\right)}^{\widetilde{} }\left(\prom A_{n+i}\right)+\left(1_E-\prom A_{n+i}\right)\sim _h$$
$$\sim _h\mu \left(\overbrace{P\left(\prom A_{n+i}\right)+\left(1_E-\prom A_{n+i}\right)}^{\widetilde{} }\right)\sim _h$$
$$\sim _h\mu \left(\overbrace{Q\left(\prom A_{n+i}\right)+\left(1_E-\prom A_{n+i}\right)}^{\widetilde{} }\right)\sim _h$$
$$\sim _h\widetilde{Q}\left(\prom A_{n+i}\right)+\left(1_E-\prom A_{n+i}\right)$$
in $Lin(n+m)$. By \lm{d} a),
$$\overbrace{P\left(\prom A_{n+i}\right)+\left(1_E-\prom A_{n+i}\right)}^{\widetilde{} }=\widetilde{P}\left(\prom A_{n+i}\right)+z\left(1_E-\prom A_{n+i}\right)\sim _h $$
$$\sim _h\widetilde{Q}\left(\prom A_{n+i}\right)+z\left(1_E-\prom A_{n+i}\right)=
\overbrace{Q\left(\prom A_{n+i}\right)+\left(1_E-\prom A_{n+i}\right)}^{\widetilde{} }$$
in $Lin(n+m)$. The assertion follows now from [R] Lemma 11.2.12 (ii).\qed

\begin{theo}\label{91}
The Bott map is bijective.
\end{theo}

\begin{center}
Step 1 Surjectivity
\end{center}

Let $a\in K_1(SF)$. There are $n\in \bn$ and $U\in \unn{\left(\check{\overbrace{SF}}\right)\!_n}$ with $a=[U]_1$. By \pr{89} a), there are $m,p\in \bn$, $p\geq n$, and $P\in \pp{F_{p+m}}$ such that
$$(z^lU)\left(\prom A_{p+i}\right)+\left(1_E-\prom A_{p+i}\right)\sim _h\widetilde{P}\quad\mbox{in}\quad \unn{\left(\check{\overbrace{SF}}\right)\!_{p+m}}\;. $$
By \lm{88} and \lm{80} c),
$$\overbrace{1_E-\prom A_{p+i}}^{\widetilde{} }=z\left(1_E-\prom A_{p+i}\right)+\left(\prom A_{p+i}\right)\sim _h$$
$$\sim _h\left(1_E-\prom A_{p+i}\right)+z^l\left(\prom A_{p+i}\right)\quad\mbox{in}\quad \unn{E_{p+m}}$$
so by \pr{951a} and \pr{984},
$$\beta _F\left([P]_0-\left[1_E-\prom A_{p+i}\right]_0\right)=[\widetilde{P} ]_1-\left[\overbrace{1_E-\prom A_{p+i}}^{\widetilde{} }\right]_1=$$
$$=\left[(z^lU)\left(\prom A_{p+i}\right)+\left(1_E-\prom A_{p+i}\right)\right]_1-$$
$$-\left[\left(1_E-\prom A_{p+i}\right)+z^l\left(\prom A_{p+i}\right)\right]_1=$$
$$=\left[\left((z^lU)\left(\prom A_{p+i}\right)+\left(1_E-\prom A_{p+i}\right)\right)\times \right.$$
$$\left.\times  \left(\left(1_E-\prom A_{p+i}\right)+z^l\left(\prom A_{p+i}\right)\right)^*\right]_1=$$
$$=\left[U\left(\prom A_{p+i}\right)+\left(1_E-\prom A_{p+i}\right)\right]_1=[U]_1=a\;.$$

\begin{center}
Step 2 Injectivity
\end{center}

Let $a\in K_0(F)$ with $\beta _Fa=0$. By \pr{916} d), there are $P,Q\in \pp{F_n}$, $PQ=0$, such that $a=[P]_0-[Q]_0$. Then $[\widetilde{P} ]_1=[\widetilde{Q} ]_1$, so $U:=\tilde{P}\tilde{Q}^*\in un_{E_n}\check{\overbrace{SF}}  $. Then
$$U=((z-1)P+1_E)((\bar{z}-1)Q+1_E ))=(z-1)P+(\bar{z}-1 )Q+1_E\,,\qquad U(1)=1_E\;.$$
By \pr{951a}, there is an $m\in \bn$ such that
$$V:=U\left(\prom A_{n+i}\right)+\left(1_E-\prom A_{n+i}\right)=\tau _{n+m,n}^FU\in \unhh{\left(\check{\overbrace{SF}}\right)\!_{n+m}}{n+m}\;.$$
Then there is a $W\in \unn{E_{n+m}}$ with $V\sim _hW$ in $\unn{\left(\check{\overbrace{SF}}\right)\!_{n+m}}$. By the above,
$$W=W(1)\sim _hV(1)=1_E\,,\qquad V\sim _h1_E \quad\mbox{in}\quad \unn{\left(\check{\overbrace{SF}}\right)\!_{n+m}}\;.$$ 
By \pr{951a},
$$\widetilde{P}\left(\prom A_{n+i}\right)+\left(1_E-\prom A_{n+i}\right)=\tau _{n+m,n}^F\tilde{P}=(\tau _{n+m,n}^FU)(\tau _{n+m,n}^F\tilde{Q})= $$
$$=V(\tau _{n+m,n}^F\tilde{Q})\sim _h\widetilde{Q}\left(\prom A_{n+i}\right)+\left(1_E-\prom A_{n+i}\right)\quad\mbox{in}\quad \unn{\left(\check{\overbrace{SF}}\right)\!_{n+m}}\,,$$
so by \pr{d} b),
$$\overbrace{P\left(\prom A_{n+i}\right)+\left(1_E-\prom A_{n+i}\right)}^{\widetilde{} }\sim _h$$
$$\sim _h\;\overbrace{Q\left(\prom A_{n+i}\right)+\left(1_E-\prom A_{n+i}\right)}^{\widetilde{}}\quad\mbox{in}\quad\unn{\left(\check{\overbrace{SF}}\right)\!_{n+m}}\;.$$
Put
$$P':=P\left(\prom A_{n+i}\right)+\left(1_E-\prom A_{n+i}\right)\,,$$
$$ Q':=Q\left(\prom A_{n+i}\right)+\left(1_E-\prom A_{n+i}\right)\;.$$
By \pr{89} b), there are $m',p'\in \bn$ such that
$$P'\left(\proo{j=1}{m'} A_{p'+j}\right)+\left(1_E-\proo{j=1}{m'} A_{p'+i}\right)\sim _h$$
$$\sim _h \;Q'\left(\proo{j=1}{m'} A_{p'+i}\right)+\left(1_E-\proo{j=1}{m'} A_{p'+i}\right)\quad\mbox{in}\quad\pp{F_{p'+m'}}\;.$$
It follows successively  
$$\left[P'\proo{j=1}{m'}A_{p'+j}\right]_0=\left[Q'\proo{j=1}{m'}A_{p'+j}\right]_0\,,$$
$$\left[P\left(\prom A_{n+i}\right)\left(\proo{j=1}{m'}A_{p'+j}\right)\right]_0=\left[Q\left(\prom A_{n+i}\right)\left(\proo{j=1}{m'}A_{p'+j}\right)\right]_0\,,$$
$$[P]_0=[Q]_0\,,\qquad\qquad a=[P]_0-[Q]_0=0\;.\qedd$$

{\it Remark.} By \h{91} and \pr{985} c), the functor $K_0$ is determined by the functor $K_1$.

\begin{co}[The six-term sequence]\label{93}
Let
$$\oc{F}{\varphi }{G}{\psi }{H}$$
be an exact sequence in \frm.
\begin{enumerate}
\item The sequence
$$\oc{SF}{S\varphi }{SG}{S\psi }{SH}$$
is exact. Let
$$\mac{\delta _2}{K_1(SH)}{K_0(SF)}$$
be its associated index map \emph{(\cor{964})} and put \emph{(\pr{985}, \h{976})}
$$\mac{\delta _0:=\theta _F^{-1}\circ \delta _2\circ \beta _H}{K_0(H)}{K_1(F)}\;.$$
We call $\delta _0$ and $\delta _1$ {\bf{the six-term index maps}}. If we denote by $\bar{\delta }_0 $ the corresponding six-term index map associated to the exact sequence in \frm (with obvious notation)
$$\oc{SF}{\varphi }{CF}{\psi }{F}$$
then $\bar{\delta }_0=\beta _F $.
\item The {\bf{six-term sequence}}
$$\begin{CD}
K_0(F)@>K_0(\varphi )>>K_0(G)@>K_0(\psi )>>K_0(H)\\
@A\delta _1AA           @.                   @VV\delta _0V\\
K_1(H)@<<K_1(\psi )<K_1(G)@<<K_1(\varphi )<K_1(F)
\end{CD}$$
is exact.
\item If $F$ (resp. $H$) is K-null (e.g. homotopic to $\{0\}$) then $\oa{K_i(G)}{K_i(\psi )}{K_i(H)}$ (resp. $\oa{K_i(F)}{K_i(\varphi) }{K_i(G)}$) is a group isomorphism for every $i\in \{0,1\}$.
\item If $G$ is K-null (e.g. homotopic to $\{0\}$) then
$$\oa{K_0(H)}{\delta _0}{K_1(F)}\,,\qquad\qquad \oa{K_1(H)}{\delta _1}{K_0(F)}$$
are group isomorphisms.
\item If $\varphi $ is K-null (e.g. factorizes through null) then the sequences
$$\oc{K_0(G)}{K_0(\psi )}{K_0(H)}{\delta _0}{K_1(F)}\,,$$
$$\oc{K_1(G)}{K_1(\psi )}{K_1(H)}{\delta _1}{K_0(F)}$$
are exact.
\item If $\psi $ is K-null (e.g. factorizes through null) then the sequences
$$\oc{K_0(H)}{\delta _0}{K_1(F)}{K_1(\varphi )}{K_1(G)}\,,$$
$$\oc{K_1(H)}{\delta _1}{K_0(F)}{K_0(\varphi )}{K_0(G)}$$
are exact.
\item The six-term index maps of a split exact sequence are equal to $0$.
\end{enumerate}
\end{co}

a) is easy to see.

b) By \h{91}, $\beta _H$ is an isomorphism. By \h{974}, the sequences
$$K_1(F)\stackrel{\scriptscriptstyle K_1(\varphi )}{\longrightarrow }K_1(G)\stackrel{\scriptscriptstyle K_1(\psi )}{\longrightarrow }K_1(H)\stackrel{\delta _1}{\longrightarrow }K_0(F)\stackrel{\scriptscriptstyle K_0(\varphi )}{\longrightarrow }K_0(G)\stackrel{\scriptscriptstyle K_0(\psi )}{\longrightarrow }K_0(H)\,,$$
$$K_1(SG)\stackrel{\scriptscriptstyle K_1(S\psi  )}{\longrightarrow }K_1(SH)\stackrel{\scriptscriptstyle \delta _2}{\longrightarrow }K_0(SF)\stackrel{K_0(S\varphi )}{\longrightarrow }K_0(SG)$$
are exact. By \pr{985} c) and \pr{978}, the diagrams

\parbox{3cm}{
$$\begin{CD}
K_0(G)@>K_0(\psi )>>K_0(H)\\
@V\beta _GVV         @VV\beta _HV\\
K_1(SG)@>>K_1(S\psi )>K_1(SH)
\end{CD}$$}
\hspace{4cm}
\parbox{3cm}{
$$\begin{CD}
K_1(F)@>K_1(\varphi )>>K_1(G)\\
@V\theta _FVV         @VV\theta _GV\\
K_0(SF)@>>K_0(S\varphi )>    K_0(SG)
\end{CD}$$}

are commutative. It follows
$$\delta _0\circ \kk{0}{\psi }=\theta _F^{-1}\circ \delta _2\circ \beta _H\circ \kk{0}{\psi }=\theta _F^{-1}\circ \delta _2\circ \kk{1}{S\psi }\circ \beta _G=0\,,$$
$Im\,\kk{0}{\psi }\subset Ker\,\delta _0$. Let $a\in Ker\,\delta _0$. Then $\delta _2\beta _Ha=\theta _F\delta _0a=0$, so there is a $b\in \kk{1}{SG}$ with $\kk{1}{S\psi }b=\beta _Ha$. It follows
$$a=\beta _H^{-1}\kk{1}{S\psi }b=\kk{0}{\psi }\beta _G^{-1}b\in Im\,\kk{0}{\psi }\,,\qquad Ker\,\delta _0\subset Im\,\kk{0}{\psi }\;.$$

c) The assertion follows immediately from b). By \pr{956} e), a null-homotopic \eo is K-null.

d) The proof is similar to the proof of c).

e) and f) follow from b) and \pr{956} f).

g) By \pr{936} and \cor{987} (with the notation of b)) $\kk{0}{\varphi }$ and $\kk{1}{\varphi }$ are injective and $\kk{0}{\psi }$ and $\kk{1}{\psi }$ are surjective and the assertion follows from b).\qed 

\begin{co}\label{29.6}
Let us consider the following commutative diagram in \frm
$$\begin{CD}
0@>>>F@>\varphi >>G@>\psi >>H@>>>0\\
@.      @V\gamma VV  @V\alpha VV @VV\beta V@.\\
0@>>>F'@>>\varphi' >G'@>>\psi' >H'@>>>0\,,
\end{CD}$$
where the horizontal lines are exact.
\begin{enumerate}
\item {\bf{(Commutativity of the six-term index maps)}} The diagrams (with obvious notation)

\parbox{3cm}{
$$\begin{CD}
K_1(H)@>\delta _1>>K_0(F)\\
@VK_1(\beta )VV    @VVK_0(\gamma )V\\
K_1(H')@>>\delta '_1>K_0(F')
\end{CD}$$}
\hspace{4cm}
\parbox{3cm}{
$$\begin{CD}
K_0(H)@>\delta _0>>K_1(F)\\
@VK_0(\beta )VV         @VVK_1(\gamma )V\\
K_0(H')@>>\delta '_0>    K_1(F')
\end{CD}$$}

are commutative. If $\kk{i}{F}=\kk{i}{F'}$, $\kk{i}{H}=\kk{i}{H'}$, and $\kk{i}{\beta }$ and $\kk{i}{\gamma }$ are the identity maps for all $i\in \z{0,1}$ then $\delta _i=\delta '_i$ for all $i\in \z{0,1}$.
\item The diagram (with obvious notation)
$$\begin{CD}
K_0(F)@=K_0(F)@>K_0(\varphi )>>K_0(G)@>K_0(\psi )>>K_0(H)@=K_0(H)\\
@V=VV@VK_0(\gamma )VV         @VK_0(\alpha )VV  @VVK_0(\beta )V@VV=V\\
K_0(F)@>K_0(\gamma )>>    K_0(F')@>K_0(\varphi ')>>K_0(G')@>K_0(\psi ')>>K_0(H')@<K_0(\beta )<<K_0(H)\\
@A\delta _1AA @A\delta'_1AA@.  @VV\delta '_0V @VV\delta _0V\\
K_1(H)@>>K_1(\beta )>K_1(H')@<<K_1(\psi ')<K_1(G')@<<K_1(\varphi ')<K_1(F')@<<K_1(\gamma )<K_1(F)\\
@A=AA @AK_1(\beta )AA @AK_1(\alpha )AA @AAK_1(\gamma )A @AA=A\\
K_1(H)@=K_1(H)@<<K_1(\psi )<K_1(G)@<<K_1(\varphi )<K_1(F)@=K_1(F)
\end{CD}$$
is commutative.
\end{enumerate}
\end{co}

a) The commutativity of the first diagram was proved in \pr{966}. By \pr{978}, the diagram
$$\begin{CD}
K_1(F)@>K_1(\gamma )>>K_1(F')\\
@V\theta _FVV    @VV\theta _{F'}V\\
K_0(SF)@>>K_0(S\gamma )>K_0(SF')
\end{CD}$$
is commutative. By \pr{966}, the diagram
$$\begin{CD}
K_1(SH)@>\delta _2>>K_0(SF)\\
@VK_1(S\beta )VV    @VVK_0(S\gamma )V\\
K_1(SH')@>>\delta '_2>K_0(SF')
\end{CD}$$ 
is commutative, where $\delta _2$ and $\delta '_2$ are defined in \cor{93} a). By \pr{985} c), the diagram
$$\begin{CD}
K_0(H)@>K_0(\beta )>>K_0(H')\\
@V\beta _HVV    @VV\beta _{H'}V\\
K_1(SH)@>>K_1(S\beta )>K_1(SH')
\end{CD}$$
is commutative. It follows, by the definition of $\delta _0$ (\cor{93} a)),
$$K_1(\gamma )\circ \delta _0=K_1(\gamma) \circ \theta _F^{-1}\circ \delta _2\circ \beta _H=\theta _{F'}^{-1}\circ K_0(S\gamma )\circ \delta _2\circ \beta _H=$$
$$=\theta _{F'}^{-1}\circ \delta '_2\circ K_1(S\beta )\circ \beta _H=\theta _{F'}^{-1}\circ \delta '_2\circ \beta _{H'}\circ K_0(\beta )=\delta '_0\circ K_0(\beta )\;.$$

b) follows from a) and \cor{93} b).\qed

\begin{center}
\chapter{Variation of the parameters}

\fbox{\parbox{8.3cm}{Throughout this chapter we endow $\{0,1\}$ with the structure of o group by identifying it with $\bzz{2}$.}}

\section{Changing $E$}
\end{center}

Let $E'$ be a commutative unital C*-algebra, $\mac{\phi }{E}{E'}$ a unital C*-homomorphism, and
$$\mae{f'}{T\times T}{\unn{E'}}{(s,t)}{\phi f(s,t)}\;.$$
Then $f'\in \f{T}{E'}$ and we may define $E_n'$ with respect to $f'$ for every \bbn\, like in \dd{a}.

Let \bbn\, and put
$$C'_n:=\si{t\in T_n}((\phi C_{n,t})\otimes id_K)V_t^{f'}\quad (\in E'_n)\;.$$
For every $s\in T_{n-1}$,
$$\si{t\in T_n}((f(s^{-1}t,t)C_{n,ts^{-1}})\otimes id_K)V_t^f=V_s^fC_n=$$
$$=C_nV_s^f=\si{t\in T_n}((f(ts^{-1},s)C_{n,ts^{-1}})\otimes id_K)V_t^f$$
so by [C2] Theorem 2.1.9 a),
$$f(s^{-1}t,t)C_{n,s^{-1}t}=f(ts^{-1},s)C_{n,ts^{-1}}$$
for every $t\in T_n$. It follows
$$f'(s^{-1}t,t)C'_{n,s^{-1}t}=f'(ts^{-1},s)C'_{n,ts^{-1}}\,,\quad V_s^{f'}C_n'=C_n'V_s^{f'}\,,\quad C'_n\in (E'_{n-1})^c\;.$$
Thus $(C_n')_{\bbn}$ satisfies the conditions of \axi{b} and we may construct a $K$-theory with respect to $T,\,E',\,f'$, and $(C'_n)_{\bbn}$, which we shall denote by $K'$.

Let $F$ be an $E'$-C*-algebra. We denote by $\bar{F} $ or by $\Phi (F)$ the \eo obtained by endowing the C*-algebra $F$ with the exterior multiplication
$$\mad{E\times F}{F}{(\alpha ,x)}{(\phi \alpha )x}\;.$$
If $\oa{F}{\varphi  }{G}$ is a morphism in $\fr{M}_{E'}$, then $\oa{\bar{F} }{\bar{\varphi } }{\bar{G} }$ is a morphism in \frm, in a natural way.

Let $F$ be an $E'$-C*-algebra and \bbn. We put for every 
$$X=\si{t\in T_n}((\alpha _t,x_t)\otimes id_K)V_t^f\in \check{\bar{F} }_n\,, $$
$$X':=\si{t\in T_n}((\phi \alpha _t,x_t)\otimes id_K)V_t^{f'}\quad (\in \check{F}_n)$$
and set 
$$\mae{\phi _{F,n}}{\check{\bar{F} }_n }{\check{F}_n }{X}{X'}\;.$$
Then $\phi _{F,n}$ is a unital C*-homomorphism (surjective or injective if $\phi $ is so ([C2] Theorem 2.1.9 a))) such that $\phi _{F,n}(Un_{E_n}\,\check{\bar{F} }_n )\subset Un_{E_n'}\,\check{F}_n $ and $\phi _{F,n}\circ \sigma _n^{\bar{F} }=\sigma _n^F\circ \phi _{F,n}$. Thus we get for every $i\in \{0,1\}$ an associated group homomorphism $\mac{\Phi _{i,F}}{K_i(\bar{F} )}{K_i'(F)}$.

Let $E''$ be a unital commutative C*-algebra, $\mac{\phi '}{E'}{E''}$ a unital C*-homomorphism, and $\phi '':=\phi '\circ \phi $. Then we may do similar constructions for $\phi '$ and $\phi ''$ as we have done for $\phi $. If $F$ is an $E''$-C*-algebra, $\Phi '(F)$ and $\Phi ''(F)$ the corresponding $E'$-C*-algebra and $E$-C*-algebra, respectively, then $\Phi ''(F)=\Phi (\Phi '(F))$. If $\Phi '_i$ and $\Phi ''_i$ are the equivalents of $\Phi _i$ with respect to $\phi '$ and $\phi ''$, respectively, then $\Phi ''_{i,F}= \Phi '_{i,F}\circ \Phi _{i,\Phi '(F)}$ for every $i\in \{0,1\}$. If $E''=E$ and $\phi ''=id_E$ then $C_n''=C_n$ for every \bbn\, and for every \eo $F$, $\Phi ''(F)=F$ and $\Phi'' _{i,F}=id_{K_i(F)}$ for every $i\in \{0,1\}$. If in addition $\phi ''':=\phi \circ \phi '=id_{E'}$ then $C'''_n=C'_n$ for every \bbn\, and for every $E'$-C*-algebra $F$, $\Phi '(\Phi (F))=F$ and $\Phi '_{i,\Phi (F)}\circ \Phi _{i,F}=id_{K'_i(F)}$ for every $i\in \{0,1\}$, i.e. the $K$-theory and the $K'$-theory "coincide".

{\it Remark.} Let $P\in \pp{E}$, $0<P<1_E$, and put
$$\mae{Pf}{T\times T}{\unn{PE}}{(s,t)}{Pf(s,t)}\;.$$
Then $Pf\in \f{T}{PE}$ and we denote by $PK$ the K-theory with respect to $T,\,PE,\,Pf,$ and $(PC_n)_{\bbn}$. Then for every $E$-C*-algebra $F$ and $i\in \{0,1\}$
$$K_i(F)\approx ((PK)_i(PF))\times (((1_E-P)K)_i((1_E-P)F))\;.$$
If $\oa{F}{\varphi }{G}$ is a morphism in \frm then
$$\mae{P\varphi }{PF}{PG}{Px}{P\varphi x}$$
is a morphism in $\fr{M}_{PE}$ and
$$K_i(\varphi )=(PK)_i(P\varphi )\times ((1_E-P)K)_i((1_E-P)\varphi )$$
for every $i\in \{0,1\}$.

\begin{p}\label{26.4}
We use the above notation and assume $i\in \{0,1\}$.
\begin{enumerate}
\item If $\oa{F}{\varphi }{G}$ is a morphism in $\fr{M}_{E'}$ then the diagram
$$\begin{CD}
K_i(\bar{F})@>K_i(\bar{\varphi } )>>K_i(\bar{G})\\
@V\Phi _{i,F}VV         @VV\Phi _{i,G}V\\
K_i'(F)@>>K'_i(\varphi )>    K'_i(G)
\end{CD}$$
is commutative.
\item For every $E'$-C*-algebra $F$ the diagram
$$\begin{CD}
K_0(\bar{F})@>\beta_{\bar{F} } >>K_1(\overline{SF})\\
@V\Phi _{0,F}VV         @VV\Phi _{1,SF}V\\
K'_0(F)@>>\beta' _F>    K'_1(SF)\,,
\end{CD}$$
is commutative, where $\beta '_F$ denotes the Bott map in the $K'$-theory.
\item If
$$\oc{F}{\varphi }{G}{\psi }{H}$$
is an exact sequence in $\fr{M}_{E'}$ then the diagram
$$\begin{CD}
K_1(\bar{H})@>\delta_1>>K_0(\bar{F} )\\
@V\Phi _{1,H}VV         @VV\Phi _{0,F}V\\
K'_1(H)@>>\delta'_1>    K'_0(F)
\end{CD}$$
is commutative, where $\delta'_1$ denotes the index maps associated to the above exact sequences in the $K'$-theory.
\end{enumerate}
\end{p}

a) For every \bbn\, and
$$X=\si{t\in T_n}((\alpha _t,x_t)\otimes id_K)V_t^{f}\in \check{\bar{F} }_n\,, $$
$$\check{\varphi }_n\phi _{F,n}X=\si{t\in T_n}(((\phi \alpha _t),\varphi x_t)\otimes id_K) V_t^{f'}=\phi _{G,n}\check{\bar{\varphi} }_nX\;.  $$

b) For every \bbn\, and $P\in Pr\,\check{\bar{F} }_n, $
$$\phi _{SF,n}\widetilde{P}=(\widetilde{P} )'=\widetilde{P'}=\widetilde{\phi _{F,n}P}\;.$$

c) Let \bbn\, and $U\in \unn{\check{\bar{H} } }_{n-1}$. By \pr{960} a), there are $V\in \unn{\check{\bar{G} } }_n$ and $P\in \pp{\check{\bar{F} } }_n$
such that
$$\check{\bar{\psi } }_nV=A_nU+B_nU^*\,,\qquad\qquad \check{\bar{\varphi } }_nP=VA_nV^*\;.  $$
Then
$$\check{\psi }_n\phi _{G,n}V=\phi _{H,n}\check{\bar{\psi } }_nV=A'_n(\phi _{H,n-1}U)+B'_n(\phi _{H,n-1}U)^*\,,  $$
$$\check{\varphi }_n\phi _{F,n}P=\phi _{G,n}\check{\bar{\varphi } }_nP=(\phi _{G,n}V)A'_n(\phi _{G,n}V)^*  $$
so by \cor{964},
$$\delta '_1\Phi _{1,H}[U]_1=\delta '_1[\phi _{H,n-1}U]_1=[\phi _{F,n}P]_0=\Phi _{0,F}[P]_0=\Phi _{0,F}\delta _1[U]_1$$
$$\delta '_1\circ \Phi _{1,H}=\Phi _{0,F}\circ \delta _1\;.\qedd$$

\begin{lem}\label{30.4}
Let $F,G$ be C*-algebras, $\mac{\varphi }{F}{G}$ a surjective C*-homomorphism, and
$$\mae{\psi }{\ccb{[0,1]}{F}}{\ccb{[0,1]}{G}}{x}{\varphi \circ x}\;.$$
\begin{enumerate}
\item $\psi $ is surjective.
\item Assume $F$ unital and let $v\in \unn{\ccb{[0,1]}{G}}$ such that there is an $x\in \unn{F}$ with $\varphi x=v(0)$. Then there is a $u\in \unn{\ccb{[0,1]}{F}}$ with $\psi u=v$ and $u(0)=x$.
\end{enumerate}
\end{lem}

a) Let $y$ be an element of $ \ccb{[0,1]}{G}$ which is piecewise linear, i.e. there is a family 
$$0=s_1<s_2<\cdots<s_{n-1}<s_n=1$$
such that for every $i\in \bnn{n-1}$ and $t\in [0,1]$,
$$y((1-t)s_i+ts_{i+1})=(1-t)y(s_i)+ty(s_{i+1})\;.$$
Since $\varphi $ is surjective, there is a family $(x_i)_{i\in \bnn{n}}$ in $F$ with $\varphi x_i=y(s_i)$ for every $i\in \bnn{n}$. Define $\mac{x}{[0,1]}{F}$ by putting
$$x((1-t)s_i+ts_{i+1}):=(1-t)x_i+tx_{i+1}$$
for every $i\in \bnn{n-1}$ and $t\in [0,1]$. For $i\in \bnn{n-1}$ and $t\in [0,1]$,
$$(\psi x)((1-t)s_i+ts_{i+1})=\varphi ((1-t)x_i+tx_{i+1})=$$
$$=(1-t)y(s_i)+ty(s_{i+1})=y((1-t)s_i+ts_{i+1})\,,$$
so $\psi x=y$, $y\in Im\,\psi $. Since the set of elements of $\ccb{[0,1]}{G}$, which are piecewise linear, is dense in $\ccb{[0,1]}{G}$ and $Im\,\psi $ is closed (as C*-homomorphism), $\psi $ is surjective.

b) Let
$$\mae{w}{[0,1]}{Un\,G}{s}{v(0)^*v(s)}\;.$$
Then $w\in \unn{\ccb{[0,1]}{G}}$ and $w(0)=1_G$. Put
$$\mae{w_t}{[0,1]}{Un\,G}{s}{w(st)}$$
for every $t\in [0,1]$. Then
$$\mad{[0,1]}{\unn{\ccb{[0,1]}{G}}}{t}{w_t}$$
is a continuous path with $w_1=w$ and $w_0=1_{\ccb{[0,1]}{G}}$. Thus 
$$w\in \unm{\ccb{[0,1]}{G}}\;.$$ 
By a), $\psi $ is surjective, so by [R] Lemma 2.1.7 (i), there is a $u'\in \unn{\ccb{[0,1]}{F}}$ with $\psi u'=w$. Put
$$\mae{u}{[0,1]}{Un\,F}{s}{xu'(0)^*u'(s)}\;.$$
Then $u\in \unn{\ccb{[0,1]}{F}}$, $u(0)=x$, and 
$$(\psi u)(s)=\varphi (u(s))=\varphi (xu'(0)^*u'(s))=\varphi (x)((\psi u')(0))^*((\psi u')(s))=$$
$$=v(0)w(0)^*w(s)=v(0)1_Gv(0)^*v(s)=v(s) $$
for every $s\in [0,1]$, i.e. $\psi u=v$.\qed

\begin{theo}\label{1.5}
$\Phi _{i,F}$ is a group isomorphism for every $i\in \{0,1\}$ and for every $E'$-C*-algebra $F$. 
\end{theo}

By \pr{26.4} b), $\Phi _{0,F}=(\beta'_F)^{-1}\circ \Phi _{1,SF}\circ \beta _{\bar{F} }$, so it suffices to prove the assertion for $\Phi _{1,F}$ only. Let \bbn\, and $U\in \unn{\check{F}_n }$. Put $V:=U(\sigma _n^FU)^*\sim _1U$. Since $\sigma _n^FV=1_{E'}$, $V$ has the form
$$V=\si{t\in T_n}((\alpha _t,x_t)\otimes id_K)V_t^{f'}$$
with $\alpha _t=\delta _{1,t}1_{E'}$ and $x_t\in F$ for every $t\in T_n$. If we put
$$W:=\si{t\in T_n}((\delta _{1,t}1_E,x_t)\otimes id_K)V_t^f$$
then $\phi _{F,n}W=V$ and we get $\Phi _{1,F}[W]_1=[V]_1=[U]_1$, so $\Phi _{1,F}$ is surjective. Thus we have to prove the injectivity of $\Phi _{1,F}$ only.

Let $a\in Ker\;\Phi _{1,F}$. We have to prove $a=0$. There are \bbn\, and 
$$U:=\si{t\in T_n}((\alpha _t,x_t)\otimes id_K)V_t^f\in \unn{\check{\bar{F} }_n }$$
with $a=[U]_1$, where $(\alpha _t,x_t)\in \check{F} $ for every $t\in T_n$. Since $[U']_1=\Phi _{1,F}[U]_1=0$, by \pr{951a}, there is an $m\in \bn$ such that
$$U'_0:=\left(\proo{i=1}{m}A'_{n+i}\right)U'+\left(1_{E'}-\proo{i=1}{m}A'_{n+i}\right)$$
is homotopic in $\unn{\check{F}_{n+m} }$ to a $U'_1\in \unn{E'_{n+m}}\;(\subset \unn{\check{F}_{n+m} })$. Thus there is a continuous path
$$\mae{U'}{[0,1]}{\unn{\check{F}_{n+m} }}{s}{U'_s}\;.$$

\begin{center}
Case 1 $\phi $ is injective
\end{center}

Put
$$W'_s:=U'_s\sigma _{n+m}^F(U'^*_sU'_0)\;(\in \unn{\check{F}_{n+m} })$$
for every $s\in [0,1]$. Then
$$\sigma _{n+m}^FW'_s=\sigma _{n+m}^FU'_0=\phi _{F,n+m}\left(\left(\proo{i=1}{m}A_{n+i}\right)(\sigma _n^{\bar{F}}U )+\left(1_E-\proo{i=1}{m}A_{n+i}\right)\right)$$
for every $s\in [0,1]$. If we put
$$W'_s=:\si{t\in T_{n+m}}((\beta _{s,t},y_{s,t})\otimes id_K)V_t^{f'}\,,$$
where $(\beta _{s,t},y_{s,t})\in \check{F} $ for all $s\in [0,1]$ and $t\in T_n$, then
$$\si{t\in T_{n+m}}((\beta _{s,t},0)\otimes id_K)V_t^{f'}=\sigma _{n+m}^FW'_s=$$
$$=\phi _{F,n+m}\left(\left(\proo{i=1}{m}A_{n+i}\right)\si{t\in T_n}((\alpha _t,0)\otimes id_K)V_t^f+\left(1_E-\proo{i=1}{m}A_{n+i}\right)\right)$$
and so by [C2] Theorem 2.1.9 a), there is a (unique) family $(\gamma _t)_{t\in T_{n+m}}$ in $E$ with $\beta _{s,t}=\phi \gamma _t$ for every $s\in [0,1]$ and $t\in T_{n+m}$. Since $\phi $ is injective, $\phi _{n+m}$ is also injective and $\phi _{n+m}(\check{\bar{F} }_{n+m} )$ may be identified with a unital C*-subalgebra of $\check{F}_{n+m} $. Thus
$$\mae{W }{[0,1]}{\unn{\check{\bar{F} }_{n+m} }}{s}{\si{t\in T_{n+m}}((\gamma _t,y_{s,t})\otimes id_K)V_t^f}$$
is a continuous path in $\unn{\check{\bar{F} }_{n+m} }$ with $\phi _{F,n+m}W_s=W'_s $ for every $s\in [0,1]$. It follows
$$\phi _{F,n+m}W_0=W'_0=U'_0=\phi _{F.
,n+m}\left(\left(\proo{i=1}{m}A_{n+i}\right)U+\left(1_E-\proo{i=1}{m}A_{n+i}\right)\right)\,, $$
$$\phi _{F,n+m}W_1=W'_1=U'_1\sigma _{n+m}^F(U'^*_1U'_0)=\sigma _{n+m}^FU'_0\in \phi _{F,n+m}(\unn{E'_{n+m}})\;. $$
Since $\phi $ is injective, $\phi _{F,n+m}$ is also injective and we get
$$\left(\proo{i=1}{m}A_{n+i}\right)U+\left(1_E-\proo{i=1}{m}A_{n+i}\right)=W_0\,,$$
$$\left(\proo{i=1}{m}A_{n+i}\right)U+\left(1_E-\proo{i=1}{m}A_{n+i}\right) \in Un_{E_{n+m}}\,\check{\bar{F}}_{n+m}\,,\qquad g=[U]_1=0\;. $$

\begin{center}
Case 2 $\phi $ is surjective
\end{center}

We put
$$\bar{U}_0:=\left(\proo{i=1}{m}A_{n+i}\right)U+\left(1_E-\proo{i=1}{m}A_{n+i}\right)\;(\in \unn{\check{\bar{F} } }_{n+m})\;.$$
Since $\phi $ is surjective, $\phi _{F,n+m}$ is also surjective ([C2] Theorem 2.1.9 a)). Since 
$$\phi _{F,n+m}\bar{U}_0=U'_0$$
it follows from \lm{30.4} b), that there is a continuous path
$$\mad{[0,1]}{\unn{\check{\bar{F} }_{n+m} }}{s}{U_s}$$
with $\phi _{F,n+m}U_s=U'_s$ for every $s\in [0,1]$ and $U_0=\bar{U}_0$ . Since $\phi _{F,n+m}U_1=U'_1\in \unn{E'_{n+m}}$, we have $\bar{U}_0\in Un_{E_{n+m}}\check{\bar{F} }_{n+m}  $ and $g=[U]_1=[\bar{U}_0 ]_1=0$.

\begin{center}
Case 3 $\phi $ is arbitrary
\end{center}

There are a unital commutative C*-algebra $E''$ and a unital C*-homomor-phisms $\mac{\phi '}{E}{E''}$ and $\mac{\phi ''}{E''}{E'}$ such that $\phi '$ is surjective, $\phi ''$ is injective, and 
$\phi =\phi ''\circ \phi '$
and the assertion follows from the first two cases and the considerations from the begin of the section.\qed

\begin{co}\label{27.6}
Let $E',E''$ be unital commutative C*-algebras such that $E=E'\times E''$ and
$$\mae{\phi '}{E}{E'}{(x',x'')}{x'}\,,$$
$$\mae{\phi ''}{E}{E''}{(x',x'')}{x''}\;.$$
If $F'$ is an $E'$-C*-algebra and $F''$ is an $E''$-C*-algebra then the map (with obvious notation)
$$\mad{K_i(\Phi '(F')\times \Phi ''(F''))}{K'_i(F')\times K''_i(F'')}{a}{(\Phi '_{i,F'}\times \Phi ''_{i,F''})(\varphi _ia)}$$ 
is a group isomorphism for every $i\in \z{0,1}$, where
$$\mac{\varphi _i}{K_i(\Phi '(F')\times \Phi ''(F''))}{K_i(\Phi '(F'))\times K_i(\Phi ''(F''))}$$
is the canonical group isomorphism \emph{(Product Theorem (\cor{938} b), \pr{975a} b)))}.\qed
\end{co}

\begin{co}\label{9.7b}
If $f(s,t)\in \bc$ for all $s,t\in T$ and $C_n\in \bc_n$ for all \bbn and if $K^{\bc}$ denotes the K-theory with respect to $T$, $\bc$, $f$, and $(C_n)_{\bbn}$ then $K_i(E)=K^{\bc}_i(\ccb{\Omega }{\bc})$ for all $i\in \z{0,1}$, where $\Omega $ denotes the spectrum of $E$.\qed
\end{co}

\begin{p}\label{28.5a}
If $F$ is an $E'$-C*-algebra then the map
$$\mae{\varphi }{E\times \Phi (F)}{\check{ \overbrace{\Phi (F)}}}{(\alpha ,x)}{(\alpha ,x-\phi \alpha )}$$
is an $E$-C*-isomorphism. 
\end{p}

For $(\alpha ,x),(\beta ,y)\in E\times \Phi (F)$ and $\gamma \in E$,
$$\varphi (\gamma (\alpha ,x))=\varphi (\gamma \alpha ,(\phi \gamma )x)=(\gamma \alpha ,(\phi \gamma )x-\phi (\gamma \alpha ))=$$
$$=(\gamma ,0)(\alpha ,x-\phi \alpha )=(\gamma ,0)\varphi (\alpha ,x)\,,$$
$$\varphi (\alpha ,x)^*=\varphi (\alpha ^*,x^*)=(\alpha ^*,x^*-\phi \alpha ^*)=(\varphi (\alpha ,x))^*\,,$$
$$\varphi (\alpha ,x)\varphi (\beta ,y)=(\alpha ,x-\phi \alpha )(\beta ,y-\phi \beta )=$$
$$=(\alpha \beta ,(\phi \alpha )y-\phi (\alpha \beta )+(\phi \beta )x-\phi (\alpha \beta )+xy-(\phi \beta )x-(\phi \alpha )y+\phi (\alpha \beta ))=$$
$$=(\alpha \beta ,xy-\phi (\alpha \beta ))=\varphi (\alpha \beta ,xy)=\varphi ((\alpha ,x)(\beta ,y))\,,$$
so $\varphi $ is an $E$-C*-homomorphism. The other assertions are easy to see.\qed

\begin{center}
\section{Changing $f$}
\end{center}

\fbox{\parbox{12cm}{In all Propositions and Corollaries of this section we use the notation and assumptions of \ee{13.4} and $F$ denotes a C*-algebra}}

\begin{lem}\label{21.1}
For every \bbn\, there is an $\varepsilon _n>0$ such that for every $m\in \bn$, $m\leq n$, and $\alpha \in \unn{\bc}$, $|\alpha -1|<\varepsilon _n$, there is a unique $\beta _\alpha \in \unn{\bc}$, $|\beta _\alpha -1|<\frac{1}{n}$, with $\beta _\alpha ^m=\alpha $; moreover the map $\alpha \mapsto \beta _\alpha $ is continuous.
\end{lem}

If $\beta ,\gamma $ are distinct elements of $\unn{\bc}$ and $\beta ^m=\gamma ^m$ then
$$|\beta -\gamma |\geq |1-e^{\frac{2\pi i}{m}}|>\frac{1}{m}\geq \frac{1}{n}$$
and the assertion follows from the continuity of the corresponding branch of the map $\alpha \mapsto \sqrt[m]{\alpha }$.\qed

\begin{de}\label{1.2a}
For every finite group $S$ we endow $\f{S}{\bc}$ with the metric
$$d_S(g,h):=\sup\me{|g(s,t)-h(s,t)|}{s,t\in S}$$
for all $g,h\in \f{S}{\bc}$. 
\end{de}

{\it Remark.} $\f{S}{\bc}$ endowed with the above metric is compact.

\begin{de}\label{705}
We put
$$\Lambda (T,E):=\me{\mac{\lambda }{T}{\unn{E}}}{\lambda (1)=1_E}$$
and
$$\mae{\delta \lambda }{T\times T}{\unn{E}}{(s,t)}{\lambda (s)\lambda (t)\lambda (st)^*}$$
for every $\lambda \in \Lambda (T,E)$.
\end{de}

\begin{lem}\label{1.2b}
Let $S$ be a finite group and $\Omega $ a compact space.  
\begin{enumerate}
\item $\me{\delta \lambda }{\lambda \in \Lambda (S,\bc)}$ is an open set of $\f{S}{\bc}$. 
\item For every $\varepsilon'>0$ there is an $\varepsilon >0$ such that for all $g,h\in \f{S}{\ccb{\Omega }{\bc}}$, if 
$$\n{g(s,t)-h(s,t)}<\varepsilon $$
for all $s,t\in S$ then there is a $\lambda \in \Lambda (S,\bc)$ such that $h=g\delta \lambda $ and $|\lambda (s)-1|<\varepsilon' $ for all $s\in S$.
\item Let $g\in \f{S}{\ccb{\Omega }{\bc}}$ and $\mac{\phi }{[0,1]\times \Omega }{\Omega }$ a continuous map. We put for every $u\in [0,1]$,
$$\mac{\phi _u:=\phi (u,\cdot )}{\Omega }{\Omega }\,,$$
$$\mae{g_u}{S\times S}{\unn{\bc}}{(s,t)}{g(s,t)\circ \phi _u}\;.$$
Then $g_u\in \f{S}{\ccb{\Omega }{\bc}}$ for every $u\in [0,1]$ and there is a $\lambda \in \Lambda (S,\bc)$ with $g_1=g_0\delta \lambda $.
\end{enumerate}
\end{lem}

a) By [K] Theorem 2.3.2 (iii), 
$$\me{\ssa{g}}{g\in \f{S}{\bc}}/\approx _{\ccc{S}}$$
is finite.
$\me{\delta \lambda }{\lambda \in \Lambda (S,\bc)}$ is obviously a closed subgroup of $\f{S}{\bc}$. By the above and [C2] Proposition 2.2.2 c), $\f{S}{\bc}$ is the union of a finite family of closed pairwise disjoint sets homeomorphic to $\me{\delta \lambda }{\lambda \in \Lambda (S,\bc)}$, so $\me{\delta \lambda }{\lambda \in \Lambda (S,\bc)}$ is open. 

b) By a), there is an $\varepsilon >0$ such that for all $g',h'\in \f{S}{\bc}$ with $d_S(g',h')<\varepsilon $ there is a $\lambda \in \Lambda (S,\bc)$ with $h'=g'\delta \lambda $. We may assume that 
$$(1+\varepsilon )^{Card\,S}-1<\varepsilon _{Card\,S}\,,$$
where $\varepsilon _{Card\,S}$ was defined in \lm{21.1}.

We put for every $\omega \in \Omega $
$$\mae{g_\omega }{S\times S}{\unn{\bc}}{(s,t)}{(g(s,t))(\omega )}\,,$$
$$\mae{h_\omega }{S\times S}{\unn{\bc}}{(s,t)}{(h(s,t))(\omega )}\;.$$

Let $\omega \in \Omega $. By the above, there is a $\lambda _\omega \in \Lambda (S,\bc)$ with $g_\omega =h_\omega \delta \lambda _\omega $. Let $s\in S$ and let \bbn\, be the least natural number with $s^n=1_S$. By [C2] Proposition 3.4.1 c),
$$\lambda _\omega (s)^n=\proo{j=1}{n-1}(g_\omega (s^j,s)^*h_\omega (s^j,s))\;.$$
For every $j\in \bn_{n-1}$,
$$\n{1_E-g(s^j,s)^*h(s^j,s)}=\n{g(s^j,s)-h(s^j,s)}<\varepsilon \,,$$
$$\n{\proo{j=1}{n-1}(g(s^j,s)^*h(s^j,s))}=\n{\proo{j=1}{n-1}(1_E-(1_E-g(s^j,s)^*h(s^j,s)))}<(1+\varepsilon )^n\,,$$
$$\n{1_E-\proo{j=1}{n-1}(g(s^j,s)^*h(s^j,s))}<(1+\varepsilon )^{n-1}-1<\varepsilon _{Card\,S}\;.$$
By \lm{21.1}, there is a unique $\gamma \in \unn{\bc}$ with
$$\gamma ^n=\proo{j=1}{n-1}(g(s^j,s)^*h(s^j,s))\,,\qquad\qquad |\gamma -1|<\frac{1}{Card\,S}\;.$$
For $\omega \in \Omega $, since $|1-\lambda _\omega (s)|<\varepsilon _{Card\,S}$, we get $\lambda _\omega (s)=\gamma (s)$. So if we put
$$\mae{\lambda (s)}{\Omega }{\bc}{\omega }{\gamma (s)}$$
we have $\lambda \in \Lambda (S,\bc)$ and $g=h\delta \lambda $. By \lm{21.1}, we may choose $\varepsilon $ in such a way that the inequality $|\lambda (s)-1|<\varepsilon '$ holds for all $s\in S$.

c) By b), there is a family $(\lambda_i )_{i\in \bnn{n}}$ in $\Lambda (S,\bc)$ and
$$0=u_0<u_1<\cdots<u_{n-1}<u_n=1$$
such that $g_{u_i}=g_{u_{i-1}}\delta \lambda _i$ for every $i\in \bnn{n}$. By induction $g_0\delta \left(\proo{i=1}{j}\lambda _i\right)=g_{u_j}$ for every $j\in \bnn{n}$. Thus if we put $\lambda :=\proo{i=1}{n}\lambda _i$ then $g_0\delta \lambda=g_1 $\qed

{\it Remark.} Let $\lambda \in \Lambda (T,E)$ and $f'=f\delta \lambda \;(\in \f{T}{E})$. For every full $E$-C*-algebra $F$ and \bbn\, we denote by $F'_n$ the equivalent of $F_n$ constructed with respect to $f'$ instead of $f$ (\dd{a}). By [C2] Proposition 2.2.2 $a_1\Rightarrow a_2$, there is for every \bbn\, a unique $E$-C*-isomorphism $\mac{\varphi _n^F}{F_n}{F'_n}$ such that for all $m,n\in \bn$, $m<n$, the diagram
$$\begin{CD}
F_m@>\varphi_m^F>>F'_m\\
@VVV         @VVV\\
F_n@>>\varphi _n^F>F'_n
\end{CD}$$
is commutative, where the vertical arrows are the canonical inclusions. We put $C'_n:=\varphi _n^EC_n$ for evrey \bbn. $(C'_n)_{\bbn}$ satisfies the conditions of \axi{b} with respect to $f'$, so we can construct a K-theory with respect to $T$, $E$, $f'$, and $(C'_n)_{\bbn}$, which we shall denote by $K^{f'}$. If $m,n\in \bn$, $m<n$, then the diagrams

\parbox{5cm}{
$$\begin{CD}
F_m@>\rho _{n,m}^F>>F_n\\
@V\varphi _m^FVV         @VV\varphi _n^FV\\  
F'_m@>>\rho _{n,m}^{F'}>F'_n
\end{CD}$$}
\hspace{2cm}
\parbox{2cm}{
$$\begin{CD}
\unn{F_m}@>\tau _{n,m}^F>>\unn{F_n}\\
@V\varphi _m^FVV         @VV\varphi _n^FV\\  
\unn{F'_m}@>>\tau _{n,m}^{F'}>\unn{F'_n}
\end{CD}$$}

\hspace{-4ex}are commutative and so we get the isomorphisms 
$$\pp{F_\rightarrow }\longrightarrow \pp{F'_\rightarrow }\,,\qquad\qquad un\;F_\leftarrow \longrightarrow un\;F'_\leftarrow\;. $$
By these considerations it can be followed that $K$ and $K^{f'}$ coincide.

\begin{de}\label{31.5}
Let $\Omega $ be the spectrum of $E$, $\Gamma $ a closed set of $\Omega $, and $F$ a C*-algebra. We denote by $\ccb{E;\,\Gamma }{F}$ the $E$-C*-algebra obtained by endowing the C*-algebra $\ccb{\Gamma }{F}$ with the structure of an $E$-C*-algebra by putting
$$\mae{\alpha x}{\Gamma }{F}{\omega }{\alpha (\omega )x(\omega )}$$
for all $(\alpha ,x)\in E\times \ccb{\Gamma }{F}$. If $\Omega '$ is an open set of $\Omega $ then the ideal and $E$-C*-subalgebra 
$$\me{x\in \ccb{E;\Omega }{F}}{x|(\Omega \setminus \Omega ')=0}$$
of $\ccb{E;\,\Omega }{F}$ will be denoted $\cbb{E;\,\Omega' }{F}$.
\end{de}

By Tietze's theorem
$$\oc{\cbb{E;\,\Omega '}{F}}{\varphi }{\ccb{E;\,\Omega }{F}}{\psi }{\ccb{E;\,\Omega \setminus \Omega '}{F}}$$
is an exact sequence in \frm, where $\varphi $ denotes the inclusion map and
$$\mae{\psi }{\ccb{E;\,\Omega }{F}}{\ccb{E;\,\Omega \setminus \Omega '}{F}}{x}{x|(\Omega \setminus \Omega ')}\;.$$

\begin{p}\label{3.5}
We denote by $\Omega $ the spectrum of $E$, by $\Gamma$ a closed set of $\Omega $, and by $\mac{\vartheta }{[0,1]\times \Omega }{\Omega }$ a continuous map such that
$$\omega \in \Omega \Longrightarrow \vartheta (0,\omega )=\omega \,,\;\vartheta (1,\omega )\in \Gamma $$
and $\vartheta (s,\omega )=\omega $ for all $s\in [0,1]$ and $\omega \in \Gamma $.
We put $E':=\ccb{\Gamma }{\bc}$, $E'':=E$, $\vartheta _s:=\vartheta (s,\cdot )$ for every $s\in [0,1]$, and
$$\mae{\phi }{E}{E'}{x}{x|\Gamma}\,,\qquad\qquad \mae{\phi '}{E'}{E''=E}{x'}{x'\circ \vartheta_1}\,,$$
$$\mae{f'}{T\times T}{\unn{E'}}{(s,t)}{\phi f(s,t)=f(s,t)|\Gamma}\,,$$
$$\mae{f''}{T\times T}{\unn{E''}}{(s,t)}{\phi 'f'(s,t)=f(s,t)\circ \vartheta_1}\;.$$
\begin{enumerate}
\item There is a $\lambda \in \Lambda (T,E)$ such that $f''=f\delta \lambda $ and the K-theories associated to $f$ and $f''$ coincide (as formulated in the above \emph{Remark}). If $\Gamma$ is a one-point set (i.e. $\Omega $ is contractible) then $f''(s,t)\in \unn{\bc}\;(\subset \unn{E})$ for all $s,t\in T$.
\item If we put 
$$\mae{\psi }{\ccb{E;\,\Omega }{F}}{\ccb{E;\,\Gamma }{F}}{x}{x|\Gamma }$$
then $K_i(\cbb{E;\,\Omega \setminus \Gamma }{F})=\{0\}$ and
$$\mac{K_i(\psi )}{K_i(\ccb{E;\Omega }{F})}{K_i(\ccb{E;\Gamma }{F})}$$
is a group isomorphism
for every $i\in \{0,1\}$.
\item If $\Gamma '$ is a compact subspace of $\Omega \setminus \Gamma $ then
$$K_i(\cbb{E;\,\Omega \setminus (\Gamma \cup \Gamma ')}{F})\approx K_{i+1}(\ccb{E;\,\Gamma '}{F})$$
for all $i\in \{0,1\}$.
\item Let $\bar{\Gamma } $ be a closed set of $\Omega $, $\mac{\bar{\varphi } }{\cbb{E;\,\Omega \setminus (\Gamma \cup \bar{\Gamma } )}{F}}{\ccb{E;\,\Omega }{F}}$ the inclusion map,
$$\mae{\bar{\psi } }{\cbb{E;\,\Omega }{F}}{\ccb{E;\,\Gamma \cup \bar{\Gamma } }{F}}{x}{x|(\Gamma \cup \bar{\Gamma } )}\,,$$
and $\delta _0,\delta _1$ the corresponding maps from the six-term sequence associated to the exact sequence in \frm
$$\oc{\cbb{E;\,\Omega \setminus (\Gamma \cup \bar{\Gamma } )}{F}}{\bar{\varphi } }{\ccb{E;\,\Omega }{F}}{\bar{\psi } }{\ccb{E;\,\Gamma \cup \bar{\Gamma } }{F}}$$
 then the sequence
$$\ob{0\longrightarrow K_i(\ccb{E;\,\Omega }{F})}{K_i(\psi )}{K_i(\ccb{E;\,\Gamma \cup \bar{\Gamma } }{F})}{\delta _i}\\$$
$${\stackrel{\delta _i}{\longrightarrow} K_{i+1}(\cbb{E;\,\Omega \setminus (\Gamma \cup \bar{\Gamma } )}{F})}\longrightarrow 0$$
is exact for every $i\in \z{0,1}$.
\end{enumerate}
\end{p}

a) By \lm{1.2b} c), for every $m\in \bn$ there is a $\lambda _m\in \Lambda (S_m,E)$ with $f''|(S_m\times S_m)=g_m\delta \lambda _m$. We put
$$\mae{\lambda }{T}{\unn{E}}{t}{\lambda _m(t)\quad \mbox{if}\quad t\in S_m}\;.$$
Then
$$f''(s,t)=\pro{m\in \bn}(g_m\delta \lambda )(s_m,t_m)=(f\delta \lambda )(s,t)$$
for all $s,t\in T$, i.e. $f''=f\delta \lambda $.

b) Let \bbn\, and $X\in \left(\check{\overbrace{\cbb{E'';\,\Omega \setminus \Gamma }{F}}} \right)_n$. Then $X$ has the form
$$X=\si{t\in T_n}((\alpha _t,x_t)\otimes id_K)V_t^{f''}\,,$$
where $\alpha _t\in E''$ and $x_t\in \cbb{E'';\,\Omega \setminus \Gamma }{F}$ for all $t\in T_n$. We put
$$X_s:=\si{t\in T_n}((\alpha _t\circ \vartheta_s,x_t\circ \vartheta_s)\otimes id_K)V_t^{f''}$$
for every $s\in [0,1]$. Then
$$\mad{[0,1]}{\left(\check{\overbrace{\cbb{E'';\,\Omega \setminus \Gamma }{F}}} \right)_n}{s}{X_s}$$
is a continuous map, $X_0=X$,
$$X_1=\si{t\in T}((\alpha _t\circ \vartheta_1,0)\otimes id_K)V_t^{f''}\,,$$
and
$$\mad{\left(\check{\overbrace{\cbb{E'':\,\Omega \setminus \Gamma }{F}}} \right)_n}{\left(\check{\overbrace{\cbb{E'';\,\Omega \setminus \Gamma }{F}}} \right)_n}{X}{X_s}$$
is an $E''$-C*-homomorphism for every $s\in [0,1]$. Thus $K_i^{f''}(\cbb{E'';\,\Omega \setminus \Gamma }{F})=\{0\}$. By a), $K_i(\cbb{E;\,\Omega \setminus \Gamma }{F})=\{0\}$. 

If $\mac{\varphi }{\cbb{E;\Omega \setminus \Gamma }{f}}{\ccb{E;\Omega }{F}}$ denotes the inclusion map then
$$\oc{\cbb{E;\,\Omega \setminus \Gamma }{F}}{\varphi }{\ccb{E;\,\Omega }{F}}{\psi }{\ccb{E;\,\Gamma }{F}}$$
is an exact sequence in \frm and the assertion follows from the six-term sequence (\cor{93} c)).

c) If we put
$$F_1:=\cbb{E;\,\Omega \setminus (\Gamma \cup \Gamma ')}{F}\,,\quad F_2:=\cbb{E;\,\Omega \setminus \Gamma }{F}\,,\quad F_3:=\ccb{E;\,\Gamma '}{F}\,,$$
$$\mae{\varphi }{F_1}{F_2}{x}{x}\,,$$
$$\mae{\psi }{F_2}{F_3}{x}{x|\Gamma '}$$
then
$$\oc{F_1}{\varphi }{F_2}{\psi }{F_3}$$
is an exact sequence in \frm and the assertion follows from b) and from the six-term sequence (\cor{93} d)).

d) $\bar{\varphi } $ factorizes through $\cbb{E;\,\Omega \setminus \Gamma }{f}$ so by b), $K_i(\bar{\varphi } )=0$ and the assertion follows from the six-term sequence \cor{93} b).\qed

\begin{co}\label{14.6}
We use the notation of \emph{\pr{3.5}}. Let $\bar{\Omega}$ be a compact space and $\mac{\bar{\vartheta} }{\Omega }{\bar{\Omega}}$ a continuous map such that the induced maps $\Omega \setminus (\Gamma \cup \Gamma ')\rightarrow \bar{\Omega}\setminus \bar{\vartheta} (\Gamma \cup \Gamma ')$, $\Gamma \rightarrow \bar{\vartheta} (\Gamma )$, and $\Gamma'\rightarrow \bar{\vartheta} (\Gamma')$ are homeomorphisms. If we put    $\bar{E}:=\ccb{\bar{\Omega}}{\bc}$ and
$$\mae{\bar{\phi} }{\bar{E}}{E}{x}{x\circ \bar{\vartheta} }$$
and take an $\bar{f} \in \f{T}{\bar{E}}$ such that $f(s,t)=\bar{\phi} \bar{f} (s,t)$ for all $s,t\in T$ and a corresponding $(\bar{C}_n )_{\bbn}\in \pro{\bbn}\bar{E}_n $ then with the notation from the beginning of \emph{section 9.1 (with $E$ and $\bar{E}$ interchanged)} 
$$\bar{K}_i\left(\cbb{\bar{E};\,\bar{\Omega}\setminus \bar{\vartheta} (\Gamma \cup \Gamma ')}{F}\right)\approx \bar{K}_{i+1}\left(\ccb{\bar{E};\,\bar{\vartheta} (\Gamma ')}{F}\right)\,,$$
for all $i\in \{0,1\}$, where $\bar{K} $ denotes the $K$-theory associated to $T$, $\bar{E}$, $\bar{f} $, and $(\bar{C}_n )_{\bbn}$. If in addition $\Gamma '$ has the same property as $\Gamma $ then 
$$\bar{K}_i\left(\ccb{\bar{E};\,\bar{\vartheta} (\Gamma) }{F}\right)\approx \bar{K}_i\left(\ccb{\bar{E};\,\bar{\vartheta} (\Gamma')}{F}\right)\;.$$
\end{co}

By our hypotheses,
$$\bar{\Phi} \left(\cbb{E;\,\Omega\setminus (\Gamma \cup \Gamma ')}{F}\right)\approx \cbb{\bar{E};\,\bar{\Omega}\setminus \bar{\vartheta} (\Gamma \cup \Gamma ')}{F}\,,$$
$$\bar{\Phi}(\ccb{E;\,\Gamma }{F})\approx \ccb{\bar{E};\,\bar{\vartheta} (\Gamma) }{F}\,,\qquad \bar{\Phi} \left(\ccb{E;\,\Gamma' }{F}\right)\approx \ccb{\bar{E};\,\bar{\vartheta} (\Gamma') }{F}\,, $$
 so by \pr{3.5} b) and \h{1.5}, 
$$\bar{K}_i\left(\cbb{\bar{E};\,\bar{\Omega}\setminus \bar{\vartheta} (\Gamma \cup \Gamma ')}{F}\right)\approx K_i\left(\cbb{E;\,\Omega\setminus (\Gamma \cup \Gamma ')}{F}\right)\approx$$
$$\approx K_{i+1}\left(\ccb{E;\,\Gamma '}{F})\approx \bar{K}_{i+1}(\ccb{\bar{E};\,\bar{\vartheta} (\Gamma ')}{F}\right)\;.$$
If the supplementary hypothesis is fulfilled then by \pr{3.5} c) and \h{1.5},
$$\bar{K}_i\left(\ccb{\bar{E};\,\bar{\vartheta} (\Gamma) }{F}\right)\approx K_i(\ccb{E;\,\Gamma}{F})\approx$$
$$\approx  K_i\left(\ccb{E;\,\Gamma') }{F}\right)\approx \bar{K}_i\left(\ccb{\bar{E};\,\bar{\vartheta} (\Gamma') }{F}\right)\;.\qedd$$

\begin{co}\label{12.6}
Assume $E=\ccb{\bt}{\bc}$.  
\begin{enumerate}
\item If $\theta _1,\,\theta _2,\,\theta _3,\,\theta _4\in \br$ such that $\theta _1\leq \theta _2<\theta _1+2\pi $, $\theta _3\leq \theta _4<\theta _3+2\pi $ then
$$K_i\left(\ccb{E;\,\me{e^{i\theta }}{\theta _1\leq \theta \leq \theta _2}}{F}\right)\approx$$
$$\approx  K_i\left(\ccb{E;\,\me{e^{i\theta }}{\theta _3\leq \theta \leq \theta _4}}{F}\right) $$
for every $i\in \z{0,1}$.
\item Let $\theta _1,\theta _2\in \br$, $\theta _1\leq \theta _2<\theta _1+2\pi $ and let $\Gamma $ be a closed set of 
$$\bt\setminus \me{e^{i\theta }}{\theta _2<\theta < \theta _1+2\pi }$$
 such that $e^{i\theta _1}\in \Gamma $ and $e^{i\theta _2}\not\in \Gamma$ if $e^{i\theta _1}\not=e^{i\theta _2}$. Then
$$K_i(\cbb{E;\,\bt\setminus \Gamma }{F})\approx K_{i+1}(\ccb{E;\;\Gamma} {F})$$
for every $i\in \z{0,1}$. Moreover
$$K_i(\cbb{E;\,\bt\setminus \Gamma }{F})\approx \ab{K_{i+1}(\ccb{E;\;\z{1}} {F})^{\Gamma }}{F\; \emph{is finite}}{\si{\bbn}K_{i+1}(\ccb{E;\,\z{1}}{F})}{F\; \emph{is infinite}}\;.$$
\item If $\Gamma _1,\,\Gamma _2$ are closed sets of $\bt$, not equal to $\bt$ and such that their cardinal numbers are equal if they are finite then
$$K_i(\ccb{E;\,\Gamma _1}{F})\approx K_i(\ccb{E;\,\Gamma _2}{F})$$
for all $i\in \z{0,1}$.
\end{enumerate}
\end{co}

a) We may assume $\theta _1\leq \theta _3<\theta _1+2\pi $. Put $\Omega ':=[\theta _1,\sup{(\theta _2,\theta _3)}]$, $E':=\ccb{\Omega '}{\bc}$,
$$\mae{\vartheta }{\Omega '}{\bt}{\alpha}{ e^{i\alpha }}\,,$$
$$\mae{\phi }{E}{E'}{x}{x\circ \vartheta }\;.$$  
Since it is possible to find an $f'\in \f{T}{E'}$ and a $(C'_n)_{\bbn}\in \pro{\bbn}E'_n$ with the desired properties, we get
$$K_i\left(\ccb{E;\,\me{e^{i\theta} }{\theta _1\leq \theta \leq \theta _2}}{F}\right)\approx K_i\left(\ccb{E;\,\z{e^{i\theta _3}}}{F}\right)\;.$$
by \cor{14.6}. Thus
$$K_i\left(\ccb{E;\,\me{e^{i\theta} }{\theta _3\leq \theta \leq \theta _4}}{F}\right)\approx K_i\left(\ccb{E;\,\z{e^{i\theta _3}}}{F}\right)\,,$$
$$K_i\left(\ccb{E;\,\me{e^{i\theta} }{\theta _1\leq \theta \leq \theta _2}}{F}\right)\approx$$
$$\approx  K_i\left(\ccb{E;\,\me{e^{i\theta} }{\theta _3\leq \theta \leq \theta _4}}{F}\right)\;.$$

b) If we put $\Omega ':=[\theta _1,\theta _1+2\pi ]$, $E':=\ccb{\Omega '}{\bc}$,
$$\mae{\vartheta }{\Omega '}{\bt}{\alpha }{e^{i\alpha }}\,,$$
$$\mae{\phi }{E}{E'}{x}{x\circ \vartheta }\,,$$
then the first assertion follows from \cor{14.6}. If $\Gamma $ is finite then the last assertion follows now from a) (and \cor{938} b) and \pr{975} b)).

Assume now $\Gamma $ infinite. Then $\Omega _0:=\bt\setminus \Gamma $ is the union of a countable set of open intervals. Let $\Xi $ be the set of  finite such intervals ordered by inclusion and for every $\Theta \in \Xi $ let $\Omega _\Theta $ be the union of the intervals of $\Theta $ and $\Gamma _\Theta :=\bt\setminus \Omega _\Theta $. By the above,
$$K_i(\cbb{E;\,\bt\setminus \Gamma _\Theta }{F})\approx K_{i+1}(\ccb{E;\,\z{1}}{F})^{\Theta }$$
for every $\Theta \in \Xi $. We get an inductive system of {E}-modules with $\cbb{E;\,\bt\setminus \Gamma }{F}$ as inductive limit. By \h{941} and \h{17.2}, $K_i(\cbb{E;\,\bt\setminus \Gamma }{F})$ is the inductive limit of $K_i(\cbb{E;\,\bt\setminus \Gamma _\Theta }{F})$ for $\Theta $ running through $\Xi $, which proves the assertion.

c) follows from b).\qed
 
 {\it Remark.} Let $\delta _0$ and $\delta _1$ be the group homomorphisms from the six-term sequence associated to the exact sequence in \frm
 $$\oc{\cbb{E;\,\bt\setminus \Gamma }{F}}{}{\ccb{E;\,\bt}{F}}{}{\ccb{E;\,\Gamma }{F}}\;.$$
 Then $\delta _0$ and $\delta _1$ do not coincide with the group isomorphism
$$K_i(\cbb{E;\,\bt\setminus \Gamma }{F})\approx K_{i+1}(\ccb{E;\;\Gamma} {F})$$
from \cor{12.6} b). 

\begin{co}\label{17.6}
If $\Omega $ is a compact space such that $E=\ccb{\Omega \times \bt}{\bc}$ then
$$K_i(\cbb{E;\,\Omega \times (\bt\setminus \z{1})}{F})\approx K_{i+1}(\ccb{E;\,\Omega \times \z{1}}{F})$$
for every $i\in \{0,1\}$.\qed
\end{co}

\begin{co}\label{22.5}
If the spectrum of $E$ is $\bbb_n$ for some \bbn\, then $K_i(\cbb{E;\,\bbb_n\setminus \{0\}}{F})=\{0\}$ and 
$$K_i(\cbb{E;\,\me{\alpha \in \br^n}{0<\n{\alpha }<1}}{F})\approx K_{i+1}(\ccb{E;\,\bs_{n-1}}{F})$$
for every $i\in \{0,1\}$.\qed 
\end{co}

\begin{co}\label{18.6}
Let $(k_j)_{j\in J}$ be a finite family in \bn, $\Omega '$ the topological sum of the family of balls $(\bbb_{k_j})_{j\in J}$, and $\Omega $ the compact space obtained from $\Omega '$ by identifying the centers of theses balls. If $\omega $ denotes the point of $\Omega $ obtained by this identification and $S$ denotes the union of $(\bs_{k_j-1})_{j\in J}$ in $\Omega $ and if $E=\ccb{\Omega }{\bc}$ then
$$K_i\left(\cbb{E;\,\Omega \setminus \z{\omega }}{F}\right)=\z{0}\,,$$
$$K_i(\cbb{E;\,(\Omega \setminus (\z{\omega }\cup S)}{F})\approx K_{i+1}(\ccb{E;\,S}{F})$$
for every $i\in \z{0,1}$.
\end{co}

If we denote by $\mac{\vartheta }{\Omega '}{\Omega }$ the quotient map, by $\Gamma $ the subset of $\Omega '$ formed by the centers of the balls $(\bbb_{k_j})_{j\in J}$, and by $\Gamma '$ the union of $(\bs_{k_j-1})_{j\in J}$ ($\Gamma '\subset \Omega '$) then the assertions follow from \pr{3.5} b), c) and \cor{14.6}.\qed

\begin{lem}\label{10.5}
Let $S$ be a finite group, $g\in \f{S}{E}$, and $\Omega $ the spectrum of $E$.
\begin{enumerate}
\item If there is an $\omega _0\in \Omega $ and a family $(\theta (s,t))_{s,t\in S}$ of selfadjoint elements of $E$ such that
$$\theta (r,s)+\theta (rs,t)=\theta (r,st)+\theta (s,t)\,,\quad g(s,t)=e^{i\theta (s,t)}(g(s,t)(\omega _0))$$ 
for all $r,s,t\in S$ then there is a $\lambda \in \Lambda (S,\bc)$ with $(g\delta \lambda )(s,t)=g(s,t)(\omega _0)$ for all $s,t\in S$.
\item If $\Omega $ is totally disconnected then there is a $\lambda \in \Lambda (S,E)$ such that 
$$((g\delta \lambda )(s,t))(\Omega )$$ 
is finite for all $s,t\in S$.
\end{enumerate}
\end{lem}

a) For every $u\in [0,1]$ put
$$\mae{g_u}{S\times S}{\unn{E}}{(s,t)}{e^{iu\theta (s,t)}(g(s,t)(\omega _0))}\;.$$
Then
$$\mad{[0,1]}{\f{S}{E}}{u}{g_u}$$
is a continuous map with $g_1=g$ and $g_0(s,t)=g(s,t)(\omega _0)$ for all $s,t\in S$. By \lm{1.2b} a),b), there are
$$0=u_0<u_1<\cdots<u_{k-1}<u_k=1$$
and a family $(\lambda _j)_{j\in \bnn{k}}$ in $\Lambda (S,\bc)$ such that $g_{u_{j-1}}=g_{u_j}\delta \lambda _j$ for every $j\in \bnn{k}$. We prove by induction that 
$$g_{u_{l-1}}=g\proo{j=l}{k}\delta \lambda _j$$
for all $l\in \bnn{k}$. This is obvious for $l=k$. Assume the identity holds for  $l\in \bnn{k}$, $l>1$. Then
$$g\proo{j=l-1}{k}\delta \lambda _j=\left(g\proo{j=l}{k}\delta \lambda _j\right)\delta \lambda _{l-1}=g_{u_{l-1}}\delta \lambda _{l-1}=g_{u_{l-2}}\,,$$
which finishes the proof by induction. If we put
$$\lambda :=\proo{j=1}{k}\lambda _j\in \Lambda (S,\bc)$$
then by the above
$$g\delta \lambda =g\proo{j=1}{k}\delta \lambda _j=g_0\;.$$

b) Let $\omega _0\in \Omega $. Since $\Omega $ is totally disconnected and $S$ is finite, by continuity, there is a clopen neighborhood $\Omega _0$ of $\omega _0$ and a family $(\theta (s,t))_{s,t\in S}$ in $Re\,\ccb{\Omega _0}{\bc}$ such that
$$\theta (r,s)+\theta (rs,t)=\theta (r,st)+\theta (s,t)\,,\qquad g(s,t)|\Omega _0=e^{i\theta (s,t)}(g(s,t)(\omega _0))$$
for all $r,s,t\in S$. By a), there is a $\lambda \in \Lambda (S,\bc)$ with
$$((g|\Omega _0)\delta \lambda)(s,t)=g(s,t)(\omega _0)$$
for all $s,t\in S$.

The assertion follows now from the fact that there is a finite partition $(\Omega _j)_{j\in J}$ of $\Omega $ with clopen sets such that $\Omega _j$ possesses the property of the above $\Omega _0$ for every $j\in J$.\qed

\begin{p}\label{10.5a}
If the spectrum of $E$ is totally disconnected then there is a $\lambda \in \Lambda (T,E)$ such that $((f\delta \lambda )(s,t))(\Omega )$ is finite for all $s,t\in T$. 
\end{p}

By \lm{10.5} b), for every $m\in \bn$ there is a $\lambda_m\in \Lambda (S_m,E) $ such that $((g_m\delta \lambda _m)(s,t))(\Omega )$ is finite for all $s,t\in S_m$. If we put
$$\mae{\lambda }{T}{\unn{E}}{t}{\lambda _m(t)\quad \mbox{if}\quad t\in S_m}$$
then $\lambda $ has the desired properties.\qed

\begin{p}\label{20.4'}
Assume that $T$, $f$, and $(C_n)_{\bbn}$ satisfy the conditions of \emph{\ee{13.4}} and of its \emph{Remark 1} and that the spectrum $\Omega $ of $E$ is simply connected.
\begin{enumerate}
\item There is a $\lambda \in \Lambda (T,E)$ such that $(f\delta \lambda )(s,t)\in \bc$ for all $s,t\in T$.
\item If $\kk{1}{\ccb{\Omega }{\bc}}=\z{0}$ for the classical $K_1$ then $\kk{1}{E}=\z{0}$ for the present theory.
\end{enumerate}
\end{p}

a) follows from \lm{10.5} a).

b) follows from a), Remark 1 of \ee{13.4}, and \pr{9.6}.\qed

\backmatter
\pagestyle{empty}

\begin{center}
REFERENCES
\end{center}

\begin{flushleft}
[C1] Corneliu Constantinescu, {\it C*-algebras.} Elsevir, 2001. \newline
[C2] Corneliu Constantinescu, {\it Projective representations of groups using Hilbert right C*-modules.} Eprint arXiv: 1111.1910 [164 pages] (11/2011) \newline
[K] Gregory Karpilowsky, {\it Projective representations of finite groups.} Marcel Dekker, Pure and Applied Mathematics 94, 1985. \newline
[R] M. R{\o}rdam, F. Larsen, N. J. Lausten {\it An Introduction to K-Theory for C*-Algebras.} London Mathematical Society, Student Texts 49, 2000. \newline 
[W] N. E. Wegge-Olsen, {\it K-theory and C*-algebras.} Oxford University Press, 1993. \newline
\end{flushleft}

\vspace{2cm}

\begin{center}
SUBJECT INDEX
\end{center}

Alexandroff K-theorem (\h{20.4})

Bott map (\pr{984}, \pr{985})

Commutativity of the index maps (\axi{27.9'e})

Commutativity of the six-term index maps (\cor{29.6} a))

Continuity axiom (\axi{5.10'})

Continuity of $K_0$ (\h{941}) and $K_1$  (\h{17.2})

$E$-C*-algebra, $E$-C*-subalgebra, $E$-ideal, $E$-linear, $E$-C*-homomorphism,  $E$-C*-isomorphism (\dd{10.3'})

Factorizes through null (\dd{28.9'c})

Full $E$-C*-algebra, full $E$-C*-subalgebra (\dd{10.3'r})

Homotopic, homotopy (\dd{939}

Homotopy axiom (\axi{27.9'b})

Homotopy invariance (\h{940}, \pr{956})

Index map (\cor{964})

Index maps (\dd{27.9'c})

Klein bottle (\dd{28.1'a})

K-null (\dd{28.9'c})

\frm -triple (\pr{28.9'e})

M$\ddot{o}$bius band(\dd{27.1'})

Null-axiom (\axi{27.9'})

Null-homotopic (\dd{939}, \dd{2.9'})

Product Theorem (\pr{14.11}, \pr{975a}, \pr{28.3'})

Projective space (\dd{4.2'a})

Schur $E$-function for $S$ (\dd{703})

Six-term axiom (\axi{27.9'd})

Split exact axiom (\axi{27.9'a})

Split Exact Theorem (\pr{936}, \cor{987})

Stability axiom (\axi{5.10'b})

Stability for $K_0$ (\h{949})

Tietze's Theorem (\cor{23.11a})

Topological six-term sequence (\pr{24.11})

Topological triple (\pr{3.12})

The triple theorem (\h{28.9'f})

Unitization (\dd{3.10'})

\vspace{2cm}

\begin{center}
SYMBOL INDEX
\end{center}

$0$,  $\frm$, $\fr{M}_{\bc}$ (\dd{10.3'})

$\pro{j\in J}F_j$ (\dd{8.6}, \dd{10.3'r})

$K_0$, $K_1$, $0$ (\dd{28.9'c})

$\delta _0$, $\delta _1$ (\dd{27.9'c})

$\Phi _{(F_j)_{j\in J},i}$, $\Psi _{(F_j)_{j\in J},i}$ (\dd{22.10'})

\frm-triple, $\varphi _{j,k}$, $\psi _{j,k}$, $\delta _{j,k,i}$ (\pr{28.9'e})

$F\otimes G$, $\varphi \otimes \psi $, $\bigotimes\limits_{j\in \emptyset }G_j $ (\dd{26.3'})

$\tilde{G} $, $\iota _G$, $\pi_G$, $\lambda_G$, $\tilde{\varphi } $ (\dd{3.10'})

$\delta _{G,i}$ (\dd{26.8'})

$\Upsilon $, $p(G)$, $q(G)$, $\Phi _{i,G,F}$, $\Upsilon $-null, $\vec{G} $  (\dd{5.7'})

$G_\Upsilon $ (\dd{25.8'})

$\bc_\Upsilon $ (\pr{5.7'a} b))

 $\Upsilon _1$, $\phi _{G,F}$ (\dd{10.10'})
 
 $\ccb{\Omega }{F}\,,\;\cbb{\Omega }{F}$ (\dd{28.3'a})
 
 $\Omega \in \Upsilon $, $p(\Omega) $, $q(\Omega) $, $\Phi _{i,\Omega ,F}$, $\Omega _\Upsilon $, $\Omega \in \Upsilon _1$, $\Omega $ is $\Upsilon $-null, (\dd{2.9'})
 
 $\bbb_n$ (\dd{28.3'b})
 
 $S\hspace{-2mm}S_n$, $\bt$ (\dd{28.3'd})
 
 $\bp_n$ (\dd{4.2'a})
 
 $\bm$, $\Gamma ^{\bm}_j$ (\dd{27.1'})
 
 $\bk$, $\Gamma _j^{\bk}$ (\dd{28.1'a})
 
 $\frc$ (\dd{10.3'r})
 
 $\check{F} $ (\dd{10.3'a})
 
 $\iota ^G,\,\pi ^G,\,\lambda ^G,\,\sigma ^G $ (\dd{924})
 
 $\check{\varphi } $ (\pr{34})
 
 $\si{j\in J}$ (\dd{18.5})
 
 $M(n)$ (\dd{9.11'a})
 
  $h$ (\axi{5.10'b})
 
$\f{S}{E}$ (\dd{703})

 $V_t$, $V_t^F$, $x\otimes id_K$, $F_n$, $\varphi _n$ (\dd{a})
 
 $A_n,\,B_n,\,C_n$ (\axi{b})
 
 $\bar{\rho }_n^F $ (\pr{911})
 
 $\rho _{n,m}^F$, \,$F_\rightarrow $,\,  $\rho _n^F$,\,  $X_\rightarrow :=X_{\rightarrow \,n}:=X_{\rightarrow \,n}^F$, \, $1_{\rightarrow \,n}:=1_{\rightarrow \,n}^F$, \, $F_{\rightarrow \,n}$, \, $Pr\,F_\rightarrow $, \, $\sim _0$, \, $\dot{P} $ (\dd{912})
 
 $\oplus $, \,$K_0(F)$, \,$[\;\cdot \;]_0$ (\pr{916},\,\,\dd{926})
 
 $\varphi _\rightarrow $ (\pr{921} a))
 
 $\bar{\tau }_n^F $ (\pr{950})
 
 $\tau _{n,m}^F,\,\tau _n^F,\,un\,F,\,un\,_E\,F,\,\unn{F_{\leftarrow \,n}},\,U_\leftarrow ,\,U_{\leftarrow \,n},\,U_{\leftarrow \,n}^F,\,1_{\leftarrow \,n},\,1_{\leftarrow \,n}^F $ (\dd{951})
 
 $\sim _1$ (\pr{950} \pr{952} c))
 
 $K_1(F),\,\oplus ,\,[\cdot ]_1$ (\dd{953})
 
 $\check{\varphi }_\leftarrow ,\,K_1(\varphi ) $ (\pr{954})
 
 $\delta _1$ (\cor{964})
 
 $CF,\,SF,\,\theta _F,\,i_F,\,j_F,\,C_\varphi ,\,S_\varphi $ (\dd{975})
 
  $\widetilde{P} $ (\dd{982})
 
 $\nu _F$ (\pr{983})
 
 $\beta _F$ (\pr{984}, \pr{985})
 
 $Trig(n),\,Pol(n,m),\,Pol(n),\,Lin(n),\,Proj(n)$ (\dd{76})
 
 $|p|,\,(p,q)_i$ (\dd{79})
 
 $\delta _0$ (\cor{93} b))
 
 $\Phi (F)\,,\;\Phi _{i,F}\,,\;K^{f'}$ (Introduction to section 9.1)
 
 $\f{S}{\bc}$ (\dd{1.2a})
 
 $\Lambda (S,E)$\,,\,$\delta \lambda $ (\dd{705})
 
 $\ccb{E;\,\Gamma }{F}\,,\;\cbb{E;\,\Omega'}{F}$ (\dd{31.5}) 
\end{document}